\documentclass[12pt]{article}
\usepackage{amsthm,amsfonts,amsmath,amssymb,latexsym,epsfig}
\input{xy}
\xyoption{all}
\input{epsf.tex}
\title{Gromov--Witten invariants\\ of symplectic quotients
       and adiabatic limits}
 
\author{A.~Rita Gaio \\
         Fac. Ci\^encias-Porto
         \and
         Dietmar~A.~Salamon\\
         ETH-Z\"urich
         }
 
\date{19 June 2001}

\newcommand{\MAT}[1]{\left[\begin{array}{#1}}
\newcommand{\RIX}{\end{array}\right]} 
\newcommand{\p}{\partial}
\newcommand{\dslash}{/\mskip-6mu/}
\newcommand{\one}{{{\mathchoice {\rm 1\mskip-4mu l} {\rm 1\mskip-4mu l}
{\rm 1\mskip-4.5mu l} {\rm 1\mskip-5mu l}}}}
\newcommand{\NN}{{\mathbb N}}
\newcommand{\Z}{{\mathbb Z}}
\newcommand{\ZZ}{{\mathbb Z}}

\newcommand{\R}{{\mathbb R}}

\newcommand{\C}{{\mathbb C}}

\newcommand{\Acal}{{\mathcal A}}   
\newcommand{\Aa}{{\mathcal A}}
\newcommand{\Bcal}{{\mathcal B}}
\newcommand{\Bb}{{\mathcal B}}

\newcommand{\Dcal}{{\mathcal D}}
\newcommand{\Dd}{{\mathcal D}}

\newcommand{\Ee}{{\mathcal E}}
\newcommand{\Fcal}{{\mathcal F}}
\newcommand{\Ff}{{\mathcal F}}
\newcommand{\Gcal}{{\mathcal G}}   
\newcommand{\Gg}{{\mathcal G}} 
\newcommand{\Hh}{{\mathcal H}}

\newcommand{\Jj}{{\mathcal J}}

\newcommand{\Ll}{{\mathcal L}}   

\newcommand{\Mcal}{{\mathcal M}}   
\newcommand{\Mm}{{\mathcal M}}
\newcommand{\Nn}{{\mathcal N}}

\newcommand{\Pp}{{\mathcal P}}

\newcommand{\Ss}{{\mathcal S}}
\newcommand{\Tcal}{{\mathcal T}}
\newcommand{\Tt}{{\mathcal T}}

\newcommand{\Uu}{{\mathcal U}}
\newcommand{\Vv}{{\mathcal V}}

\newcommand{\Xx}{{\mathcal X}}

\newcommand{\eps}{{\varepsilon}}
\renewcommand{\i}{{\iota}}

\renewcommand{\phi}{{\varphi}}

\newcommand{\Om}{{\Omega}}
\newcommand{\om}{{\omega}}
\newcommand{\im}{{\rm im }}        
\newcommand{\id}{{\rm id}}         

\newcommand{\supp}{{\rm supp}}     
 
\newcommand{\Lie}{{\rm Lie}}          
\newcommand{\Vect}{{\rm Vect}}        
\newcommand{\Vol}{{\rm Vol}}          
\newcommand{\End}{{\rm End}}          
\renewcommand{\d}{{\rm d}}

\newcommand{\w}{{\rm w}}

\newcommand{\ad}{{\rm ad}}
\newcommand{\Ad}{{\rm Ad}}

\newcommand{\ev}{{\rm ev}}
\newcommand{\bev}{{\overline{\rm ev}}}

\newcommand{\dvol}{{\rm dvol}}
\newcommand{\loc}{{\rm loc}}

\newcommand{\G}{{\rm G}}

\newcommand{\U}{{\rm U}}

\newcommand{\g}{{\mathfrak g}}         
\newcommand{\Cinf}{C^{\infty}}
\newcommand{\BG}{{\mathrm{BG}}}
\newcommand{\EG}{{\mathrm{EG}}}

\newcommand{\QH}{{\mathrm{QH}}}

\newcommand{\GW}{{\mathrm GW}}
\newcommand{\maslov}{{\mathrm m}}

\newcommand{\inner}[2]{\langle #1, #2\rangle}   
\newcommand{\INNER}[2]{\left\langle #1, #2\right\rangle}

\def\abracket#1{\langle#1\rangle}
\def\NABLA#1{{\mathop{\nabla\kern-.5ex\lower1ex\hbox{$#1$}}}}
\def\Nabla#1{\nabla\kern-.5ex{}_{#1}}
\def\Tabla#1{\tilde\nabla\kern-.5ex{}_{#1}}
\renewcommand{\Tilde}{\widetilde}
\renewcommand{\p}{{\partial}}

\newcommand{\IMP}{\Longrightarrow}
\newcommand{\IFF}{\Longleftrightarrow}
\newcommand{\INTO}{\hookrightarrow}
\newcommand{\TO}{\longrightarrow}
\newtheorem{theorem}{Theorem}[section]
\newtheorem{corollary}[theorem]{Corollary}

\newtheorem{lemma}[theorem]{Lemma}
\newtheorem{proposition}[theorem]{Proposition}

\newtheorem{remark}[theorem]{Remark}
\newtheorem{example}[theorem]{Example}

\begin{document}
 
\maketitle
 
 
\begin{abstract}
We study pseudoholomorphic curves in symplectic quotients
as adia\-batic limits of solutions of a system of nonlinear
first order elliptic partial differential equations 
in the ambient symplectic manifold.  
The symplectic manifold carries a Hamiltonian 
group action. The equations involve the
Cauchy-Riemann operator over a Riemann surface,
twisted by a connection,
and couple the curvature of the connection 
with the moment map.  Our main theorem asserts that the
genus zero invariants of Hamiltonian group actions defined 
by these equations are related to the genus zero Gromov--Witten
invariants of the symplectic quotient (in the monotone case)
via a natural ring homomorphism from the equivariant cohomology
of the ambient space to the quantum cohomology of the quotient. 
\end{abstract}
 

 
\section{Introduction}\label{sec:intro}

The main theorem of this paper asserts that 
under certain hypotheses there is a 
ring homomorphism from the equivariant cohomology of 
a symplectic manifold $M$ with a Hamiltonian
$\G$-action to the quantum cohomology 
of the symplectic quotient $\bar M$ such that the following 
diagram commutes
$$
\xymatrix
{
{\rm H}^*(M_\G)\ar[rr]^\phi\ar[dr]_{\Phi_B} & 
& \QH^*(\bar M)\ar[dl]^{\GW_{\bar B}} \\
& \Z &
}.
$$
Here $\GW_{\bar B}$ denotes the genus zero 
Gromov--Witten invariants of $\bar M$ 
with fixed marked points associated to 
a homology class $\bar B\in H_2(\bar M;\Z)$,
and $\Phi_B$ denotes the genus zero invariants of 
Hamiltonian group actions associated to 
the equivariant homology class $B=\kappa(\bar B)\in H_2(M_\G;\Z)$.
The latter invariant was introduced in~\cite{CGS,CGMS,Mu}.
The homomorphism $\phi$ is defined indirectly as a 
consequence of a comparison theorem for the two 
invariants.  A more direct definition in terms of 
vortices over the complex plane with values in $M$
will be given elsewhere. The proof of the comparison 
theorem is based on an adiabatic limit analysis
which relates the solutions of the equations used in 
the definition of the invariants $\Phi$ to pseudoholomorphic
curves in the symplectic quotient.  
Our hypotheses are that the moment map is proper,
that $M$ is convex at infinity, and that the quotient $\bar M$
is smooth. These hypotheses are needed to even state the 
result. In addition we assume that there are no 
holomorphic spheres in the ambient manifold (and hence $M$ 
is necessarily noncompact) and that $\bar M$ is monotone.  
These hypotheses are of technical
nature and it might be possible to remove them.
But this would require more analysis than is carried
out in the present paper. Before stating the main results
more precisely (Theorem~A and Corollary~A') we begin with
a brief discussion of the invariants introduced 
in~\cite{CGS,CGMS,Mu}.

\subsection*{Invariants of Hamiltonian group actions}

Let $(M,\om)$ be a symplectic manifold 
(not necessarily compact) and $\G$ 
be a compact connected Lie group with 
Lie algebra $\g$. We fix an invariant 
inner product $\inner{\cdot}{\cdot}$ 
on $\g$ and identify $\g$ with its dual $\g^*$. 
We assume that $\G$ acts on $M$ by Hamiltonian 
symplectomorphisms and that the action is generated 
by an equivariant moment map $\mu:M\to\g$. 
This means that, for every $\eta\in\g$,
the vector field $X_\eta\in\Vect(M)$ that generates
the action is determined by 
$
     \i(X_\eta)\om = d\inner{\mu}{\eta}.
$
Let $\pi:P\to\Sigma$ be a principal $\G$-bundle
over a compact oriented Riemann surface 
$(\Sigma,j_\Sigma,\dvol_\Sigma)$.
We fix a smooth family $\Sigma\to\Jj_\G(M,\om):z\mapsto J_z$
of $\G$-invariant and $\om$-compatible almost complex 
structures on $M$. This determines a family of metrics 
$
     \inner{\cdot}{\cdot}_z := \om(\cdot,J_z\cdot).
$
The invariants are derived from the equations 
\begin{equation}\label{eq:1}
     \bar\partial_{J,A}(u)=0,\qquad
     *F_A+\mu(u)=0,
\end{equation}
for a pair $(u,A)$, where $u:P\to M$ is a $\G$-equivariant map
and $A$ is a connection on~$P$.
Here $\bar\partial_{J,A}$ is the 
nonlinear Cauchy-Riemann operator twisted by $A$
and $F_A$ is the curvature of $A$. 
Both terms in the second identity in~(\ref{eq:1}) 
are sections of the Lie algebra bundle $\g_P:=P\times_\Ad\g$. 
In contrast to the theory
of pseudoholomorphic curves, equations~(\ref{eq:1}) 
involve the volume form $\dvol_\Sigma$
(via the Hodge $*$-operator in the second equation)
and not just the complex structure $j_\Sigma$.
Equations~(\ref{eq:1}) are invariant under 
the action of the gauge group~$\Gg=\Gg(P)$ 
(of equivariant maps from $P$ to $\G$) by
$$
     g^*(u,A) = (g^{-1}u,g^{-1}dg+g^{-1}Ag).
$$
>From a geometric point of view, the solutions of~(\ref{eq:1})
correspond to the absolute minima of the energy functional
$$
     E(u,A)
     =
     \frac{1}{2}\int_{\Sigma}
     \left(|d_Au|^2+|F_A|^2+|\mu(u)|^2\right)
     \,\dvol_{\Sigma}
$$
in a given homotopy class.
If the pair $(u,A)$ is a solution of~(\ref{eq:1}) 
then it has energy 
$$
     E(u,A) = \int_\Sigma\bigl(
     u^*\om - d\inner{\mu(u)}{A}
     \bigr) =: \inner{[\om-\mu]}{[u]}
$$
and this number is an invariant of the equivariant 
homology class represented by the map $u$. 
We impose the following hypothesis throughout this paper.
\begin{description}
{\it
\item[(H1)]
The moment map $\mu$ is proper,
zero is a regular value of $\mu$,
and $\G$ acts freely on $\mu^{-1}(0)$.
}
\end{description}
Under this hypothesis the quotient
$$
     \bar M:=M\dslash\G:=\mu^{-1}(0)/\G
$$
is a compact symplectic manifold. The induced symplectic
form will be denoted by $\bar\om$.  The equivariant
homology class $[u]\in H_2(M_\G;\Z)$ is defined by 
the following diagram, which also shows how it is related
to the characteristic class $[P]\in H_2(\BG;\Z)$
and to the class $[\bar u]\in H_2(\bar M;\Z)$
in the case $\mu\circ u\equiv 0$. Note that, since $\G$
is connected, the equivariant homology class $[u]$
determines the isomorphism class of the bundle $P$. 
We denote $M_\G:=M\times_\G\EG$.
$$
\xymatrix
{
[\Sigma] \;\; \in \hspace{-20pt} & 
{\rm H}_2(\Sigma;\Z) \ar[r]^{\bar{u}_*}  &
{\rm H}_2(\bar M;\Z)\ar[d]^\kappa &
\hspace{-20pt} \ni \;\;[\bar u]   \\
& {\rm H}_2(P_\G;\Z) \ar[u]^\simeq
\ar[r]^{u_*} \ar[d] &
{\rm H}_2(M_\G;\Z) \ar[dl] &
\hspace{-20pt} \ni\;\; [u]  \\
 [P] \;\; \in \hspace{-20pt} &
{\rm H}_2(\BG;\Z)  & & 
}
$$

Fix a homology class $\bar B\in H_2(\bar M;\Z)$, 
let $B:=\kappa(\bar B)\in H_2(M_\G;\Z)$, and denote
the space of solutions of~(\ref{eq:1}) that represent
this homology class by 
$$
     \Tilde{\Mm}_{B,\Sigma}
     := \left\{(u,A)\in\Cinf_\G(P,M)\times\Aa(P)\,|\,
        [u] = B,\,
        u\mbox{ and }A\mbox{ satisfy }(\ref{eq:1})
        \right\}.
$$
Here $P\to\Sigma$ denotes a principal $\G$-bundle 
whose characteristic class $[P]\in H_2(\BG;\Z)$ 
is determined by $B$ as above.
The quotient by the action of the gauge group
will be denoted by 
$$
     \Mm_{B,\Sigma} := \Tilde{\Mm}_{B,\Sigma}/\Gg.
$$
We impose another hypothesis which guarantees
compactness~\cite{CGMS}.  
\begin{description}
\item[(H2)]
{\it There exists a $\G$-invariant almost complex structure
$J\in\Jj_\G(M,\om)$, a proper $\G$-invariant function
$f:M\to[0,\infty)$, and a constant $c>0$ such that
$$
     f(x)\ge c\qquad\IMP\qquad
     \inner{\Nabla{\xi}\nabla f(x)}{\xi}>0
$$
for every nonzero vector $\xi\in T_xM$ and
$$
     f(x)\ge c\qquad\IMP\qquad
     df(x)JX_{\mu(x)}(x)\ge 0.
$$
Moreover, 
$
     \int_{S^2} v^*\om = 0
$ 
for every smooth map $v:S^2\to M$.}
\end{description}
This hypothesis implies that $\sup_P(f\circ u)\le c$
for every solution $(u,A)$ of~(\ref{eq:1}) 
over any Riemann surface and in any homology class
(see~\cite{CGMS}).  In~\cite{CGS} it is shown that~(H2) follows
from~(H1) in the case of linear actions on $\C^n$. 
In~\cite{CGMS} it is shown that the moduli space $\Mm_{B,\Sigma}$ 
is a smooth compact manifold of dimension
$$
     \dim\,\Mm_{B,\Sigma}
     = \left(\frac12\dim\,M-\dim\,\G\right)\chi(\Sigma)
       + 2\inner{c_1^\G(TM)}{B}
$$
for a generic $J$, provided that~$(H1-2)$ are satisfied, 
$B$ is a nontorsion homology class, and the area of $\Sigma$
is sufficiently large.  The latter condition, 
together with the energy identity, guarantees that every 
solution of~(\ref{eq:1}) is somewhere close to the zero set
of the moment map.  The class $c_1^\G(TM)\in H^2(M_\G;\Z)$
in the dimension formula denotes the equivariant first Chern class 
of the complex vector bundle $(TM,J)$.

Consider the evaluation map 
$
     \ev_\G:\Mm_{B,\Sigma}\to M_\G,
$
defined by 
$$
     \ev_\G([u,A]):=[u(p_0),\Theta_0(u,A)],
$$ 
where $p_0\in P$ is fixed and $\Theta_0:\Tilde{\Mm}_{B,\Sigma}\to\EG$
is a smooth map such that
$$
     \Theta_0(g^{-1}u,g^*A) = g(p_0)^{-1}\Theta_0(u,A).
$$
This means that $\Theta_0$ is a classifying map for the 
principal $\G$-bundle $\Pp_{B,\Sigma}\to\Mm_{B,\Sigma}$
obtained as the quotient of $\Tilde{\Mm}_{B,\Sigma}$
by the based gauge group
$
     \Gg_0 := \left\{g\in\Gg\,|\,g(p_0)=\one\right\}.
$
Let $\alpha\in H^*(M_\G;\Z)$ be a class of degree
$
     \deg(\alpha) = \dim\,\Mm_{B,\Sigma}
$
and define 
$$
     \Phi_{B,\Sigma}(\alpha)
     := \int_{\Mm_{B,\Sigma}}\ev_\G^*\alpha.       
$$
In~\cite{CGMS} it is shown that this integer is independent 
of the almost complex structure $J$, the metric on $\Sigma$, 
and the point $p_0$ used to define it. 

Now let 
$
     D:=\{z\in\C\,|\,|z|\le1\}
$
and consider the space of maps 
$v:D\to M$ that map the boundary $\p D$ 
to a $\G$-orbit in $\mu^{-1}(0)$:
$$
     \Vv
     := \left\{v:D\to M\,|\,\exists g:\R/2\pi\Z\to\G\;
        \exists x\in\mu^{-1}(0)\;\forall\theta\in\R\;\;
        v(e^{i\theta})=g(\theta)x
        \right\}.
$$
Let $\maslov:\Vv\to\Z$ denote the function which assigns to each 
element $v\in\Vv$ the Maslov index 
of the loop of symplectic matrices obtained 
from the linear maps $g(\theta):T_xM\to T_{g(\theta)x}M$ 
in a trivialization along $v$.
Every smooth map $\bar v:S^2\to\bar M$ lifts to 
a map $v\in\Vv$ and in this case the Maslov index
$\maslov(v)$ is equal to the first Chern number
$\inner{c_1(T\bar M)}{\bar v_*[S^2]}$. 
The minimal Maslov number will be denoted by
$$
     N := \inf_{v\in\Vv,\,\,\maslov(v)>0}\maslov(v).
$$
This is a lower bound for the minimal Chern number of $\bar M$.
We impose a third hypothesis.
\begin{description}
\item[(H3)]
{\it There exists a constant $\tau>0$ such that
$$
     \int_Dv^*\om = \tau\maslov(v)
$$
for every $v\in\Vv$.}
\end{description}
This hypothesis implies that the quotient 
$\bar M$ is a monotone symplectic manifold
and that the energy of every holomorphic sphere in $\bar M$
is an integer multiple of $\hbar:=\tau N$. 
The main result of this paper asserts that
under hypotheses~(H1-3) the invariant $\Phi_{B,S^2}$ 
agrees with the corresponding genus zero Gromov--Witten invariant 
of $\bar M$, provided that the cohomology classes 
$\alpha_i$ have degrees less than $2N$. 

\subsection*{The main theorem}

\noindent{\bf Theorem~A}.
{\it
Assume~{\rm (H1-3)} and let $\bar B\in H_2(\bar M;\Z)$
and $\alpha_1,\dots,\alpha_k\in H^*_\G(M;\Z)$ be given 
such that
$$
     \deg(\alpha_i)<2N
$$
for $i=1,\dots,k$ and
$$
     \sum_{i=1}^k\deg(\alpha_i)
     := \left(\frac12\dim\,M-\dim\,\G\right)\chi(\Sigma)
        + 2\inner{c_1^\G(TM)}{B},
$$
where $B:=\kappa(\bar B)\in H^\G_2(M;\Z)$. Then
$$
     \Phi_{B,S^2}(\alpha_1\smile\cdots\smile\alpha_k)
     = \GW_{\bar B,S^2}(\bar\alpha_1,\dots,\bar\alpha_k),
$$
where $\bar\alpha_i:=\kappa(\alpha_i)\in H^*(\bar M;\Z)$.
}

\medskip
\noindent
{\bf Remarks.}
{\bf (i)}
In the definition of $\Phi_{B,\Sigma}$ the 
point $p_0\in P$ at which the map $u$ 
is evaluated is fixed and the cohomology class
$\ev_\G^*\alpha\in H^*(\Mm_{B,\Sigma};\Z)$
is independent of the choice of the point $p_0$
used in the definition of $\ev_\G$.
The Gromov--Witten invariants in Theorem~A are also to be 
understood with fixed marked points on $S^2$ in the definitions 
of the evaluation maps, and with almost complex structures that
are allowed to depend on the base point $z\in S^2$. 

\smallskip
\noindent{\bf (ii)}
If, in addition to~$(H3)$, we assume $[\bar\om]=\tau c_1(T\bar M)$
then the proof of Theorem~A goes through word by word for the higher
genus case (with fixed marked points). 
In general, the extension to general Riemann 
surfaces requires a refined version of the 
compactness theorem in Section~\ref{sec:bubbling}
which takes account of the preservation of the homotopy
class in the limit, as in Gromov compactness. 
With similar refined arguments one should be able 
to deal with the case of varying marked points 
or of varying complex structures on $\Sigma$. 

\smallskip
\noindent{\bf (iii)}
The assertion of Theorem~A does not continue
to hold in the case 
$
     \deg(\alpha_i)\ge 2N.
$
For example, consider the standard action of $S^1$
on $\C^n$, let $P\to S^2$ be an $S^1$-bundle
of degree $d\ge 0$, and denote by 
$c\in H^2({\mathrm B}S^1;\Z)=H^2_{S^1}(\C^n;\Z)$
the positive generator. Then the minimal Chern number 
is $N=n$, the dimension of $\Mm_{d,S^2}$ is $2nd+2n-2$,
and we have
$
     \Phi_{d,S^2}(c^m) = 1
$
whenever 
$
     m=nd+n-1.
$ 
The corresponding Gromov--Witten invariant 
(for a $k$-tuple of classes $c^{m_1},\dots, c^{m_k}$
with $m_1+\cdots+m_k=m$) counts holomorphic
spheres of degree~$d$ in $\C P^{n-1}$ passing at $k$
given points $z_1,\dots,z_k\in S^2$
through generic copies of $\C P^{n-1-m_i}$ for 
$i=1,\dots,k$. Thus the Gromov--Witten invariant is zero 
whenever $\deg(c^{m_i})=2m_i\ge 2n$ for some~$i$.

\subsection*{Equivariant and quantum cohomology}

Consider the monotone case. 
The quantum cohomology $\QH^*(\bar M)$ of $\bar M$ is the ring 
of all formal sums of the form
$$
     \bar\alpha 
     = \sum_{\bar B\in H_2(\bar M;\Z)}
       \bar\alpha_{\bar B}e^{\bar B},
$$
where $\bar\alpha_{\bar B}\in H^*(\bar M;\Z)/{\rm torsion}$, 
such that
$$
     \#\left\{\bar B\in H_2(\bar M;\Z)\,|\,\bar\alpha_{\bar B}\ne 0,\,
     \inner{[\bar\om]}{\bar B}\le c\right\} <  \infty 
$$
for all $c>0$. The degree convention is 
$
     \deg(e^{\bar B}) := 2\inner{c_1(T\bar M)}{\bar B}.
$
Choose an integral basis $\bar e_0,\dots,\bar e_n$ of
$H^*(\bar M;\Z)/{\rm torsion}$ and let $\bar e_i^*$ denote the dual
basis in the sense that 
$$
     \int_{\bar M}\bar e_i\smile\bar e_j^* = \delta_{ij}.
$$
Then the product structure on $\QH^*(\bar M)$ is defined by
$$
     \bar\alpha_1 * \bar\alpha_2
     := \sum_{\bar B_1,\bar B_2,\bar B}\sum_{i=0}^n
        \GW_{\bar B-\bar B_1-\bar B_2,S^2}
        (\bar\alpha_{1\bar B_1},\bar\alpha_{2\bar B_2},\bar e_i^*)
        \bar e_ie^{\bar B}.
$$
The sum is over all quadruples $i,\bar B_1,\bar B_2,\bar B$ 
such that 
$$
     \deg(\bar\alpha_1)+\deg(\bar\alpha_2)
     =\deg(\bar e_i)+2c_1(\bar B), 
$$
where we abbreviate $c_1(\bar B):=\inner{c_1(T\bar M)}{\bar B}$. 

The Gromov--Witten invariant associated to a
Riemann surface $\Sigma$, with a fixed complex structure $j_\Sigma$ 
and fixed marked points $z_1,\dots,z_k$, can be extended to a map
$
    \GW_{\bar B,\Sigma}:\QH^*(\bar M)\otimes\cdots\otimes\QH^*(\bar M)\to\Z
$
by the formula
$$
    \GW_{\bar B,\Sigma}(\bar\alpha_1,\dots,\bar\alpha_k)
    := \sum_{\bar B_i}\GW_{\bar B-\bar B_1-\cdots-\bar B_k,\Sigma}
       (\bar\alpha_{1\bar B_1},\dots,\bar\alpha_{k\bar B_k}).
$$
With this convention the gluing formula for the Gromov--Witten
invariants~\cite{MS1,RT} can be expressed in the form
\begin{equation}\label{eq:gluing}
    \GW_{\bar B,\Sigma}(\bar\alpha_1,\dots,\bar\alpha_k)
    = \GW_{\bar B,\Sigma}(\bar\alpha_1*\cdots*\bar\alpha_k).
\end{equation}
We abbreviate $H^*(X):=H^*(X;\Z)/{\rm torsion}$.

\medskip
\noindent{\bf Corollary~A'}.
{\it Assume~$(H1-3)$ and suppose that 
$H^*(M_\G)$ is generated by classes 
of degree less than $2N$.  Then there exists a 
unique (surjective) ring homomorphism 
$
     \phi:H^*(M_\G)\TO\QH^*(\bar M)
$
such that, for every $\alpha\in H^*(M_\G)$ 
and every $\bar B\in H_2(\bar M;\Z)$,}
$$
     \deg(\alpha)<2N\qquad\IMP\qquad
     \phi(\alpha)=\kappa(\alpha),
$$
$$
    \Phi_{\kappa(\bar B),S^2}(\alpha) 
    = \GW_{\bar B,S^2}(\phi(\alpha)).
$$

\begin{proof}
Let $\alpha\in H^*(M_\G)$ and choose $\alpha_{ij}\in H^*(M_\G)$ 
such that $\deg(\alpha_{ij})<2N$ and
\begin{equation}\label{eq:alpha}
     \alpha = \sum_{i=1}^k\alpha_{i1}\smile\cdots\smile\alpha_{i\ell}.
\end{equation}
Define 
\begin{equation}\label{eq:phialpha}
     \phi(\alpha) := \sum_{i=1}^k
     \kappa(\alpha_{i1}) * \cdots * \kappa(\alpha_{i\ell}).
\end{equation}
We prove that $\phi(\alpha)$ is independent of the choice
of $\alpha_{ij}$.  To see this, note that,
since the cohomology of $\bar M$ is generated by classes 
of degree less than $2N$, so is the quantum cohomology. 
This means that a quantum cohomology class 
$\bar\alpha\in\QH^*(\bar M)$ is zero if and only if 
$
     \GW_{\bar B,S^2}(\bar\alpha,\bar\beta_1,\dots,\bar\beta_m)=0
$
for every $\bar B\in H_2(\bar M;\Z)$ and all
$\bar\beta_1,\dots\bar\beta_m\in H^*(\bar M)$
such that $\deg(\bar\beta_j)<2N$ for all~$j$.
Now suppose that the expression
on the right of~(\ref{eq:phialpha}) is nonzero.
Then, by what we have just observed, there exist 
cohomology classes $\bar\beta_1,\dots,\bar\beta_m$ of degrees
less than $2N$ and a homology class $\bar B\in H_2(\bar M;\Z)$ 
such that 
$$
     \sum_{i=1}^k\GW_{\bar B,S^2}
     (\kappa(\alpha_{i1}),\dots,\kappa(\alpha_{i\ell}),
     \bar\beta_1,\dots,\bar\beta_m) \ne 0.
$$
Since the homomorphism 
$
     \kappa:H^*(M_\G)\to H^*(\bar M)
$
is surjective (cf.~\cite{Kw}), there exist classes 
$
     \beta_j\in H^*(M_\G)
$
(of degrees less than $2N$) such that $\kappa(\beta_j)=\bar\beta_j$
for every $j$. Hence, by Theorem~A,
$$
     \sum_{i=1}^k\Phi_{\kappa(\bar B),S^2}
     (\alpha_{i1}\smile\cdots\smile\alpha_{i\ell}\smile
     \beta_1\smile\dots\smile\beta_m) \ne 0,
$$
and hence $\alpha\ne 0$.  
This shows that $\phi$ is well defined. 
The map $\phi$ is obviously a ring homomorphism.
The formula $\Phi_{B,S^2}(\alpha)=\GW_{\bar B,S^2}(\phi(\alpha))$
follows immediately from Theorem~A and the gluing 
formula~(\ref{eq:gluing}) for the Gromov--Witten invariants.    
\end{proof}

The homomorphism
$
     \phi:H^*(M_\G)\to\QH^*(\bar M)
$
can be defined geometrically in terms of the
vortex equations over $\C$:
\begin{equation}\label{eq:vortex-C}
     \p_su+L_u\Phi+J(\p_tu+L_u\Psi)=0,\qquad
     \p_s\Psi-\p_t\Phi+[\Phi,\Psi]+\mu(u)=0.
\end{equation}
For every finite energy solutions of~(\ref{eq:vortex-C})
in radial gauge there exist a loop 
$g:S^1\to\G$ and a point $x_0\in\mu^{-1}(0)$
such that 
\begin{equation}\label{eq:limit-C}
    \lim_{r\to\infty}u(re^{i\theta}) = g(e^{i\theta})x_0
\end{equation}
(see Section~\ref{sec:vortex}). 
Every map $u:\C\to M$ that satisfies~(\ref{eq:limit-C})
determines an equivariant homology class $B=[u]\in H_2(M_\G;\Z)$.
Now the moduli space 
$
    \Mm_B(J) 
$
of gauge equivalence classes of solutions 
of~(\ref{eq:vortex-C}) and~(\ref{eq:limit-C})
that represent the class $B$ has two evaluation maps 
$
    \ev_0:\Mm_B\to M_\G
$
and
$
    \ev_\infty:\Mm_B\to\bar M.
$
The map $\phi$ can be defined by 
$$
    \phi(\alpha) = \sum_{i=0}^n\sum_{\bar B}
    \left(\int_{\Mm_{\kappa(\bar B)}}
    \ev_0^*\alpha\smile\ev_\infty^*\bar e_i^*\right)
    \bar e_ie^{\bar B}.
$$
The details of this construction will be carried out elsewhere.

\subsection*{Outline of the proof of Theorem~A}

The proof of Theorem~A is based on an adiabatic limit
argument in which the metric on the Riemann surface 
is scaled by a large factor $\eps^{-2}$.
Then equations~(\ref{eq:1}) have the form
\begin{equation}\label{eq:eps}
     \bar\partial_{J,A}(u)=0,\qquad
      *F_A+\eps^{-2}\mu(u)=0.
\end{equation}
The solutions of~(\ref{eq:eps}) minimize the $\eps$-dependent 
energy
$$
    E^\eps(u,A)
    =
    \frac{1}{2}\int_{\Sigma}
    \left(|d_Au|^2+\eps^2|F_A|^2+\eps^{-2}|\mu(u)|^2\right)
    \,\dvol_{\Sigma},
$$
and the value of this functional at a solution of~(\ref{eq:eps})
is independent of $\eps$ in a given equivariant homology class. 
In this paper we examine the limit behaviour of the 
solutions of~(\ref{eq:eps}) as $\eps$ tends to zero
for Riemann surfaces of any genus. 
The limit equations have the form
\begin{equation}\label{eq:jhol}
     \bar\p_{J,A}(u)=0,\qquad \mu(u)=0.
\end{equation}
The solutions of~(\ref{eq:jhol}) can be 
interpreted as pseudoholomorphic curves in the
symplectic quotient $\bar M=\mu^{-1}(0)/\G$
with respect to the induced family of almost complex
structures $\bar J_z$ (see Section~\ref{sec:hol}). 
We impose a further hypothesis that is satisfied
for a generic family of $\G$-invariant
almost complex structures on $M$:
\begin{description}
{\it
\item[(H4)] 
Every nonconstant $\bar J$-holomorphic curve
$\bar u:\Sigma\to\bar M$ is regular in the sense that 
the linearized Cauchy-Riemann operator along $\bar u$
is surjective.
}
\end{description}
This hypothesis guarantees that the moduli space 
of holomorphic curves in $\bar M$ is smooth.  

The proof of Theorem~A requires three preliminary
theorems which are of interest in their own rights.
Theorem~B constructs a $\Gcal(P)$-equi\-va\-ri\-ant map
$$
     (u_0,A_0)\mapsto (u_\eps,A_\eps) =:\Tilde{\Tcal}^{\eps}(u_0,A_0)
$$
which assigns to every regular solution of~(\ref{eq:jhol})
a nearby solution of~(\ref{eq:eps}) 
for $\eps>0$ sufficiently small.  
How small $\eps$ must be chosen depends
(continuously) on the given pair $(u_0,A_0)$. 
Theorem~C shows that the map $\Tilde{\Tt}^\eps$
constructed in Theorem~B is {\it locally surjective}
in the sense that every solution of~(\ref{eq:eps}) 
that is sufficiently close to a solution $(u_0,A_0)$ 
of~(\ref{eq:jhol}) must be in the image of~$\Tilde{\Tcal}^{\eps}$. 
The neighbourhood in which surjectivity holds depends
on $\eps$: it becomes smaller as $\eps$ tends to zero.
Theorem~D strengthens the local surjectivity 
result of Theorem~C.  We remove the 
assumption that the solution of~(\ref{eq:eps}) 
is close to some given solution of~(\ref{eq:jhol}). 
However, we consider only solutions of~(\ref{eq:eps})
that satisfy a suitable $L^\infty$-bound
on the first derivatives and prove that every 
solution of~(\ref{eq:eps}) that satisfies this bound
lies in the image of $\Tilde{\Tcal}^{\eps}$ 
for $\eps$ small.  The proof of Theorem~A is then
based on a bubbling argument in the small $\eps$
limit which establishes a one-to-one correspondence 
between the solutions of~(\ref{eq:jhol}) and those
of~(\ref{eq:eps}) in a zero dimensional setting, where
additional conditions have been imposed.

In Section~\ref{sec:hol} we review standard results
about the moduli space of pseudoholomorphic curves 
in the symplectic quotient $M\dslash\G$ 
and rephrase them in terms of solutions 
of~(\ref{eq:jhol}). Theorems~B, C, and~D
will be stated in Section~\ref{sec:adiabatic}. 
The remaining sections are devoted to the proofs
of the four main theorems.  

While the general outline of the proof of Theorem~A 
is analogous to the proof of the Atiyah--Floer conjecture 
in~\cite{DS1,DS2} there are several new ingredients in the 
present paper.  

In~\cite{DS2} the moduli space 
of Floer connecting orbits is a finite set, while the 
moduli space $\Mm^0$ of pseudoholomorphic curves is, 
in general, a (noncompact) manifold of positive dimension.
Hence, in constructing the map $\Tt^\eps$ 
from (a compact subset of) $\Mm^0$ to $\Mm^\eps$, 
care must be taken to establish that the constants 
depend continuously on the point in $\Mm^0$. This refers
to the linear and quadratic estimates needed in the proof 
of Theorem~B (Sections~\ref{sec:linear} and~\ref{sec:quad}).
Secondly, we extend the estimates of~\cite{DS1}
for the $(1,p,\eps)$-norms to the $(2,p,\eps)$-norms
(see Lemmata~\ref{le:difference}-\ref{le:onto1} 
and Proposition~\ref{prop:quadra}).  
For the linear estimates the extension to the higher 
derivatives is quite subtle because of the 
$\eps$-dependent norms.  Another new ingredient 
arises from the presence of nonlinearities in the 
highest order terms of the Cauchy--Riemann equations.
This requires more delicate quadratic estimates 
(Proposition~\ref{prop:quadra}) for the proof of 
local uniqueness (Theorem~\ref{thm:uniqbaby}). 

The proof of the $\eps$-local slice theorem in
Section~\ref{sec:gauge} (which can be viewed as a 
simpler analogue of Theorem~B) is considerably harder
than the analogous result in~\cite{DS2},
since a) we must establish estimates for the $(2,p,\eps)$-norms,
b) we must prove that the constants depend continuously 
on the point in $\Mm^0$, c) the manifold $M$ is not an 
affine space so we must deal with additional lower order terms,
and d) we give a proof of the linear estimate in Lemma~\ref{le:gauge}
for $p>2$ (the analogue in~\cite{DS2} was only established for $p=2$).

The proof of local surjectivity (Theorem~C)
requires a subtle tubular neighbourhood theorem 
for the moduli spaces $\Mm^\eps$.
The result is quantitative with constants independent 
of $\eps$. In particular, the proof involves an estimate
for the derivative of the map
$(u_0,A_0)\to (\xi_0,\alpha_0,g)$
given by the $\eps$-local slice theorem.
This is where the estimates for the
$(2,p,\eps)$-norms are needed. As a result
the entire adiabatic limit argument has to be 
carried out for these higher norms.
In comparison, the analogous result in~\cite[Proposition 6.3]{DS2}
can be disposed of with a simple time shift argument 
and only requires estimates in the $(1,p,\eps)$-norm.

Another new ingredient in the present paper is the 
apriori estimate in Lemma~\ref{le:mu}. 
It asserts that every solution of~(\ref{eq:eps}) 
which satisfies a certain $L^\infty$ bound must be 
$\eps^{3/2}$-close to the zero set of the moment map.
As a consequence we obtain in Theorem~\ref{thm:D}
a much stronger surjectivity result for the map~$\Tcal^\eps$;
compare with~\cite[Theorem 8.1]{DS2}.
This strenghtened form of~\cite[Theorem 8.1]{DS2} is needed
to close a gap in the proof of~\cite[Theorem 9.1]{DS2}, 
namely to prove that the holomorphic sphere 
in the symplectic quotient appearing in the bubbling 
argument on page~634 in~\cite{DS2} is nonconstant.
At the same time the bubbling argument in the 
proof of~\cite[Theorem 9.1]{DS2} can be simplified:
it suffices to work with the sequence 
$c_\nu:=\sup(|d_{A\nu}u_\nu|+\eps_\nu^{-1}|F_{A_\nu}|)$ 
instead of
$c_\nu:=\sup(|d_{A\nu}u_\nu|+\eps_\nu^{-1}|F_{A_\nu}|^{1/2})$.
The modified bubbling argument is carried out 
in the present context in Section~\ref{sec:bubbling}.

In Section~\ref{sec:vortex} we establish the asymptotic 
behaviour and the quantization of the energy for
solutions of the nonlinear vortex equations on
the complex plane. (An analogous result for 
anti-self-dual instantons is used without proof 
in~\cite{DS2}.)
In Section~\ref{sec:bubbling} we construct a classifying map 
on an open set in $\Cinf_\G(P,M)\times\Aa(P)$,
which contains the moduli spaces $\Mm^\eps$ for 
all $\eps\in[0,\eps_0]$, with values in a finite
dimensional approximation of $\EG$; and we prove 
$C^1$-convergence for the composition of the 
resulting evaluation map with $\Tcal^\eps$. 
All these results are needed in the proof of our main theorem.

 
\section{Pseudoholomorphic curves}\label{sec:hol}

For $z\in\Sigma$ let $\bar J_z$ denote the 
almost complex structure on $\bar M$ induced 
by $J_z$, let $\bar P\TO\bar M$ 
denote the principal $\G$-bundle 
$
     \bar P:=\mu^{-1}(0)\subset M,
$
and let $\bar A$ denote the connection on $\bar P$
determined by $\om$ and~$J$. 
If $(u,A)$ is a solution of~(\ref{eq:jhol}) 
then $u$ descends to a $\bar J$-holomorphic curve 
$\bar u:\Sigma\to\bar M$ and 
$
     A
$
is the pullback of $\bar A$ under $\bar u$. 
Two gauge equivalent solutions descend to the same map $\bar u$ 
and every $\bar J$-holomorphic curve $\bar u:\Sigma\to\bar M$
lifts to a solution of~(\ref{eq:jhol}) 
for some principal $\G$-bundle $P$
(isomorphic to the pullback of $\bar P$ 
under~$\bar u$).  
 
Fix a homology class $\bar B\in H_2(\bar M;\Z)$, let 
$B:=\kappa(\bar B)\in H_2^\G(M;\Z)$, and consider the 
space 
$$
    \Tilde{\Mcal}_{B,\Sigma}^0 
    := \left\{(u,A)\in\Cinf_\G(P,M)\times\Acal(P)\,|\,
      [u]=B\mbox{ and }(\ref{eq:jhol})\mbox{ holds}\right\}.
$$
This space is invariant under the action of the gauge
group~$\Gg(P)$. Under our standing hypothesis~(H4) 
the quotient
$$
    \Mcal_{B,\Sigma}^0 
    := \Tilde{\Mcal}_{B,\Sigma}^0/\Gg(P)
$$
is a smooth manifold of dimension
$$
    \dim\,\Mcal_{B,\Sigma}^0
    = \left(\frac12\dim \, M-\dim\, \G\right)\chi(\Sigma)
    + 2\abracket{c_1^\G(TM),B}
$$
(see~\cite[Theorem~3.3.4]{MS1}).
Note that $\Mcal_{B,\Sigma}$ and $\Mcal_{B,\Sigma}^0$
have the same dimension.
 
For later reference we now introduce explicit 
notation for a local pa\-ra\-me\-tri\-zation of 
$\Mcal_{B,\Sigma}^0$ by the kernel of the 
linearized operator. Linearizing equations~(\ref{eq:jhol}) 
at a solution $(u,A)$ gives rise
to the Cauchy--Riemann operator
$$
     \Dd^0 := \Dd^0_{(u,A)}:\Om^0(\Sigma,H_u)
     \to \Om^{0,1}(\Sigma,H_u)
$$
given by 
$$
     \Dd^0\xi_0 := \pi_uD_{(u,A)}\xi_0,
$$
where $D_{(u,A)}:\Om^0(\Sigma,u^*TM/\G)\to\Om^{0,1}(u^*TM/\G)$
is the operator~(\ref{eq:Du}) 
in Appendix~\ref{app:CR}.  The bundle $H_u\to\Sigma$
and the projection $\pi_u:u^*TM\to H_u$ are defined
as follows. Consider the bundle $H\TO\Sigma\times\mu^{-1}(0)$
with fibres
$$
     H_{z,x} 
     := \ker\,d\mu(x)\cap 
        \ker\,d\mu(x)J_z.
$$
There is an orthogonal decomposition
$$
     T_xM
     = \im\, L_x\oplus H_{z,x}\oplus \im\,J_zL_x
$$
for every $(z,x)\in\Sigma\times\mu^{-1}(0)$,
where $L_x:\g\to T_xM$ the infinitesimal action, i.e.
$$
     L_x\eta:=X_\eta(x).
$$
Its dual operator with respect to the metric determined by $J_z$
is given by 
$$
     L_x^* = L_x^{*_z} = d\mu(x)J_z(x).
$$
Now let $u:P\to\mu^{-1}(0)$ be an equivariant map 
and consider the pullback of $H$ under 
the map $\tilde u:P\to\Sigma\times\mu^{-1}(0)$,
given by 
$
     \tilde u(p):=(\pi(p),u(p)).
$
This pullback is a $\G$-equivariant vector bundle over 
$P$ and its quotient 
$$
     H_u := {\tilde u}^*H/\G\TO\Sigma
$$
is naturally isomorphic to the pullback
of the tangent bundle $T\bar M$ under the 
induced map $\bar u:\Sigma\to\bar M$. 
Let $\pi_u:u^*TM/\G\to H_u$ denote the orthogonal
projection onto the harmonic part.
Thus $\pi_u[\xi]:=[\pi_u\xi]$ where the 
lifted projection $u^*TM\to u^*H$
(also denoted by $\pi_u$) is given by 
\begin{equation}\label{eq:pi-u}
     \pi_u\xi
     := \xi - L_u(L_u^*L_u)^{-1}L_u^*\xi
        + JL_u(L_u^*L_u)^{-1}L_u^*J\xi
\end{equation}
for a $\G$-equivariant section $\xi:P\to u^*TM$. 

\begin{theorem}\label{thm:M0}
Assume~{\rm (H1)} and~{\rm (H4)} 
and fix a constant $p>2$.  
For ev\-ery $(\bar u_0,\bar A_0)\in\Tilde{\Mm}_{B,\Sigma}^0$ 
there exist a sequence of positive constants 
$\delta,c,c_1,c_2,\dots$ and a map
$$
     \Ff^0 := \Ff^0_{(\bar u_0,\bar A_0)}:
     B_\delta^0\to\Tilde{\Mm}^0_{B,\Sigma},\qquad
     B_\delta^0
     := \{\xi_0\in\ker\,\Dd^0_{(\bar u_0,\bar A_0)}\,|\,
        \left\|\xi_0\right\|_{L^p}<\delta\},
$$
such that the following holds.

\smallskip
\noindent{\bf (i)}
If $\xi_0\in\Bb_\delta^0$ then
there exists a unique pair of sections 
$\xi_1\in\Om^0(\Sigma,H_{\bar u_0})$ 
and $\xi_2\in\Om^0(\Sigma,\im\,JL_{\bar u_0}/\G)$
such that 
$$
     \xi_1\in\im\,(\Dd^0_{(\bar u_0,\bar A_0)})^*,\qquad
     \left\|\xi_1\right\|_{W^{1,p}}
     + \left\|\xi_2\right\|_{W^{1,p}}
     \le c\left\|\xi_0\right\|_{W^{1,p}}.
$$
and the pair $(u_0,A_0)$, given by
$$
     u_0 := \exp_{\bar u_0}^{\mu^{-1}(0)}(\xi_0+\xi_1)
          = \exp_{\bar u_0}(\xi_0+\xi_1+\xi_2),\quad
     A_0 := -(L_{u_0}^*L_{u_0})^{-1}L_{u_0}^*du_0,
$$
satisfies~(\ref{eq:jhol}). The pair $(u_0,A_0)$
is the image of $\xi_0$ under $\Ff^0$. 

\smallskip
\noindent{\bf (ii)}
For every integer $k\ge 1$ and every $\xi_0\in B_\delta^0$
we have 
$$
     \left\|\xi_1\right\|_{W^{k,p}}
     + \left\|\xi_2\right\|_{W^{k,p}}
     \le c_k\left\|\xi_0\right\|_{W^{k,p}}^2,\qquad
     \left\|A_0-\bar A_0\right\|_{W^{k,p}}
     \le c_k\left\|\xi_0\right\|_{W^{k,p}},
$$
where $\xi_1$, $\xi_2$, and $A_0$ are as in~(i).

\smallskip
\noindent{\bf (iii)}
The map $\Ff^0$ is smooth and
$
     d\Fcal^0(0)\xi_0=(\xi_0,\alpha_0),
$
where $\alpha_0\in\Om^1(\Sigma,\g_P)$
is uniquely determined by the equation
$$
     D\bar\p_{J,\bar A_0}(\bar u_0)\xi_0
     + X_{\alpha_0}(\bar u_0)^{0,1} = 0.
$$
\end{theorem}
 
Theorem~\ref{thm:M0} is a standard result in the theory 
of holomorphic curves (cf.~\cite{MS1}).  
It follows from Fredholm theory and an 
infinite dimensional version of the implicit function 
theorem. In most applications the moduli space 
$\Mcal_{B,\Sigma}^0$ is not compact.  
However, it can be exhausted by the compact subsets 
$$
   \Mcal_{B,\Sigma}^0(c_0)
   :=\Tilde{\Mcal}_{B,\Sigma}^0(c_0)/\Gcal(P),
$$   
where $c_0>0$ and 
$$   
    \Tilde{\Mcal}_{B,\Sigma}^0(c_0)
    := \left\{(u,A)\in\Tilde{\Mcal}_{B,\Sigma}^0\,|\,
       \left\|d_Au\right\|_{L^\infty} 
       + \left\|F_A\right\|_{L^\infty}
       \le c_0
       \right\}.
$$   
Note that $\Tilde{\Mcal}_{B,\Sigma}^0(c_0)$ 
is invariant under the action of $\Gcal(P)$.
For later reference we prove the following lemma.
     
\begin{lemma}\label{le:c0}
Fix a reference connection $\hat A\in\Aa(P)$.
Then, for every $c_0>0$ and every integer $\ell\in\NN$,
there exists a constant $c=c(c_0,\ell)>0$ such that,
for every $(u_0,A_0)\in\Tilde{\Mm}_{B,\Sigma}^0(c_0)$,
we have
$$
     \inf_{g\in\Gg(P)}
     \left(
     \|g^{-1}u\|_{C^\ell}
     + \|g^*A_0-\hat A\|_{C^\ell}
     \right)
     \le c.
$$
\end{lemma}

\begin{proof} 
Suppose, by contradiction, that there 
exists a sequence $(u_\nu,A_\nu)\in\Tilde{\Mm}_{B,\Sigma}^0(c_0)$
such that 
$
     \|g^{-1}u_\nu\|_{C^\ell}
     + \|g^*A_\nu-\hat A\|_{C^\ell}\ge\nu
$
for every $\nu$ and every $g\in\Gg(P)$. By~\cite[Theorem~B.4.2]{MS1}
there exists a subsequence, still denoted by $(u_\nu,A_\nu)$,
such that the induced maps $\bar u_\nu:\Sigma\to M\dslash\G$
converge in the $\Cinf$-topology to a 
smooth $\bar J$-holomorphic curve.  The limit
curve represents the same homotopy class 
as the approximating curves and hence
can be represented by a pair 
$(u,A)\in\Tilde{\Mm}_{B,\Sigma}^0(c_0)$.
Since the sequence $\bar u_\nu$ converges to 
$\bar u:\Sigma\to M\dslash\G$ in the $C^{\ell+1}$-topology,
there exists a constant $\nu_0\ge 0$
such that, for every $\nu\ge\nu_0$,
there exist a gauge transformation $g_\nu$ 
and a section $\xi_\nu\in\Om^0(\Sigma,H_u)$ such that
$$
    g_\nu^{-1}u_\nu = \exp_u(\xi_\nu),\qquad
    \lim_{\nu\to\infty}\left\|\xi_\nu\right\|_{C^{\ell+1}}=0.
$$
The formulae
$$
    g_\nu^*A_\nu = -(L_{g_\nu^{-1}u_\nu}^*
    L_{g_\nu^{-1}u_\nu})^{-1}L_{g_\nu^{-1}u_\nu}d(g_\nu^{-1}u_\nu),\qquad
    A = -(L_u^*L_u)^{-1}L_udu,
$$
show that $g_\nu^*(u_\nu,A_\nu)$ converges to $(u,A)$ 
in the $C^\ell$ topology.  
This contradicts the choice of the sequence
$(u_\nu,A_\nu)$ and hence proves the lemma. 
\end{proof}

\begin{theorem}\label{thm:jhol}
Assume~$(H1)$ and~$(H4)$, let $\bar B\in H_2(M;\Z)$ be a nontorsion
homology class, and let $(\Sigma,\dvol_\Sigma,j_\Sigma)$
be a compact Riemann surface. Then, for every $c_0>0$ 
and every $p>2$, there exist positive constants $c$ and $\delta$ 
such that the following holds. If $\bar u:\Sigma\to\bar M$
is a smooth map such that $[\bar u]=\bar B$ and
$$
     \|d\bar u\|_{L^\infty}\le c_0,\qquad
     \|\bar\p_{\bar J}(\bar u)\|_{L^p}\le\delta
$$
then there exists a section 
$\bar\xi\in\Om^0(\Sigma,\bar u^*T\bar M)$ such that
$$
     \bar\p_{\bar J}(\exp_{\bar u}(\bar\xi))=0,\qquad
     \|\bar\xi\|_{W^{1,p}}\le c
     \|\bar\p_{\bar J}(\bar u)\|_{L^p}.
$$
\end{theorem}

\begin{proof}
This is again a standard result for 
pseudoholomorphic curves and the proof is almost
word by word the same as that of~\cite[Theorem~2.5]{DS2}.
Here is a sketch. One argues by contradiction. 
If the result were false,
there would be a sequence of smooth maps 
$\bar u_i:\Sigma\to\bar M$ that satisfies
$$
    \sup_i\|d\bar u_i\|_{L^\infty}<\infty,\qquad
    \lim_{i\to\infty}\|\bar\p_{\bar J}(\bar u_i)\|_{L^p} = 0,
$$
but which does not satisfy the conclusion of the theorem
for any constant~$c$. This means that the $W^{1,p}$-distance 
of $\bar u_i$ to the space of $\bar J$-holomorphic 
curves is not controlled uniformly by the $L^p$-norm
of $\bar\p_{\bar J}(\bar u_i)$.  Now, 
by the Arz\'ela--Ascoli and Banach--Alaoglu theorems, 
a suitable subsequence of $\bar u_i$ converges,
strongly with respect to the sup-norm and weakly
in $W^{1,p}$, to a $\bar J$-holomorphic 
curve $\bar u$. It follows from standard elliptic
regularity for $\bar J$-holomorphic curves
that $\bar u_i$ then converges 
strongly with respect to the $W^{1,p}$-norm.
To see this, write $\bar u_i=\exp_{\bar u}(\bar\xi_i)$
and observe that
\begin{eqnarray*}
     \|\bar\xi_i\|_{W^{1,p}}
&\le &
     c_1\left(
     \|D_{\bar u}\bar\xi_i\|_{L^p}
     + \|\bar\xi_i\|_{L^p}
     \right)  \\
&\le &
     c_2\left(
     \|\bar\p_{\bar J}(\bar u_i)\|_{L^p}
     + \|\bar\xi_i\|_{W^{1,p}}\|\bar\xi_i\|_{L^\infty}
     \right)
     + c_1\|\bar\xi_i\|_{L^p}.
\end{eqnarray*}
Here the first inequality is the elliptic estimate
for the Cauchy--Riemann operator $D_{\bar u}$
and the second is the quadratic estimate for 
$\bar\p_{\bar J}$. 
With this established it follows from hypothesis~(H4)
and the implicit function theorem for the operator 
$\bar\p_{\bar J}$ that there exists a sequence of 
$\bar J$-holomorphic curves $\bar u_{0i}$ whose
$W^{1,p}$-distance to $\bar u_i$ is bounded above by 
a fixed constant times the $L^p$-norm of
$\bar\p_{\bar J}(\bar u_i)$ (see~\cite[Theorem~2.1]{DS2}).
This shows that the sequence $\bar u_i$ does after
all satisfy the conclusion of the theorem, in contradiction
to our assumption. 
\end{proof}


\section{Adiabatic limits}\label{sec:adiabatic}
     
Before stating our main results we introduce some notation.
Fix an equivariant homology class $\bar B\in H_2(\bar M;\ZZ)$, 
let $B:=\kappa(\bar B)$, and denote
$$
     \Bb := \left\{
     (u,A)\in\Cinf_\G(P,M)\times\Aa(P)\,|\,
     [u]=B\right\}.
$$
This space is an infinite dimensional Fr\^echet manifold 
with tangent space 
$$
     T_{(u,A)}\Bb := \Om^0(\Sigma,u^*TM/\G)\times\Om^1(\Sigma,\g_P).
$$
It carries an action of the gauge group
$\Gg=\Gg(P)$ by 
$
     g^*(u,A)=(g^{-1}u,g^*A).
$
Consider the vector bundle $\Ee\to\Bb$ with fibres
$$
     \Ee_{(u,A)}:=\Om^{0,1}(\Sigma,u^*TM/\G)\oplus\Om^0(\Sigma,\g_P)
$$
and the $\Gg$-equivariant section $\Ff^\eps:\Bb\to\Ee$ given by 
$$
     \Ff^\eps(u,A):=(\bar\p_{J,A}(u),*F_A+\eps^{-2}\mu(u)).
$$
The zero set of this section is the space 
$$
     \Tilde{\Mm}^\eps_{B,\Sigma} 
     := \left\{(u,A)\in\Bb\,|\,
        u\mbox{ and }A\mbox{ satisfy }(\ref{eq:eps})\right\}.
$$
Its quotient by the action of the gauge group will be denoted by
$$
     \Mm^\eps_{B,\Sigma}
     := \Tilde{\Mm}^\eps_{B,\Sigma}/\Gg(P).
$$
The following theorem asserts the existence of solutions
of~(\ref{eq:eps}) for sufficiently small $\eps$ near every
regular solution of~(\ref{eq:jhol}).  The result is 
quantitative and the estimates are expressed in terms 
of suitable $\eps$-dependent norms.  Moreover, an operator
$
     \Dd^\eps:T_{(u,A)}\Bb\to\Ee_{(u,A)}\oplus\Om^0(\Sigma,\g_P)
$
appears.  This operator is the augmented vertical differential 
of $\Ff^\eps$.  The operator and the norms 
will be defined in Section~\ref{sec:linear}.  

\bigskip
\noindent{\bf Theorem~B}.
{\it
Assume~$(H1)$ and~$(H4)$ and let $\bar B\in H_2(\bar M;\Z)$
be a nontorsion homology class.  Then, for every $c_0>0$ and every $p>2$,
there exist positive constants $\eps_0$, $c$, and $\delta$ such that 
for every $\eps \in (0,\eps_0]$ there exists a 
$\Gcal(P)$-equivariant map
$$
     \Tilde{\Tt}^{\eps}:
     \Tilde{\Mm}_{B,\Sigma}^0(c_0)
     \to \Tilde{\Mm}_{B,\Sigma}^{\eps}
$$
that satisfies the following conditions. 

\smallskip
\noindent{\bf (a)}
If $(u_0,A_0)\in\Tilde{\Mcal}_{B,\Sigma}^0(c_0)$ then
$$
      \Tilde{\Tcal}^{\eps}(u_0,A_0)
      = (\exp_{u_0}(\xi_\eps),A_0+\alpha_\eps),
$$
where 
$
     \zeta_\eps=(\xi_\eps,\alpha_\eps)\in T_{(u_0,A_0)}\Bb
$
satisfies
$$
     \left\|\zeta_\eps\right\|_{2,p,\eps;(u_0,A_0)}
     \le c\eps^2,
$$
\begin{equation}\label{eq:Teps2}
     -d_{A_0}^*\alpha_\eps + \eps^{-2}L_{u_0}^*\xi_\eps=0,\qquad
     \zeta_\eps\in\im\,({\Dcal_{(u_0,A_0)}^\eps})^*.
\end{equation}

\smallskip
\noindent{\bf (b)}
If $(u_0,A_0)\in\Tilde{\Mcal}_{B,\Sigma}^0(c_0)$ and 
$
     (u,A) = (\exp_{u_0}(\xi),A_0+\alpha)
     \in\Tilde{\Mm}_{B,\Sigma}^{\eps}
$
where 
$
     \zeta=(\xi,\alpha)\in T_{(u_0,A_0)}\Bb
$
satisfies~(\ref{eq:Teps2}) and
$$
     \left\|\zeta\right\|_{1,p,\eps;(u_0,A_0)}
     \le \delta\eps^{2/p+1/2},
$$
then $(u,A) = \Tilde\Tcal^\eps(u_0,A_0)$.
}

\medbreak

The map $\Tilde{\Tcal}^{\eps}$ of Theorem~A descends 
to a map between the quotient 
spaces which we denote~by 
$$
   \Tt^\eps:\Mm_{B,\Sigma}^0(c_0)
     \TO\Mm_{B,\Sigma}^\eps.
$$
Assertion~(a) is proved by a Newton type iteration 
(see Section~\ref{sec:exist}).
It requires linear and quadratic estimates for the 
$\eps$-dependent norms with constants that are independent
of $\eps$.  These estimates are proved in Sections~\ref{sec:linear}
and~\ref{sec:quad}. Assertion~(b) is a strengthened form of the 
corresponding uniqueness statement. Here the neighbourhood 
in which uniqueness holds is larger than in the existence
result (namely it is of radius $c\eps^{2/p+1/2}$ 
instead of $c\eps^2$). The uniqueness statement
shows that the maps $\Tt^\eps$ are independent 
of~$c_0$ in the sense that
two such maps corresponding to different values of~$c_0$
(but the same value of~$\eps$) agree on the intersection
of their domains. The next theorem shows that
$\Tt^\eps$ is locally surjective.

\bigskip
\noindent{\bf Theorem~C}.
{\it
Assume~$(H1)$ and~$(H4)$ and let $\bar B\in H_2(\bar M;\Z)$
be a nontorsion homology class. 
Then, for every $c_0>0$ and every $p>2$,
there exist positive constants $\eps_0$ and $\delta$ 
such that the following holds for every $\eps\in(0,\eps_0]$. 
If 
$$
     (\bar u_0,\bar A_0)\in\Tilde{\Mcal}_{B,\Sigma}^0(c_0-1),\qquad
     (u,A) = (\exp_{\bar u_0}(\bar\xi),\bar A_0+\bar\alpha)
     \in \Tilde\Mcal_{B,\Sigma}^\eps,
$$
where 
$
     \bar\zeta=(\bar\xi,\bar\alpha) \in T_{(\bar u_0,\bar A_0)}\Bb
$
satisfies
$$
     \left\|\bar\zeta\right\|_{1,p,\eps;(\bar u_0,\bar A_0)}
     \le \delta\eps^{2/p+1/2},
$$
then $(u,A)\in\Tilde{\Tt}^\eps(\Tilde{\Mm}_{B,\Sigma}^0(c_0))$.
}

\bigskip

This result is restated more precisely in 
Theorem~\ref{thm:loc-onto} in Section~\ref{sec:tub}. 
There it is proved that 
$$
     g^*(u,A) = \Tilde\Tcal^\eps(u_0,A_0)
$$
for some gauge transformation $g$ and some pair
$(u_0,A_0)$ in the image of the map $\Ff^0$ 
of Theorem~\ref{thm:M0}.  Moreover, it is shown that
the distances of $g$ to $\one$ (in the $(2,p,\eps)$-norm)
and of $(u_0,A_0)$ to $(\bar u_0,\bar A_0)$
(in any norm) are controlled by the $(1,p,\eps)$-norm
of $\bar\zeta$. 

Theorem~C strengthens the local uniqueness result of Theorem~B~(b) 
in that condition~(\ref{eq:Teps2}) is no longer required.  
The proof relies on an $\eps$-dependent
local slice theorem (Section~\ref{sec:gauge})
and on the construction of a tubular neighbourhood
of the moduli space $\Mm_{B,\Sigma}^0(c_0)$ in which the 
normal bundle is the intersection of the $\eps$-dependent
local slice with the image of the adjoint operator
${\Dd^\eps}^*$ (Section~\ref{sec:tub}). 

The next theorem strengthens the local surjectivity 
result of Theorem~C.  It does not require the solution
$(u,A)$ of~(\ref{eq:eps}) to be close to any solution 
of~(\ref{eq:jhol}).  However, it only applies to solutions
that satisfy a uniform $L^\infty$-bound 
on $d_Au$ and for which $u$ takes values 
in the compact set 
$$
     M^C
     := \left\{x\in M\,|\,|\mu(x)|\le C,\,
        |\eta|\le C|L_x\eta|\,\,\forall\,\,\eta\in\g\right\}.
$$

\noindent{\bf Theorem~D}.
{\it
Assume~$(H1)$ and~$(H4)$ and let $\bar B\in H_2(\bar M;\Z)$
be a nontorsion homology class. Then, for every $C>0$, 
there exist positive constants $\eps_0$ and $c_1$ such that 
the following holds for every $\eps\in(0,\eps_0]$.
If $(u,A)\in\Tilde{\Mcal}_{B,\Sigma}^\eps$
such that 
\begin{equation}\label{eq:infty}
     \|d_Au\|_{L^\infty}\le C,\qquad
     u(P)\subset M^C
\end{equation}
then 
$
     (u,A)\in\Tilde{\Tcal}^\eps
     (\Tilde{\Mcal}_{B,\Sigma}^0(c_1)).
$
}

\bigskip

Under hypotheses~$(H1-2)$ the 
moduli space $\Mm^\eps_{B,\Sigma}$ is compact~\cite{CGMS}.  
In this case all solutions of~(\ref{eq:eps}) satisfy 
$$
     \|d_Au\|_{L^\infty} + \|\mu(u)\|_{L^\infty}\le C_\eps
$$
for some $\eps$-dependent constant $C_\eps$.  
However, this does not guarantee surjectivity 
because, on the one hand, the constant $C_\eps$ may
diverge to infinity as $\eps$ tends to zero
and, on the other hand, the solutions of~(\ref{eq:eps})
may not all satisfy the second condition in the 
definiton of $M^C$, namely that the image of $u$ 
belongs to the set of regular points of $\mu$. 
There may be sequences $(\eps_i,u_i,A_i)$
of solutions of~(\ref{eq:eps}) with $\eps_i\to0$ 
such that either $u_i(P)$ intersects the set 
of singular points of $\mu$ or 
$d_{A_i}u_i$ does not stay bounded,
and then bubbling occurs in the small $\eps$ limit.  
Under the hypotheses of Theorem~A we shall prove that
such bubbling cannot occur and establish a bijection
between suitable zero dimensional moduli spaces. 


\section{Linear estimates}\label{sec:linear}

The estimates in this section follow the 
ones in~\cite[Section~4]{DS2}.
In adapting the proofs to the present context 
we encounter additional zeroth order terms.
These arise from the Levi-Civita connection and the 
almost complex structure on $M$;
they are not present in~\cite{DS2} 
where $M$ is replaced by the space of
connections over a Riemann surface
and the almost complex structure 
by the Hodge $*$-operator.
We extend the results of~\cite{DS2}
by including estimates for the second derivatives. 
Moreover, in the present case it is crucial that
the constants depend continuously on the pair 
$(u,A)$.  In~\cite{DS2} the moduli space is 
a finite set and so the question of continuous 
dependence does not arise. 

For $u\in\Cinf_\G(P,M)$ we introduce the spaces
\begin{eqnarray*}
     \Xx_u 
     &:= &\Om^0(\Sigma,u^*TM/\G) 
        \times \Om^1(\Sigma,\g_P), \\
     \Xx'_u 
     &:= &\Om^{0,1}(\Sigma,u^*TM/\G) 
        \times \Om^0(\Sigma,\g_P)
        \times\Om^0(\Sigma,\g_P).
\end{eqnarray*}
Thus $\Xx_u = T_{(u,A)}\Bb$ and 
$\Xx'_u = \Ee_{(u,A)}\times\Om^0(\Sigma,\g_P)$
for every $A\in\Aa(P)$. 
If the map $u$ is understood from the context 
then we shall omit the subscript $u$. 
It is convenient to introduce the norms
\begin{eqnarray*}
     \left\|\xi\right\|_{1,p,\eps;(u,A)}  
&:= &
     \left\|\xi\right\|_{L^p}
     + \eps\left\|\Tabla{A}\xi\right\|_{L^p}, \\
     \left\|\xi\right\|_{2,p,\eps;(u,A)}  
&:= &
     \left\|\xi\right\|_{L^p}
     + \eps\left\|\Tabla{A}\xi\right\|_{L^p}
     + \eps^2\left\|{\Tabla{A}}^*\Tabla{A}\xi\right\|_{L^p}, \\
     \left\|\alpha\right\|_{1,p,\eps;A}
&:= &
     \left\|\alpha\right\|_{L^p}
     + \eps\left\|d_A\alpha\right\|_{L^p}
     + \eps\left\|d_A^*\alpha\right\|_{L^p},  \\
     \left\|\alpha\right\|_{2,p,\eps;A}
&:= &
     \left\|\alpha\right\|_{L^p}
     + \eps\left\|d_A\alpha\right\|_{L^p}
     + \eps\left\|d_A^*\alpha\right\|_{L^p}
     + \eps^2\left\|d_A^*d_A\alpha+d_Ad_A^*\alpha\right\|_{L^p}
\end{eqnarray*}
for $\xi\in\Om^{0,k}(\Sigma,u^*TM/\G)$, 
$\alpha\in\Om^k(\Sigma,\g_P)$, 
$k=0,1$, and $1\le p\le\infty$. 
Here $\Tabla{A}$ denotes the Hermitian
connection on $u^*TM/\G$ defined by~(\ref{eq:Tabla}).
For $\ell=0,1,2$, $1\le p\le\infty$, and
$$
     \zeta=(\xi,\alpha)\in \Xx_u,\qquad
     \zeta' = (\xi',\phi',\psi')\in\Xx'_u
$$
we consider the norms
\begin{eqnarray*}
     \left\|\zeta\right\|_{\ell,p,\eps;(u,A)}  
&:= & 
     \left\|\xi\right\|_{\ell,p,\eps;(u,A)}
     + \eps\left\|\alpha\right\|_{\ell,p,\eps;(u,A)}, \\
     \left\|\zeta'\right\|_{\ell,p,\eps;(u,A)}  
&:= & 
     \left\|\xi'\right\|_{\ell,p,\eps;(u,A)}  
     + \eps\left\|\phi'\right\|_{\ell,p,\eps;(u,A)}  
     + \eps\left\|\psi'\right\|_{\ell,p,\eps;(u,A)},
\end{eqnarray*}
where 
$
     \left\|\xi\right\|_{0,p,\eps;(u,A)}
     := \left\|\xi\right\|_{L^p}.
$
These norms are gauge invariant, e.g.
$$
     \left\|(g^{-1}\xi,g^{-1}\alpha g)\right\|_{\ell,p,\eps;(g^{-1}u,g^*A)}
     = \left\|(\xi,\alpha)\right\|_{\ell,p,\eps;(u,A)}.
$$
If the pair $(u,A)$ is understood from the context 
we shall drop it to simplify the notation.
In particular, we abbreviate
$$
     \left\|\zeta\right\|_{\infty,\eps}
     := \left\|\zeta\right\|_{0,\infty,\eps;(u,A)}.
$$

The augmented vertical differential of $\Ff^\eps$ 
at a zero $(u,A)\in\Bb$ is the operator 
$$
     \Dd^\eps=\Dd_{(u,A)}^{\eps}:
     \Xx_u \to \Xx_u'
$$
given by 
\begin{equation}\label{eq:Deps}
     \Dd^{\eps}
     \left(\begin{array}{c}\xi \\ \alpha\end{array}\right)
     = \left(
       \begin{array}{c}
         D\xi + (L_u\alpha)^{0,1} \\
         \eps^{-2}L_u^*\xi - d_A^*\alpha  \\
         \eps^{-2}d\mu(u)\xi + *d_A\alpha
       \end{array}
       \right),
\end{equation}
where
$
     D = D_{(u,A)}:
     \Om^0(\Sigma,u^*TM/\G)
     \to \Om^{0,1}(\Sigma,u^*TM/\G)
$
is the Cauchy--Riemann operator defined by~(\ref{eq:Du}). 
The second coordinate in the definition  of $\Dd^\eps$
corresponds to the local slice condition
for the $\Gg$-action. For the definition of the adjoint operator 
it is convenient to use the $\eps$-dependent inner products
associated to the $(0,2,\eps)$-norms.  In addition we use twice
the standard inner product on the space $\Om^{0,1}(\Sigma,u^*TM/\G)$.
Then the adjoint of $\Dd^\eps$ is given by
$$
     {\Dd^{\eps}}^*
     \left(\begin{array}{c}\xi' \\ \phi' \\ \psi'
     \end{array}\right)
     = \left(
       \begin{array}{c}
         2D^*\xi' + L_u\phi' + JL_u\psi' \\ 
         2\eps^{-2}L_u^*\xi' - d_A\phi' - *d_A\psi'
       \end{array}
       \right)
$$
for 
$
     (\xi',\phi',\psi')
     \in\Xx'.
$
The sole purpose of the factor $2$
is to render the off-diagonal terms in the operator
${\Dd^\eps}^*\Dd^\eps$ of zeroth order.

\begin{remark}[Local coordinates]\label{rmk:local}\rm
Let $\upsilon:U\to\Sigma$ be a holomorphic coordinate chart 
defined on an open set $U\subset\C$ 
and let $\tilde\upsilon:U\to P$ 
be a lift of~$\upsilon$.   
In this trivialization the map $u$, the connection $A$, 
the vector field~$\xi$ along $u$,
and the $1$-form $\alpha$
are represented by 
$$
\begin{array}{rclcrcl}
     u^\loc &:= & u\circ\tilde\upsilon, &\qquad &
     A^\loc &:= & {\tilde\upsilon\,}^*A = \Phi\,ds+\Psi\,dt, \\
     \xi^\loc &:= & \xi\circ\tilde\upsilon, &\qquad &
     \alpha^\loc &:= & {\tilde\upsilon\,}^*\alpha=\phi\,ds+\psi\,dt,
\end{array}
$$
where $\Phi,\Psi,\phi,\psi$ are Lie algebra valued
functions on $U$.  The volume form on $U$ is given by
$$
    \lambda^2\,ds\wedge dt := \upsilon^*\dvol_\Sigma
$$
for some function $\lambda:U\to(0,\infty)$
and the metric has the form $\lambda^2(ds^2+dt^2)$.
{From} now on we shall drop the superscript ``$\loc$''
and introduce the notation
$$
\begin{array}{rlrl}
     v_s\,\,:=&\p_su+X_\Phi(u), &
     v_t\,\,:=&\p_tu+X_\Psi(u), \\
     \Nabla{A,s}\xi\,\,:=&\Nabla{s}\xi+\Nabla\xi X_{\Phi}(u), &
     \Nabla{A,t}\xi\,\,:=&\Nabla{t}\xi+\Nabla\xi X_{\Psi}(u), \\
     \Tabla{A,s}\xi\,\,:=&\Nabla{A,s}\xi 
       - \frac12 J(\Nabla{v_s}J+\p_sJ)\xi, &
     \Tabla{A,t}\xi\,\,:=&\Nabla{A,t}\xi 
       - \frac12 J(\Nabla{v_t}J+\p_tJ)\xi, \\
     \Nabla{A,s}\eta\,\,:=&\p_s\eta+[\Phi,\eta], &
     \Nabla{A,t}\eta\,\,:=&\p_t\eta+[\Psi,\eta], 
\end{array}
$$
for $\eta:U\to\g$ and a vec\-tor field $\xi:U\to u^*TM$ 
along $u$. Then 
$$
     d_Au 
     = v_s\,ds + v_t\,dt, \qquad
     \Nabla{A}\xi 
     = \Nabla{A,s}\xi\,ds + \Nabla{A,t}\xi\,dt,
$$
and
$$
\begin{array}{rcl}
     *F_A &=& \lambda^{-2}
     \left(\p_s\Psi-\p_t\Phi+[\Phi,\Psi]\right),  \\
     *d_A\alpha &=& \lambda^{-2}
     \left(\Nabla{A,s}\psi-\Nabla{A,t}\phi\right), \\
     d_A^*\alpha
     &=&
     -\lambda^{-2}
     (\Nabla{A,s}\phi+\Nabla{A,t}\psi).
\end{array}
$$
In local coordinates a $(0,1)$-form on $\Sigma$ with values in 
$u^*TM/\G$ has the form $\frac12(\xi'ds - J\xi'dt)$,
where $\xi'(s,t)\in T_{u(s,t)}M$. In particular,
\begin{eqnarray*}
     \bar\p_{J,A}(u) 
&= &
     \frac12(v_s+Jv_t)\,ds + \frac12(v_t-Jv_s)\,dt, \\
     (\Nabla{A}\xi)^{0,1} 
&= &
     \frac12\left(\Nabla{A,s}\xi+J\Nabla{A,t}\xi\right)ds 
     + \frac12\left(\Nabla{A,t}\xi-J\Nabla{A,s}\xi\right)dt,
\end{eqnarray*}
We represent a $(0,1)$-form by twice the coefficient of
$ds$.  Then 
$$
     D\xi = 
     \Tabla{A,s}\xi + J\Tabla{A,t}\xi 
     + \frac14N(\xi,v_s-Jv_t) + \frac12(J\p_sJ-\p_tJ)\xi,
$$
where 
$
     N(\xi_1,\xi_2)=2J((\Nabla{\xi_2}J)\xi_1-(\Nabla{\xi_1}J)\xi_2)
$
denotes the Nijenhuis tensor, and
$$
     2D^*\xi' = 
     \frac{1}{\lambda^2}\left(-\Tabla{A,s}\xi' + J\Tabla{A,t}\xi'
       - \frac12J(\Nabla{\xi'}J)(v_s-Jv_t)
       + \frac12(J\p_sJ-\p_tJ)\xi\right).
$$
The Weitzenb\"ock formula has the form
$
     D^*D\xi = \frac12{\Tabla{A}}^*\Tabla{A}\xi
     + \mbox{l.o.t.}
$
In the K\"ahler case we have $\Nabla{A}=\Tabla{A}$,
$\nabla J=0$, and $\p_sJ=\p_tJ=0$.  Hence in this case
\begin{eqnarray*}
     D^*D\xi
&= &
     -\frac{1}{2\lambda^2}\left(
       \Nabla{A,s}\Nabla{A,s}\xi+\Nabla{A,t}\Nabla{A,t}\xi
       \right)  \\
&&
     - \frac{1}{2\lambda^2}JR(v_s,v_t)\xi
     - \frac{1}{2\lambda^2}J\Nabla{\xi}
     X_{\p_s\Psi-\p_t\Phi+[\Phi,\Psi]}(u).
\end{eqnarray*}
In local coordinates the operators $\Dd^\eps$ and
${\Dd^\eps}^*$ have the form
$$
     \Dd^\eps\zeta
     = \left(\begin{array}{c}
       D\xi + L_u\phi + JL_u\psi \\
       \lambda^{-2}(\Nabla{A,s}\phi + \Nabla{A,t}\psi)
       + \eps^{-2}L_u^*\xi \\
       \lambda^{-2}(\Nabla{A,s}\psi - \Nabla{A,t}\phi)
       + \eps^{-2}d\mu(u)\xi
       \end{array}\right),
$$
$$
     {\Dd^\eps}^*\zeta'
     = \left(\begin{array}{c}
       2D^* + L_u\phi' + JL_u\psi' \\
       -\Nabla{A,s}\phi' + \Nabla{A,t}\psi'
       + \eps^{-2}L_u^*\xi' \\
       -\Nabla{A,s}\psi' - \Nabla{A,t}\phi'
       + \eps^{-2}d\mu(u)\xi'
       \end{array}\right).
$$
\end{remark}
 
\begin{proposition}\label{prop:estimate}
For every $p\ge 2$ and every $c_0>0$ there exist 
positive constants $\eps_0$ and $c$ such that
\begin{eqnarray*}
     \|\zeta\|_{1,p,\eps;(u,A)}
&\le & 
     c\left(
     \eps\|\Dd^\eps\zeta\|_{0,p,\eps}
     +\|\pi_u\xi\|_{L^p}
     \right), \\
     \|\zeta-\pi_u\zeta\|_{1,p,\eps;(u,A)}
&\le &
     c\eps\left(
     \|\Dd^\eps\zeta\|_{0,p,\eps}
     +\|\pi_u\xi\|_{L^p}
     \right), \\
     \|\zeta'\|_{1,p,\eps;(u,A)}
&\le & 
     c\left(
     \eps\|{\Dd^\eps}^*\zeta'\|_{0,p,\eps}
     +\|\pi_u\xi'\|_{L^p}
     \right), \\
     \|\zeta'-\pi_u\zeta'\|_{1,p,\eps;(u,A)}
&\le &
     c\eps\left(
     \|{\Dd^\eps}^*\zeta'\|_{0,p,\eps}
     +\|\pi_u\xi'\|_{L^p}
     \right),
\end{eqnarray*}
for all $(u,A)\in\Tilde{\Mm}_{B,\Sigma}^0(c_0)$,
$
     \zeta=(\xi,\alpha)\in\Xx_u,
$
$
     \zeta'=(\xi',\phi',\psi')\in\Xx'_u,
$
and $\eps\in(0,\eps_0]$.  Here we abbreviate
$\Dd^\eps:=\Dd^\eps_{(u,A)}$ and
$\pi_u\zeta:=(\pi_u\xi,0)$ and 
$\pi_u\zeta':=(\pi_u\xi',0,0)$, 
where $\pi_u$ is defined by~(\ref{eq:pi-u}). 
\end{proposition}

In this paper we prove Proposition~\ref{prop:estimate}  
only in the case $p=2$.  The proof for $p>2$ is similar to 
the proof of an analogous result in~\cite{Sa}. 
The details for the present case will 
be carried out elsewhere. 

\begin{lemma}\label{le:adjoint}
{\bf (i)}
If $\bar\p_{J,A}(u)=0$ then
$$
     \Dd^\eps{\Dd^{\eps}}^*\zeta'
     = \left(
       \begin{array}{c}
         2DD^*\xi' + 2\eps^{-2}(L_uL_u^*\xi')^{0,1} 
          + (DJ-JD)L_u\psi' \\ 
         \Delta_\eps\phi' + [*F_A+\eps^{-2}\mu(u),\psi']  \\
         \Delta_\eps\psi' - [*F_A+\eps^{-2}\mu(u),\phi']
          + 2\eps^{-2}L_u^*(DJ-JD)^*\xi'
       \end{array}
       \right)
$$
for $\zeta':=(\xi',\phi',\psi')\in\Xx'_u$, where  
$\Delta_\eps:=d_A^*d_A+\eps^{-2}L_u^*L_u$.

\smallskip
\noindent{\bf (ii)}
If $\bar\p_{J,A}(u)=0$ and $\mu(u)=0$ then
$$
     {\Dd^{\eps}}^*\Dd^\eps\zeta
     = \left(
       \begin{array}{c}
         2D^*D\xi + \eps^{-2}L_uL_u^*\xi + \eps^{-2}JL_uL_u^*J^*\xi
         + Q^*\alpha \\ 
         d_A^*d_A\alpha+d_Ad_A^*\alpha+\eps^{-2}L_u^*L_u\alpha 
         + \eps^{-2}Q\xi
       \end{array}
       \right),
$$
for $\zeta=(\xi,\alpha)\in\Xx_u$, where
$Q:\Om^0(\Sigma,u^*TM/\G)\to\Om^1(\Sigma,\g_P)$ denotes
the zeroth order operator
\begin{eqnarray*}
     Q\xi 
&:= &
     2L_u^*D\xi - d_AL_u^*\xi - *d_Ad\mu(u)\xi \\   
&= &
     \rho(\xi,d_Au) - *\rho(J\xi,d_Au) 
     + *L_u^*\dot J\xi + \frac12L_u^*N(\xi,\p_{J,A}(u)).
\end{eqnarray*}
\end{lemma}

\begin{proof}
We shall repeatedly use the identities
$$
\begin{array}{rclcrcl}
     d_A^*\alpha &= & -*d_A*\alpha, &\qquad &
     *d_Ad_A\phi &= & [*F_A,\phi],  \\
           L_u^* &= & d\mu(u)J,  &\qquad &
  d\mu(u)L_u\phi &= & - [\mu(u),\phi]
\end{array}
$$
for $\alpha\in\Om^1(\Sigma,\g_P)$ and 
$\phi\in\Om^0(\Sigma,\g_P)$.  To prove~(i) note 
that the triple 
$
     (\tilde\xi',\tilde\phi',\tilde\psi')
     := \Dd^\eps{\Dd^\eps}^*(\xi',\phi',\psi')
$
is given by 
\begin{eqnarray*}
     \tilde\xi' 
&= &
     D(2D^*\xi'+L_u\phi'+JL_u\psi')
     + (L_u(2\eps^{-2}L_u^*\xi'-d_A\phi'-*d_A\psi'))^{0,1} \\
&= &
     2DD^*\xi' + 2\eps^{-2}(L_uL_u^*\xi')^{0,1}  \\
&&
       +\, DL_u\phi'-(L_ud_A\phi')^{0,1} 
       + DJL_u\psi'-(L_u*d_A\psi')^{0,1}, \\
     \tilde\phi'
&= &
     \eps^{-2}L_u^*(2D^*\xi'+L_u\phi+JL_u\psi)
     -d_A^*(2\eps^{-2}L_u^*\xi'-d_A\phi-*d_A\psi) \\
&= &
     d_A^*d_A\phi' + \eps^{-2}L_u^*L_u\phi' 
       + 2\eps^{-2}(DL_u-L_ud_A)^*\xi'
       + [*F_A+\eps^{-2}\mu(u),\psi'],  \\
     \tilde\psi'
&= &
     \eps^{-2}d\mu(u)(2D^*\xi'+L_u\phi'+JL_u\psi')
     + *d_A(2\eps^{-2}L_u^*\xi'-d_A\phi'-*d_A\psi') \\
&= &
     d_A^*d_A\psi' + \eps^{-2}L_u^*L_u\psi' 
       - [*F_A+\eps^{-2}\mu(u),\phi']  \\
&&
       +\, 2\eps^{-2}(L_u^*J^*D^*\xi' + d_A^**L_u^*\xi').
\end{eqnarray*}
The assertion now follows from the fact that
\begin{equation}\label{eq:JLu}
     J(L_u\alpha)^{0,1} = (L_u*\alpha)^{0,1},\qquad
     L_u^*J^*\xi' = -*L_u^*\xi',
\end{equation}
for $\alpha\in\Om^1(\Sigma,\g_P)$ and
$\xi'\in\Om^{0,1}(\Sigma,u^*TM/\G)$,
and
\begin{equation}\label{eq:DLu}
     \bar\p_{J,A}(u)=0\qquad\IMP\qquad
     DL_u\phi = (L_ud_A\phi)^{0,1}
\end{equation}
for $\phi\in\Om^0(\Sigma,\g_P)$. 
The first equation in~(\ref{eq:JLu}) follows from 
the fact that $*\alpha = -\alpha\circ J_\Sigma$
for every $1$-form $\alpha$ on $\Sigma$
(with values in any vector bundle).
The second equation in~(\ref{eq:JLu}) follows
from the first by duality. 
Equation~(\ref{eq:DLu}) follows from the 
fact that the section 
$(u,A)\mapsto\bar\p_{J,A}(u)$ 
of the vector bundle over $\Bb$ with fibres
$\Om^{0,1}(\Sigma,u^*TM/\G)$
is $\Gg(P)$-equivariant. 

To prove~(ii) note that the pair 
$
     (\tilde\xi,\tilde\alpha)
     := {\Dd^\eps}^*\Dd^\eps(\xi,\alpha)
$
is given by 
\begin{eqnarray*}
     \tilde\xi
&= &
     2D^*(D\xi+(L_u\alpha)^{0,1}) \\
&&
     +\, L_u(\eps^{-2}L_u^*\xi-d_A^*\alpha)
     + JL_u(\eps^{-2}d\mu(u)\xi+*d_A\alpha) \\
&= &
     2D^*D\xi + 2\eps^{-2}(L_uL_u^*\xi')^{0,1} 
     + \left(2D^*(L_u\alpha)^{0,1} - L_ud_A^*\alpha + JL_u*d_A\alpha
       \right),  \\
     \tilde\alpha
&= &
     2\eps^{-2}L_u^*(D\xi+(L_u\alpha)^{0,1})  \\
&&
     -\, d_A(\eps^{-2}L_u^*\xi-d_A^*\alpha) 
     - *d_A(\eps^{-2}d\mu(u)\xi+*d_A\alpha) \\
&= &
     d_A^*d_A\alpha 
     + d_Ad_A^*\alpha 
     + \eps^{-2}L_u^*L_u\alpha 
     + \eps^{-2}\left(2L_u^*D\xi - d_AL_u^*\xi - *d_Ad\mu(u)\xi\right).
\end{eqnarray*}
Here we have used the fact that $\mu(u)=0$ and hence 
$2L_u^*(L_u\alpha)^{0,1}=L_u^*L_u\alpha$. 
The formula for the operator 
$
     Q := 2L_u^*D - d_AL_u^* - *d_Ad\mu(u)
$
follows by computing in local coordinates. 
\end{proof}

\begin{proof}[Proof of Proposition~\ref{prop:estimate} for $p=2$]
Let $\zeta':=(\xi',\phi',\psi'):=\Dd^\eps\zeta$.  Then, by 
Lemma~\ref{le:adjoint}, the formula 
$$
     {\Dd^\eps}^*\Dd^\eps\zeta = {\Dd^\eps}^*\zeta'
$$
is equivalent to
\begin{eqnarray*}
     2D^*D\xi + \eps^{-2}L_uL_u^*\xi + \eps^{-2}JL_uL_u^*J^*\xi
     + Q^*\alpha
&= &
     2D^*\xi' + L_u\phi'+JL_u\psi',  \\
     d_A^*d_A\alpha+d_Ad_A^*\alpha+\eps^{-2}L_u^*L_u\alpha
     + \eps^{-2}Q\xi
&= &
     2\eps^{-2}L_u^*\xi' - d_A\phi' - *d_A\psi'.
\end{eqnarray*}
Take the $L^2$-inner product of the first equation with 
$\xi$ and of the second equation with $\eps^2\alpha$.
The sum of the resulting identities gives
\begin{eqnarray*}
&&
     \eps^{-2}\left\|L_u^*\xi\right\|^2
     + \eps^{-2}\left\|L_u^*J\xi\right\|^2
     + 2\left\|D\xi\right\|^2  
     + \left\|L_u\alpha\right\|^2
     + \eps^2\left\|d_A\alpha\right\|^2
     + \eps^2\left\|d_A^*\alpha\right\|^2  \\
&&\quad
     = 2\inner{\xi'}{D\xi} 
       + 2\inner{\xi'}{L_u\alpha} 
       - 2 \inner{\alpha}{Q\xi}  \\
&&\qquad
       +\, \inner{\phi'}{L_u^*\xi}
       - \inner{\psi'}{L_u^*J\xi}  
       - \eps^2\inner{\phi'}{d_A^*\alpha}
       + \eps^2\inner{\psi'}{*d_A\alpha}  \\
&&\quad
     \le 3\left\|\xi'\right\|^2 
       + \left\|D\xi\right\|^2
       + 2^{-1}\left\|L_u\alpha\right\|^2 
       + \delta\left\|\alpha\right\|^2
       + \delta^{-1}\left\|Q\xi\right\|^2   \\
&&\qquad
       +\, \eps^2\left\|\phi'\right\|^2
       + \eps^2\left\|\psi'\right\|^2
       + 2^{-1}\eps^{-2} \left\|L_u^*\xi\right\|^2 
       + 2^{-1}\eps^{-2}\left\|L_u^*J\xi\right\|^2 \\
&&\qquad
       +\, 2^{-1}\eps^2\left\|d_A^*\alpha\right\|^2
       + 2^{-1}\eps^2\left\|d_A\alpha\right\|^2. 
\end{eqnarray*}
Here all norms are $L^2$-norms
and all inner products are $L^2$-inner products. 
Choose $\delta>0$ so small that 
$
     \delta\left\|\alpha\right\|^2
     \le 4^{-1}\left\|L_u\alpha\right\|^2
$
for all $\alpha$.  Then 
\begin{eqnarray*}
&&
     \eps^{-2}\left\|L_u^*\xi\right\|^2
     + \eps^{-2}\left\|L_u^*J\xi\right\|^2
     + \left\|D\xi\right\|^2  
     + \left\|L_u\alpha\right\|^2
     + \eps^2\left\|d_A\alpha\right\|^2
     + \eps^2\left\|d_A^*\alpha\right\|^2  \\
&&\quad\le
     12\left\|\xi'\right\|^2 + 4\eps^2\left\|\phi'\right\|^2
     + 4\eps^2\left\|\psi'\right\|^2
     + 4\delta^{-1}\left\|Q\right\|_{L^\infty}^2
       \left\|\xi\right\|^2   \\
&&\quad\le
     12\left\|\Dd^\eps(\xi,\alpha)\right\|_{0,2,\eps}^2
     + 4\delta^{-1}\left\|Q\right\|_{L^\infty}^2
       \left\|\xi\right\|_{L^2}^2.
\end{eqnarray*}
Now the required estimates follow from
the inequalities
\begin{eqnarray*}
     \|\Tabla{A}\pi_u\xi\|_{L^2}
&\le &
     c\left(\|\Tabla{A}\xi\|_{L^2}
      + \left\|\xi\right\|_{L^2}
     \right), \\
     \left\|\xi-\pi_u\xi\right\|_{L^2}
&\le &
     c'\left(\left\|L_u^*\xi\right\|_{L^2}
        + \left\|L_u^*J\xi\right\|_{L^2}
        \right), \\
     \|\Tabla{A}\xi\|_{L^2}
&\le &
     c''\left(\|D\xi\|_{L^2}+ \|\xi\|_{L^2}\right).
\end{eqnarray*}
The first inequality follows from~(\ref{eq:difference})
below.  In the second inequality the constant
$c'$ can be chosen as an upper bound for the norms of 
the linear maps $L_x(L_x^*L_x)^{-1}$ 
over all $x\in\mu^{-1}(0)$. The third inequality 
is the $L^2$-estimate for the Cauchy--Riemann operator
and it follows from the Weitzenb\"ock formula. 
The constant $c''$ is gauge invariant and depends continuously 
on the pair $(u,A)$ with respect to the $C^1$-norm 
and hence can be chosen independent of 
$(u,A)\in\Tilde{\Mm}^0_{B,\Sigma}(c_0)$.  
This proves the proposition in the case $p=2$.
\end{proof}

The next lemma expresses the Sobolev inequalities
in terms of the $\eps$-dependent norms. 

\begin{lemma}\label{lemma:p}
For every $p>2$ and every $c_0>0$ there exists
a constant $c>0$ such that
$$
        \|\zeta\|_{\infty,\eps} 
        \le c{\eps}^{-2/p}\|\zeta\|_{1,p,\eps;(u,A)},\qquad
        \|\zeta\|_{1,\infty,\eps} 
        \le c{\eps}^{-2/p}
        \|\zeta\|_{2,p,\eps;(u,A)}
$$
for all $(u,A)\in\Tilde{\Mm}_{B,\Sigma}^0(c_0)$,
$
     \zeta\in\Xx_u,
$
and $\eps\in(0,1)$.
\end{lemma}

\begin{proof}
Multiply the metric on $\Sigma$ by $\eps^{-2}$.
Then the $W^{k,p}$-norm of $(\xi,\alpha)$ 
with respect to the rescaled metric is equal 
to $\eps^{-2/p}$ times the $(k,p,\eps)$-norm 
of $\zeta$, and the $L^\infty$-norm
with respect to the rescaled metric
is equal to the $(\infty,\eps)$-norm.
Hence the estimates follows from the Sobolev
embedding theorem for the rescaled metric. 
The constant is gauge invariant and it depends continuously 
on $u$ (with respect to the $C^1$-norm) 
and $A$ (with respect to the $C^0$-norm).
By Lemma~\ref{le:c0}, the estimate holds
with a uniform constant~$c$. 
\end{proof}

\begin{lemma}\label{le:difference}
For every $p\ge2$ and every $c_0>0$ 
there exist positive constants $\eps_0$ and $c$ such that
\begin{eqnarray*}
      \|\pi_u\Dd^\eps\zeta - \Dd^0\pi_u\zeta\|_{k,p,\eps}
&\le &
      c\|\xi-\pi_u\xi\|_{k,p,\eps}, \\
      \|\pi_u\Dd^\eps\zeta - \Dd^\eps\pi_u\zeta\|_{k,p,\eps}
&\le &
      c\|\xi\|_{k,p,\eps}, \\
      \|\pi_u{\Dd^\eps}^*\zeta' - {\Dd^0}^*\pi_u\zeta'\|_{k,p,\eps}
&\le &
      c\|\xi'-\pi_u\xi'\|_{k,p,\eps}, \\
      \|\pi_u{\Dd^\eps}^*\zeta' - {\Dd^\eps}^*\pi_u\zeta'\|_{k,p,\eps}
&\le &
      c\|\xi'\|_{k,p,\eps}
\end{eqnarray*}
for every $(u,A)\in\Tilde\Mcal_{B,\Sigma}^0(c_0)$,
$ 
     \zeta=(\xi,\alpha)\in\Xx_u,
$
$
     \zeta'=(\xi',\phi',\psi')\in\Xx'_u,
$
$\eps\in (0,1]$, and $k=0,1$.  Here we abbreviate 
$\Dd^\eps:=\Dd_{(u,A)}^\eps$ for $\eps\ge 0$. 
\end{lemma}

\begin{proof}
We prove first that, for every vector field $v\in\Vect(\Sigma)$,
there exists a constant $c=c(p,c_0,v)>0$ such that
\begin{equation}\label{eq:difference}
\begin{array}{rcl}
     \|\pi_u\Tabla{A,v}\xi-\Tabla{A,v}\pi_u\xi
     \|_{W^{k,p},A}
&\le &
     c\|\xi\|_{W^{k,p},A},\\
     \|\pi_u\Tabla{A,v}\xi-\pi_u\Tabla{A,v}\pi_u\xi
     \|_{W^{k,p},A}
&\le &
     c\|\xi-\pi_u\xi\|_{W^{k,p},A},\\
\end{array}
\end{equation}
for $(u,A)\in\Tilde{\Mm}^0_{B,\Sigma}(c_0)$,
$\xi\in\Om^0(\Sigma,u^*TM/\G)$ and $k=0,1$.
Here the $W^{1,p}$-norm labelled by $A$ is understood
as the (gauge invariant) $(1,p,\eps)$-norm for 
$\eps=1$.  
To prove~(\ref{eq:difference}) we choose 
local holomorphic coordinates $s+it$
on $\Sigma$. Thus $\xi(s,t)\in T_{(u(s,t)}M$, and 
$v_s$, $v_t$, $\Nabla{A,s}\xi$, and $\Nabla{A,t}\xi$ 
are as in Remark~\ref{rmk:local}. Write  
$$ 
     \xi = \pi_u\xi + L_u\eta_1 + JL_u\eta_2,
$$
where $\eta_i(s,t)\in\g$. 
Define $B_s(s,t):\g\to T_{u(s,t)}M$
and $B_t(s,t):\g\to T_{u(s,t)}M$ by
$$
     B_s\eta:=\Nabla{v_s}X_\eta(u),\qquad
     B_t\eta:=\Nabla{v_t}X_\eta(u).
$$
Then 
$$
     \Nabla{A,s}L_u\eta - L_u\Nabla{A,s}\eta = B_s\eta,\qquad
     \Nabla{A,t}L_u\eta - L_u\Nabla{A,t}\eta = B_t\eta
$$
and hence 
$$
     \Nabla{A,s}\pi_u\xi - \pi_u\Nabla{A,s}\xi 
     = \pi_u(B_s\eta_1+JB_s\eta_2+(\Nabla{v_s}J+\p_sJ)L_u\eta_2).
$$
Since
$
    \eta_1 = (L_u^*L_u)^{-1}L_u^*(\xi-\pi_u\xi)
$
and
$
    \eta_2 = -(L_u^*L_u)^{-1}L_u^*J(\xi-\pi_u\xi),
$
we have
$$
     \|\Nabla{A,s}\pi_u\xi - \pi_u\Nabla{A,s}\xi\|_{L^p}
     \le c\|\xi-\pi_u\xi\|_{L^p}.
$$
This proves~(\ref{eq:difference}) for the local
vector field $\p/\p s$.  For $\p/\p t$ the proof
is analogous.  Hence the result follows for any 
linear combination of these vector fields supported
in the given coordinate chart, and hence for every
vector field on $\Sigma$. For 
$\xi'\in\Om^{0,1}(\Sigma,u^*TM/\G)$
there are similar inequalities.  

By~(\ref{eq:difference}), there exists 
a constant $c'=c'(p,c_0)>0$ such that
$$
\begin{array}{rcl}
     \|\pi_u(\Tabla{A}\xi)^{0,1}-(\Tabla{A,v}\pi_u\xi)^{0,1}
     \|_{W^{k,p},A} 
&\le &
     c'\|\xi\|_{W^{k,p},A},\\
     \|\pi_u(\Tabla{A}\xi)^{0,1}-\pi_u(\Tabla{A,v}\pi_u\xi)^{0,1}
     \|_{W^{k,p},A} 
&\le &
     c'\|\xi-\pi_u\xi\|_{W^{k,p},A},
\end{array}
$$
for $(u,A)\in\Tilde{\Mm}^0_{B,\Sigma}(c_0)$,
$\xi\in\Om^0(\Sigma,u^*TM/\G)$, and $k=0,1$.  
Now the operator $D=D_{(u,A)}$ is given by 
$$
     D\xi = (\Nabla{A}\xi)^{0,1} - J(\Nabla{\xi}J)\p_{J,A}(u)
$$
and hence 
\begin{equation}\label{eq:Dpi}
\begin{array}{rcl}
     \|\pi_uD\xi-D\pi_u\xi\|_{W^{k,p},A} 
&\le &
     c''\|\xi\|_{W^{k,p},A},  \\
     \|\pi_uD(\xi-\pi_u\xi)\|_{W^{k,p},A} 
&\le &
     c''\|\xi-\pi_u\xi\|_{W^{k,p},A},
\end{array}
\end{equation}
for $(u,A)\in\Tilde{\Mm}^0_{B,\Sigma}(c_0)$,
$\xi\in\Om^0(\Sigma,u^*TM/\G)$, and $k=0,1$.
Since
$$
    \pi_u\Dd^\eps\zeta = \pi_uD\xi,\qquad
    \Dd^0\pi_u\zeta = \pi_uD\pi_u\xi,\qquad
    \Dd^\eps\pi_u\zeta = (D\pi_u\xi,0,0),
$$
the required estimates for the operator $\Dd^\eps$
follow from~(\ref{eq:Dpi}).  The proof for the 
adjoint operator is analogous. 
\end{proof}

In the following we use the notation 
$$
     \Tabla{A,v}\zeta 
     := (\Tabla{A,v}\xi,\Nabla{A,v}\alpha),\qquad
     \Tabla{A,v}\xi'
     := (\Tabla{A,v}\xi',\Nabla{A,v}\phi',\Nabla{A,v}\psi')
$$
for $v\in\Vect(\Sigma)$, $\zeta=(\xi,\alpha)\in\Xx_u$,
and $\zeta'=(\xi',\phi',\psi')\in\Xx'_u$,
where $\Tabla{A}$ is the connection
on $u^*TM/\G$ defined by~(\ref{eq:Tabla}) and
$\Nabla{A,v}\alpha\in\Om^1(\Sigma,\g_P)$
is the covariant derivative induced by the 
connection $A$ and the Levi-Civita connection 
on $\Sigma$.

\begin{lemma}\label{le:com}
For every $p\ge2$, every $c_0>0$, and every 
equivariant vector field $v\in\Vect(\Sigma)$
there exists a constant $c>0$ such that
\begin{eqnarray}
     \|
     \Dd^\eps\Tabla{A,v}\zeta - \Tabla{A,v}\Dd^\eps\zeta
     \|_{k,p,\eps}
&\le & 
     c\eps^{-1}\|\zeta\|_{k+1,p,\eps}, 
     \label{eq:com1} \\
     \|
     {\Dd^\eps}^*\Tabla{A,v}\zeta' - \Tabla{A,v}{\Dd^\eps}^*\zeta'
     \|_{k,p,\eps}
&\le & 
     c\eps^{-1}\|\zeta'\|_{k+1,p,\eps}
     \label{eq:com2}
\end{eqnarray}
for all $(u,A)\in\Tilde\Mcal_{B,\Sigma}^0(c_0)$,
$\zeta\in\Xx_u$, $\zeta'\in\Xx_u$, $\eps\in(0,1]$,
and $k=0,1$.
\end{lemma}

\begin{proof}
We compute in local coordinates.  
Let $\zeta'=(\xi',\phi',\psi'):=\Dd^\eps\zeta$.
Then 
\begin{eqnarray}\label{eq:local}
     \xi'
&= &
     \Tabla{A,s}\xi + J\Tabla{A,t}\xi
     + \frac14N(\xi,v_s-Jv_t)
     + \frac12(J\p_sJ-\p_tJ)\xi
     + L_{u}\phi + JL_{u}\psi, 
     \nonumber\\
     \phi'
&= &   
     \lambda^{-2}\left(\Nabla{A,s}\phi + \Nabla{A,t}\psi\right)
     + \eps^{-2}L_{u}^*\xi, \\
     \psi'
&= &
     \lambda^{-2}\left(\Nabla{A,s}\psi - \Nabla{A,t}\phi\right) 
     + \eps^{-2}d\mu(u)\xi.
     \nonumber
\end{eqnarray}
Here 
$
     A=\Phi ds+\Psi dt,
$
$
    \alpha=\phi ds+\psi dt,
$
and $v_s$, $v_t$, $\Nabla{A,s}\phi$, $\Tabla{A,s}\xi$,
$\Nabla{A,t}\phi$, and $\Tabla{A,t}\xi$
are as in Remark~\ref{rmk:local}.
It suffices to prove the estimate
for the local operators $\Tabla{A,s}$ and $\Tabla{A,t}$.  Let
$\zeta'_s=(\xi'_s,\phi'_s,\psi'_s)$
be defined by~(\ref{eq:local}) with $(\xi,\phi,\psi)$
replaced by $(\Tabla{A,s}\xi,\Nabla{A,s}\phi,\Nabla{A,s}\psi)$.
Since $\Tabla{A}J=0$ we obtain
\begin{eqnarray*}
     \Tabla{A,s}\xi' - \xi'_s
&= &
     J(\Tabla{A,s}\Tabla{A,t}\xi - \Tabla{A,t}\Tabla{A,s}\xi) 
     + \Nabla{v_s}X_\phi(u) + J\Nabla{v_s}X_\psi(u) \\
&&
     +\,\frac14\Tabla{A,s}(N(\xi,v_s-Jv_t))
     - \frac14N(\Tabla{A,s}\xi,v_s-Jv_t) \\
&&
     +\,\frac12\Tabla{A,s}((J\p_sJ-\p_tJ)\xi)
     - \frac12(J\p_sJ-\p_tJ)\Tabla{A,s}\xi \\
&&
     -\, \frac12J(\Nabla{v_s}J+\p_sJ)(L_u\phi-JL_u\psi),  \\
     \Nabla{A,s}\phi' - \phi'_s
&= &
     \lambda^{-2}\left(
     \Nabla{A,s}\Nabla{A,t}\psi - \Nabla{A,t}\Nabla{A,s}\psi 
     \right)
     + (\p_s\lambda^{-2})
     \left(\Nabla{A,s}\phi + \Nabla{A,t}\psi\right) \\
&&
     +\, \eps^{-2}\rho(v_s,\xi)
     - \frac12\eps^{-2}d\mu(u)(\Nabla{v_s}J+\p_sJ)\xi,  \\
     \Nabla{A,s}\psi' - \psi'_s
&= &
     -\,\lambda^{-2}\left(
     \Nabla{A,s}\Nabla{A,t}\phi + \Nabla{A,t}\Nabla{A,s}\phi
     \right) 
     + (\p_s\lambda^{-2})
     \left(\Nabla{A,s}\psi - \Nabla{A,t}\phi\right) \\
&&
     -\, \eps^{-2}\rho(v_s,J\xi)
     - \frac12\eps^{-2}L_u^*(\Nabla{v_s}J+\p_sJ)\xi.
\end{eqnarray*}
Here we have used Lemma~\ref{le:rho}.
For the vector field $\p/\p s$, multiplied by any 
cutoff function, the estimates~(\ref{eq:com1}) 
and~(\ref{eq:com2}) follow from these three identities.  
The proof for $\p/\p t$ is similar, and so is the 
proof for the adjoint operator.
\end{proof}

\begin{lemma}\label{le:onto}
Let $p\ge2$ and $c_0>0$.
Suppose that $\Dd^0:=\Dcal_{(u,A)}^0$ is onto
for every $(u,A)\in\Tilde\Mcal_{B,\Sigma}^0(c_0)$.
Then there exist positive constants $\eps_0$ and $c$ 
such that the operator $\Dd^\eps:=\Dcal_{(u,A)}^\eps$ is onto
for every $(u,A)\in\Tilde\Mcal_{B,\Sigma}^0(c_0)$
and every $\eps\in(0,\eps_0]$ and 
\begin{eqnarray}
     \left\|\zeta'\right\|_{k+1,p,\eps}
&\le &
     c\Bigl(\eps
     \left\|{\Dcal^\eps}^*\zeta'\right\|_{k,p,\eps}
     + \left\|\pi_u{\Dcal^\eps}^*\zeta'\right\|_{k,p,\eps}
     \Bigr),
     \label{eq01} \\
     \left\|\zeta'-\pi_u\zeta'\right\|_{k+1,p,\eps}
&\le &
     c\eps\left\|{\Dcal^\eps}^*\zeta'\right\|_{k,p,\eps},
     \label{eq02}
\end{eqnarray}
for  
$
     \zeta'\in\Xx'_u
$
and $k=0,1$. 
\end{lemma}

\begin{proof}
By elliptic regularity, there exists a constant 
$C_0>0$ such that
\begin{equation}
\label{eqlemma:0}
       \|\xi_0\|_{L^p} \le C_0 \|{\Dd^0}^*\xi_0\|_{L^p}
\end{equation}
for every $(u,A)\in\Tilde\Mcal_{B,\Sigma}^0(c_0)$
and every $\xi_0'\in\Om^{0,1}(\Sigma,H_u)$.
Hence
\begin{eqnarray*}
     \|\pi_u\zeta'\|_{L^p}
&\le &
     C_0\|{\Dd^0}^*\pi_u\zeta'\|_{L^p}  \\
&\le &
     C_0\|\pi_u{\Dd^\eps}^*\zeta'\|_{L^p} 
     + C_0\|\pi_u{\Dd^\eps}^*\zeta'
       - {\Dd^0}^*\pi_u\zeta'\|_{L^p} \\
&\le &
     C_0\|\pi_u{\Dd^\eps}^*\zeta'\|_{L^p} 
     + C_0c_1\|\zeta'-\pi_u\zeta'\|_{0,p,\eps}   \\
&\le &
     C_0\|\pi_u{\Dd^\eps}^*\zeta'\|_{L^p} 
     + C_0c_1c_2\eps
     \left(
       \|{\Dd^\eps}^*\zeta'\|_{0,p,\eps}
       + \|\pi_u\zeta'\|_{L^p}
     \right).
\end{eqnarray*}
Here $c_1$ is the constant of Lemma~\ref{le:difference}
and $c_2$ is the constant of Proposition~\ref{prop:estimate}.
With $C_0c_1c_2\eps\le1/2$ we obtain 
\begin{equation}\label{eq05}
     \|\pi_u\zeta'\|_{L^p} 
     \le c_3
     \bigl(
       \eps\|{\Dd^\eps}^*\zeta'\|_{0,p,\eps}
       + \|\pi_u{\Dd^\eps}^*\zeta'\|_{L^p}
     \bigr),
\end{equation}
where $c_3:=2C_0c_1c_2$. The inequality~(\ref{eq01}) 
for $k=0$ now follows from~(\ref{eq05}) 
and Proposition~\ref{prop:estimate}.
To prove~(\ref{eq02}) for $k=0$ we use 
Proposition~\ref{prop:estimate} and~(\ref{eq05})
again to obtain
\begin{eqnarray*}
     \|\zeta'-\pi_u\zeta'\|_{1,p,\eps}
&\le &
     c_2\eps
     \left(
     \|{\Dd^\eps}^*\zeta'\|_{0,p,\eps}
     + \|\pi_u\zeta'\|_{L^p}
     \right)  \\
&\le &
     c_2\eps(1+c_3\eps)\|{\Dd^\eps}^*\zeta'\|_{0,p,\eps}
     + c_2c_3\eps\|\pi_u{\Dd^\eps}^*\zeta'\|_{L^p} \\
&\le &
     c_4\eps\|{\Dd^\eps}^*\zeta'\|_{0,p,\eps},
\end{eqnarray*}
where $c_4:=c_2(1+2c_3)$. 

Now let $v\in\Vect(\Sigma)$.  By definition of the 
$(1,p,\eps)$-norm and~(\ref{eq:difference}),
there exists a constant $c_5=c_5(v,p,c_0)$ such that 
$$
     \|\Tabla{A,v}\zeta\|_{0,p,\eps}
     \le c_5\eps^{-1}\|\zeta\|_{1,p,\eps},\qquad
     \|\pi_u\Tabla{A,v}\zeta-\Tabla{A,v}\pi_u\zeta\|_{L^p}
     \le c_5\|\zeta\|_{0,p,\eps},
$$
for $(u,A)\in\Tilde\Mcal_{B,\Sigma}^0(c_0)$,
$\zeta\in\Xx_u$ and $\eps \in (0,1]$.
Let $c_6=c_6(v,p,c_0)$ be the constant of Lemma~\ref{le:com}.
Then, by~(\ref{eq01}) with $k=0$ and Lemma~\ref{le:com}, we have
\begin{eqnarray*}
     \eps\|\Tabla{A,v}\zeta'\|_{1,p,\eps}
&\le &
     c\eps\left(
     \eps\|{\Dd^\eps}^*\Tabla{A,v}\zeta'\|_{0,p,\eps}
     + \|\pi_u{\Dd^\eps}^*\Tabla{A,v}\zeta'\|_{L^p}
     \right)  \\
&\le &
     c\eps^2\|
     {\Dd^\eps}^*\Tabla{A,v}\zeta' - \Tabla{A,v}{\Dd^\eps}^*\zeta'
     \|_{0,p,\eps} 
     + c\eps^2\|\Tabla{A,v}{\Dd^\eps}^*\zeta'\|_{0,p,\eps} \\ 
&&
     +\, c\eps\|
       \pi_u({\Dd^\eps}^*\Tabla{A,v}\zeta' 
       - \Tabla{A,v}{\Dd^\eps}^*\zeta')
       \|_{L^p} \\     
&&
     +\, c\eps\|
       \pi_u\Tabla{A,v}{\Dd^\eps}^*\zeta' 
       - \Tabla{A,v}\pi_u{\Dd^\eps}^*\zeta'
       \|_{L^p}
     + c\eps\|\Tabla{A,v}\pi_u{\Dd^\eps}^*\zeta'\|_{L^p}  \\
&\le &
     2cc_5\left(
     \eps\|{\Dd^\eps}^*\zeta'\|_{1,p,\eps}
     + \|\pi_u{\Dd^\eps}^*\zeta'\|_{1,p,\eps} 
     \right)
     + 2cc_6\|\zeta'\|_{1,p,\eps}  \\
&\le &
     2c(c_5+cc_6)
     \left(
     \eps\|{\Dd^\eps}^*\zeta'\|_{1,p,\eps}
     + \|\pi_u{\Dd^\eps}^*\zeta'\|_{1,p,\eps} 
     \right).
\end{eqnarray*}
The last inequality follows again from~(\ref{eq01})
with $k=0$. The estimate~(\ref{eq01}) for $k=1$ 
now follows by taking the sum over finitely 
many suitably chosen vector fields~$v$. 

To prove~(\ref{eq02}) for $k=1$ we observe that 
$\pi_u{\Dd^\eps}^*\pi_u\zeta'={\Dd^0}^*\pi_u\zeta'$ 
and choose $c_7$ such that 
$
     \|\pi_u\zeta\|_{1,p,\eps}
     \le c_7\|\zeta\|_{1,p,\eps}
$
for every $\zeta\in\Xx_u$.  
Let $c_8$ be the constant of Lemma~\ref{le:difference}.  
Then, by~(\ref{eq01}) with $k=1$ 
and Lemma~\ref{le:difference}, 
we have 
\begin{eqnarray*}
    \|\zeta'-\pi_u\zeta'\|_{2,p,\eps}
&\le &
    c\left(\eps\|{\Dd^\eps}^*(\zeta'-\pi_u\zeta')\|_{1,p,\eps}
    + \|\pi_u{\Dd^\eps}^*(\zeta'-\pi_u\zeta')\|_{1,p,\eps}
    \right)  \\
&\le &
    c\eps\left(
    \|{\Dd^\eps}^*\zeta'\|_{1,p,\eps}
    + \|\pi_u{\Dd^\eps}^*\zeta'\|_{1,p,\eps}
    \right) \\
&&
    +\,c\eps \|
    {\Dd^\eps}^*\pi_u\zeta'-\pi_u{\Dd^\eps}^*\zeta'
    \|_{1,p,\eps} \\
&&
    +\, c\|
    \pi_u{\Dd^\eps}^*\zeta'-{\Dd^0}^*\pi_u\zeta'
    \|_{1,p,\eps} \\
&\le &
    c(1+c_7)\eps\|{\Dd^\eps}^*\zeta'\|_{1,p,\eps}
    + cc_8\eps\|\zeta'\|_{1,p,\eps}
    + cc_8\|\zeta'-\pi_u\zeta'\|_{1,p,\eps}  \\
&\le &
    c_9\eps\|{\Dd^\eps}^*\zeta'\|_{1,p,\eps}.
\end{eqnarray*}
The last inequality follows from~(\ref{eq01})
and~(\ref{eq02}) with $k=0$.
\end{proof}

\begin{lemma}\label{le:onto1}
Let $p\ge2$ and $c_0>0$.
Suppose that $\Dd^0:=\Dd_{(u,A)}^0$ is onto
for every $(u,A)\in\Tilde\Mcal_{B,\Sigma}^0(c_0)$.
Then there exist positive constants
$c$ and $\eps_0$ such that
\begin{eqnarray}
       \left\|{\Dd^\eps}^*\zeta'\right\|_{k+1,p,\eps}
&\le &
       c\Bigl(\eps
       \left\|\Dd^\eps{\Dd^\eps}^*\zeta'\right\|_{k,p,\eps}
       +\left\|\pi_u\Dd^\eps{\Dd^\eps}^*\zeta'\right\|_{k,p,\eps}
       \Bigr),
       \label{eq:onto1} \\
       \left\|{\Dd^\eps}^*\zeta'
         -\pi_u{\Dd^\eps}^*\zeta'\right\|_{k+1,p,\eps}
&\le &
       c\eps
       \left\|\Dd^\eps{\Dd^\eps}^*\zeta'\right\|_{k,p,\eps}
       \label{eq:onto2}
\end{eqnarray}
for every $(u,A)\in\Tilde\Mcal_{B,\Sigma}^0(c_0)$,
$\zeta'\in\Xx_u'$, $\eps\in(0,\eps_0]$, and $k=0,1$. 
\end{lemma}

\begin{proof}
The proof has nine steps. 

\medskip
\noindent{\bf Step~1.}
{\it Let $q>1$ such that $1/p+1/q=1$.
Then there exists a constant $c_0>0$ such that
$$
    \left\|\xi_0\right\|_{L^p}+\left\|\xi_0\right\|_{L^q}
    \le c_0\left\|\xi_0\right\|_{L^2},\qquad
    \left\|\xi_0'\right\|_{L^q}
    \le c_0\left\|{\Dd^0}^*\xi_0'\right\|_{L^q}
$$
for every $(u,A)\in\Tilde\Mcal_{B,\Sigma}^0(c_0)$,
every $\xi_0\in\ker\Dd^0$, 
and every $\xi_0'\in \Om^{0,1}(\Sigma,H_u)$.
}

\medskip
\noindent
These are standard estimates for elliptic pdes.
The first estimate uses $L^2$ regularity
for the operator $\Dd^0$,
the Sobolev embedding $W^{1,2}\INTO L^p$,
and the H\"older inequality. 
The second estimate uses $L^q$ regularity
for ${\Dd^0}^*$ and the fact that 
${\Dd^0}^*$ is injective. 

\medskip
\noindent{\bf Step~2.}
{\it
There exists a constant $c_1>0$ such that
$$
     \left\|{\Dd^0}^*\xi_0'\right\|_{L^p}
     \le c_1 \sup_{\xi_0''\ne 0}
     \frac{\inner{{\Dd^0}^*\xi_0'}{{\Dd^0}^*\xi_0''}}
     {\left\|{\Dd^0}^*\xi_0''\right\|_{L^q}}
$$
for every $(u,A)\in\Tilde\Mcal_{B,\Sigma}^0(c_0)$
and every $\xi_0'\in \Om^{0,1}(\Sigma,H_u)$.
}

\medskip
\noindent
Let $\xi_1,\ldots,\xi_m$ be an
$L^2$-orthonormal basis of $\ker\Dcal^0$.
Given $\xi_0'$ choose $\xi\in L^q(\Sigma,H_u)$ such that
$$
     \inner{\xi}{{\Dcal^0}^*\xi_0'}
     = \left\|{\Dcal^0}^*\xi_0'\right\|_{L^p},\qquad
     \left\|\xi\right\|_{L^q} = 1.
$$
Let $\xi_0''\in W^{1,q}(\Sigma,\Lambda^{0,1}T^*\Sigma\otimes H_u)$
be the unique section such that 
$$
     \xi = {\Dcal^0}^*\xi_0'' 
         + \sum_{j=1}^m \inner{\xi}{\xi_j}\xi_j.
$$
Then
\begin{eqnarray*}
     \left\|{\Dcal^0}^*\xi_0'\right\|_{L^p}
&=&
     \inner{{\Dcal^0}^*\xi_0''}{{\Dcal^0}^*\xi_0'} \\
&=&
     \left\|\xi - \sum_{j=1}^m \inner{\xi}{\xi_j}\xi_j\right\|_{L^q}
     \frac{\inner{{\Dcal^0}^*\xi_0''}{{\Dcal^0}^*\xi_0'}}
          {\left\|{\Dcal^0}^*\xi_0''\right\|_{L^q}}  \\
&\le &
     \left(
     1 + \sum_{j=1}^m\left\|\xi_j\right\|_{L^p}
                     \left\|\xi_j\right\|_{L^q}
     \right)
     \frac{\inner{{\Dcal^0}^*\xi_0''}{{\Dcal^0}^*\xi_0'}}
          {\left\|{\Dcal^0}^*\xi_0''\right\|_{L^q}}  \\
&\le &
     (1 + m{c_0}^2)
     \frac{\inner{{\Dcal^0}^*\xi_0''}{{\Dcal^0}^*\xi_0'}}
          {\left\|{\Dcal^0}^*\xi_0''\right\|_{L^q}}.
\end{eqnarray*}

\medskip
\noindent{\bf Step~3.}
{\it
There exists a constant $c_2>0$ such that
$$
     \left\|\pi_u{\Dd^\eps}^*\zeta'\right\|_{L^p}
\le 
     c_2 \biggl(
     \left\|
       \pi_u{\Dd^\eps}^*\zeta'
       - {\Dd^0}^*\pi_u\zeta'
     \right\|_{L^p} 
     + \left\|
       \Dd^0\pi_u{\Dd^\eps}^*\zeta'
     \right\|_{L^p}
     \biggr)
$$
for every $(u,A)\in\Tilde\Mcal_{B,\Sigma}^0(c_0)$,
every
$
    \zeta'\in \Xx'_u, 
$
and every $\eps\in(0,1]$.
}

\medskip
\noindent
For every
$
    \xi_0'\in \Om^{0,1}(\Sigma,H_u)
$ 
we have
\begin{eqnarray*}
     \frac{\inner{{\Dcal^0}^*\xi_0'}{{\Dcal^0}^*\pi_u\zeta'}}
          {\left\|{\Dcal^0}^*\xi_0'\right\|_{L^q}} 
&= &
     \frac{\inner{{\Dcal^0}^*\xi_0'}{{\Dcal^0}^*\pi_u\zeta'
           - \pi_u{\Dcal^\eps}^*\zeta'}}
          {\left\|{\Dcal^0}^*\xi_0'\right\|_{L^q}}
     +
     \frac{\inner{\xi_0'}{\Dcal^0\pi_u{\Dcal^\eps}^*\zeta'}}
          {\left\|{\Dcal^0}^*\xi_0'\right\|_{L^q}}  \\
&\le &
     \left\|{\Dcal^0}^*\pi_u\zeta'-\pi_u{\Dcal^\eps}^*\zeta'
     \right\|_{L^p}   
     +
     \left\|\Dcal^0\pi_u{\Dcal^\eps}^*\zeta'\right\|_{L^p}
     \frac{\left\|\xi_0'\right\|_{L^q}}
          {\left\|{\Dcal^0}^*\xi_0'\right\|_{L^q}} \\
&\le &
     \left\|
     {\Dcal^0}^*\pi_u\zeta'-\pi_u{\Dcal^\eps}^*\zeta'
     \right\|_{L^p}
     + c_0\left\|\Dcal^0\pi_u{\Dcal^\eps}^*\zeta'
          \right\|_{L^p}.
\end{eqnarray*}
Here the last inequality follows from Step~1. 
Now, by Step~2, 
\begin{eqnarray*}
     \left\|\pi_u{\Dd^\eps}^*\zeta'\right\|_{L^p}
&\le &
     \left\|
       {\Dcal^0}^*\pi_u\zeta'-\pi_u{\Dcal^\eps}^*\zeta'
     \right\|_{L^p}
     + 
     \left\|{\Dcal^0}^*\pi_u\zeta'\right\|_{L^p}  \\
&\le &
     \left\|
       {\Dcal^0}^*\pi_u\zeta'-\pi_u{\Dcal^\eps}^*\zeta'
     \right\|_{L^p}
     + 
     c_1 \sup_{\xi_0'\ne 0}
     \frac{\inner{{\Dd^0}^*\pi_u\zeta'}{{\Dd^0}^*\xi_0'}}
          {\left\|{\Dd^0}^*\xi_0'\right\|_{L^q}}  \\
&\le &
     (1+c_1)\left\|
       {\Dcal^0}^*\pi_u\zeta'-\pi_u{\Dcal^\eps}^*\zeta'
     \right\|_{L^p}
     + c_0c_1\left\|
        \Dcal^0\pi_u{\Dcal^\eps}^*\zeta'
       \right\|_{L^p}.
\end{eqnarray*}

\medskip
\noindent{\bf Step~4.}
{\it
There exist positive constants $\eps_0$ and $c_3$ such that 
$$
     \left\|\pi_u{\Dd^\eps}^*\zeta'\right\|_{L^p}
\le 
     c_3\Bigl(
     \eps\left\|{\Dd^\eps}^*\zeta'\right\|_{0,p,\eps} 
     + \eps\left\|\Dd^\eps{\Dd^\eps}^*\zeta'\right\|_{0,p,\eps} 
     + \left\|\pi_u\Dd^\eps{\Dd^\eps}^*\zeta'\right\|_{L^p} 
     \Bigr)
$$
for every $(u,A)\in\Tilde\Mcal_{B,\Sigma}^0(c_0)$,
every
$
    \zeta'\in \Xx'_u, 
$
and every $\eps\in(0,\eps_0]$.
}

\medskip
\noindent
We apply Lemma~\ref{le:difference}
to both operators $\Dd^\eps$ and ${\Dd^\eps}^*$.
Then, by Step~3, 
\begin{eqnarray*}
     \left\|\pi_u{\Dd^\eps}^*\zeta'\right\|_{L^p}
&\le &
     c_2 \Bigl(
     \left\|
       \pi_u{\Dd^\eps}^*\zeta'
       - {\Dd^0}^*\pi_u\zeta'
     \right\|_{L^p}  \\
&&\quad
     +\, \left\|
       \Dd^0\pi_u{\Dd^\eps}^*\zeta' - \pi_u\Dd^\eps{\Dd^\eps}^*\zeta'
     \right\|_{0,p,\eps}
     +\,
     \left\|
       \pi_u\Dd^\eps{\Dd^\eps}^*\zeta'
     \right\|_{L^p}
     \Bigr) \\
&\le &
     c_2\Bigl(
     c\left\|
       \zeta'-\pi_u\zeta'
     \right\|_{0,p,\eps} 
     + c\left\|
       {\Dd^\eps}^*\zeta' - \pi_u{\Dd^\eps}^*\zeta'
     \right\|_{0,p,\eps}  \\
&&\quad
     +\,
     \left\|
       \pi_u\Dd^\eps{\Dd^\eps}^*\zeta'
     \right\|_{L^p}
     \Bigr) \\
&\le &
     c_4\Bigl(
     \eps
     \left\|{\Dd^\eps}^*\zeta'\right\|_{0,p,\eps} 
     + \eps\left\|\Dd^\eps{\Dd^\eps}^*\zeta'\right\|_{0,p,\eps}
     + \left\|\pi_u\Dd^\eps{\Dd^\eps}^*\zeta'\right\|_{L^p}
     \Bigr).
\end{eqnarray*}
The last inequality follows from Lemma~\ref{le:onto}
and Proposition~\ref{prop:estimate}. 

\medskip
\noindent{\bf Step~5.}
{\it
We prove~(\ref{eq:onto1}) for $k=0$.
}

\medskip
\noindent
By Proposition~\ref{prop:estimate} and Step~4, 
\begin{eqnarray*}
     \left\|{\Dd^\eps}^*\zeta'\right\|_{1,p,\eps}
&\le &
     c\Bigl(
     \eps\left\|\Dd^\eps{\Dd^\eps}^*\zeta'\right\|_{1,p,\eps}
     + \left\|\pi_u{\Dd^\eps}^*\zeta'\right\|_{L^p}  
     \Bigr) \\
&\le &
     c\Bigl(
     c_3\eps\left\|{\Dd^\eps}^*\zeta'\right\|_{0,p,\eps} 
     + (1+c_3)\eps\left\|\Dd^\eps{\Dd^\eps}^*\zeta'\right\|_{0,p,\eps} \\
&&\quad
     +\, c_3\left\|\pi_u\Dd^\eps{\Dd^\eps}^*\zeta'\right\|_{L^p} 
     \Bigr)
\end{eqnarray*}
for all $(u,A)\in\Tilde\Mcal_{B,\Sigma}^0(c_0)$
and 
$
     \zeta'\in\Xx_u'
$
and $\eps\in(0,\eps_0]$.
With $cc_3\eps\le1/2$ we obtain~(\ref{eq:onto1}) for $k=0$.

\medskip
\noindent{\bf Step~6.}
{\it
We prove~(\ref{eq:onto2}) for $k=0$.
}

\medskip
\noindent
By Proposition~\ref{prop:estimate} and Step~4,
\begin{eqnarray*}
     \left\|{\Dd^\eps}^*\zeta'-\pi_u{\Dd^\eps}^*\zeta'\right\|_{1,p,\eps}
&\le &
     c\eps\Bigl(
     \left\|\Dd^\eps{\Dd^\eps}^*\zeta'\right\|_{1,p,\eps}
     + \left\|\pi_u{\Dd^\eps}^*\zeta'\right\|_{L^p}
     \Bigr)  \\
&\le &
     c\eps\Bigl(
     c_3\eps\left\|{\Dd^\eps}^*\zeta'\right\|_{0,p,\eps} 
     + (1+c_3)\eps\left\|\Dd^\eps{\Dd^\eps}^*\zeta'\right\|_{0,p,\eps} \\
&&\quad
     +\, c_3\left\|\pi_u\Dd^\eps{\Dd^\eps}^*\zeta'\right\|_{L^p} 
     \Bigr)  \\
&\le &
     c_5\eps\left\|\Dd^\eps{\Dd^\eps}^*\zeta'\right\|_{0,p,\eps}.
\end{eqnarray*}
Here the last inequality follows from Step~5.

\medskip
\noindent{\bf Step~7.}
{\it
There exist positive constants $\eps_0$ and $c_6$ such that 
\begin{eqnarray*}
     \eps^2\|
      [\Dd^\eps{\Dd^\eps}^*,\Tabla{A,v}]\zeta' 
     \|_{L^p} 
&\le &
     c_6\|\zeta'\|_{2,p,\eps}, \\
     \eps\|
      [\pi_u\Dd^\eps{\Dd^\eps}^*,\Tabla{A,v}]\zeta' 
     \|_{L^p} 
&\le &
     c_6\left(
     \eps\|\Dd^\eps{\Dd^\eps}^*\zeta'\|_{0,p,\eps}
     + \|\pi_u\Dd^\eps{\Dd^\eps}^*\zeta'\|_{L^p}
     + \|\zeta'\|_{1,p,\eps}
     \right)
\end{eqnarray*}
for every $(u,A)\in\Tilde\Mcal_{B,\Sigma}^0(c_0)$,
every
$
    \zeta'\in\Xx_u', 
$
and every $\eps\in(0,\eps_0]$.
}

\medskip
\noindent
The first estimate follows immediately 
from Lemma~\ref{le:com}.  To prove the second 
estimate, recall from Lemma~\ref{le:adjoint} that
$$
     \pi_u\Dd^\eps{\Dd^\eps}^*\zeta' 
     = \pi_u2DD^*\xi' + \pi_u(DJ-JD)L_u\psi'
$$
where 
$
     D:\Om^0(\Sigma,u^*TM/\G)\to\Om^{0,1}(\Sigma,u^*TM/\G)
$
is the Cauchy--Riemann operator defined by~(\ref{eq:Du})
and $R:=(DJ-JD)L_u$ is a zeroth order operator
(Lemma~\ref{le:D}).  Hence 
$$
     2DD^*\pi_u\xi'
     = - 2[\pi_u,DD^*]\xi' 
       + \pi_u\Dd^\eps{\Dd^\eps}^*\zeta' 
       - \pi_uR\psi'.
$$
By~(\ref{eq:Dpi}) in the proof of Lemma~\ref{le:difference}, 
the commutator $[\pi_u,DD^*]$ is a first order operator 
in $\xi'$. Hence there exists a constant
$c_7=c_7(p,c_0)$ such that
$$
    \|{\Tabla{A}}^*\Tabla{A}\pi_u\xi'\|_{L^p}
    \le c_7\left(
        \|\pi_u\Dd^\eps{\Dd^\eps}^*\zeta'\|_{L^p} 
        + \eps^{-1}\|\zeta'\|_{1,p,\eps}
        \right)
$$
Moreover, by Lemma~\ref{le:onto}, 
\begin{eqnarray*}
     \|{\Tabla{A}}^*\Tabla{A}(\xi'-\pi_u\xi')\|_{L^p}
&\le &
     \eps^{-2}\|\xi'-\pi_u\xi'\|_{2,p,\eps}  \\
&\le &
     c_8\eps^{-1}\|{\Dd^\eps}^*\zeta'\|_{1,p,\eps}  \\
&\le &
     c_9\left(
     \|\Dd^\eps{\Dd^\eps}^*\zeta'\|_{0,p,\eps}
     + \eps^{-1}\|\pi_u\Dd^\eps{\Dd^\eps}^*\zeta'\|_{0,p,\eps}
     \right).
\end{eqnarray*}
The last inequality follows from Step~5. 
Now the commutator
$$
     [\Tabla{A,v},\pi_u\Dd^\eps{\Dd^\eps}^*]\zeta'
     = 2[\Tabla{A,v},\pi_uDD^*]\xi' 
       + [\Tabla{A,v},\pi_uR]\psi'
$$
is a second order operator in $\xi'$ 
and a zeroth order operator in $\psi'$.
Hence the assertion follows from the last two 
inequalities.

\medskip
\noindent{\bf Step~8.}
{\it
We prove~(\ref{eq:onto1}) for $k=1$.
}

\medskip
\noindent
Let $c_{10}$ be the constant in~(\ref{eq:onto1}) for $k=0$
and $c_{11}$ be the constant of Lemma~\ref{le:com}  
Then, for every $v\in\Vect(\Sigma)$, we have 
\begin{eqnarray*}
&&
     \eps\|\Tabla{A,v}{\Dd^\eps}^*\zeta'\|_{1,p,\eps}  \\
&&\le 
     \eps\|{\Dd^\eps}^*\Tabla{A,v}\zeta'\|_{1,p,\eps}
     + \eps\|
       \Tabla{A,v}{\Dd^\eps}^*\zeta'-{\Dd^\eps}^*\Tabla{A,v}\zeta'
       \|_{1,p,\eps} \\
&&\le 
     c_{10}\eps^2\|\Dd^\eps{\Dd^\eps}^*\Tabla{A,v}\zeta'\|_{0,p,\eps}
     + c_{10}\eps\|\pi_u\Dd^\eps{\Dd^\eps}^*\Tabla{A,v}\zeta'
       \|_{L^p}  
     + c_{11}\|\zeta'\|_{2,p,\eps} \\
&&\le 
     c_{10}\eps^2\|
     \Tabla{A,v}\Dd^\eps{\Dd^\eps}^*\zeta'
     \|_{0,p,\eps}
     + c_{10}\eps\|
       \Tabla{A,v}\pi_u\Dd^\eps{\Dd^\eps}^*\zeta'
       \|_{L^p}  \\
&&\quad
     +\, c_{10}\eps^2\|
         [\Dd^\eps{\Dd^\eps}^*,\Tabla{A,v}]\zeta'
         \|_{0,p,\eps} 
     + c_{10}\eps\|
         [\pi_u\Dd^\eps{\Dd^\eps}^*,\Tabla{A,v}]\zeta'
         \|_{L^p} 
     + c_{11}\|\zeta'\|_{2,p,\eps} \\
&&\le 
     c_{10}\eps^2\|
     \Tabla{A,v}\Dd^\eps{\Dd^\eps}^*\zeta'
     \|_{0,p,\eps}
     + c_{10}\eps\|
       \Tabla{A,v}\pi_u\Dd^\eps{\Dd^\eps}^*\zeta'
       \|_{L^p}  \\
&&\quad
     +\,c_6c_{10}\left(
     \eps\|\Dd^\eps{\Dd^\eps}^*\zeta'\|_{0,p,\eps} 
     + \|\pi_u\Dd^\eps{\Dd^\eps}^*\zeta'\|_{L^p} 
     \right)
     + (2c_6c_{10}+c_{11})\|\zeta'\|_{2,p,\eps}  \\
&&\le 
     c_{12}
     \left(
     \eps\|\Dd^\eps{\Dd^\eps}^*\zeta'\|_{1,p,\eps}
     + \|\pi_u\Dd^\eps{\Dd^\eps}^*\zeta'\|_{1,p,\eps}
     \right).
\end{eqnarray*}
The penultimate inequality follows from Step~7,
and the last step from~(\ref{eq01}) and~(\ref{eq:onto1})
and the definition of the $(1,p,\eps)$-norm.
Now~(\ref{eq:onto1}) for $k=1$ follows by taking the sum over 
finitely many vector fields $v\in\Vect(\Sigma)$.

\medskip
\noindent{\bf Step~9.}
{\it
We prove~(\ref{eq:onto2}) for $k=1$.
}

\medskip
\noindent
By Step~8, suppose that~(\ref{eq:onto1})
holds with $k=1$ and $c=c_{13}$,
choose $c_{14}$ such that 
$
     \|\pi_u\zeta\|_{1,p,\eps}
     \le c_{14}\|\zeta\|_{1,p,\eps}
$
for every $\zeta\in\Xx_u$, and let $c_{15}$
be the constant of Lemma~\ref{le:difference}.
Then 
\begin{eqnarray*}
    \|{\Dd^\eps}^*\zeta'-\pi_u{\Dd^\eps}^*\zeta'\|_{2,p,\eps}
&\le &
    c_{13}\eps\|
    \Dd^\eps({\Dd^\eps}^*\zeta'-\pi_u{\Dd^\eps}^*\zeta')
    \|_{1,p,\eps} \\
&&
    + \,c_{13}\|
    \pi_u\Dd^\eps({\Dd^\eps}^*\zeta'-\pi_u{\Dd^\eps}^*\zeta')
    \|_{1,p,\eps} \\
&\le &
    c_{13}\eps\left(
    \|\Dd^\eps{\Dd^\eps}^*\zeta'\|_{1,p,\eps}
    + \|\pi_u\Dd^\eps{\Dd^\eps}^*\zeta'\|_{1,p,\eps}
    \right) \\
&&
    +\,c_{13}\eps\|\Dd^\eps\pi_u{\Dd^\eps}^*\zeta'
      - \pi_u\Dd^\eps{\Dd^\eps}^*\zeta'\|_{1,p,\eps} \\
&&
    +\,c_{13}\|
    (\pi_u\Dd^\eps-\Dd^0\pi_u){\Dd^\eps}^*\zeta'
    \|_{1,p,\eps} \\
&\le &
    (1+c_{14})c_{13}\eps
    \|\Dd^\eps{\Dd^\eps}^*\zeta'\|_{1,p,\eps}  \\
&&
    +\,c_{15}c_{13}\left(
    \eps\|{\Dd^\eps}^*\zeta'\|_{1,p,\eps} 
    + \|{\Dd^\eps}^*\zeta'-\pi_u{\Dd^\eps}^*\zeta'
    \|_{1,p,\eps}\right) \\
&\le &
    c_{16}\eps\|\Dd^\eps{\Dd^\eps}^*\zeta'\|_{1,p,\eps}.
\end{eqnarray*}
The last inequality follows from Steps~5
and~6. 
\end{proof}


\section{Quadratic estimates}\label{sec:quad}

Fix $p>2$, $c_0>0$ and
$(u,A)\in\Tilde\Mcal_{B,\Sigma}^0(c_0)$,
and consider the map
$$
   \Ff^\eps = \Ff_{(u,A)}^\eps:\Xx_u\to\Xx'_u
$$
given by
\begin{equation}\label{eq:Feps}
   \Fcal^\eps(\xi,\alpha)
   =
   \begin{pmatrix}
   \rho(\xi)
   \bigl(\bar\partial_{J,A+\alpha}(\exp_u(\xi))\bigr) \\
   \eps^{-2}L_u^*\xi-d_A^*\alpha \\
   \eps^{-2}\mu(\exp_u(\xi))+*F_{A+\alpha}
\end{pmatrix}.
\end{equation}
Here $\rho(\xi):T_{\exp_u(\xi)}M\to T_uM$
denotes parallel transport along the geodesic
$r\mapsto\exp_u(r\xi)$ with respect to the 
Hermitian connection $\tilde\nabla:=\nabla-\frac12J\nabla J$ 
on $TM$. The differential of $\Fcal^\eps$ 
at zero is given by
$$
      d\Fcal^\eps(0)=\Dcal_{(u,A)}^\eps.
$$
Let $\Fcal^\eps_i$ denote the $i$th component of $\Ff^\eps$. 
Since $\Ff^\eps_2$ is a linear map,
the following proposition only deals with the
first and third components of $\Ff^\eps$.

\begin{proposition}\label{prop:quadra}
For every $p>2$ and every $c_0>0$
there exists a cons\-tant $c>0$ 
such that the following holds
for every $(u,A)\in\Tilde\Mcal_{B,\Sigma}^0(c_0)$,
any two pairs
$
    \zeta=(\xi,\alpha),\hat\zeta=(\hat\xi,\hat\alpha) \in \Xx_u,
$
and every $\eps\in(0,1]$.

\smallskip
\noindent{\bf (i)}
If $\|\xi\|_{L^\infty}+\|\hat\xi\|_{L^\infty}\le 1$ then
\begin{eqnarray*}
&&
     \|
     \Ff^\eps_1(\zeta+\hat\zeta) -\Fcal^\eps_1(\zeta)
     - d\Fcal^\eps_1(\zeta)\hat\zeta
     \|_{L^p} \\
&&
     \le c
     \|\hat\xi\|_{L^\infty}
     \Bigl(
     \|\hat\xi\|_{L^p}
     + \|\Tabla{A}\hat\xi\|_{L^p}
     + \|\hat\alpha\|_{L^p} 
     \Bigr) \\
&&\quad
     +\,
     c\|\hat\xi\|_{L^\infty}^2
     \Bigl(
     \|\Tabla{A}\xi\|_{L^p}
      + \|\alpha\|_{L^p}
      + \|\Tabla{A}\hat\xi\|_{L^p}
      + \|\hat\alpha\|_{L^p}
     \Bigr).
\end{eqnarray*}
If, in addition, 
$
     \|\Tabla A\xi\|_{L^p}+\|\Tabla{A}\hat\xi\|_{L^p}
     +\|\alpha\|_{L^p}+\|\hat\alpha\|_{L^p} 
     \le \eps^{2/p-1}
$
then
\begin{eqnarray*}
&&
     \|
     \Ff^\eps_1(\zeta+\hat\zeta) -\Fcal^\eps_1(\zeta)
     - d\Fcal^\eps_1(\zeta)\hat\zeta
     \|_{L^p} \\
&&
     \le c \|\hat\xi\|_{L^\infty}
     \Bigl(
     \eps^{-1}\|\hat\xi\|_{L^p}
     + \|\Tabla A\hat\xi\|_{L^p}
     + \|\hat\alpha\|_{L^p}
     \Bigr).
\end{eqnarray*}

\smallskip
\noindent{\bf (ii)}
If $\|\xi\|_{L^\infty}+\|\hat\xi\|_{L^\infty}\le 1$
then
$$
     \|
     \Ff^\eps_3(\zeta+\hat\zeta) - \Ff^\eps_3(\zeta) 
     - d\Ff^\eps_3(\zeta)\hat\zeta
     \|_{L^p}
     \le c\Bigl(
     \|\hat\alpha\|_{L^\infty}\|\hat\alpha\|_{L^p}
     + \eps^{-2}\|\hat\xi\|_{L^\infty}\|\hat\xi\|_{L^p}
     \Bigr).
$$

\smallskip
\noindent{\bf (iii)}
If $\|\xi\|_{L^\infty}\le1$ then 
\begin{eqnarray*}
     \|d\Ff^\eps_1(\zeta)\hat\zeta - d\Ff^\eps_1(0)\hat\zeta\|_{L^p}  
&\le &
     c \|\hat\xi\|_{L^\infty}
     \Bigl(
     \|\xi\|_{L^p} 
     + \|\Tabla{A}\xi\|_{L^p}
     + \|\alpha\|_{L^p}
     \Bigr)   \\
&&
     +\,c \|\xi\|_{L^\infty}
     \Bigl(
     \|\hat\xi\|_{L^p} 
     + \|\Tabla{A}\hat\xi\|_{L^p}
     + \|\hat\alpha\|_{L^p}
     \Bigr).
\end{eqnarray*}

\smallskip
\noindent{\bf (iv)}
If $\|\xi\|_{L^\infty} \le 1$ then
$$
     \|d\Ff^\eps_3(\zeta)\hat\zeta - d\Ff^\eps_3(0)\hat\zeta\|_{L^p}  
     \le c \Bigl(
     \eps^{-2}\|\xi\|_{L^p}\|\hat\xi\|_{L^\infty}
     + \|\alpha\|_{L^p}\|\hat\alpha\|_{L^\infty}
     \Bigr).
$$
\end{proposition}

The estimates in Proposition~\ref{prop:quadra}
differ from the ones in~\cite{DS2} in that 
the first derivatives of $\xi$ appear on the
right hand sides of the inequalities.
This is because the nonlinearities in the 
Cauchy-Riemann equations appear in the first
order terms whereas the nonlinearities 
in the anti-self-duality equations only 
appear in the zeroth order terms. 
In our equations the nonlinear terms 
involving $\alpha$ are of zeroth order. 
Hence no derivatives of $\alpha$ appear 
in the quadratic estimates. This is crucial 
for our adiabatic limit argument.

\begin{proof}[Proof of Proposition~\ref{prop:quadra}]
In local holomorphic coordinates $s+it$ on $\Sigma$ 
the map~$\Fcal^\eps$ is given by
$$
     \Fcal^\eps(\zeta)
     = \begin{pmatrix}
          \rho(\xi)
          \Bigl(
          \partial_s u_\xi+X_{\Phi+\phi}(u_\xi)
          +J\bigl(
          \partial_t u_\xi+X_{\Psi+\psi}(u_\xi)
          \bigr)
          \Bigr)
          \\
          \eps^{-2}L_u^*\xi
          + \lambda^{-2}
          \left(\Nabla{A,s}\phi + \Nabla{A,t}\psi\right) 
          \\
          \eps^{-2}\mu(u_\xi)
          + \lambda^{-2}\left(\p_s(\Psi+\psi)-\p_t(\Phi+\phi)
            + [\Phi+\phi, \Psi+\psi]
            \right)
       \end{pmatrix},
$$
where $u_\xi:=\exp_u(\xi)$,
$\zeta:=(\xi,\phi,\psi)$
and $\alpha:=\phi\,ds+\psi\,dt$.
Suppose that 
$
     \left\|\xi\right\|_{L^\infty}\le 1.
$
The second derivatives of $\Ff^\eps_1$ and
$\Ff^\eps_3$ satisfy the following {\em pointwise} estimates
for suitable constants $c_1=c_1(u,A,v_s,v_t)$
and $c_2=c_2(u)$ (i.e. $c_2$ does not depend on 
the derivatives of $u$):
\begin{eqnarray}
      |d^2\Fcal^\eps_1(\zeta)(\zeta_1,\zeta_2)|
&\le &
      c_1
      \Bigl(
      \bigl(1+\left|\alpha\right|+|\Tabla{A}\xi|\bigr)
      \left|\xi_1\right|\left|\xi_2\right|   
      \nonumber \\
&&
      +\, \left|\xi_1\right|
      \bigl(|\Tabla{A}\xi_2| + \left|\alpha_2\right|\bigr)
      + \left|\xi_2\right|
      \bigl(|\Tabla{A}\xi_1| + \left|\alpha_1\right|\bigr)
      \Bigr)
      \label{eq:d2F1}  \\
      |d^2\Fcal^\eps_3(\zeta)(\zeta_1,\zeta_2)|
&\le &
      c_2\bigl(
      |\alpha_1||\alpha_2|
      +\eps^{-2}|\xi_1||\xi_2|
      \bigr).
      \label{eq:d2F3}
\end{eqnarray}
The estimate~(\ref{eq:d2F3}) is obvious 
and~(\ref{eq:d2F1}) will be proved below.
Now consider the identities
\begin{eqnarray}
     \Fcal^\eps(\zeta+\hat\zeta) -\Fcal^\eps(\zeta)
       - d\Fcal^\eps(\zeta)\hat\zeta
&=&
     \int_0^1 (1-r) d^2\Fcal^\eps(\zeta+r\hat\zeta)(\hat\zeta,\hat\zeta)\,dr
     \label{eq:quad1} \\         
     d\Fcal^\eps(\zeta)\hat\zeta - d\Fcal^\eps(0)\hat\zeta
&=&
     \int_0^1 d^2\Fcal^\eps(r\zeta)(\zeta,\hat\zeta)\,dr.
     \label{eq:quad2} 
\end{eqnarray}
To prove assertions~(i) and~(ii) replace
$(\zeta,\zeta_1,\zeta_2)$
by $(\zeta+r\hat\zeta,\hat\zeta,\hat\zeta)$
in~(\ref{eq:d2F1}) and~(\ref{eq:d2F3}),
insert the resulting inequalities in~(\ref{eq:quad1}),
and integrate over $\Sigma$.
Moreover, to derive the second assertion in~(i) 
from the first we use the inequality
$
     \eps^{2/p-1}\|\hat\xi\|_{L^\infty}
     \le c(\eps^{-1}\|\hat\xi\|_{L^p}
         + \|\Tabla{A}\hat\xi\|_{L^p})
$
of Lemma~\ref{lemma:p}.
To prove assertions~(iii) and~(iv) replace
$(\zeta,\zeta_1,\zeta_2)$
by
$(r\zeta,\zeta,\hat\zeta)$
in~(\ref{eq:d2F1}) and~(\ref{eq:d2F3}),
insert the resulting inequalities in~(\ref{eq:quad2}),
and integrate over $\Sigma$.

To prove~(\ref{eq:d2F1})
we give an explicit formula for the second 
derivative of $\Fcal^\eps_1$
in local coordinates on $M$:
$$
     d^2\Fcal^\eps_1(\zeta)(\zeta_1,\zeta_2)
     =
     d^2\rho(\xi)(\xi_1,\xi_2)H
     +\bigl(d\rho(\xi)\xi_2\bigr)H_1
     +\bigl(d\rho(\xi)\xi_1\bigr)H_2
     +\rho(\xi)H_{12}
$$     
Here
$E_j:=E_j(u,\xi)$ for $j=1,2$ 
(see Appendix~\ref{app:invariant}),
and $H$, $H_1$, $H_2$ and $H_{12}$ 
are defined as follows.
The section $H=H(\zeta)$ is given by
\begin{eqnarray*}
         H
         &=&
         E_1(v_s+X_\phi(u))
         +
         E_2(\Nabla{A,s}\xi
         +
         \Nabla{\xi}X_\phi(u))
         \\
         &&
         +
         J(u_\xi)
         \Bigl(
         E_1(v_t+X_\psi(u))
         +
         E_2(\Nabla{A,t}\xi
         +
         \Nabla{\xi}X_\psi(u))
         \Bigr).
\end{eqnarray*}
Now we use the following notation for $i=1,2$:
\begin{eqnarray*}
        \p E_i(v,w)
        &:= &
        \bigl(
        \p_2E_i(u,\xi)v
        \bigr)w,
        \\
        \p^2E_i(v_1,v_2,w)
        &:= &
        \bigl(
        \p_2\p_2E_i(u,\xi)(v_1,v_2)
        \bigr)w.
\end{eqnarray*}
The section $H_j=H_j(\zeta;\zeta_j)$
is linear in $\zeta_j$ and is defined by
\begin{eqnarray*}
         H_j
&=&
         \p E_1(\xi_j,v_s+X_\phi(u))
         + E_1X_{\phi_j}(u)
         + \p E_2\left(\xi_j,\Nabla{A,s}\xi
           + \Nabla{\xi}X_\phi(u)\right) \\
&&
         +\, E_2\left(\Nabla{A,s}\xi_j
         + \Nabla{\xi_j}X_\phi(u)
         + \Nabla{\xi}X_{\phi_j}(u)\right)\\
&&
         +\, \Bigl(dJ(u_\xi)E_2\xi_j\Bigr)
         \Bigl(
         E_1(v_t+X_\psi(u))
         + E_2(\Nabla{A,t}\xi
         + \Nabla{\xi}X_\psi(u))
         \Bigr) \\
&&
         +\, J(u_\xi)\Bigl(
           \p E_1(\xi_j,v_t+X_\psi(u))
           + \p E_2\left(\xi_j,\Nabla{A,t}\xi+\Nabla{\xi}X_\psi(u)
         \right)  \\
&&
         \qquad
           +\, E_1X_{\psi_j}(u)
           + E_2\left(\Nabla{A,t}\xi_j
             + \Nabla{\xi_j}X_\psi(u)
             + \Nabla{\xi}X_{\psi_j}(u)\right)
         \Bigr).
\end{eqnarray*}
The section $H_{12}=H_{12}(\zeta;\zeta_1,\zeta_2)$
is bilinear in $\zeta_1$ and $\zeta_2$,
and is defined by
\begin{eqnarray*}
         H_{12}
&=&
         \p^2 E_1
         (\xi_1,\xi_2,v_s+X_\phi(u))
         + \p E_1(\xi_2,X_{\phi_1}(u))
         + \p E_1(\xi_1,X_{\phi_2}(u)) \\
&&
         +\, E_2\left(\Nabla{\xi_2}X_{\phi_1}(u)
             + \Nabla{\xi_1}X_{\phi_2}(u)\right)
         +   \p^2 E_2\left(\xi_1,\xi_2,\Nabla{A,s}\xi
             + \Nabla{\xi}X_\phi(u)\right) \\
&&
         +\, \p E_2\left(\xi_2,\Nabla{A,s}\xi_1
            + \Nabla{\xi_1}X_\phi(u)
            + \Nabla{\xi}X_{\phi_1}(u)\right) \\
&&
         +\, \p E_2\left(\xi_1,\Nabla {A,s}\xi_2
            + \Nabla{\xi_2}X_\phi(u)
            + \Nabla{\xi}X_{\phi_2}(u)\right) \\
&&
         +\, \bigl(d^2J(u_\xi)(E_2\xi_1,E_2\xi_2)\bigr) 
         \Bigl(
         E_1\left(v_t+X_\psi(u)\right)
         +  E_2\left(\Nabla{A,t}\xi+\Nabla{\xi}X_{\psi}(u)\right)
         \Bigr)  \\
&&
         +\,\bigl(dJ(u_\xi)\p E_2(\xi_1,\xi_2)\bigr)
         \Bigl(
         E_1\left(v_t+X_\psi(u)\right)
         + E_2\left(\Nabla{A,t}\xi+\Nabla{\xi}X_\psi(u)\right)
         \Bigr) \\
&&
         + \bigl(dJ(u_\xi)E_2\xi_1\bigr)
         \Bigl(
         \p E_1(\xi_2,v_t+X_\psi(u))  
         + \p E_2\left(
             \xi_2,\Nabla {A,t}\xi+\Nabla{\xi}X_\psi(u)
           \right) \\
&&
         \qquad
         +\, E_1X_{\psi_2}(u)
         + E_2\left(\Nabla{A,t}\xi_2
           + \Nabla{\xi_2}X_\psi(u)
           + \Nabla{\xi}X_{\psi_2}(u)\right)
         \Bigr) \\
&&
         +\,\bigl(dJ(u_\xi)E_2\xi_2\bigr)
         \Bigl(
         \p E_1(\xi_1,v_t+X_\psi(u))
         + \p E_2(\xi_1,\Nabla{A,t}\xi
           + \Nabla{\xi}X_\psi(u))  \\
&&
         \qquad
         +\, E_1X_{\psi_1}(u) 
         +   E_2\left(\Nabla{A,t}\xi_1
             + \Nabla{\xi_1}X_\psi(u)
             + \Nabla{\xi}X_{\psi_1}(u)\right)
         \Bigr) \\
&&
         +\, J(u_\xi)\Bigl(
         \p^2 E_1\left(\xi_1,\xi_2,v_t+X_\psi(u)\right) \\
&&
         \qquad
         +\, \p E_1(\xi_2,X_{\psi_1}(u))
         + \p E_1(\xi_1,X_{\psi_2}(u))  \\
&&
         \qquad
         +\, \p^2 E_2\left(\xi_1,\xi_2,
           \Nabla{A,t}\xi+\Nabla{\xi}X_\psi(u)\right)  \\
&&
         \qquad
         +\, \p E_2\left(\xi_2,\Nabla{A,t}\xi_1
         +\Nabla{\xi_1}X_\psi(u)
         +\Nabla{\xi}X_{\psi_1}(u)\right) \\
&&
         \qquad
         +\, \p E_2\left(\xi_1,
         \Nabla{A,t}\xi_2
         +\Nabla{\xi_2}X_\psi(u)
         +\Nabla{\xi}X_{\psi_2}(u)
         \right) \\
&&
         \qquad
         +\, E_2\left(\Nabla{\xi_2}X_{\psi_1}(u)
           + \Nabla{\xi_1}X_{\psi_2}(u)\right)
         \Bigr).
\end{eqnarray*}
The inequalities~(\ref{eq:d2F1}) and~(\ref{eq:d2F3})
now follow by a term by term inspection of 
$H$, $H_1$, $H_2$, and $H_{12}$, assuming $\|\xi\|_{L^\infty}\le 1$.
\end{proof} 

\begin{proposition}\label{prop:quadra1}
For every $p>2$ and every $c_0>0$
there exists a cons\-tants $c>0$ 
such that the following holds
for every $(u,A)\in\Tilde\Mcal_{B,\Sigma}^0(c_0)$,
any two pairs
$
    \zeta=(\xi,\alpha),\hat\zeta=(\hat\xi,\hat\alpha)\in\Xx_u,
$
and every $\eps\in(0,1]$.

\smallskip
\noindent{\bf (i)}
If $\|\xi\|_{L^\infty}+\|\hat\xi\|_{L^\infty}\le1$ and 
$\|\zeta\|_{1,p,\eps}+\|\hat\zeta\|_{1,p,\eps}\le\eps^{2/p}$
then
$$
     \|
     \Ff^\eps(\zeta+\hat\zeta) - \Ff^\eps(\zeta) - d\Ff^\eps(\zeta)\hat\zeta
     \|_{0,p,\eps}
     \le c\eps^{-1-2/p}
     \|\hat\zeta\|_{1,p,\eps}^2.
$$

\smallskip
\noindent{\bf (ii)}
If $\|\xi\|_{L^\infty}+\|\hat\xi\|_{L^\infty}\le1$ and 
$\|\zeta\|_{2,p,\eps}+\|\hat\zeta\|_{2,p,\eps}\le\eps^{2/p}$
then
$$
     \|
     \Ff^\eps(\zeta+\hat\zeta) - \Ff^\eps(\zeta) - d\Ff^\eps(\zeta)\hat\zeta
     \|_{1,p,\eps}
     \le c\eps^{-1-2/p}
     \|\hat\zeta\|_{2,p,\eps}^2.
$$

\smallskip
\noindent{\bf (iii)}
If $\|\xi\|_{L^\infty} \le 1$ then
$$
     \|d\Ff^\eps(\zeta)\hat\zeta- d\Ff^\eps(0)\hat\zeta\|_{0,p,\eps}
     \le c\eps^{-1-2/p}
     \|\zeta\|_{1,p,\eps}\|\hat\zeta\|_{1,p,\eps}. 
$$

\smallskip
\noindent{\bf (iv)}
If $\|\xi\|_{L^\infty}\le1$ and  
$\|\zeta\|_{1,p,\eps}\le\eps^{2/p}$
then
$$
     \|d\Ff^\eps(\zeta)\hat\zeta- d\Ff^\eps(0)\hat\zeta\|_{1,p,\eps}
     \le c\eps^{-1-2/p}\|\zeta\|_{2,p,\eps}\|\hat\zeta\|_{2,p,\eps}.
$$
\end{proposition}

\begin{proof}
Assertions~(i) and~(iii) follow
immediately from Proposition~\ref{prop:quadra}. 
To prove~(ii) we observe that 
in estimating the quadratic terms in $\p\Ff^\eps_1$ 
we encounter products of the following forms
\begin{itemize}
\item
$\p^2\xi\cdot\hat\xi\cdot\hat\xi$ and $\p\phi\cdot\hat\xi\cdot\hat\xi$.
Here the $L^p$-norms of $\p^2\xi$ and $\p\phi$
can be estimated by $\eps^{2/p-2}$ and the 
$L^\infty$-norm of $\hat\xi\cdot\hat\xi$
by $\eps^{-4/p}\|\hat\zeta\|_{1,p,\eps}^2$. 
\item 
$\xi\cdot\p\hat\xi\cdot\p\hat\xi$, $\xi\cdot\hat\xi\cdot\p^2\hat\xi$, 
$\xi\cdot\p\hat\xi\cdot\hat\phi$, and $\xi\cdot\p\hat\phi\cdot\hat\xi$.  
The $L^p$-norms of these products can be estimated 
by 
$
     \eps^{-2-2/p}\|\hat\zeta\|_{1,p,\eps}\|\hat\zeta\|_{2,p,\eps}.
$
\item
$\p\xi\cdot\p\hat\xi\cdot\hat\xi$, 
$\phi\cdot\p\hat\xi\cdot\hat\xi$, 
and $\p\xi\cdot\hat\phi\cdot\hat\xi$.
In these cases the $L^p$-norm of $\p\xi$ is bounded 
by $\eps^{2/p-1}$ and the $L^\infty$-norms
of $\p\hat\xi\cdot\hat\xi$ and $\hat\phi\cdot\hat\xi$
are bounded by 
$
     \eps^{-1-4/p}\|\hat\zeta\|_{1,p,\eps}\|\hat\zeta\|_{2,p,\eps}. 
$
\end{itemize}
Similarly, in estimating the quadratic terms in $\p\Ff^\eps_3$ 
we encounter products of the following forms
\begin{itemize}
\item
$\eps^{-2}\hat\xi\cdot\p\hat\xi$ 
and $\p\hat\phi\cdot\hat\psi$.
The $L^p$-norms of these products can be 
estimated by 
$
     \eps^{-3-2/p}\|\hat\zeta\|_{1,p,\eps}^2. 
$
\item
$\eps^{-2}\p\xi\cdot\hat\xi\cdot\hat\xi$.
Here the $L^p$-norm of $\p\xi$ bounded by 
$\eps^{2/p-1}$ and the $L^\infty$-norm of 
$\eps^{-2}\hat\xi\cdot\hat\xi$ is bounded by 
$
     \eps^{-2-4/p}\|\hat\zeta\|_{1,p,\eps}^2. 
$
\end{itemize}
This proves~(ii).  The proof of~(iv) is similar.
\end{proof}

Assertions~(i) and~(iii) in Proposition~\ref{prop:quadra1} 
are weaker than Proposition~\ref{prop:quadra};
in the former the first derivatives of $\alpha$
appear on the right hand sides of the estimates.
The full strength of Proposition~\ref{prop:quadra}
will be required in the proof of 
Theorem~\ref{thm:uniqbaby} below.
 
 
\section{Proof of Theorem~B}\label{sec:exist}

>From now on we assume~(H1) and~(H4). 
In this section we establish the existence of a 
$\Gg(P)$-equivariant map 
$
     \Tilde{\Tcal}^{\eps}:
     \Tilde{\Mcal}_{B,\Sigma}^0(c_0)
     \to\Tilde{\Mcal}_{B,\Sigma}^{\eps}
$
that satisfies the requirements of Theorem~B. 

\begin{theorem}\label{thm:exist}
For every $c_0>0$ and every $p>2$
there exist positive cons\-tants 
$\eps_0$ and $c$ 
such that for every 
$\eps \in (0,\eps_0]$ the following holds. 
For every $(u_0,A_0)\in\Tilde{\Mcal}_{B,\Sigma}^0(c_0)$
there exists a unique pair 
$
     \zeta_\eps=(\xi_\eps,\alpha_\eps)\in\Xx_{u_0}
$
such that 
\begin{equation}\label{eq:ueps}
      (u_\eps,A_\eps)
      := (\exp_{u_0}(\xi_\eps),A_0+\alpha_\eps)
      \in \Tilde{\Mcal}_{B,\Sigma}^{\eps},
\end{equation}
\begin{equation}\label{eq:exist}
     -d_{A_0}^*\alpha_\eps + \eps^{-2}L_{u_0}^*\xi_\eps=0,\qquad
     \zeta_\eps\in \im\,(\Dcal_{(u,A)}^\eps)^*,
\end{equation}
\begin{equation}\label{eq:eps2}
     \left\|\zeta_\eps\right\|_{2,p,\eps;(u_0,A_0)}
     \le c\eps^2.
\end{equation}
The map $(u_0,A_0)\mapsto(u_\eps,A_\eps)$ is $\Gg(P)$-equivariant
and will be denoted by $\Tilde{\Tcal}^{\eps}$.
\end{theorem}

The next theorem shows that uniqueness holds
under a slightly weaker hypothesis, namely in a 
larger neighbourhood of $(u,A)$. 

\begin{theorem}\label{thm:uniqbaby}
For every $c_0>0$ and every $p>2$
there exist po\-si\-ti\-ve cons\-tants
$\delta$ and $\eps_0$
such that for every $\eps \in (0,\eps_0]$
the following holds.
Suppose that $(u_0,A_0)\in\Tilde\Mm_{B,\Sigma}^0(c_0)$ and
$
    (\xi,\alpha) \in \Xx_{u_0}
$
satisfy~(\ref{eq:ueps}), (\ref{eq:exist}), and
\begin{equation}\label{eq:unique}
     \left\|\xi\right\|_{1,p,\eps}
     + \eps \left\|\alpha\right\|_{L^p}
     + \eps^{2/p}\left\|(\xi,\alpha)\right\|_{\infty,\eps}
     \le \delta\eps^{2/p+1/2}.
\end{equation}
Then
$
      (\exp_{u_0}(\xi),A_0+\alpha) = \Tilde\Tt^\eps(u,A).
$
\end{theorem}

\begin{corollary}\label{cor:uniqbaby}
For every $c_0>0$ and every $p>2$
there exist po\-si\-ti\-ve cons\-tants
$\delta$ and $\eps_0$
such that for every $\eps\in(0,\eps_0]$
the following holds.
Suppose that $(u_0,A_0)\in\Tilde\Mm_{B,\Sigma}^0(c_0)$ and
$
     \zeta=(\xi,\alpha) \in \Xx_{u_0}
$
satisfy~(\ref{eq:ueps}), (\ref{eq:exist}), and
\begin{equation}\label{eq:uniqbaby}
     \left\|\zeta\right\|_{1,p,\eps;(u_0,A_0)}
     \le \delta\eps^{2/p+1/2}.
\end{equation}
Then
$
      (\exp_{u_0}(\xi),A_0+\alpha)
      = \Tilde\Tcal^\eps(u,A).
$
\end{corollary}

\begin{proof}
Theorem~\ref{thm:uniqbaby} and Lemma~\ref{lemma:p}.
\end{proof}

\begin{proof}[Proof of Theorem~\ref{thm:exist}]
The proof is similar to that of Theorem~5.1 in~\cite{DS2}.
However, in the present case the nonlinearities (in the
quadratic estimates) appear in the highest order terms,
and we establish estimates for the $(2,p,\eps)$-norms 
and not just the $(1,p,\eps)$-norms (as in~\cite{DS2}).  
We assume throughout that the exponential map at each 
point in $\mu^{-1}(0)$ is defined in a ball of radius one. 

Abbreviate $\Dcal^\eps:=\Dcal_{(u_0,A_0)}^\eps$ and
let $\Fcal^\eps:\Xx_{u_0}\to\Xx'_{u_0}$
be defined by~(\ref{eq:Feps}). 
Then 
$$
     d\Fcal^\eps(0)=\Dcal^\eps,\qquad
     \Fcal^\eps(0)=(0,0,*F_{A_0}).
$$  
Hence, by Lemma~\ref{le:c0}, there exists a constant $C_0>0$,
depending only on $c_0$ and $p$, such that 
$$
     \left\|\Fcal^\eps(0)\right\|_{1,p,\eps}
     = \eps\left\|F_{A_0}\right\|_{L^p}
       + \eps^2\left\|d_{A_0}*F_{A_0}\right\|_{L^p}
     \le C_0\eps.
$$
We use Newton iteration to obtain a zero of~$\Fcal^\eps$,
and hence a solution of~(\ref{eq:eps}).
Let
$
     \zeta_\nu = (\xi_\nu,\alpha_\nu)\in\Xx_{u_0}
$
be the sequence defined re\-cur\-sively 
by~$\zeta_0:=0$ and
\begin{equation}\label{eq:recursive}
        \zeta_{\nu+1}
        := \zeta_\nu + \hat\zeta_\nu,\qquad
        \hat\zeta_\nu\in\im\,{\Dd^\eps}^*,\qquad
        \Dcal^\eps\hat\zeta_\nu
        = -\Fcal^\eps(\zeta_\nu).
\end{equation}
We prove by induction over $\nu$ that there exist positive 
cons\-tants~$\eps_0,c_1,C$, depending only on $c_0$ and $p$, 
such that
\begin{eqnarray}
       \label{eq:exist1}
       \|\hat\zeta_\nu\|_{2,p,\eps}
       &\le &
       c_1\|\Fcal^\eps(\zeta_\nu)\|_{1,p,\eps},
       \\
       \label{eq:exist2}
       \|\hat\zeta_\nu\|_{2,p,\eps}
       &\le &
       2^{-\nu}C_0c_1
       \eps^2,
       \\
       \label{eq:exist3}
       \|\Fcal^\eps(\zeta_{\nu+1})\|_{1,p,\eps}
       &\le&
       C \eps^{1-2/p}\|\hat\zeta_\nu\|_{2,p,\eps}
\end{eqnarray}
for $\eps\in(0,\eps_0]$ and $\nu\ge 0$. 
The constants are chosen such that
the linear estimates of Lemma~\ref{le:onto1} hold 
for $0<\eps\le\eps_0$ with $c=c_1\ge1$, the quadratic estimates
of Proposition~\ref{prop:quadra1} hold for $0<\eps\le1$ with $c=c_2\ge1$, 
the $L^\infty$ estimates of Lemma~\ref{lemma:p} hold
for $0<\eps\le 1$ with $c=c_\infty\ge1$, and 
$$
     C=3C_0c_1c_2,\qquad
     Cc_1\eps_0^{1-2/p}\le1/2,\qquad
     3C_0c_1c_\infty\eps_0^{2-2/p} \le 1.
$$
For $\nu=0$ the estimates~(\ref{eq:exist1}) and~(\ref{eq:exist2})  
follow from Lemma~\ref{le:onto1}.  Namely, 
by~(\ref{eq:onto1}) with $k=1$, we have
$$
      \|\hat\zeta_0\|_{2,p,\eps}
      \le c_1\eps\|\Ff^\eps(0)\|_{1,p,\eps}
      \le C_0c_1\eps^2.
$$
The estimate~(\ref{eq:exist3}) for $\nu=0$ follows from 
the identity
$
      d\Fcal^\eps(0)\hat\zeta_0
      =-\Fcal^\eps(0)
$
and Proposition~\ref{prop:quadra1}~(ii).  Namely, since
$C_0c_1c_\infty\eps^{2-2/p} \le 1$ we have 
$$
     \|\hat\zeta_0\|_{2,p,\eps}
     \le C_0c_1\eps^2 \le\eps^{2/p},\qquad
     \|\hat\zeta_0\|_{\infty,\eps}
     \le C_0c_1c_\infty\eps^{2-2/p} 
     \le 1.
$$
Hence the hypotheses of Proposition~\ref{prop:quadra1}~(ii)
are satisfied with $\zeta=0$ and $\hat\zeta=\hat\zeta_0$,
and hence
\begin{eqnarray*}
     \|\Ff^\eps(\zeta_1)\|_{1,p,\eps}
&= &
     \|\Ff^\eps(\hat\zeta_0) 
     - \Ff^\eps(0) - d\Ff^\eps(0)\hat\zeta_0
     \|_{1,p,\eps} \\
&\le &
     c_2\eps^{-1-2/p}\|\hat\zeta_0\|_{2,p,\eps}^2  \\
&\le &
     C_0c_1c_2\eps^{1-2/p}\|\hat\zeta_0\|_{2,p,\eps}.
\end{eqnarray*}
Since $C\ge C_0c_1c_2$ this proves~(\ref{eq:exist3})
for $\nu=0$.  Now assume that the sequences 
$\zeta_0,\dots,\zeta_\nu$ and $\hat\zeta_0,\dots,\hat\zeta_{\nu-1}$
have been constructed up to some integer $\nu\ge 1$ and
that the estimates~(\ref{eq:exist1}-\ref{eq:exist3}) 
have been established for all integers up to $\nu-1$. 
Then, by~(\ref{eq:exist2}),
\begin{equation}\label{eq:exist4}
     \left\|\zeta_\nu\right\|_{2,p,\eps}
     \le \sum_{j=0}^{\nu-1}
         \|\zeta_{j+1}-\zeta_j\|_{2,p,\eps}      
     = \sum_{j=0}^{\nu-1}\|\hat\zeta_j\|_{2,p,\eps}
     \le 2C_0c_1\eps^2,
\end{equation}
and hence 
$$
     \left\|\zeta_\nu\right\|_{\infty,\eps}   
     \le 2C_0c_1c_\infty\eps^{2-2/p}
     \le 1. 
$$
This shows that $\xi_\nu(p)$ lies in the domain
of the exponential map at $u_0(p)$ for every $p\in P$
and so $\zeta_\nu$ lies in the domain of $\Ff^\eps$. 
Let $\hat\zeta_\nu$ and $\zeta_{\nu+1}$ be defined 
by~(\ref{eq:recursive}). Then, by Lemma~\ref{le:onto1},
$\hat\zeta_\nu$ satisfies the estimate~(\ref{eq:exist1}).
To prove~(\ref{eq:exist2}) we observe that,
by the induction hypothesis,  
\begin{eqnarray*}
     \|\hat\zeta_\nu\|_{2,p,\eps}
&\le &
     c_1 \|\Fcal^\eps(\zeta_\nu)\|_{1,p,\eps}  \\
&\le &
     Cc_1\eps^{1-2/p}\|\hat\zeta_{\nu-1}\|_{2,p,\eps}  \\ 
&\le &
     2^{-1}\|\hat\zeta_{\nu-1}\|_{2,p,\eps}  \\
&\le &
     2^{-\nu}C_0c_1\eps^2.
\end{eqnarray*}
To prove~(\ref{eq:exist3}) we observe that, 
by~(\ref{eq:exist4}),
$$
\begin{array}{rcccl}
     \|\zeta_\nu\|_{2,p,\eps}
     + \|\hat\zeta_\nu\|_{2,p,\eps}
     &\le & 3C_0c_1\eps^2
     &\le & \eps^{2/p},  \\
     \|\zeta_\nu\|_{\infty,\eps}
     + \|\hat\zeta_\nu\|_{\infty,\eps}
     &\le &3C_0c_1c_\infty\eps^{2-2/p} 
     &\le &1.
\end{array}
$$
Thus the hypotheses of Proposition~\ref{prop:quadra1}~(ii)
and~(iv) are satisfied with $\zeta=\zeta_\nu$ 
and $\hat\zeta=\hat\zeta_\nu$. Hence
\begin{eqnarray*}
      \|\Fcal^\eps(\zeta_{\nu+1})\|_{1,p,\eps}
&\le &
     \|\Fcal^\eps(\zeta_\nu+\hat\zeta_\nu)
     - \Fcal^\eps(\zeta_\nu)
     - d\Fcal^\eps(\zeta_\nu)\hat\zeta_\nu
     \|_{1,p,\eps}  \\
&&
     +\,\|d\Fcal^\eps(\zeta_\nu)\hat\zeta_\nu
       - d\Fcal^\eps(0)\hat\zeta_\nu\|_{1,p,\eps}  \\
&\le &
      c_2\eps^{-1-2/p}
      \bigl(
      \|\hat\zeta_\nu\|_{2,p,\eps}
      + \|\zeta_\nu\|_{2,p,\eps}
      \bigr)
      \|\hat\zeta_\nu\|_{2,p,\eps} \\
&\le &
      3C_0c_1c_2\eps^{1-2/p}
      \|\hat\zeta_\nu \|_{2,p,\eps} \\
&= &
      C\eps^{1-2/p}\|\hat\zeta_\nu\|_{2,p,\eps}.
\end{eqnarray*}
This completes the induction. 

By~(\ref{eq:exist2}), the sequence $\zeta_\nu$
is Cauchy in the $W^{2,p}$-norm.
Moreover, by examining the second component of $\Ff^\eps$
we find that $\zeta_\nu$ satisfies~(\ref{eq:exist}) for every
$\nu$ and hence so does its limit
$$
     \zeta_\eps
     := (\xi_\eps,\alpha_\eps) 
     := \lim_{\nu\to\infty}(\xi_\nu,\alpha_\nu).
$$
By~(\ref{eq:exist4}), this limit also satisfies~(\ref{eq:eps2})
with $c:=2C_0c_1$. Moreover, by~(\ref{eq:exist2})
and~(\ref{eq:exist3}), the sequence 
$\Ff^\eps(\zeta_\nu)$
converges to zero in the $(1,p,\eps)$-norm and hence 
$\Ff^\eps(\zeta_\eps)=0$. Hence 
$\zeta_\eps$ satisfies~(\ref{eq:ueps}) 
and it follows from elliptic regularity that 
$\zeta_\eps$ is smooth. This proves existence. 

We prove uniqueness.  Suppose
$
     \zeta=(\xi,\alpha)\in\Xx_{u_0}
$ 
satisfies~(\ref{eq:ueps}), (\ref{eq:exist}),
and $\|\zeta\|_{1,p,\eps}\le c\eps^2$.  
Then, by Lemma~\ref{le:onto1} and 
Proposition~\ref{prop:quadra1}~(i) and~(iii), 
\begin{eqnarray*}
     \|\zeta-\zeta_\eps\|_{1,p,\eps}
&\le & 
     c_1
     \|\Dcal^\eps(\zeta-\zeta_\eps)\|_{0,p,\eps} \\
&\le &
     c_1 
     \|\Fcal^\eps(\zeta) - \Fcal^\eps(\zeta_\eps)
       - d\Fcal^\eps(\zeta_\eps)(\zeta-\zeta_\eps)
     \|_{0,p,\eps}  \\
&&
     +\, c_1 
     \|d\Fcal^\eps(0)(\zeta-\zeta_\eps)
      - d\Fcal^\eps(\zeta_\eps)(\zeta-\zeta_\eps)
     \|_{0,p,\eps}  \\
&\le &
     c_1c_2\eps^{-1-2/p}
     \bigl(
     \|\zeta -\zeta_\eps\|_{1,p,\eps}
     + \|\zeta_\eps\|_{1,p,\eps}
     \bigr)
     \|\zeta -\zeta_\eps\|_{1,p,\eps} \\
& \le &
     3cc_1c_2\eps^{1-2/p}
     \|\zeta-\zeta_\eps\|_{1,p,\eps}.
\end{eqnarray*}
If $3cc_1c_2\eps^{1-2/p}<1$ then
$\zeta=\zeta_\eps$.  This proves uniqueness.

Since the conditions~(\ref{eq:ueps}), (\ref{eq:exist}), 
and~(\ref{eq:eps2}) are gauge invariant, it follows
that the map $(u,A)\mapsto(u_\eps,A_\eps)$
is $\Gg(P)$-equivariant. 
\end{proof}

\begin{proof}[Proof of Theorem~\ref{thm:uniqbaby}]
In this proof we drop the subscript $0$. 
Fix two pairs 
$
     (u,A)\in\Tilde{\Mm}_{B,\Sigma}^0(c_0)
$
and
$
     \zeta=(\xi,\alpha)\in\Xx_u
$
that satisfy~(\ref{eq:exist}), (\ref{eq:unique}),
and
$$
     (\exp_u(\xi),A+\alpha)\in\Tilde{\Mm}_{B,\Sigma}^\eps.
$$
We prove that $\zeta$ satisfies~(\ref{eq:eps2}),
provided that $\delta$ and $\eps$ are sufficiently small. 
By ellipticity of the operator $\Dd^0:=\Dd_{(u,A)}^0$,
there exists a constant $c_1=c_1(p,c_0)>0$ such that
$$
     \|\Tabla{A}\pi_{u}\xi\|_{L^p}
     \le c_1\left(\|\Dd^0\pi_u\xi\|_{L^p}
         + \|\pi_u\xi\|_{L^p}
         \right). 
$$
Now let $c_2=c_2(p,c_0)$ be the constant 
of Lemma~\ref{le:difference}
and $c_3=c_3(p,c_0)$ be the constant of 
Proposition~\ref{prop:estimate}.  Then   
\begin{eqnarray*}
     \|\Tabla{A}\xi\|_{L^p}
&\le &
     \eps^{-1}\|\xi-\pi_u\xi\|_{1,p,\eps}
     +
     \|\Tabla{A}\pi_u\xi\|_{L^p} \\
&\le &
     \eps^{-1}\|\xi-\pi_u\xi\|_{1,p,\eps}
     +
     c_1\left(\|\Dcal^0\pi_u\xi\|_{L^p}
         + \|\pi_u\xi\|_{L^p}
         \right) \\
&\le &
     \eps^{-1}\|\xi-\pi_u\xi\|_{1,p,\eps}  \\
&&
     +\,
     c_1
     \left(\|(\Dcal^0\pi_u-\pi_u\Dcal^\eps)\xi\|_{L^p}
     +
     \|\pi_u\Dcal^\eps\zeta\|_{L^p}
     +
     \|\pi_u\xi\|_{L^p}\right) \\
&\le &
     (\eps^{-1}+c_1c_2)\|\xi-\pi_u\xi\|_{1,p,\eps} 
     + c_1\left(
     \|\pi_u\Dcal^\eps\zeta\|_{L^p}
     +
     \|\pi_u\xi\|_{L^p}\right) \\
&\le &
     c_4
     \left(\|\Dcal^\eps\zeta\|_{0,p,\eps}
     +
     \|\xi\|_{L^p}\right),
\end{eqnarray*}
where $c_4:=c_3(1+c_1c_2)+c_1$. 
Hence 
\begin{eqnarray*}
      \|\Tabla{A}\xi\|_{L^p} + \|\alpha\|_{L^p}
&\le &
      c_4\|\Dcal^\eps\zeta\|_{0,p,\eps}
      + c_4\|\xi\|_{L^p}
      + \eps^{-1}\|\zeta - \pi_u \zeta\|_{1,p,\eps}  \\
&\le &
      (c_4+c_5)\|\Dcal^\eps\zeta\|_{0,p,\eps}
      + c_4\|\xi\|_{L^p},
\end{eqnarray*}
where $c_5=c_5(p,c_0)$ is the constant of Lemma~\ref{le:onto1}. 
Since 
$$
     \Fcal^\eps(\zeta)=0,\qquad
     \Fcal^\eps(0)=(0,0,*F_A),\qquad
     \Dcal^\eps=d\Fcal^\eps(0)
$$ 
we obtain
\begin{eqnarray}\label{eq:unique0}
     \|\Tabla{A}\xi\|_{L^p} + \|\alpha\|_{L^p}
&\le &
     c_6\|\Fcal^\eps(\zeta)-\Fcal^\eps(0) - d\Fcal^\eps(0)\zeta
     \|_{0,p,\eps}  
     \nonumber \\
&&
     +\, c_6\eps\|F_A\|_{L^p}
     + c_4\|\xi\|_{L^p},
\end{eqnarray}
where $c_6:=c_4+c_5$. Now we use the refined quadratic estimate 
of Proposition~\ref{prop:quadra} with $c=c_7$. 
By~(\ref{eq:unique}), we have
$$
     \|\Tabla{A}\xi\|_{L^p}+\|\alpha\|_{L^p}
     \le \delta\eps^{2/p-1/2}
     \le \eps^{2/p-1}
$$
(provided that $\delta\le 1$).  
Thus the hypotheses of Proposition~\ref{prop:quadra}~(i) 
and~(ii) are satisfied with $\zeta=0$
and $\hat\zeta$ replaced by $\zeta$.
Hence, by~(\ref{eq:unique0}), 
\begin{eqnarray*}
      \|\Tabla{A}\xi\|_{L^p} + \|\alpha\|_{L^p}
&\le &
      c_6c_7
      \|\xi\|_{L^\infty}
      \left(
        \eps^{-1}\|\xi\|_{L^p}+\|\Tabla{A}\xi\|_{L^p}
        + \|\alpha\|_{L^p}
      \right)   \\
&&
      +\,c_6c_7\eps
      \left(
        \|\alpha\|_{L^\infty}\|\alpha\|_{L^p}
        + \eps^{-2}\|\xi\|_{L^\infty}\|\xi\|_{L^p}
      \right)   \\
&&
     +\, c_6\eps \|F_A\|_{L^p} + c_4\|\xi\|_{L^p} \\
&\le &
      3c_6c_7
      \|\zeta\|_{\infty,\eps}
      \left(
        \eps^{-1}\|\xi\|_{L^p}
        + \|\Tabla{A}\xi\|_{L^p}
        + \|\alpha\|_{L^p}
      \right)   \\
&&
     +\, c_6\eps \|F_A\|_{L^p} + c_4\|\xi\|_{L^p}  \\
&\le & 
     3c_6c_7\delta \eps^{1/2}
     \left(\eps^{-1}\|\xi\|_{L^p} 
       + \|\Tabla A \xi\|_{L^p}
       + \|\alpha\|_{L^p}
     \right)  \\
&&
     +\, c_0c_6\eps + c_4\|\xi\|_{L^p}  \\
&\le & 
     3c_6c_7\delta\eps^{1/2}\left(
      \|\Tabla{A}\xi\|_{L^p} + \|\alpha\|_{L^p}
     \right)  \\
&&
     +\, c_0c_6\eps + 
     c_4\delta\eps^{2/p+1/2}
     + 3c_6c_7\delta^2\eps^{2/p}.
\end{eqnarray*}
Here the last two inequalities follow
from~(\ref{eq:unique}).  With
$3c_6c_7\delta\eps^{1/2}\le1/2$ we have
\begin{equation}\label{eq:unique1}
     \|\Tabla{A}\xi\|_{L^p} + \|\alpha\|_{L^p}
     \le c_8(\delta+\eps^{1-2/p})\eps^{2/p},
\end{equation}
where $c_8:=2c_0c_6+2c_4+6c_6c_7$. 
Since $\zeta$ satisfies~(\ref{eq:exist}) we can apply 
Lemma~\ref{le:onto1} (with $c=c_5$) to obtain
\begin{eqnarray}\label{eq:unique2}
   \|\zeta\|_{1,p,\eps}
&\le &
   c_5
   \left(\eps\|\Dd^\eps\zeta\|_{L^p}
   + \|\pi_u\Dd^\eps\zeta\|_{L^p}
   \right) 
   \nonumber \\
&\le &
   c_5
   \left(2\|d\Fcal^\eps_1(0)\zeta\|_{L^p}
   + \eps^2 \|d\Fcal^\eps_3(0)\zeta\|_{L^p}\right).
\end{eqnarray}
By Proposition~\ref{prop:quadra}~(i) with $c=c_7$
and~(\ref{eq:unique1}), we now have
\begin{eqnarray*}
     \left\|d\Fcal^\eps_1(0)\zeta\right\|_{L^p}
&= &
     \left\|\Fcal^\eps_1(\zeta)-\Fcal^\eps_1(0)
     - d\Fcal^\eps_1(0)\zeta\right\|_{L^p} \\
&\le &
     c_7
     \|\xi\|_{L^\infty}
     \Bigl(
     \|\xi\|_{L^p} + \|\Tabla{A}\xi\|_{L^p}
     + \|\alpha\|_{L^p} 
     \Bigr)  \\
&& 
     +\, c_7 \|\xi\|_{L^\infty}^2
     \Bigl(\|\Tabla{A}\xi\|_{L^p} + \|\alpha\|_{L^p}\Bigr)  \\
&\le &
     2c_7
     \|\xi\|_{L^\infty}
     \Bigl(
     \|\xi\|_{L^p} + \|\Tabla{A}\xi\|_{L^p} + \|\alpha\|_{L^p}
     \Bigr) \\
&\le &
     2c_7
     \|\zeta\|_{\infty,\eps} \,
     \left(
     \delta\eps^{2/p+1/2}
     + c_8(\delta+\eps^{1-2/p})\eps^{2/p}
     \right) \\
&\le &
     c_9 
     (\delta+\eps^{1-2/p})\|\zeta\|_{1,p,\eps}.
\end{eqnarray*}
Here we have used the fact that
$\|\xi\|_{L^\infty}\le\delta\eps^{1/2}\le1$.
Moreover, the penultimate inequality follows 
from~(\ref{eq:unique}) and~(\ref{eq:unique1})
and the last inequality, with a suitable constant 
$c_9=c_9(p,c_0)$, follows from Lemma~\ref{lemma:p}. 
By Proposition~\ref{prop:quadra}~(ii)
with $c=c_7$, we have
\begin{eqnarray*}
     \|d\Fcal^\eps_3(0)\zeta\|_{L^p}
&\le &
   \|\Fcal^\eps_3(0)
   -
   \Fcal^\eps_3(\zeta)
   -
   d\Fcal^\eps_3(0)\zeta\|_{L^p}
   +\|F_A\|_{L^p}
   \\
&\le &
   c_7\left(
   \|\alpha\|_{L^\infty}\|\alpha\|_{L^p}
   + \eps^{-2}\|\xi\|_{L^\infty}\|\xi\|_{L^p}\right)
   + \|F_A\|_{L^p} \\
&\le &
   2c_7 
   \left(\eps^{-1}\|\alpha\|_{L^p}
   + \eps^{-2}\|\xi\|_{L^p}\right)
   \|\zeta\|_{\infty,\eps}
   + \|F_A\|_{L^p} \\
&\le &
   2c_7 
   \delta\eps^{2/p}\eps^{-3/2}
   \|\zeta\|_{\infty,\eps}
   + \|F_A\|_{L^p} \\
&\le &
   c_{10}
   \delta \eps^{-3/2}
   \|\zeta\|_{1,p,\eps}
   + \|F_A\|_{L^p}.
\end{eqnarray*}
Here the penultimate inequality follows 
from~(\ref{eq:unique}) and the last follows from 
Lemma~\ref{lemma:p}.  Combining these two estimates
with~(\ref{eq:unique2}) we obtain
$$
      \|\zeta\|_{1,p,\eps}
      \le c_5 \left(
          2c_9(\delta+\eps^{1-2/p})
          + c_{10}\delta \eps^{1/2}
          \right)
          \|\zeta\|_{1,p,\eps}
      + c_5\eps^2\|F_A\|_{L^p}.
$$
If $\delta$ and $\eps$ are sufficiently small, we obtain
$$
     \|\zeta\|_{1,p,\eps}
     \le 2c_5\eps^2\|F_A\|_{L^p}
     \le 2c_0c_5\eps^2.
$$  
Hence the result follows from the uniqueness argument
at the end of the proof of Theorem~\ref{thm:exist}.
\end{proof}

Corollary~\ref{cor:uniqbaby} has a slightly 
stronger hypothesis than Theorem~\ref{thm:uniqbaby}, 
however, it does not seem to have a simpler proof.  
In order to significantly simplify the proof 
we would have to further strengthen the hypthesis
and assume 
$$
     \left\|\zeta\right\|_{1,p,\eps}
     \le\delta\eps^{2/p+1}
$$
with a small constant $\delta$
(instead of 
$
     \left\|\zeta\right\|_{1,p,\eps}
     \le\delta\eps^{2/p+1/2}
$
as in Corollary~\ref{cor:uniqbaby}). 
Under this hypothesis uniqueness can be established 
with the same straight forward argument that 
is used at the end of the proof of Theorem~\ref{thm:exist}.
However, such a weaker result just fails 
to suffice for the proof of Theorem~D.
Namely, in Section~\ref{sec:D} we shall establish 
an inequality of the form
$$
     \left\|\zeta\right\|_{1,p,\eps}
     \le c\eps^{2/p+1}
$$
under the hypotheses of Theorem~D.
In this inequality the constant $c$ 
is not small and so the argument in the 
proof of Theorem~\ref{thm:exist} does 
not suffice to give uniqueness. 
However, if $\eps$ is chosen so small that
$c\eps^{1/2}\le\delta$ then 
we can use Corollary~\ref{cor:uniqbaby}
to obtain uniqueness. 

 
\section{Relative Coulomb gauge}\label{sec:gauge}

This section is of preparatory nature. 
We prove a local slice theorem
for the action of the gauge group 
$\Gcal=\Gcal(P)$ on
$
    \Bcal=\Cinf_\G(P,M)\times\Acal(P).
$
The infinitesimal action is the operator
$
    d_{(u,A)}:\Om^0(\Sigma,\g_P)
    \to T_{(u,A)}\Bb
$
given by
$$
    d_{(u,A)}\eta := (-L_u\eta,d_A\eta).
$$
Denote by $d_{(u,A)}^{*_\eps}$
its formal adjoint with respect to the
$\eps$-inner pro\-duct, i.e. 
$$
    d_{(u,A)}^{*_\eps}(\xi,\alpha)
    := \eps^2 d_A^*\alpha - L_u^*\xi
$$
for $(\xi,\alpha)\in T_{(u,A)}\Bcal$.
The next proposition restates the local slice theorem 
for the $\Gg$-action on $\Bb$ with $\eps$-dependent 
norms for elements $(u_0,A_0)$ of the moduli space 
$\Tilde{\Mm}_{B,\Sigma}^0$. The result continues to hold 
for every element $(u_0,A_0)\in\Bb$ with $\mu(u_0)=0$. 
However, in this generality, more care must be taken in 
determining the norm on $\Bb$ with respect to which 
the constants $c$ and $\delta$ depend continuously
on $(u_0,A_0)$. In the case of $J$-holomorphic curves the 
$W^{1,p}$-norm controls all higher derivatives and 
therefore the choice of the norm is immaterial. 


\begin{proposition}\label{prop:slice}
Assume~$(H1)$.  For every $p>2$ and every $c_0>0$ 
there exist positive constants
$\delta$ and $c$ such that for every 
$\eps\in(0,1]$ the following holds.
Let $(u_0,A_0)\in\Tilde\Mcal_{B,\Sigma}^0(c_0)$
and 
$
    \zeta=(\xi,\alpha) \in T_{(u_0,A_0)}\Bb
$ 
such that 
\begin{equation}\label{eq:pre-slice}
    \left\|\zeta\right\|_{1,p,\eps;(u_0,A_0)}
    \le \delta \eps^{2/p}.
\end{equation}
Denote
$
    (u,A):=(\exp_{u_0}(\xi),A_0+\alpha).
$
Then there exist a unique pair
$
    \zeta_0=(\xi_0,\alpha_0) \in T_{(u_0,A_0)}\Bb
$
and a unique section $\eta_0\in\Om^0(\Sigma,\g_P)$
such that
$$
    d_{(u_0,A_0)}^{*_\eps}\zeta_0=0,\qquad
    g^*(u,A)=(\exp_{u_0}(\xi_0),A_0+\alpha_0),\qquad
    g:=e^{\eta_0},
$$
and
\begin{equation}\label{eq:propslice1}
   \|\eta_0\|_{2,p,\eps;A_0} 
   + \left\|\zeta_0\right\|_{1,p,\eps;(u_0,A_0)}
   \le c\left\|\zeta\right\|_{1,p,\eps;(u_0,A_0)}.
\end{equation}
\end{proposition}

Proposition~\ref{prop:slice} can be understood as a
quantitative version of the implicit function theorem
with $\eps$-dependent norms and constants independent
of $\eps$.  As in the case of Theorem~\ref{thm:exist} 
we shall prove it with a Newton type iteration.  
The relevant linear estimates are established in
Lemma~\ref{le:laplace} below, and the quadratic estimate
in Lemma~\ref{le:gauge}. 

\begin{lemma}\label{le:laplace}
For every $p\ge 2$ and every $c_0>0$ there exists a 
constant $c>0$ such that the following holds
for every
$
     (u_0,A_0)\in\Tilde{\Mm}_{B,\Sigma}^0(c_0)
$
and every $\eps\in(0,1]$. For every 
$
     \zeta=(\xi,\alpha)
     \in W^{1,p}(\Sigma,u_0^*TM/\G)\oplus T^*\Sigma\otimes\g_P)
$
there exists a unique $\eta\in W^{2,p}(\Sigma,\g_P)$ 
such that
\begin{equation}\label{eq:z-eta}
    d_{(u_0,A_0)}^{*_\eps} d_{(u_0,A_0)}\eta
    = d_{(u_0,A_0)}^{*_\eps}\zeta.
\end{equation}
Moreover, $\eta$ satisfies the estimates
\begin{equation}\label{eq:eta} 
    \left\|\eta\right\|_{1,p,\eps}
    \le c \left\|\zeta\right\|_{0,p,\eps},\qquad
    \left\|\eta\right\|_{2,p,\eps}
    \le c \left\|d_{(u_0,A_0)}^{*_\eps}\zeta\right\|_{L^p}.
\end{equation}
\end{lemma}

\begin{lemma}\label{le:gauge}
For every $p>2$ and every $c_0>0$ there exist
positive constants $\delta$ and $c$ such that 
the following holds for every
$
     (u_0,A_0)\in\Tilde{\Mm}_{B,\Sigma}^0(c_0)
$
and every $\eps\in(0,1]$. Assume that 
$
     \zeta_0=(\xi_0,\alpha_0)\in T_{(u_0,A_0)}\Bb
$
and $\eta\in\Om^0(\Sigma,\g_P)$ satisfy the inequality
\begin{equation}\label{eq:gauge}
     \left\|\eta\right\|_{2,p,\eps}
     + \left\|\zeta_0\right\|_{1,p,\eps}
     \le \delta\eps^{2/p}.
\end{equation}
Then there exists a unique pair
$\zeta_1=(\xi_1,\alpha_1)\in T_{(u_0,A_0)}\Bb$ 
such that
\begin{equation}\label{eq:zeta'}
     (\exp_{u_0}(\xi_1),A_0+\alpha_1) 
     = g^*(\exp_{u_0}(\xi_0),A_0+\alpha_0),
\end{equation}
where $g:=e^\eta$ and
\begin{equation}\label{eq:gauge1}
     \left\|\zeta_1-\zeta_0\right\|_{0,p,\eps}
     \le c \left\|\eta\right\|_{1,p,\eps},\qquad
     \left\|\zeta_1-\zeta_0\right\|_{1,p,\eps}
     \le  c \left\|\eta\right\|_{2,p,\eps}.
\end{equation}
Moreover,
\begin{equation}\label{eq:gauge2}
    \left\|d_{(u_0, A_0)}^{*_\eps}
     \left(\zeta_1 - \zeta_0 - d_{(u_0,A_0)}\eta\right)
    \right\|_{L^p} 
    \le c\eps^{-2/p}\left(
         \left\|\zeta_0\right\|_{1,p,\eps}
         + \left\|\eta\right\|_{2,p,\eps}
         \right)
    \left\|\eta\right\|_{1,p,\eps}.
\end{equation}
\end{lemma}

\begin{lemma}\label{le:eta}
For every $p\ge 2$ and every $c_0>0$ there exist positive
constant $\delta$ and $c$ such that the following holds
for every
$
     (u_0,A_0)\in\Tilde{\Mm}_{B,\Sigma}^0(c_0)
$
and every $\eps\in(0,1]$.
If $\eta_1,\eta_2\in\Om^0(\Sigma,\g_P)$ satisfy 
$\left\|\eta_1\right\|_{L^\infty}\le\delta$ 
and $\left\|\eta_2\right\|_{L^\infty}\le\delta$
then there exists a unique element $\eta\in\Om^0(\Sigma,\g_P)$ 
such that
$$
     e^{\eta} = e^{\eta_1}e^{\eta_2},\qquad
     2^{-1}\left\|\eta\right\|_{L^\infty}
     \le \left\|\eta_1+\eta_2\right\|_{L^\infty}
     \le 2\left\|\eta\right\|_{L^\infty}.
$$
Moreover, $\eta$ satisfies the estimate
$$
     c^{-1}\left\|\eta\right\|_{2,p,\eps;A_0}
     \le \left\|\eta_1+\eta_2\right\|_{2,p,\eps;A_0}
     \le c\left\|\eta\right\|_{2,p,\eps;A_0}. 
$$
\end{lemma}

\begin{proof}
For a fixed connection $A_0$ and $\eps=1$ 
the result is obvious.  Choose $c_1,c_2,c_3$
such that
$$
     c_k^{-1}\left\|\eta\right\|_{W^{k,p}}
     \le \left\|\eta_1+\eta_2\right\|_{W^{k,p}}
     \le c_k\left\|\eta\right\|_{W^{k,p}}
$$
for $k=1,2,3$, whenever $\eta_1,\eta_2,\eta$ are
sufficiently small in the $C^0$-norm and satisfy
$e^{\eta} = e^{\eta_1}e^{\eta_2}$.  Here 
the $W^{k,p}$ norms are understood with respect
to the connection $A_0$. It follows that
$$
     \eps\left\|d_{A_0}(\eta_1+\eta_2)\right\|_{L^p}
     \le \eps c_1\left(
         \left\|d_{A_0}\eta\right\|_{L^p}
         + \left\|\eta\right\|_{L^p}
         \right)
     \le c_1\left\|\eta\right\|_{1,p,\eps;A_0}
$$
and hence 
$$
     \left\|\eta_1+\eta_2\right\|_{1,p,\eps;A_0}
     \le (c_0+c_1)\left\|\eta\right\|_{1,p,\eps;A_0}
$$
for $0\le\eps\le 1$.  The other three inequalities 
follow by similar arguments.  This proves the lemma 
for a fixed connection $A_0$.  Moreover, the constant
$c$ depends continuously on $A_0$ with respect to
the $C^1$-norm, and is gauge invariant (with respect
to the action of $\Gg$ on $\Om^0(\Sigma,\g_P)$ by conjugation).
Hence, by Lemma~\ref{le:c0}, it can be chosen independent
of $A_0$ as long as $(u_0,A_0)\in\Tilde{\Mm}_{B,\Sigma}^0(c_0)$
for some $u_0$. 
\end{proof}

\begin{proof}[Proof of Lemma~\ref{le:laplace}]
The operator
$
     d_{(u_0,A_0)}^{*_{\eps}}d_{(u_0,A_0)}:
     W^{2,p}(\Sigma,\g_P)\to L^p(\Sigma,\g_P)
$
is given by
$$
     d_{(u_0,A_0)}^{*_\eps}d_{(u_0,A_0)}\eta
     = \eps^2d_{A_0}^*d_{A_0}\eta
       + L_{u_0}^*L_{u_0} \eta.
$$
By our standing hypotheses, $\mu^{-1}(0)$ is compact and 
$L_x:\g\to T_xM$ is injective for every $x\in\mu^{-1}(0)$.  
Hence there exists a constant $c_1>0$ such that
\begin{equation}\label{eq:Lx}
     c_1^{-1}\left|\eta\right|
     \le \left|L_x\eta\right|_z 
     \le c_1\left|\eta\right|
\end{equation}
for every $x\in\mu^{-1}(0)$, every $\eta\in\g$,
and every $z\in\Sigma$. (Here $|\cdot|_z$ denotes
the metric on $M$ induced by $J_z$ and $\om$.)
Hence the operator $d_{(u_0,A_0)}^{*_{\eps}}d_{(u_0,A_0)}$
is injective and hence, by elliptic regularity,
it is bijective. 

Next we prove that there exists a constant 
$c_2=c_2(p,c_0)>0$ such that, for every 
pair $(u_0,A_0)\in\Tilde{\Mm}_{B,\Sigma}^0(c_0)$,
every $\eta\in\Om^0(\Sigma,\g_P)$, and every
$\eps\in(0,1]$, we have 
\begin{equation}\label{eq:interpolate}
     \left\|d_{A_0}\eta\right\|_{L^p}
     \le \eps\left\|d_{A_0}^*d_{A_0}\eta\right\|_{L^p}
         + c_2\eps^{-1}\left\|\eta\right\|_{L^p}.
\end{equation}
For a fixed connection $A_0\in\Aa(P)$ this follows
directly from the interpolation inequality 
in~\cite[Theorem~7.27]{GT} and the $L^p$-estimate
for the operator $d_{A_0}^*d_{A_0}$.  Now the identity
\begin{eqnarray*}
     d_A^*d_A\eta - d_{A_0}^*d_{A_0}\eta
&= &
     [A-A_0\wedge d_{A_0}\eta]
     + *[*(A-A_0)\wedge d_{A_0}\eta] \\
&&
     - \, *[d_{A_0}*(A-A_0),\eta]
     + *[*(A-A_0)\wedge [A-A_0,\eta]]
\end{eqnarray*}
shows that the constant in~(\ref{eq:interpolate}) 
depends continuously on $A$ with respect to the 
$C^1$-norm. Moreover, the inequality~(\ref{eq:interpolate}) 
is gauge invariant. Hence it follows from Lemma~\ref{le:c0}
(with $\ell=2$) and the Arz\'ela-Ascoli theorem
that the estimate~(\ref{eq:interpolate}) holds with a uniform 
constant $c_2$ for all $A_0$ such that
$(u_0,A_0)\in\Tilde{\Mm}_{B,\Sigma}^0(c_0)$ for some $u_0$.

Using the identity
$$
     d\left|\eta\right|^{p-2}
     = (p-2)\left|\eta\right|^{p-4}\inner{\eta}{d_{A_0}\eta}
     \in\Om^1(\Sigma)
$$
for $\eta\in\Om^0(\Sigma,\g_P)$
and integration by parts we obtain
$$
    \int_\Sigma
    \left|\eta\right|^{p-2}\left|d_{A_0}\eta\right|^2
    = \int_\Sigma
      \left|\eta\right|^{p-2}\inner{\eta}{d_{A_0}^*d_{A_0}\eta}
    - (p-2)\int_\Sigma
      \left|\eta\right|^{p-4}
      \left|\inner{\eta}{d_{A_0}\eta}\right|^2.
$$
The last term on the right is negative. 
Now~(\ref{eq:z-eta}) is equivalent to
$$
     \eps^2d_{A_0}^*\alpha - L_{u_0}^*\xi
     = \eps^2d_{A_0}^*d_{A_0}\eta + L_{u_0}^*L_{u_0}\eta.
$$
Hence, by the previous identity and~(\ref{eq:Lx}), we have
\begin{eqnarray*}
&&
     \int_\Sigma
     \left(c_1^{-2}\left|\eta\right|^p
     + \eps^2\left|\eta\right|^{p-2}\left|d_{A_0}\eta\right|^2
     \right)  \\
&&\le 
     \int_\Sigma
     \left|\eta\right|^{p-2}
     \left(
     \left|L_{u_0}\eta\right|^2 
      + \eps^2\left|d_{A_0}\eta\right|^2
     \right) \\
&&\le 
     \int_\Sigma
     \left|\eta\right|^{p-2}
     \inner{\eta}{\eps^2d_{A_0}^*d_{A_0}\eta+L_{u_0}^*L_{u_0}\eta} \\
&&= 
     \int_\Sigma
     \left|\eta\right|^{p-2}
     \inner{\eta}{\eps^2 d_{A_0}^*\alpha-L_{u_0}^*\xi} \\
&&= 
     \int_\Sigma 
     \left|\eta\right|^{p-2}\left(
     \eps^2\inner{d_{A_0}\eta}{\alpha}
     + \eps^2(p-2)|\eta|^{-2}
       \inner{\inner{\eta}{d_{A_0}\eta}}
             {\inner{\eta}{\alpha}}
     - \inner{L_{u_0}\eta}{\xi}  
     \right) \\
&&\le 
     c_1\int_\Sigma 
     \left|\eta\right|^{p-1}\left|\xi\right| 
     + \eps^2(p-1)\int_\Sigma
     \left|\eta\right|^{p-2}
     \left|d_{A_0}\eta\right|\left|\alpha\right|  \\
&&\le 
     c_1\int_\Sigma 
     \left|\eta\right|^{p-1}\left|\xi\right| 
     + \frac{\eps^2(p-1)^2}{2}\int_\Sigma
     \left|\eta\right|^{p-2}\left|\alpha\right|^2
     + \frac{\eps^2}{2}\int_\Sigma
     \left|\eta\right|^{p-2}
     \left|d_{A_0}\eta\right|^2.
\end{eqnarray*}
Therefore, by H\"older's inequality,
$$
     c_1^{-2}\left\|\eta\right\|_{L^p}^p
     \le 
     c_1\left\|\eta\right\|_{L^p}^{p-1}
     \left\|\xi\right\|_{L^p}
     + \frac{\eps^2(p-1)^2}{2}
       \left\|\eta\right\|_{L^p}^{p-2}
       \left\|\alpha\right\|_{L^p}^2,
$$
$$
     c_1^{-2}\left\|\eta\right\|_{L^p}^p
     \le\left\|\eta\right\|_{L^p}^{p-1}
        \left\|\tilde\eta\right\|_{L^p},\qquad
     \tilde\eta
     := \eps^2d_{A_0}^*d_{A_0}\eta+L_{u_0}^*L_{u_0}\eta.
$$
Hence 
$$
     c_1^{-2}\left\|\eta\right\|_{L^p}^2
     \le 
     c_1\left\|\eta\right\|_{L^p}
     \left\|\xi\right\|_{L^p}
     + \frac{\eps^2(p-1)^2}{2}
       \left\|\alpha\right\|_{L^p}^2,
$$
and hence
$$
     \left\|\eta\right\|_{L^p}^2
     \le 
     c_1^6\left\|\xi\right\|_{L^p}^2
     + c_1^2\eps^2(p-1)^2
       \left\|\alpha\right\|_{L^p}^2.
$$
Thus we have proved the inequalities
\begin{equation}\label{eq:eta1}
    \left\|\eta\right\|_{L^p}
    \le c_1\max\{p-1,c_1^2\}
        \left\|\zeta\right\|_{0,p,\eps},\qquad
    \left\|\eta\right\|_{L^p}
    \le c_1^2\left\|\tilde\eta\right\|_{L^p}.
\end{equation}
By~(\ref{eq:Lx}), (\ref{eq:interpolate}),
and~(\ref{eq:eta1}), 
\begin{eqnarray*}
    \left\|\eta\right\|_{2,p,\eps;A_0} 
&= &
    \left\|\eta\right\|_{L^p}
      + \eps\left\|d_{A_0}\eta\right\|_{L^p}
      + \eps^2\left\|d_{A_0}^*d_{A_0}\eta\right\|_{L^p}  \\
&\le &
    (1+c_2)\left\|\eta\right\|_{L^p}
    + 2\eps^2\left\|d_{A_0}^*d_{A_0}\eta\right\|_{L^p} \\
&\le &
    (1+c_2+2c_1^2)\left\|\eta\right\|_{L^p}
    + 2\left\|\tilde\eta\right\|_{L^p} \\
&\le &
    \left(2+c_1^2(1+c_2+2c_1^2)\right)
    \left\|\tilde\eta\right\|_{L^p}.
\end{eqnarray*}
This proves the second estimate in~(\ref{eq:eta}).

\smallbreak

To prove the first estimate in~(\ref{eq:eta}) 
we use a rescaling argument in local holomorphic
coordinates on $\Sigma$. Cover $\Sigma$ by finitely many 
open sets, each of which is holomorphically
diffeomorphic to the unit square in $\C$,
suppose that the coordinate charts extend 
to a closed square of side length two,
and choose trivializations of the bundle $P$ 
over each of these (extended) open sets. 
In these coordinates we write the
metric in the form $\lambda^2(ds^2+dt^2)$, 
and we write $A:=A_0=\Phi\,ds+\Psi\,dt$,
$\alpha=\phi\,ds+\psi\,dt$.  
Moreover, $u:=u_0:[0,2]^2\to M$, $\xi:[0,2]^2\to TM$
is a vector field along $u$, and $\eta:[0,2]^2\to\g$. 
In this notation equation~(\ref{eq:z-eta}) has the form 
\begin{equation}\label{eq:z-eta-loc}
    \Nabla{s}\Nabla{s}\eta + \Nabla{t}\Nabla{t}\eta
    = \Nabla{s}\phi + \Nabla{t}\psi
      + (\lambda/\eps)^2 L_u^*(L_u\eta+\xi),
\end{equation}
where 
$
     \Nabla{s}\eta:=\Nabla{A,s}^\G\eta=\p_s\eta+[\Phi,\eta]
$
and 
$
     \Nabla{t}\eta:=\Nabla{A,t}^\G\eta=\p_t\eta+[\Psi,\eta].
$ 
Now we introduce new functions, defined on 
the square $[0,2/\eps]^2$, by 
$$
\begin{array}{rclcrcl}
   \tilde\eta(s,t) &:= &\eta(\eps s,\eps t),&\qquad &
   \tilde\lambda(s,t) &:= &\lambda(\eps s,\eps t), \\
   \tilde\xi(s,t) &:= &\xi(\eps s,\eps t),&\qquad &
   \tilde u(s,t) &:= &u(\eps s,\eps t), \\
   \tilde\phi(s,t) &:= &\eps\phi(\eps s,\eps t),&\qquad &
   \tilde\Phi(s,t) &:= &\eps\Phi(\eps s,\eps t), \\
   \tilde\psi(s,t) &:= &\eps\psi(\eps s,\eps t),&\qquad &
   \tilde\Psi(s,t) &:= &\eps\Psi(\eps s,\eps t).
\end{array}
$$
Then~(\ref{eq:z-eta-loc}) is equivalent to
$$
     \Tabla{s}\Tabla{s}\tilde\eta + \Tabla{t}\Tabla{t}\tilde\eta
     = \Tabla{s}\tilde\phi + \Tabla{t}\tilde\psi
       + \tilde\lambda^2 L_{\tilde u}^*
         (L_{\tilde u}\tilde\eta+\tilde\xi),
$$
where
$
     \Tabla{s}\tilde\eta:=\p_s\tilde\eta+[\tilde\Phi,\tilde\eta]
$
and 
$
     \Tabla{t}\tilde\eta:=\p_t\tilde\eta+[\tilde\Psi,\tilde\eta].
$ 
This equation can be written in the form
$$
     \Delta\tilde\eta 
     = \p_s\tilde f + \p_t\tilde g + \tilde h
$$
where $\Delta:=\p_s\p_s+\p_t\p_t$ and
$\tilde f,\tilde g,\tilde h:[0,2/\eps]^2\to\g$
are given by 
$$
     \tilde f:= \tilde\phi-2[\tilde\Phi,\tilde\eta],\qquad
     \tilde g:= \tilde\psi-2[\tilde\Psi,\tilde\eta],
$$
$$
     \tilde h
     := \tilde\lambda^2 L_{\tilde u}^*
        (L_{\tilde u}\tilde\eta+\tilde\xi)
        + [\tilde\Phi,\tilde\phi-[\tilde\Phi,\tilde\eta]]
        + [\tilde\Psi,\tilde\psi-[\tilde\Psi,\tilde\eta]]
        + [\p_s\tilde\Phi+\p_t\tilde\Psi,\tilde\eta].
$$
Hence there exists a constant $c_3>0$ such that,
for all real numbers $a,b$ such that 
$1/2\le a<b\le 2/\eps-1/2$, we have
$$
     \int_{[a,b]^2}\left(
     |\Tabla{s}\tilde\eta|^p
     + |\Tabla{t}\tilde\eta|^p
     \right)
     \le c_3
     \int_{[a-1/2,b+1/2]^2}\left(
     |\tilde f|^p
     + |\tilde g|^p
     + |\tilde h|^p
     + |\tilde\eta|^p
     \right).
$$
Here the constant $c_3$ is independent of $a$ and $b$. 
It follows that
\begin{eqnarray*}
     \int_{[a,b]^2}\left(
     |\Tabla{s}\tilde\eta|^p
     + |\Tabla{t}\tilde\eta|^p
     \right)\tilde\lambda^{2-p}
&\le &
     c_4\int_{[a-1/2,b+1/2]^2}
     \left(|\tilde \phi|^p
     + |\tilde \psi|^p\right)\tilde\lambda^{2-p} \\
&&
     +\,c_4\int_{[a-1/2,b+1/2]^2}
     \left(|\tilde\xi|^p
     + |\tilde\eta|^p\right)\tilde\lambda^2,
\end{eqnarray*}
where the constant $c_4$ depends on the metric and 
on the $C^1$-norms of~$\tilde\Phi$ and~$\tilde\Psi$. 
With $a=1/2\eps$, $b=3/2\eps$, and $0<\eps\le1$
we obtain
\begin{eqnarray*}
     \eps^p\int_{[1/2,3/2]^2}\left(
     |\Nabla{s}\eta|^p
     + |\Nabla{t}\eta|^p
     \right)\lambda^{2-p}
&\le &
     c_4\eps^p\int_{[0,2]^2}
     \left(|\phi|^p
     + |\psi|^p\right)\lambda^{2-p} \\
&&
     +\,c_4\int_{[0,2]^2}
     \left(|\xi|^p + |\eta|^p\right)\lambda^2.
\end{eqnarray*}
Hence, by taking the sum over the coordinate charts,
$$
     \eps\left\|d_{A_0}\eta\right\|_{L^p}
     \le N^{1/p}c_4\left(
     \eps\left\|\alpha\right\|_{L^p}
     + \left\|\xi\right\|_{L^p}
     + \left\|\eta\right\|_{L^p}
     \right).
$$
Here $N$ is the number of open sets in the cover
and the constant $c_4$ depends continuously
on $A_0$ with respect to the $C^1$-norm.
Hence, by Lemma~\ref{le:c0}, $c_4$ 
can be chosen independent of the pair 
$(u_0,A_0)\in\Tilde{\Mm}_{B,\Sigma}^0(c_0)$.
Combining the last inequality with~(\ref{eq:eta1})
we obtain the first estimate in~(\ref{eq:eta})
as claimed. 
\end{proof}

In the following proof we use the identity
\begin{equation}\label{eq:series}
      (e^\eta)^*A-A-d_A\eta
      = \sum_{k=1}^\infty\frac{(-1)^k}{(k+1)!} 
        \ad(\eta)^k d_A\eta
\end{equation}
for $A\in\Aa(P)$ and $\eta\in\Om^0(\Sigma,\g_P)$, where
$\ad(\eta)\alpha:=[\eta,\alpha]$ for $\alpha\in\Om^1(\Sigma,\g_P)$. 

\begin{proof}[Proof of Lemma~\ref{le:gauge}]
Throughout the proof we denote by $c_1,c_2,c_3,\dots$ 
positive constants depending only on $p$ and $c_0$
(and not on the pair $(u_0,A_0)$).  
Fix a pair $(u_0,A_0)\in\Tilde{\Mm}_{B,\Sigma}^0(c_0)$
and choose a positive constant $\delta_0$ that is smaller
than the injectivity radius of $M$ on the compact set $u_0(P)$.  
Suppose that $\eta\in\Om^0(\Sigma,\g_P)$ and 
$\zeta_0=(\xi_0,\alpha_0)\in T_{(u_0,A_0)}\Bb$ satisfy 
the hypotheses of Lemma~\ref{le:gauge} with a sufficiently 
small constant $\delta>0$. 
Let $c_1$ be the constant of Lemma~\ref{lemma:p}. 
Then, by~(\ref{eq:gauge}),
$$
     \left\|\xi_0\right\|_{L^\infty}
     + \left\|\eta\right\|_{L^\infty} 
     \le c_1\eps^{-2/p}\left(\left\|\xi_0\right\|_{1,p,\eps} 
         + \left\|\eta\right\|_{2,p,\eps}\right)
     \le c_1\delta.
$$
Hence, if $\delta$ is sufficiently small, it follows
that the $C^0$-distance between $u_0$ and 
$e^{-r\eta}\exp_{u_0}(\xi_0)$ is smaller than
$\delta_0$ for every $r\in[0,1]$.  Hence there 
exists a unique smooth path 
$
     [0,1]\to T_{(u_0,A_0)}\Bb:
     r\mapsto \zeta_r = (\xi_r,\alpha_r)
$
starting at $\zeta_0$ such that 
$$
     (\exp_{u_0}(\xi_r),A_0+\alpha_r) 
     = g^*(u,A),
$$
where
\begin{equation}\label{eq:Z}
     (u,A) := (\exp_{u_0}(\xi_0),A_0+\alpha_0),\qquad
     g := e^{r\eta}.
\end{equation}
The endpoint $\zeta_1$ of this
path obviously satisfies~(\ref{eq:zeta'}).
We prove the inequalities
\begin{equation}\label{eq:Z1}
     \left\|\p_r\zeta_r\right\|_{0,p,\eps}
     \le c\left\|\eta\right\|_{1,p,\eps},\qquad  
     \left\|\p_r\zeta_r\right\|_{1,p,\eps}
     \le c\left\|\eta\right\|_{2,p,\eps},
\end{equation}
\begin{equation}\label{eq:Z2}
     \left\|
     d_{(u_0,A_0)}^{*_\eps}
     \left(\p_r\zeta_r
     - d_{(u_0,A_0)}\eta\right)
     \right\|_{L^p}
     \le c\eps^{-2/p}\left(
         \left\|\zeta_0\right\|_{1,p,\eps}
         + \left\|\eta\right\|_{2,p,\eps}
         \right)
       \left\|\eta\right\|_{1,p,\eps}
\end{equation}
for $0\le r\le 1$, 
where the constant $c$ depends only on $c_0$ and $p$.  
Then the inequalities~(\ref{eq:gauge1}) and (\ref{eq:gauge2})
follow by integrating the function 
$r\mapsto\p_r\zeta_r$ over the interval $0\le r\le 1$. 

For every $u\in\Cinf_\G(P,M)$ whose $C^0$-distance to $u_0$
is less than $\delta_0$ we define the linear operator
$
     Z(u):\Om^0(\Sigma,\g_P)
     \to \Om^0(\Sigma,u_0^*TM/\G)
$
by 
$$
     Z(u)\hat\eta := d\exp_{u_0}^{-1}(u)L_u\hat\eta
$$
for $\hat\eta\in\Om^0(\Sigma,\g_P)$. Then 
\begin{equation}\label{eq:Zdot}
     \p_r\zeta_r
     = (-Z(g^{-1}u)\eta,d_{g^*A}\eta_0),
\end{equation}
where $(u,A)$ and $g$ are as in~(\ref{eq:Z}).
We prove the first inequality in~(\ref{eq:Z1}).  
Since $A=A_0+\alpha_0$ we have
\begin{equation}\label{eq4.20}
    d_{g^*A}\eta
    = d_{A_0}\eta
      + [g^*A_0-A_0,\eta]
      + [g^{-1}\alpha_0g,\eta]
\end{equation}
and we must estimate the three terms on the right 
with $g:=e^{r\eta}$. Since
$
     \left\|\eta\right\|_{L^\infty}\le c_1\delta
$
it follows from~(\ref{eq:series}), with $\eta$ 
replaced by $r\eta$, that
$$
     \left\|[g^*A_0-A_0,\eta]\right\|_{L^p}
     \le c_2\left\|d_{A_0}\eta\right\|_{L^p}
            \left\|\eta\right\|_{L^\infty}
     \le c_3 \eps^{-1-2/p} \left\|\eta\right\|_{1,p,\eps}^2.
$$
Moreover,
$$
     \left\|[g^{-1}\alpha_0g,\eta]\right\|_{L^p}
     \le \eps^{-1-2/p}
     \left\|\zeta_0\right\|_{0,p,\eps}
     \left\|\eta\right\|_{1,p,\eps}.
$$
Hence, by~(\ref{eq:gauge}) and~(\ref{eq4.20}), 
$
   \left\|d_{g^*A}\eta\right\|_{L^p}
   \le c_4\eps^{-1}\left\|\eta\right\|_{1,p,\eps}
$
and hence the first inequality in~(\ref{eq:Z1})
follows from~(\ref{eq:Zdot}). 

Next we prove the second inequality in~(\ref{eq:Z1}) 
and~(\ref{eq:Z2}). Using the identity
$$
    d_{A_0}[g^{-1}\alpha_0g,\eta]
    = [g^{-1}(d_{A_0}\alpha_0)g,\eta]
      + [[(A_0-g^*A_0)\wedge g^{-1}\alpha_0g],\eta]
      - [g^{-1}\alpha_0g\wedge d_{A_0}\eta]
$$
and~(\ref{eq:series}) we obtain
\begin{eqnarray*}
     \left\|[g^{-1}\alpha_0g,\eta]\right\|_{1,p,\eps}
&= &
     \left\|[g^{-1}\alpha_0g,\eta]\right\|_{L^p}  \\
&&
     +\,\eps\left\|d_{A_0}[g^{-1}\alpha_0g,\eta]\right\|_{L^p}
     + \eps\left\|d_{A_0}[*g^{-1}\alpha_0g,\eta]\right\|_{L^p} \\
&\le &
     c_5\eps^{-1-2/p}\left\|\zeta_0\right\|_{1,p,\eps}
     \left\|\eta\right\|_{1,p,\eps}.
\end{eqnarray*}
Similarly, using the identity
$$
    d_{A_0}[(A_0-g^*A_0),\eta]
    = [d_{A_0}(A_0-g^*A_0),\eta]
      - [(A_0-g^*A_0)\wedge d_{A_0}\eta]
$$
and~(\ref{eq:series}) we obtain
$$
    \left\|[(A_0-g^*A_0),\eta]\right\|_{1,p,\eps}
    \le c_6\eps^{-1-2/p}\left\|\eta\right\|_{2,p,\eps}
        \left\|\eta\right\|_{1,p,\eps}.
$$
Hence, by~(\ref{eq4.20}),
$$
    \left\|d_{g^*A}\eta-d_{A_0}\eta\right\|_{1,p,\eps}
    \le c_7\eps^{-1-2/p}\left(
          \left\|\zeta_0\right\|_{1,p,\eps}
          + \left\|\eta\right\|_{2,p,\eps}
          \right)
          \left\|\eta\right\|_{1,p,\eps}.
$$
Moreover, since $Z(u_0)=L_{u_0}$, we have
\begin{eqnarray*}
&&
     \left\|
     Z(g^{-1}u)\eta - L_{u_0}\eta
     \right\|_{1,p,\eps}  \\
&&\le 
     \left\|Z(g^{-1}u)\eta-Z(u)\eta\right\|_{1,p,\eps}  
     + \left\|Z(u)\eta-Z(u_0)\eta\right\|_{1,p,\eps}  \\
&&\le 
     c_8\left(
     \left\|\eta\right\|_{1,p,\eps}\left\|\eta\right\|_{L^\infty}
     + \left\|\xi_0\right\|_{1,p,\eps}\left\|\eta\right\|_{L^\infty}
     + \left\|\xi_0\right\|_{L^\infty}\left\|\eta\right\|_{1,p,\eps}
     \right)  \\
&&\le 
     c_9\eps^{-2/p}\left(
     \left\|\zeta_0\right\|_{1,p,\eps}
       + \left\|\eta\right\|_{2,p,\eps}
     \right)
     \left\|\eta\right\|_{1,p,\eps}.
\end{eqnarray*}
Here we have used the inequality 
$\left\|\eta\right\|_{L^\infty}\le c_1\delta$
from~(\ref{eq:gauge}) and Lemma~\ref{lemma:p}.  
The constants $c_7$ and $c_9$
in the last two estimates depend continuously on the pair
$(u_0,A_0)$ with respect to the $C^1$-norm and 
are gauge invariant.  Hence, by Lemma~\ref{le:c0}, 
they can be chosen independent of 
$(u_0,A_0)\in\Tilde{\Mm}_{B,\Sigma}^0(c_0)$. 
Hence the second inequality in~(\ref{eq:Z1})
follows from~(\ref{eq:Zdot}). 
To prove~(\ref{eq:Z2}) we observe that
\begin{eqnarray*}
&&
     d_{(u_0,A_0)}^{*_\eps}\p_r\zeta_r
     - d_{(u_0,A_0)}^{*_\eps}d_{(u_0,A_0)}\eta  \\
&&= 
     \eps^2d_{A_0}^*\left(d_{g^*A}\eta-d_{A_0}\eta\right)
     + L_{u_0}^*\left(Z(g^{-1}u)\eta - L_{u_0}\eta\right),
\end{eqnarray*}
where $(u,A)$ and $g:=e^{r\eta}$ are as in~(\ref{eq:Z}). 
The terms on the right have been estimated 
above and this proves~(\ref{eq:Z2}). 
Thus we have proved the existence of $\zeta_1$.  
The inequality~(\ref{eq:gauge1}) 
with $\delta$ sufficiently small guarantees 
that the $C^0$ distance between $u_0$ and 
$\exp_{u_0}(\xi_1)$ is smaller than the injectivity radius.
This proves uniqueness. 
\end{proof}

\begin{proof}[Proof of Proposition~\ref{prop:slice}]
The proof is based on a Newton type iteration. Let 
$$
      (u_1,A_1):=(u,A)=(\exp_{u_0}(\xi),A_0+\alpha),\qquad
      \zeta_1:=\zeta.
$$
For $\nu\ge 2$ we define 
$\zeta_\nu=(\xi_\nu,\alpha_\nu)\in T_{(u_0,A_0)}\Bb$
inductively by 
$$
     (\exp_{u_0}(\xi_{\nu+1}),A_0+\alpha_{\nu+1})
     := (u_{\nu+1},A_{\nu+1})
     := g_\nu^*(u_\nu,A_\nu),
$$
where
$
     (u_\nu,A_\nu) 
     := (\exp_{u_0}(\xi_\nu),A_0+\alpha_\nu),
$ 
$g_\nu:= e^{\hat\eta_\nu}$,
and $\hat\eta_\nu\in\Om^0(\Sigma,\g_P)$ 
is the unique solution of the equation
$$
     d_{(u_0,A_0)}^{*_\eps}
     d_{(u_0,A_0)}\hat\eta_\nu
     + d_{(u_0,A_0)}^{*_\eps}\zeta_\nu
     = 0.
$$
To construct these sequences we must ensure that in each
step $\zeta_\nu$ and $\hat\eta_\nu$ satisfy the hypotheses
of Lemma~\ref{le:gauge} so that $\zeta_{\nu+1}$ can be chosen
as in  the assertion of Lemma~\ref{le:gauge}. 
We shall prove this below.  And we shall also prove
that these sequences satisfy the following estimates.
\begin{eqnarray}
     \left\|\zeta_\nu\right\|_{1,p,\eps}
&\le & 
     C\left\|\zeta\right\|_{1,p,\eps},  
     \label{eq6.3} \\     
     \left\|d_{(u_0,A_0)}^{*_\eps}\zeta_\nu\right\|_{L^p}
&\le &
     C\eps^{-2/p}
     \left\|\zeta_{\nu-1}\right\|_{1,p,\eps}
     \left\|d_{(u_0,A_0)}^{*_\eps}\zeta_{\nu-1}\right\|_{L^p},  
     \label{eq6.4} \\
     \left\|d_{(u_0,A_0)}^{*_\eps}\zeta_\nu\right\|_{L^p}
&\le & 
     2^{1-\nu}
     \left\|d_{(u_0,A_0)}^{*_{\eps}}\zeta\right\|_{L^p},
     \label{eq6.5} \\
     \left\|\hat\eta_\nu\right\|_{2,p,\eps}
&\le &
     C2^{-\nu}\left\|\zeta\right\|_{1,p,\eps},
     \label{eq6.6}
\end{eqnarray}
The constants $C$ and $\delta$ are chosen as follows. 
Suppose that the constants $c_1,c_2,c_3,c_4,c_5\ge 1$ 
and $\delta_0,\delta_3,\delta_4\in(0,1]$ 
satisfy the following conditions.
\begin{itemize}
\item
The injectivity radius of $M$ on $u_0(P)$ is bigger than
$\delta_0$.
\item
The inequality~(\ref{eq:Lx}) holds with $c_1$ 
for every $x\in\mu^{-1}(0)$.
\item
The assertion of Lemma~\ref{le:laplace}
holds for $0<\eps\le 1$ with $c$ replaced by $c_2$.
\item
The assertion of Lemma~\ref{le:gauge}
holds for $0<\eps\le 1$ with $c$ replaced by $c_3$
and $\delta$ replaced by $\delta_3$.
\item
The assertion of Lemma~\ref{le:eta} 
holds for $0<\eps\le 1$ with $c$ replaced by $c_4$
and $\delta$ replaced by $\delta_4$. 
\item
The assertion of Lemma~\ref{lemma:p}
holds for $0<\eps\le 1$ with $c$ replaced by $c_5$.
\end{itemize}
Now choose positive constants $C$ and $\delta$ such that 
$$
     c_1c_2(1+2c_2)c_3\le C,\quad 
     2c_2(1+c_2)c_3C\delta\le 1,\quad
     2C\delta \le\delta_3,\quad
     4c_4c_5C\delta \le\delta_4. 
$$
We prove that the estimates~(\ref{eq6.3}-\ref{eq6.6}) 
hold for $\nu=1$. Since $C/2\ge c_1c_2$, the inequality~(\ref{eq6.6})
with $\nu=1$ follows from Lemma~\ref{le:laplace}. 
Since $C\ge 1$, the inequality~(\ref{eq6.3}) 
holds for $\nu=1$. The inequality~(\ref{eq6.4}) 
is vacuous for $\nu=1$ and~(\ref{eq6.5}) is obvious. 

Now suppose that the sequences have been constructed
and the inequalities~(\ref{eq6.3}-\ref{eq6.6})
have been established up to some integer $\nu\ge1$.  
Then
$$
     \left\|\hat\eta_\nu\right\|_{2,p,\eps}
     + \left\|\zeta_\nu\right\|_{1,p,\eps}
     \le C(1+2^{-\nu})
         \left\|\zeta\right\|_{1,p,\eps}
     \le 2C\delta\eps^{2/p}
     \le\delta_3\eps^{2/p}.
$$
Hence the hypotheses of Lemma~\ref{le:gauge}
are satisfied with $\zeta_0$ replaced by $\zeta_\nu$ 
and $\eta$ replaced by $\hat\eta_\nu$.  Choose 
$\zeta_{\nu+1}=(\xi_{\nu+1},\alpha_{\nu+1})$
as in the assertion of Lemma~\ref{le:gauge}. 
By Lemma~\ref{le:laplace}, we have
$$
     \left\|\hat\eta_\nu\right\|_{1,p,\eps}
     \le c_2\left\|\zeta_\nu\right\|_{0,p,\eps},\qquad
     \left\|\hat\eta_\nu\right\|_{2,p,\eps}
     \le c_2\left\|d_{(u_0,A_0)}^{*_\eps}\zeta_\nu\right\|_{L^p}.
$$
Moreover, $d_{(u_0,A_0)}^{*_\eps}(\zeta_\nu+d_{(u_0,A_0)}\eta_\nu)=0$,
and hence, by~(\ref{eq:gauge2}),
\begin{eqnarray*}
     \left\|d_{(u_0,A_0)}^{*_\eps}\zeta_{\nu+1}\right\|_{L^p}
&\le &
     c_3\eps^{-2/p}\left(
         \left\|\zeta_\nu\right\|_{1,p,\eps}
         + \left\|\hat\eta_\nu\right\|_{2,p,\eps}
         \right)
     \left\|\hat\eta_\nu\right\|_{1,p,\eps}  \\
&\le &
     c_2(1+c_2)c_3\eps^{-2/p}\left\|\zeta_\nu\right\|_{1,p,\eps}
     \left\|d_{(u_0,A_0)}^{*_\eps}\zeta_\nu\right\|_{L^p}.
\end{eqnarray*}
Since $c_2(1+c_2)c_3\le C$, this proves~(\ref{eq6.4})
with $\nu$ replaced by $\nu+1$. Moreover, by~(\ref{eq6.3}),
\begin{eqnarray*}
     \left\|d_{(u_0,A_0)}^{*_\eps}\zeta_{\nu+1}\right\|_{L^p}
&\le &
     c_2(1+c_2)c_3C\eps^{-2/p}\left\|\zeta\right\|_{1,p,\eps}
     \left\|d_{(u_0,A_0)}^{*_\eps}\zeta_\nu\right\|_{L^p}  \\
&\le &
     c_2(1+c_2)c_3C\delta
     \left\|d_{(u_0,A_0)}^{*_\eps}\zeta_\nu\right\|_{L^p}.
\end{eqnarray*}
Since $2c_2(1+c_2)c_3C\delta\le 1$, this proves~(\ref{eq6.5}) 
with $\nu$ replaced by $\nu+1$. Now let $\hat\eta_{\nu+1}$
be the unique solution of
$
     d_{(u_0,A_0)}^*d_{(u_0,A_0)}\hat\eta_{\nu+1}
     + d_{(u_0,A_0)}^*\zeta_{\nu+1} = 0.
$
Then, by Lemma~\ref{le:laplace} and~(\ref{eq6.5}), 
$$
     \left\|\hat\eta_{\nu+1}\right\|_{2,p,\eps}
     \le c_2
         \left\|d_{(u_0,A_0)}^{*_\eps}\zeta_{\nu+1}\right\|_{L^p}
     \le c_22^{-\nu}
         \left\|d_{(u_0,A_0)}^{*_\eps}\zeta\right\|_{L^p}
     \le c_1c_22^{-\nu}\left\|\zeta\right\|_{1,p,\eps}.
$$
Since $2c_1c_2\le C$ this implies~(\ref{eq6.6}) 
with $\nu$ replaced by $\nu+1$.  
It remains to prove~(\ref{eq6.3}) with $\nu$ replaced by $\nu+1$.
By~(\ref{eq:gauge1}) and~(\ref{eq6.5}), we have
\begin{equation}\label{eq:cauchy}
     \left\|\zeta_{j+1}-\zeta_j\right\|_{1,p,\eps}
\le 
     c_3\left\|\hat\eta_j\right\|_{2,p,\eps}
\le 
     c_2c_32^{1-j}\left\|d_{(u_0,A_0)}^{*_\eps}\zeta\right\|_{L^p}
\le 
     c_1c_2c_32^{1-j}\left\|\zeta\right\|_{1,p,\eps}
\end{equation}
for $j=1,\dots,\nu$.  Hence
$$
     \left\|\zeta_{\nu+1}\right\|_{1,p,\eps} 
     \le \left\|\zeta\right\|_{1,p,\eps} 
         + \sum_{j=1}^\nu
         \left\|\zeta_{j+1}-\zeta_j\right\|_{1,p,\eps} 
     \le (1+2c_1c_2c_3)\left\|\zeta\right\|_{1,p,\eps}.
$$
Since $1+2c_1c_2c_3\le C$ this proves~(\ref{eq6.3}) 
with $\nu$ replaced by $\nu+1$.  This completes
the induction. 

By~(\ref{eq:cauchy}), $\zeta_\nu$ is a Cauchy sequence 
in the $W^{1,p}$-norm.  Moreover, 
$$
     (u_\nu,A_\nu)
     = \bigl(\exp_{u_0}(\xi_\nu),A_0+\alpha_\nu\bigr)
     = h_\nu^*(u,A),
$$
where $h_\nu:=g_1g_2\cdots g_{\nu-1}$. 
We prove by induction that there exists a sequence 
$\eta_\nu\in\Om^0(\Sigma,\g_P)$ such that
\begin{equation}\label{eq:h}
     h_\nu = e^{\eta_\nu},\qquad
     \left\|\eta_{\nu+1}-\eta_{\nu}\right\|_{2,p,\eps}
     \le c_4C2^{-\nu}\left\|\zeta\right\|_{1,p,\eps}.
\end{equation}
For $\nu=1$ we set $h_1:=\one$ and $\eta_1:=0$.
Suppose that the sequence has been constructed 
for all integers up to $\nu\ge1$. Then
\begin{equation}\label{eq:eta'}
     \left\|\eta_\nu\right\|_{2,p,\eps}
     \le \sum_{j=1}^{\nu-1}
         \left\|\eta_{j+1}-\eta_j\right\|_{2,p,\eps}
     \le c_4C\left\|\zeta\right\|_{1,p,\eps}.
\end{equation}
Hence, by Lemma~\ref{lemma:p}, (\ref{eq:pre-slice}), and~(\ref{eq6.6}), 
$$
      \left\|\eta_\nu\right\|_{L^\infty}
      \le c_5\eps^{-2/p}\left\|\eta_\nu\right\|_{2,p,\eps}
      \le c_4c_5C\eps^{-2/p}\left\|\zeta\right\|_{1,p,\eps}
      \le c_4c_5C\delta
      \le \delta_4/4,
$$
$$
     \left\|\hat\eta_\nu\right\|_{L^\infty}
     \le c_5\eps^{-2/p}\left\|\hat\eta_\nu\right\|_{2,p,\eps}
     \le c_5C\eps^{-2/p}\left\|\zeta\right\|_{1,p,\eps}
     \le c_5C\delta
     \le \delta_4/4.
$$
By Lemma~\ref{le:eta}, there exists a section
$\eta_{\nu+1}\in\Om^0(\Sigma,\g_P)$ such that
$$
     e^{\eta_{\nu+1}} 
     = e^{\eta_\nu}e^{\hat\eta_\nu}
     = h_\nu g_\nu
     = h_{\nu+1},\qquad
     \left\|\eta_{\nu+1}\right\|_{L^\infty}
     \le 2\left\|\eta_\nu+\hat\eta_\nu\right\|_{L^\infty}
     \le \delta_4.
$$
Applying Lemma~\ref{le:eta} to 
$-\eta_\nu$ and $\eta_{\nu+1}$ we find
$$
     \left\|\eta_{\nu+1}-\eta_\nu\right\|_{2,p,\eps}
     \le c_4\left\|\hat\eta_\nu\right\|_{2,p,\eps}
     \le c_4C2^{-\nu}\left\|\zeta\right\|_{1,p,\eps}.
$$
The last inequality follows from~(\ref{eq6.6}). 
This completes the induction.
Thus we have proved that $h_\nu$ satisfies~(\ref{eq:h})
and hence is a Cauchy sequence in $\Gg^{2,p}(P)$. 
Denote 
$$
     \zeta := \lim_{\nu\to\infty}\zeta_\nu,\qquad
     h := \lim_{\nu\to\infty}h_\nu,\qquad
     \eta := \lim_{\nu\to\infty}\eta_\nu.
$$
Then
$$
     e^\eta=h,\qquad
     h^*(u,A)=(\exp_{u_0}(\xi),A_0+\alpha),\qquad
     d_{(u_0,A_0)}^{*_\eps}\zeta=0.
$$
The last equation follows from~(\ref{eq6.5}).
Moreover, by~(\ref{eq6.3}) and~(\ref{eq:eta'}), 
we have
$
     \left\|\eta\right\|_{2,p,\eps}
     + \left\|\zeta\right\|_{1,p,\eps}
     \le C(1+c_4)\left\|\zeta\right\|_{1,p,\eps}.
$
Hence~(\ref{eq:propslice1}) holds with $c:=C(1+c_4)$. 

To complete the existence proof we must show that 
$\eta$ and $\zeta$ are smooth.  We shall prove that 
the sequence $\zeta_\nu$ is bounded on $W^{k,p}$ 
for every $k$.  Here it suffices to obtain  rather
crude estimates with constants which depend on $\eps$ 
and are allowed to diverge as $\eps$ tends to zero. 
We fix a constant $\eps>0$ and prove by induction
that for every integer $k\ge 1$ there exists a constant 
$c_k=c_k(p,\eps,u_0,A_0,u,A)$ such that, for every $\nu$,
\begin{equation}\label{eq:k-1}
     \left\|\zeta_\nu\right\|_{W^{k,p}}\le c_k,\qquad
     \left\|\hat\eta_\nu\right\|_{W^{k+1,p}}\le c_k2^{-\nu}.
\end{equation}
For $k=1$ this follows from~(\ref{eq6.3}) and~(\ref{eq6.6}).
Now let $k\ge 2$ and assume that these estimates have 
been established with $k$ replaced by $k-1$. 
Observe that there exists a constant $C_k\ge1$ 
such that, for every $\nu$, 
\begin{eqnarray*}
     \left\|\zeta_{\nu+1}-\zeta_\nu\right\|_{W^{k,p}}
&\le & 
     C_k\left\|\hat\eta_\nu\right\|_{W^{k+1,p}},  \\
     \left\|\hat\eta_\nu\right\|_{W^{k+1,p}}
&\le & 
     C_k\left\|d_{(u_0,A_0)}^{*_\eps}\zeta_\nu\right\|_{W^{k-1,p}}, \\
     \left\|d_{(u_0,A_0)}^{*_\eps}\zeta_{\nu+1}\right\|_{W^{k-1,p}}
&\le & 
     C_k\left(
     \left\|\zeta_\nu\right\|_{W^{k,p}}
      + \left\|\hat\eta_\nu\right\|_{W^{k+1,p}}
     \right)
     \left\|\hat\eta_\nu\right\|_{W^{k,p}}.
\end{eqnarray*}
The first two inequalities are obvious, and the last follows
by inspecting the formula~(\ref{eq:Zdot}) in the proof 
of Lemma~\ref{le:gauge}.  Combining these inequalities
with the induction hypothesis, we obtain 
\begin{eqnarray*}
     \left\|\zeta_{\nu+1}\right\|_{W^{k,p}}
&\le &
     \left\|\zeta_\nu\right\|_{W^{k,p}}
         + C_k\left\|\hat\eta_\nu\right\|_{W^{k+1,p}},  \\
     \left\|\hat\eta_{\nu+1}\right\|_{W^{k+1,p}}
&\le &
     C_k^2c_{k-1}\left(
     \left\|\zeta_\nu\right\|_{W^{k,p}}
     + \left\|\hat\eta_\nu\right\|_{W^{k+1,p}}
     \right)2^{-\nu}.
\end{eqnarray*}
Abbreviate 
$$
      a_\nu := \left\|\zeta_{\nu+\nu_0}\right\|_{W^{k,p}}
      + C_k\left\|\hat\eta_{\nu+\nu_0}\right\|_{W^{k+1,p}}
$$
and choose $\nu_0$ so large that $C_k^3c_{k-1}2^{-\nu_0}\le 1$.
Then 
\begin{eqnarray*}
      a_{\nu+1} 
&\le &
      \left\|\zeta_{\nu+\nu_0}\right\|_{W^{k,p}}
      + C_k\left\|\hat\eta_{\nu+\nu_0}\right\|_{W^{k+1,p}}
      + C_k\left\|\hat\eta_{\nu+\nu_0+1}\right\|_{W^{k+1,p}} \\
&= &
      a_\nu
      + C_k\left\|\hat\eta_{\nu+\nu_0+1}\right\|_{W^{k+1,p}} \\
&\le &
      a_\nu
      + C_k^3c_{k-1}\left(
     \left\|\zeta_\nu\right\|_{W^{k,p}}
     + \left\|\hat\eta_\nu\right\|_{W^{k+1,p}}
     \right)2^{-\nu-\nu_0}   \\
&\le &
      (1+2^{-\nu})a_\nu
\end{eqnarray*}
for all $\nu$ and hence the sequence $a_\nu$ is bounded.  
It follows that the sequences $\left\|\zeta_\nu\right\|_{W^{k,p}}$
and $2^\nu\left\|\hat\eta_\nu\right\|_{W^{k+1,p}}$ are bounded.
Thus we have proved that $\hat\eta_\nu$ 
and $\zeta_\nu$ satisfy~(\ref{eq:k-1}).
This completes the induction.   
It follows that $\zeta$ is smooth and hence, 
so is $\eta$.  This proves existence. 

We prove uniqueness.  Choose $\delta>0$ so small that 
$$
     c_5c\delta \le\delta_0,\qquad
     2c_4c\delta \le\delta_3,\qquad
     c_5c\delta\le\delta_4,\qquad
     2c_2c_3c_4c\delta<1.
$$
Assume that $\zeta_0,\zeta_1\in T_{(u_0,A_0)}\Bb$ 
and $\eta_0,\eta_1\in\Om^0(\Sigma,\g_P)$
satisfy the requirements of the proposition.
Then
$$
     d_{(u_0,A_0)}^{*_\eps}\zeta_i=0,\qquad 
     g_i^*(u,A) = (\exp_{u_0}(\xi_i),A_0+\alpha_i),
$$
for $i=0,1$, where $g_i:=e^{\eta_i}$. 
By Lemma~\ref{lemma:p}, we have 
$$
     \left\|\eta_i\right\|_{L^\infty}
     \le c_5\eps^{-2/p}\left\|\eta_i\right\|_{2,p,\eps}
     \le c_5c\eps^{-2/p}\left\|\zeta\right\|_{1,p,\eps}
     \le c_5c\delta
     \le \delta_4
$$
for $i=0,1$. 
Hence, by Lemma~\ref{le:eta}, there exists a unique element
$\eta\in\Om^0(\Sigma,\g_P)$ such that 
$$
     g := e^\eta = g_0^{-1}g_1,\qquad
     c_4^{-1}\left\|\eta_1-\eta_0\right\|_{2,p,\eps}
     \le \left\|\eta\right\|_{2,p,\eps}
     \le c_4\left\|\eta_1-\eta_0\right\|_{2,p,\eps}.
$$
The gauge transformation $g$ satisfies
$$
     g^*(\exp_{u_0}(\xi_0),A_0+\alpha_0)
     = (\exp_{u_0}(\xi_1),A_0+\alpha_1).
$$
Moreover, 
$$
     \left\|\zeta_0\right\|_{1,p,\eps} 
     + \left\|\eta\right\|_{2,p,\eps} 
     \le 2c_4c\left\|\zeta\right\|_{1,p,\eps} 
     \le 2c_4c\delta\eps^{2/p}
     \le \delta_3\eps^{2/p}.
$$
Hence $\zeta_0$ and $\eta$ satisfy the hypotheses 
of  Lemma~\ref{le:gauge}. We use Lemma~\ref{le:laplace}
and the estimate~(\ref{eq:gauge2}) of Lemma~\ref{le:gauge}
to obtain 
\begin{eqnarray*}
     \left\|\eta\right\|_{2,p,\eps} 
&\le &
     c_2\left\|
       d_{(u_0,A_0)}^{*_\eps}d_{(u_0,A_0)}\eta
     \right\|_{L^p} \\
&\le &
     c_2c_3\eps^{-2/p}
     \left(
     \left\|\zeta_1-\zeta_0\right\|_{1,p,\eps}
     + \left\|\eta\right\|_{2,p,\eps}
     \right)
     \left\|\eta\right\|_{1,p,\eps}  \\
&\le &
     c_2c_3c_4\eps^{-2/p}
     \left(
     \left\|\zeta_1-\zeta_0\right\|_{1,p,\eps}
     + \left\|\eta_1-\eta_0\right\|_{2,p,\eps}
     \right)
     \left\|\eta\right\|_{1,p,\eps}  \\
&\le &
     2c_2c_3c_4c\eps^{-2/p}\left\|\zeta\right\|_{1,p,\eps}
     \left\|\eta\right\|_{1,p,\eps}  \\
&\le &
     2c_2c_3c_4c\delta
     \left\|\eta\right\|_{1,p,\eps}.
\end{eqnarray*}
Since $2c_2c_3c_4c\delta<1$ we have $\eta=0$
and hence $\eta_1=\eta_0$.
Hence $\alpha_0=\alpha_1$ and $\exp_{u_0}(\xi_0)=\exp_{u_0}(\xi_1)$.
By~(\ref{eq:pre-slice}), (\ref{eq:propslice1}), 
and Lemma~\ref{lemma:p}, we have 
$$
     \left\|\xi_i\right\|_{L^\infty}
     \le c_5\eps^{-2/p}\left\|\zeta_i\right\|_{1,p,\eps}
     \le c_5c\eps^{-2/p}\left\|\zeta\right\|_{1,p,\eps}
     \le c_5c\delta
     \le \delta_0
$$
for $i=0,1$. Hence $\xi_0=\xi_1$. 
\end{proof}


\section{Proof of Theorem~C}\label{sec:tub}

In this section we prove that the map 
$
    \Tilde{\Tcal}^{\eps}:
    \Tilde{\Mcal}_{B,\Sigma}^0(c_0)
    \to
    \Tilde{\Mcal}_{B,\Sigma}^{\eps}
$
introduced in Theorem~\ref{thm:exist}
is locally surjective.  This is the content 
of Theorem~C and is restated more precisely 
as follows. 

\begin{theorem}\label{thm:loc-onto}
Assume~$(H1)$ and~$(H4)$ and let $\bar B\in  H_2(\bar M;\Z)$
be a nontorsion homology class.
Then, for every $c_0>0$ and every $p>2$,
there exist positive constants $\eps_0$ and $\delta$ 
such that the following holds for every $\eps\in(0,\eps_0]$. 
If $(\bar u_0,\bar A_0)\in\Tilde{\Mcal}_{B,\Sigma}^0(c_0-1)$
and
$
     (u,A) = (\exp_{\bar u_0}(\bar\xi),\bar A_0+\bar\alpha)
     \in \Tilde\Mcal_{B,\Sigma}^\eps
$
where 
$
     \bar\zeta=(\bar\xi,\bar\alpha) \in T_{(\bar u_0,\bar A_0)}\Bb
$
satisfies
$$
     \left\|\bar\zeta\right\|_{1,p,\eps;(\bar u_0,\bar A_0)}
     \le \delta\eps^{2/p+1/2},
$$
then there exist 
$
     \bar\xi_0 \in\ker\,\Dd^0_{(\bar u_0,\bar A_0)}
$
and $\eta_0\in\Om^0(\Sigma,\g_P)$ such that
$$
     g^*(u,A) = \Tilde\Tcal^\eps(u_0,A_0),\qquad
     g := e^{\eta_0},\qquad
     (u_0,A_0) := \Ff^0_{(\bar u_0,\bar A_0)}(\bar\xi_0),
$$
$$
    \left\|\bar\zeta_0\right\|_{W^{1,p}}
    + \left\|\eta_0\right\|_{2,p,\eps;A_0}
    \le c\left\|\bar\zeta
    \right\|_{1,p,\eps;(\bar u_0,\bar A_0)}.
$$
\end{theorem}

Here $\Ff^0_{(\bar u_0,\bar A_0)}$ 
is the map of Theorem~\ref{thm:M0}.
The proof of Theorem~\ref{thm:loc-onto}
is based on Corollary~\ref{cor:uniqbaby}
and on the construction of a tubular 
neighbourhood of the moduli space 
$\Mm_{B,\Sigma}^0(c_0)$ in the quotient $\Bcal/\Gcal$.

\begin{proposition}\label{prop:tubular}
Assume~$(H1)$ and~$(H4)$ and let $\bar B\in  H_2(\bar M;\Z)$
be a nontorsion homology class.
For every $p>2$ and every $c_0>0$ there exist positive 
constants $\delta$, $\eps_0$, and $c$ such that, 
for every $\eps\in(0,\eps_0]$, the following holds.
Let 
$
     (\bar u_0,\bar A_0) \in 
     \Tilde\Mcal_{B,\Sigma}^0(c_0-1)
$
and
$
     (u,A) = (\exp_{\bar u_0}(\bar\xi),\bar A_0+\bar\alpha),
$
where the pair
$
     \bar\zeta:=(\bar\xi,\bar\alpha)
     \in T_{(\bar u_0,\bar A_0)}\Bb
$ 
satisfies
\begin{equation}\label{eq:tub}
     \left\|\bar\zeta\right\|_{1,p,\eps;(\bar u_0,\bar A_0)}
     \le \delta\eps^{2/p}.
\end{equation}
Then there exist
$
     \bar\xi_0\in\ker\Dd_{(\bar u_0,\bar A_0)}^0,
$
$
     \eta_0\in\Om^0(\Sigma,\g_P),
$
$
     (u_0,A_0)\in\Tilde\Mcal_{B,\Sigma}^0(c_0),
$ 
and  
$
     \zeta_0=(\xi_0,\alpha_0)\in T_{(u_0,A_0)}\Bb,
$
such that
\begin{equation}\label{eq:tub1}
     g^*(u,A)=(\exp_{u_0}(\xi_0),A_0+\alpha_0),\qquad
     (u_0,A_0) = \Ff^0_{(\bar u_0,\bar A_0)}(\bar\xi_0),
\end{equation}
where $g:=e^{\eta_0}$, and
\begin{equation}\label{eq:tub2}
     d_{(u_0,A_0)}^{*_\eps}\zeta_0=0,\qquad
     \zeta_0\in\im\,(\Dd_{(u_0,A_0)}^\eps)^*,
\end{equation}
\begin{equation}\label{eq:tub3}
     \left\|\bar\xi_0\right\|_{W^{1,p}}
     + \left\|\eta_0\right\|_{2,p,\eps;\bar A_0}
     + \left\|\zeta_0\right\|_{1,p,\eps;(u_0,A_0)}
     \le c\left\|\bar\zeta\right\|_{1,p,\eps;(\bar u_0,\bar A_0)}.
\end{equation}
\end{proposition}

\begin{figure}[htp]
   \centerline
   {\hbox{\epsfysize=160pt\epsffile{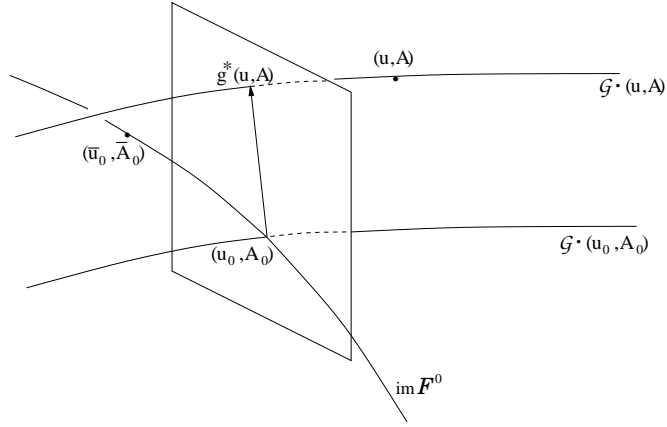}}}
\caption{{A tubular neighbourhood of $\Mm^0$.}}
\label{fig:morse}
\end{figure}

The proof of Proposition~\ref{prop:tubular} is based
on Proposition~\ref{prop:slice}. The latter can
be restated as follows.  Let $\Hh^\eps\subset T\Bb$ 
denote the horizontal subbundle with fibres
$$
     \Hh^\eps_{(u,A)}
     := \ker\,d_{(u,A)}^{*_\eps}
     \subset T_{(u,A)}\Bb. 
$$
Given a pair $(u,A)\in\Bb$ and constants 
$p$, $c_0$, $\delta$, $\eps$ denote by 
$
     \Uu^0 = \Uu^0(\delta,\eps)
     \subset\Tilde{\Mm}^0_{B,\Sigma}(c_0)
$
the open set
$$
     \Uu^0 := 
     \Bigl\{(u_0,A_0)\in\Tilde{\Mm}^0_{B,\Sigma}(c_0)\Big|
     (u,A)=(\exp_{u_0}(\xi),A_0+\alpha),
     \left\|\zeta\right\|_{1,p,\eps}
     \le\delta\eps^{2/p}
     \Bigr\}.  
$$
If $\delta$ and $\eps$ are sufficiently small
then Proposition~\ref{prop:slice} defines
two maps 
$$
     \Ss^\eps:\Uu^0\to\Hh^\eps,\qquad
     \Nn^\eps:\Uu^0\to\Om^0(\Sigma,\g_P)
$$
such that $\Ss^\eps$ is a section of $\Hh^\eps$
over $\Uu^0$ and, for every $(u_0,A_0)\in\Uu^0$, 
the pair $(\xi_0,\alpha_0):=\Ss^\eps(u_0,A_0)$ and the 
gauge transformation $g:=e^\eta$, 
where $\eta:=\Nn^\eps(u_0,A_0)$,
satisfy $g^*(u,A)=(\exp_{u_0}(\xi_0),A_0+\alpha_0)$
and~(\ref{eq:propslice1}).  In particular,
$$
     \|\Ss^\eps(u_0,A_0)\|_{1,p,\eps}
     \le c\|\zeta\|_{1,p,\eps},
$$
where $(u,A)=(\exp_{u_0}(\xi),A_0+\alpha)$
and $\left\|\zeta\right\|_{1,p,\eps}\le\delta\eps^{2/p}$. 
In this notation Proposition~\ref{prop:tubular} 
asserts that for every $(u,A)\in\Bb$, whose distance to 
$\Tilde{\Mm}_{B,\Sigma}^0$ in the $(1,p,\eps)$-norm
is less than $\delta\eps^{2/p}$ for a sufficiently
small constant $\delta$, there exists 
a pair $(u_0,A_0)\in\Uu^0$ such that $\Ss^\eps(u_0,A_0)$
lies in the image of $(\Dd_{(u_0,A_0)}^\eps)^*$.

\begin{lemma}\label{le:slice}
For every $p>2$ and every $c_0>0$ there exist
positive constants $\delta$, $\eps_0$, and $c$ such that
the following holds for every $\eps\in(0,\eps_0]$.  
Let $(u,A)\in\Bb$ and $I\subset\R$ be an interval.
Suppose that 
$
     I\to\Uu^0(\delta,\eps):r\mapsto(u_0(r),A_0(r))
$
is a smooth path, and let 
$
     \zeta(r)=(\xi(r),\alpha(r))\in T_{(u_0(r),A_0(r))}\Bb
$
be the corresponding vector field along this path 
that satisfies
$$
     (u,A)=(\exp_{u_0}(\xi),A_0+\alpha),\qquad
     \left\|\zeta\right\|_{1,p,\eps}\le\delta\eps^{2/p}.
$$
Then the function $r\mapsto\zeta_\eps(r):=\Ss^\eps(u_0(r),A_0(r))$
satisfies the estimate
\begin{equation}\label{eq:slice}
     \left\|(\p_ru_0,\p_rA_0) + \Tabla{r}\zeta_\eps
     \right\|_{1,p,\eps}
     \le c\left(
         \left\|\zeta\right\|_{1,p,\eps} 
         + \left\|\eps^2d_{A_0}^*\p_rA_0-L_{u_0}^*\p_ru_0\right\|_{L^p}
         \right),
\end{equation}
where 
$
     \Tabla{r}\zeta_\eps := (\Tabla{r}\xi_\eps,\p_r\alpha_\eps)
$
and
$ 
     \Tabla{r}\xi_\eps
     := \Nabla{r}\xi_\eps - \frac12J(\Nabla{\p_ru_0}J)\xi_\eps.  
$ 
\end{lemma}

\begin{proof}
Let 
$
     g(r):=e^{\eta_\eps(r)},
$
where
$
     \eta_\eps(r) := \Nn^\eps(u_0(r),A_0(r)) \in\Om^0(\Sigma,\g_P),
$
and denote $\dot u_0:=\p_ru_0$, $\dot A_0:=\p_rA_0$, $\dot g:=\p_rg$. 
Let $\delta_1$ and $c_1$ be the constants
of Proposition~\ref{prop:slice}.
Then 
\begin{equation}\label{eq:slice-est}
     \left\|\eta_\eps(r)\right\|_{2,p,\eps}
     + \left\|\zeta_\eps(r)\right\|_{1,p,\eps}
     \le c_1\left\|\zeta(r)\right\|_{1,p,\eps}
     \le c_1\delta_1\eps^{2/p},
\end{equation}
$$
     g^{-1}u = \exp_{u_0}(\xi_\eps),\qquad
     g^*A = A_0+\alpha_\eps,\qquad
     \eps^2d_{A_0}^*\alpha_\eps - L_{u_0}^*\xi_\eps = 0.
$$
Differentiating these identities we obtain
\begin{equation}\label{eq:slice2}
     - L_{g^{-1}u}(g^{-1}\dot g) = E_1\dot u_0 + E_2\Nabla{r}\xi_\eps,\qquad
     d_{g^*A}(g^{-1}\dot g) = \dot A_0 + \dot\alpha_\eps,
\end{equation}
where $E_1:=E_1(u_0,\xi_\eps)$ and $E_2:=E_2(u_0,\xi_\eps)$
(see Appendix~\ref{app:invariant}),
and
\begin{equation}\label{eq:slice3}
     \eps^2d_{A_0}^*\dot\alpha_\eps 
      - \eps^2*[\dot A_0\wedge *\alpha_\eps]
      - L_{u_0}^*\Nabla{r}\xi_\eps - \rho(\dot u_0,\xi_\eps) 
     = 0,
\end{equation}
where $\rho\in\Om^2(M,\g)$ is given by 
$
     \inner{\eta}{\rho(\xi_1,\xi_2)}
     := \inner{\Nabla{\xi_1}X_\eta}{\xi_2}
$
(see Lemma~\ref{le:rho}).  Inserting the expressions 
for $\Nabla{r}\xi_\eps$ and $\dot\alpha_\eps$
in~(\ref{eq:slice2}) into~(\ref{eq:slice3}) gives
\begin{eqnarray*}
&&
     \eps^2d_{A_0}^*d_{g^*A}(g^{-1}\dot g)
     + L_{u_0}^*E_2^{-1}L_{g^{-1}u}(g^{-1}\dot g) \\
&&
     = \eps^2d_{A_0}^*\dot A_0 + \eps^2*[\dot A_0\wedge *\alpha_\eps]
       - L_{u_0}^*E_2^{-1}E_1\dot u_0 + \rho(\dot u_0,\xi_\eps).
\end{eqnarray*}
Since $g^{-1}u=\exp_{u_0}(\xi_\eps)$ and
$g^*A=A_0+\alpha_\eps$ we have, by Lemma~\ref{le:group}, 
$$
\begin{array}{rcl}
     E_2^{-1}L_{g^{-1}u}(g^{-1}\dot g) 
&= &
     L_{u_0}(g^{-1}\dot g)
       + \left(E_2^{-1}E_1-\one\right) L_{u_0}(g^{-1}\dot g)
       + \Nabla{\xi_\eps}X_{g^{-1}\dot g}(u_0), \\
     d_{A_0}^*d_{g^*A}(g^{-1}\dot g)
&= &
     d_{A_0}^*d_{A_0}(g^{-1}\dot g)
       + [d_{A_0}^*\alpha_\eps,g^{-1}\dot g]
       + *[*\alpha_\eps\wedge d_{A_0}(g^{-1}\dot g)].
\end{array}
$$
Hence
\begin{eqnarray*}
&&
     \eps^2d_{A_0}^*d_{A_0}(g^{-1}\dot g)
     + L_{u_0}^*L_{u_0}(g^{-1}\dot g) \\
&&
     = -\, \eps^2[d_{A_0}^*\alpha_\eps,g^{-1}\dot g]
       - \eps^2 *[*\alpha_\eps\wedge d_{A_0}(g^{-1}\dot g)]  \\
&&\quad
     -\, L_{u_0}^*\left(E_2^{-1}E_1-\one\right) L_{u_0}(g^{-1}\dot g)
     - L_{u_0}^*\Nabla{\xi_\eps}X_{g^{-1}\dot g}(u_0)  \\
&&\quad
     +\, \eps^2*[\dot A_0\wedge *\alpha_\eps]
     + L_{u_0}^*(\one-E_2^{-1}E_1)\dot u_0 + \rho(\dot u_0,\xi_\eps) \\
&&\quad
     +\,\eps^2d_{A_0}^*\dot A_0 -  L_{u_0}^*\dot u_0.
\end{eqnarray*}
By Lemma~\ref{le:laplace} and Lemma~\ref{lemma:p}, 
there exists a constant $c_2>0$ such that
\begin{eqnarray*}
     \left\|g^{-1}\dot g\right\|_{2,p,\eps}
&\le &
     c_2\eps^{-2/p}\left\|\zeta_\eps\right\|_{1,p,\eps}
         \left\|g^{-1}\dot g\right\|_{2,p,\eps}  \\
&&
     +\, c_2\left(
         \left\|\zeta_\eps\right\|_{0,p,\eps} 
         + \left\|\eps^2d_{A_0}^*\dot A_0-L_{u_0}^*\dot u_0\right\|_{L^p}
         \right) \\
&\le &
     c_1c_2\delta_1
         \left\|g^{-1}\dot g\right\|_{2,p,\eps}  \\
&&
     +\, c_2\left(
         \left\|\zeta_\eps\right\|_{0,p,\eps} 
         + \left\|\eps^2d_{A_0}^*\dot A_0-L_{u_0}^*\dot u_0\right\|_{L^p}
         \right).
\end{eqnarray*}
The last inequality follows from~(\ref{eq:slice-est}).
With $c_1c_2\delta_1\le1/2$ it follows that
$$
     \left\|g^{-1}\dot g\right\|_{2,p,\eps}
     \le 2c_2\left(
         \left\|\zeta_\eps\right\|_{0,p,\eps} 
         + \left\|\eps^2d_{A_0}^*\dot A_0-L_{u_0}^*\dot u_0\right\|_{L^p}
         \right).
$$
Hence~(\ref{eq:slice}) follows from~(\ref{eq:slice2})
and~(\ref{eq:propslice1}). 
\end{proof}

Consider the vector bundle 
$$
     \Vv^0 \TO \Tilde{\Mm}_{B,\Sigma}^0
$$
whose fibre over $(u_0,A_0)\in\Tilde{\Mm}_{B,\Sigma}^0$ 
is the finite dimensional vector space $\Vv_{(u_0,A_0)}^0$ 
of all pairs
$$
     (\xi_0,\alpha_0)
     \in\Om^0(\Sigma,H_{u_0})\times\Om^1(\Sigma,\g_P)
$$
that satisfy the equation
\begin{equation}\label{eq:kerD0}
     D\bar\p_{J,A_0}(u_0)\xi_0
     + X_{\alpha_0}(u_0)^{0,1}
     = 0.
\end{equation}
This space can be identified
with the kernel of the operator $\Dd^0_{(u_0,A_0)}$.
Namely, the kernel of $\Dd^0_{(u_0,A_0)}$ consists 
of all sections $\xi_0\in\Om^0(\Sigma,H_{u_0})$ that 
satisfy~(\ref{eq:kerD0}) for some $1$-form 
$\alpha_0\in\Om^1(\Sigma,\g_P)$, and the $1$-form
$\alpha_0$ is uniquely determined by~$\xi_0$. 
Thus 
$
     \Vv^0 \TO \Tilde{\Mm}_{B,\Sigma}^0
$
is a vector bundle of rank $m:=\dim\Mm_{B,\Sigma}^0$.

\begin{lemma}\label{le:commute}
For every $p>2$ and every $c_0>0$ there exist
positive constants $\eps_0$ and $c$ such that
the following holds for every $\eps\in(0,\eps_0]$.  
Let $I\subset\R$ be an interval and
$$
     I\to\Tilde{\Mm}_{B,\Sigma}^0(c_0):
     r\mapsto(u_0(r),A_0(r))
$$
be a smooth path such that, for every $r\in I$, 
\begin{equation}\label{eq:path}
     \left\|\p_ru_0(r)\right\|_{C^2(\Sigma)}
     + \left\|\p_rA_0(r)\right\|_{C^1(\Sigma)}
     \le c_0.
\end{equation}
Then every smooth vector field
$r\mapsto\zeta(r)=(\xi(r),\alpha(r))\in T_{(u_0(r),A_0(r))}\Bb$
satisfies the inequality
$$
     \left\|\Tabla{r}\Dd^\eps\zeta - \Dd^\eps\Tabla{r}\zeta
     \right\|_{k,p,\eps}  
     \le c\eps^{-1}\left\|\zeta\right\|_{k+1,p,\eps}
$$
for $k=0,1$, where
$
     \Dd^\eps := \Dd^\eps_{(u_0(r),A_0(r))}
$
and 
$
     \Tabla{r}\zeta
     := (\Nabla{r}\xi-\frac12J(\Nabla{\p_ru_0}J)\xi,\p_r\alpha).
$
\end{lemma}

\begin{proof}
We denote 
$$
     (\xi,\alpha) := \zeta,\qquad
     \zeta'
     := (\xi',\phi',\psi')
     := \Dd^\eps\zeta,\qquad
     \zeta'_r
     := (\xi'_r,\phi'_r,\psi'_r)
     := \Dd^\eps\Tabla{r}\zeta.
$$
Moreover, we drop the subscript $0$ and write $(u,A):=(u_0,A_0)$. 
Then, in local holomorphic coordinates on $\Sigma$ and a local
frame of $P$, $\zeta'$ is given by 
\begin{eqnarray*}
     \xi'
&= &
     \Tabla{A,s}\xi + J\Tabla{A,t}\xi
     + \frac14N(\xi,v_s-Jv_t)
     + \frac12(J\p_sJ-\p_tJ)\xi
     + L_u\phi + JL_u\psi, \\
     \phi'
&= &   
     \lambda^{-2}\left(\Nabla{A,s}\phi + \Nabla{A,t}\psi\right) 
     + \eps^{-2}L_u^*\xi, \\
     \psi'
&= &
     \lambda^{-2}\left(\Nabla{A,s}\psi - \Nabla{A,t}\phi\right) 
     + \eps^{-2}d\mu(u)\xi.
\end{eqnarray*}
Here we use the notation of Remark~\ref{rmk:local}.
Differentiating these formulae with respect to $r$
we obtain
\begin{eqnarray*}
     \Tabla{r}\xi' - \xi'_r
&= &
     \Tabla{r}\Tabla{A,s}\xi - \Tabla{A,s}\Tabla{r}\xi 
     + J(\Tabla{r}\Tabla{A,t}\xi - \Tabla{A,t}\Tabla{r}\xi)  \\
&&
     +\,\frac14\Tabla{r}N(\xi,v_s-Jv_t)
     - \frac14N(\Tabla{r}\xi,v_s-Jv_t)  \\
&&
     +\,\frac12\Tabla{r}((J\p_sJ-\p_tJ)\xi)
     - \frac12(J\p_sJ-\p_tJ)\Tabla{r}\xi  \\
&&
     +\,\Nabla{\p_ru}X_\phi(u)
     + J\Nabla{\p_ru}X_\psi(u)
     - \frac12J(\Nabla{\p_ru}J)(L_u\phi-JL_u\psi), \\
     \p_r\phi'-\phi'_r
&= &
     \lambda^{-2}\left(
     [\p_r\Phi,\phi]+[\p_r\Psi,\psi]
     \right) \\
&&
     +\, \eps^{-2}\left(\rho(\p_ru,\xi)
     - \frac12d\mu(u)(\Nabla{\p_ru}J)\xi\right), \\
     \p_r\psi'-\psi'_r
&= &
     \lambda^{-2}
     \left([\p_r\Phi,\psi] - [\p_r\Psi,\phi]\right) \\
&&
     -\, \eps^{-2}\left(\rho(\p_ru,J\xi)
       + \frac12L_u^*(\Nabla{\p_ru}J)\xi\right).
\end{eqnarray*}
Here $\rho=\rho_{s,t}\in\Om^2(M,\g)$ is defined by Lemma~\ref{le:rho}.
The required estimates follow from these three 
identities via a term by term inspection. 
\end{proof}

\begin{lemma}\label{le:kernel}
For every $p>2$ and every $c_0>0$ there exist
positive constants $\eps_0$ and $c$ such that
the following holds for every $\eps\in(0,\eps_0]$.  
Let 
$
     r\mapsto(u_0(r),A_0(r))
$
be as in Lemma~\ref{le:commute} and 
suppose that
$
     r\mapsto\zeta_0(r)
$ 
is a smooth section of $\Vv^0$ along this path. 
Abbreviate $\Dd^\eps:=\Dd^\eps_{(u_0(r),A_0(r))}$ and let
$\zeta_\eps(r)\in\ker\Dd^\eps$ be given by 
$$
     \zeta_\eps(r) := \zeta_0(r) - {\Dd^\eps}^*
     \left(\Dd^\eps{\Dd^\eps}^*\right)^{-1}\Dd^\eps\zeta_0(r). 
$$
Then
\begin{eqnarray}
     \|\zeta_\eps - \zeta_0\|_{0,2,\eps}
&\le &
     c\eps^2\|\xi_0\|_{L^2},
     \label{eq:kernel0}  \\
     \|\Tabla{r}\zeta_\eps-\Tabla{r}\zeta_0\|_{0,2,\eps}
&\le &
     c\left(
     \|\xi_0\|_{L^2}
     + \|\Tabla{r}\xi_0\|_{L^2}
     \right).
     \label{eq:kernel1}
\end{eqnarray}
\end{lemma}

\begin{proof}
Let
$
     (\xi_\eps,\alpha_\eps):=\zeta_\eps 
$
for $\eps\ge 0$ and $\zeta:=(\xi,\alpha):=\zeta_\eps-\zeta_0$. 
Then 
$$
    \zeta \in \im\,{\Dd^\eps}^*,\qquad
    \Dd^\eps\zeta
    = \left(\begin{array}{c}
      0 \\ d_{A_0}^*\alpha_0 \\ -*d_{A_0}\alpha_0
      \end{array}\right).
$$
Hence, by Lemma~\ref{le:onto1}, there exist constants 
$c_1,c_2>0$ (depending only on $c_0$)
such that
$$
     \left\|\zeta\right\|_{1,2,\eps}
\le 
     c_1\eps\left\|\Dd^\eps\zeta\right\|_{0,2,\eps}  
= 
     c_1\eps^2\left(
        \left\|d_{A_0}\alpha_0\right\|_{L^2}
        + \left\|d_{A_0}^*\alpha_0\right\|_{L^2}
        \right)  
\le 
     c_2\eps^2\left\|\xi_0\right\|_{L^2}.
$$
The last inequality follows from~(\ref{eq:kerD0}) and the 
basic elliptic estimates for the operator~$\Dd^0$. 
Thus we have proved~(\ref{eq:kernel0}).
To prove~(\ref{eq:kernel1}) let
$$
     \zeta' := \zeta'(r) 
     := - (\Dd^\eps{\Dd^\eps}^*)^{-1}\Dd^\eps\zeta_0(r)
     \in\Xx'_{u_0(r)}
$$
so that
$
     \zeta = {\Dd^\eps}^*\zeta'.
$
Then, by Lemma~\ref{le:onto} with $p=2$,
$$
     \left\|\zeta'\right\|_{2,2,\eps}
     \le c_3\left\|{\Dd^\eps}^*\zeta'\right\|_{1,2,\eps}
      =  c_3\left\|\zeta\right\|_{1,2,\eps}
     \le c_2c_3\eps^2\left\|\xi_0\right\|_{L^2},
$$
and hence, by Lemma~\ref{le:commute} (with $c=c_4$ and $\Dd^\eps$
replaced by ${\Dd^\eps}^*$), 
$$
     \|{\Dd^\eps}^*\Tabla{r}\zeta'
       - \Tabla{r}{\Dd^\eps}^*\zeta'
     \|_{1,2,\eps}
     \le c_4\eps^{-1}\left\|\zeta'\right\|_{2,2,\eps}
     \le c_2c_3c_4\eps\left\|\xi_0\right\|_{L^2}.
$$
Now it follows from Lemmata~\ref{le:onto1} (with $c=c_5$)
and~\ref{le:commute} (with $c=c_4$) that
\begin{eqnarray*}
     \|{\Dd^\eps}^*\Tabla{r}\zeta'\|_{1,2,\eps}
&\le &
     c_5\|\Dd^\eps{\Dd^\eps}^*\Tabla{r}\zeta'
     \|_{0,2,\eps} \\
&\le &
     c_5\|
     \Dd^\eps({\Dd^\eps}^*\Tabla{r}\zeta'
        - \Tabla{r}{\Dd^\eps}^*\zeta')
     \|_{0,2,\eps} \\
&&
     +\, c_5\|
     \Dd^\eps\Tabla{r}\zeta - \Tabla{r}\Dd^\eps\zeta
     \|_{0,2,\eps} 
     + c_5\|\Tabla{r}\Dd^\eps\zeta\|_{0,2,\eps} \\
&\le & 
     c_6\eps^{-1}\|
        {\Dd^\eps}^*\Tabla{r}\zeta'
        - \Tabla{r}{\Dd^\eps}^*\zeta'
     \|_{1,2,\eps} \\
&&
     +\, c_5\left(
     c_4\eps^{-1}\|\zeta\|_{1,2,\eps}
     + \|\Tabla{r}\Dd^\eps\zeta\|_{0,2,\eps}
     \right) \\
&\le &
     c_7\left(
     \|\xi_0\|_{L^2}
     + \|\Tabla{r}\xi_0\|_{L^2}
     \right).
\end{eqnarray*}
Hence 
\begin{eqnarray*}
     \|\Tabla{r}\zeta\|_{0,2,\eps}
&\le &
     \|
      \Tabla{r}{\Dd^\eps}^*\zeta'
      - {\Dd^\eps}^*\Tabla{r}\zeta'
     \|_{0,2,\eps}
     + \|{\Dd^\eps}^*\Tabla{r}\zeta'\|_{0,2,\eps}  \\
&\le &
     (c_7+c_2c_3c_4\eps)
     \left(\|\xi_0\|_{L^2}+\|\Tabla{r}\xi_0\|_{L^2}\right).
\end{eqnarray*}
This proves~(\ref{eq:kernel1}). 
\end{proof}

The estimate~(\ref{eq:kernel1}) is fairly crude.
More careful considerations give an additional
factor $\eps$. However, we shall not use this fact. 

\begin{lemma}\label{le:tubular}
For every $p>2$ and every $c_0>0$ there exist
positive constants $\delta$, $\eps_0$, and $c$ such that
the following holds for every $\eps\in(0,\eps_0]$.  
Let $(u,A)\in\Bb$ and $I\subset\R$ be an interval.
Suppose that the path
$
     I\to\Uu^0(\delta,\eps):r\mapsto(u_0(r),A_0(r))
$
and the vector field $r\mapsto\zeta(r)\in T_{(u_0(r),A_0(r))}\Bb$ 
are as in Lemma~\ref{le:slice}.  Moreover, assume~(\ref{eq:path})
and let $r\mapsto\zeta_0(r)$ and $r\mapsto\zeta_\eps(r)$  
be as in Lemma~\ref{le:kernel}. 
Then 
\begin{eqnarray*}
&&
     \left|
     \frac{d}{dr} \INNER{\zeta_\eps}{\Ss^\eps(u_0,A_0)}_\eps
     + \inner{\xi_0}{\dot u_0}
     \right|   \\
&&\quad
     \le 
     c\left(\|\xi_0\|_{L^2}
      + \|\Tabla{r}\xi_0\|_{L^2}\right)
     \left(\eps^2
     + \|\zeta\|_{1,p,\eps}
     + \|L_{u_0}^*\p_ru_0\|_{L^p}
     \right).
\end{eqnarray*}
\end{lemma}

\begin{proof}
Abbreviate $\Ss^\eps:=\Ss^\eps(u_0,A_0)$.
Consider the identity
\begin{eqnarray*}
    \frac{d}{dr} \INNER{\zeta_\eps}{\Ss^\eps}_\eps
     + \INNER{\xi_0}{\p_ru_0}_\eps 
&= &
    \INNER{\Tabla{r}\zeta_\eps}{\Ss^\eps}_\eps
    + \INNER{\zeta_\eps} 
    {(\p_ru_0,\p_rA_0)+\Tabla{r}\Ss^\eps}_\eps  \\
&&
    +\, \INNER{\zeta_0-\zeta_\eps}
    {(\p_ru_0,\p_rA_0)}_\eps
    - \eps^2\inner{\alpha_0}{\p_rA_0}.
\end{eqnarray*}
By Lemma~\ref{le:kernel}, the $(0,2,\eps)$-norm of 
$\Tabla{r}\zeta_\eps$ is bounded above by a constant 
times $\|\xi_0\|_{L^2}+\|\Tabla{r}\xi_0\|_{L^2}$.
By Proposition~\ref{prop:slice}, the $(0,2,\eps)$-norm
of $\Ss^\eps$ is bounded above by a constant times
$\|\zeta\|_{1,p,\eps}$.  Hence the first term satisfies
the required bound. For the second term the estimate 
follows from Lemma~\ref{le:slice} and the fact that
the $(0,2,\eps)$-norm of $\zeta_\eps$ is bounded above
by a constant times $\|\xi_0\|_{L^2}$.  For the third term 
we use~(\ref{eq:path}) and~(\ref{eq:kernel0})
and for the last the estimate follows from~(\ref{eq:path}). 
\end{proof}

\begin{proof}[Proof of Proposition~\ref{prop:tubular}]
Let $U\subset\R^m$ be an open set containing zero and 
$$
     U\TO\Tilde{\Mm}_{B,\Sigma}^0(c_0):x\mapsto (u_0(x),A_0(x))
$$
be the composition of the map $\Ff^0_{(\bar u_0,\bar A_0)}$ 
defined in Theorem~\ref{thm:M0} with a Hilbert space 
isomorphism $\R^m\to\ker\Dd_{(\bar u_0,\bar A_0)}^0$.
Then 
$$
     (u_0(0),A_0(0)) = (\bar u_0,\bar A_0),\qquad
     (\p_iu_0(0),\p_iA_0(0))\in\Vv^0_{(\bar u_0,\bar A_0)}
$$
for $i=1,\dots,m$; in particular, 
$$
     L_{\bar u_0}^*\p_iu_0(0) = 0,\qquad 
     \INNER{\p_iu_0(0)}{\p_ju_0(0)}_{L^2(\Sigma)} = \delta_{ij}.  
$$
Now choose $m$ smooth sections
$\zeta_{10},\dots,\zeta_{m0}:U\to\Vv^0$ so that
$$
     \zeta_{j0}(x) = (\xi_{j0}(x),\alpha_{j0}(x)) 
     \in \Vv^0_{(u_0(x),A_0(x))}
$$
and
$$
     \zeta_{j0}(0) = - (\p_ju_0(0),\p_jA_0(0)),\qquad
     \INNER{\xi_{i0}(x)}{\xi_{j0}(x)}_{L^2(\Sigma)} = \delta_{ij}
$$
for $x\in U$ and $j=1,\dots,m$.
Given $x\in U$ we abbreviate 
$
     \Dd^\eps:=\Dd_{(u_0(x),A_0(x))}^\eps.
$
If $\eps$ is sufficiently small then, by Lemma~\ref{le:onto1},
this operator is surjective for every $x\in U$.  In this
case we define 
$
     \zeta_{j\eps}(x) \in \ker\Dd^\eps
$
by
$$
     \zeta_{j\eps}(x) := \zeta_{j0}(x) - {\Dd^\eps}^*
     \left(\Dd^\eps{\Dd^\eps}^*\right)^{-1}\Dd^\eps\zeta_{j0}(x)
$$
for $j=1,\dots,m$. By Lemma~\ref{le:kernel},
these vectors form a basis of $\ker\,\Dd^\eps$ 
for $\eps$ sufficiently small.
Now let $\delta_1$ and $c_1$ be the 
constants of Proposition~\ref{prop:slice}.  Choose
$\delta_0>0$ so small that
$$
     |x|<\delta_0\eps^{2/p}\qquad\IMP\qquad
     (u_0(x),A_0(x))\in\Uu^0(\delta_1,\eps)
$$
for $x\in U$ and $0<\eps\le1$. 
Let $\Ss^\eps:\Uu^0\to\Hh^\eps$ be the map of 
Proposition~\ref{prop:slice} as introduced above.
Define $\theta=(\theta_1,\dots,\theta_m):U\to\R^m$ by 
$$
     \theta_j(x) 
     := \INNER{\zeta_{j\eps}(x)}{\Ss^\eps(u_0(x),A_0(x))}_\eps,
$$
where $\INNER{\cdot}{\cdot}_\eps$
denotes the $(0,2,\eps)$-inner product 
on $T_{(u_0(x),A_0(x))}\Bb$. Then 
$$
     \theta(x)=0 \qquad\IFF\qquad
     \Ss^\eps(u_0(x),A_0(x))\in\im\,{\Dd^\eps}^*.
$$
We shall establish the existence of a zero of 
$\theta$ with the inverse function theorem.  
We must prove that $|d\theta(x)-\one|\le1/2$
on a ball of radius $r$ and that $\theta(0)$ 
is less than $r/2$. 

To see this, we first observe that 
$$
     \left|\theta_j(0)\right|
     \le \left\|\zeta_{j0}(0)\right\|_{0,2,\eps}
         \left\|\Ss^\eps(\bar u_0,\bar A_0)\right\|_{0,2,\eps}
     \le c_2\left\|\bar\zeta\right\|_{1,p,\eps}
     \le c_2\delta\eps^{2/p}.
$$
Here we have used the fact that the $(0,2,\eps)$-norm 
of $\zeta_{j\eps}(0)$ is less than or equal to 
the $(0,2,\eps)$-norm of $\zeta_{j0}(0)= - (\p_ju_0(0),\p_jA_0(0))$,
that the $L^2$-norm of $\p_jA_0(0)$ is controlled by 
$\|\p_ju_0(0)\|_{L^2}=1$, 
that the $(0,2,\eps)$-norm of $\Ss^\eps(\bar u_0,\bar A_0)$
is controlled by its $(1,p,\eps)$-norm,
and that, by Proposition~\ref{prop:slice} and~(\ref{eq:tub}), 
the latter is bounded above by 
$
     c_1\|\bar\zeta\|_{1,p,\eps}
     \le c_1\delta\eps^{2/p}. 
$
Thus we have proved that
\begin{equation}\label{eq:tub0}
     \left|\theta(0)\right|
     \le \sqrt{m}c_2\left\|\bar\zeta\right\|_{1,p,\eps}
     \le \sqrt{m}c_2\delta\eps^{2/p}. 
\end{equation}

Now let $\zeta(x)=(\xi(x),\alpha(x))\in T_{(u_0(x),A_0(x))}\Bb$ 
be the unique smooth section defined by 
$$
     (u,A) = (\exp_{u_0(x)}(\xi(x)),A_0(x)+\alpha(x)),\qquad
     \zeta(0) = \bar\zeta,
$$
for $x$ sufficiently small.
Then there exists a constant $c_3>0$ such that 
$$
     \left\|\zeta(x)\right\|_{1,p,\eps}
     \le \left\|\bar\zeta\right\|_{1,p,\eps} + c_3|x|
$$
for $x$ sufficiently small and $0<\eps\le1$.  
Hence, by Proposition~\ref{prop:slice}, we have
that, for $|x|<\delta_0\eps^{2/p}$ and $0<\eps\le1$, 
\begin{eqnarray}\label{eq:tub4}
&&
     \left\|\Nn^\eps(u_0(x),A_0(x))\right\|_{2,p,\eps}
     + \left\|\Ss^\eps(u_0(x),A_0(x))\right\|_{1,p,\eps} 
     \nonumber \\
&&
     \le c_1\left\|\zeta(x)\right\|_{1,p,\eps} 
     \le c_1\left(
         \left\|\bar\zeta\right\|_{1,p,\eps} + c_3|x|
         \right).
\end{eqnarray}
Moreover, there exists a constant $c_4>0$ such that 
$$
     |\delta_{ij}+\INNER{\xi_{j0}(x)}{\p_iu_0(x)}_{L^2}|
     + \|L_{u_0(x)}^*\p_iu_0(x)\|_{L^p}
     \le  c_4|x|, 
$$
$$
     \|\xi_{j0}(x)\|_{L^2}
     + \|\Tabla{i}\xi_{j0}(x)\|_{L^2}
     \le  c_4
$$
for $x$ sufficiently small.
Now suppose that $\delta_1$ and $\eps_0$ have been chosen so
small that the assertion of Lemma~\ref{le:tubular}
holds, with $c$ replaced by $c_5$, 
for the paths $x_i\mapsto(u_0(x),A_0(x))$, 
$x_i\mapsto\zeta(x)$,
$x_i\mapsto\zeta_{j0}(x)$, 
and $x_i\mapsto\zeta_{j\eps}(x)$. 
Then 
\begin{eqnarray*}
     \left|\p_i\theta_j(x)-\delta_{ij}\right|
&\le &
     \left|\delta_{ij}+\INNER{\xi_{j0}(x)}{\p_iu_0(x)}\right|
     + \left|\p_i\theta_j(x)+\INNER{\xi_{j0}(x)}{\p_iu_0(x)}\right|  \\
&\le &
     c_4|x| +
     c_5\bigl(\|\xi_{j0}(x)\|_{L^2}
      + \|\Tabla{i}\xi_{j0}(x)\|_{L^2}\bigr)\cdot \\
&&\qquad 
     \cdot\bigl(\eps^2
     + \|\zeta(x)\|_{1,p,\eps}
     + \|L_{u_0(x)}^*\p_iu_0(x)\|_{L^p}
     \bigr) \\
&\le &
     c_4|x| 
     + c_4c_5\bigl(\eps^2+\|\bar\zeta\|_{1,p,\eps}
     + (c_3 + c_4)|x|\bigr)  \\
&\le &
     c_6\bigl(
     \eps^2 + \|\bar\zeta\|_{1,p,\eps} + |x|
     \bigr)
\end{eqnarray*}
for $|x|<\delta_0\eps^{2/p}$ and $0<\eps\le\eps_0$. 
Thus the Jacobian $d\theta(x)\in\R^{m\times m}$ satisfies
$$
     \left|d\theta(x)-\one\right|
     \le c_7\left(\eps^2+\delta\eps^{2/p}+|x|\right).
$$
Choose $\delta_0$ and $\eps_0$ so small that
$c_7(\eps_0^2+2\delta_0\eps_0^{2/p})\le 1/2$.
Then 
$$
     |x|\le \delta_0\eps^{2/p},\quad 
     0<\eps\le\eps_0,\quad
     0<\delta\le\delta_0
     \qquad\IMP\qquad 
     \left|d\theta(x)-\one\right|\le 1/2.
$$
Hence the inverse function theorem asserts that 
$\theta(B_r(0))\supset B_{r/2}(\theta(0))$
whenever $r<\delta_0\eps^{2/p}$.
Now suppose that $\sqrt{m}c_2\delta<\delta_0/2$.
Then, by~(\ref{eq:tub0}), we have
$
     2\left|\theta(0)\right| < \delta_0\eps^{2/p}
$
and hence we can apply the inverse function theorem
with $r=2\left|\theta(0)\right|$. 
Then $B_{r/2}(\theta(0))$ contains zero and,
by the inverse function theorem, there exists
a point $x_0\in\R^m$ such that 
$$
      \theta(x_0) = 0,\qquad 
      |x_0| \le 2 \left|\theta(0)\right| 
      \le 2\sqrt{m}c_2\left\|\bar\zeta\right\|_{1,p,\eps}.
$$
The last inequality follows from~(\ref{eq:tub0}). 
Now define 
$$
      (u_0,A_0) := (u_0(x_0),A_0(x_0)),\qquad
      \zeta_0 := \Ss^\eps(u_0,A_0),\qquad
      \eta_0 := \Nn^\eps(u_0,A_0). 
$$
Then~(\ref{eq:tub1}) and~(\ref{eq:tub2}) are 
satisfied by definition.  The estimate~(\ref{eq:tub3})
follows from Proposition~\ref{prop:slice}:
$$
      \left\|\eta_0\right\|_{2,p,\eps}
      + \left\|\zeta_0\right\|_{1,p,\eps}
      \le c_1\left\|\zeta(x_0)\right\|_{1,p,\eps}  
      \le c_1\left(\left\|\bar\zeta\right\|_{1,p,\eps} + c_3|x_0|\right) 
      \le c_8\left\|\bar\zeta\right\|_{1,p,\eps}.
$$
Moreover, the vector $\bar\xi_0$
in the assertion of Proposition~\ref{prop:tubular} 
is the image of $x_0$ under our Hilbert space 
isomorphism $\R^m\to\ker\Dd_{(\bar u_0,\bar A_0)}^0$.
Hence, by elliptic regularity for the Cauchy--Riemann 
operator, its $W^{1,p}$-norm is bounded by $|x_0|$ 
and hence by the $(1,p,\eps)$-norm of~$\bar\zeta$.
\end{proof}

\begin{proof}[Proof of Theorem~\ref{thm:loc-onto}]
Let $p>2$ and $c_0>0$ be given. Choose positive constants 
$\eps_0$, $\delta_1$, $\delta_2$, and $c_1$ 
such that Proposition~\ref{prop:tubular}
holds with $\delta$ replaced by $\delta_1$ 
and $c$ replaced by $c_1$, 
Corollary~\ref{cor:uniqbaby} holds with 
$\delta$ replaced by $\delta_2$,
and both results hold for $0<\eps\le\eps_0$. 
Now choose $\delta>0$ so small that
$$
    \delta{\eps_0}^{1/2}\le\delta_1,\qquad
    c_1\delta \le \delta_2.
$$
Let $\eps\in(0,\eps_0]$ and suppose that 
$(u,A)$ and $(\bar u_0,\bar A_0)$
satisfy the hypotheses of Theorem~B, namely
$$
     (\bar u_0,\bar A_0)\in\Tilde{\Mm}_{B,\Sigma}^0(c_0-1)
$$
and 
$$
     (u,A) = (\exp_{\bar u_0}(\bar\xi),\bar A_0+\bar\alpha)
     \in\Tilde{\Mm}_{B,\Sigma}^\eps,
$$
where 
$
     \bar\zeta:=(\bar\xi,\bar\alpha)
     \in T_{(\bar u_0,\bar A_0)}\Bb
$
satisfies 
$$
     \left\|\bar\zeta\right\|_{1,p,\eps;(\bar u_0,\bar A_0)}
     \le \delta\eps^{2/p+1/2}
     \le \delta_1\eps^{2/p}.
$$
By Proposition~\ref{prop:tubular}, there exist
$$
     \bar\xi_0\in\ker\Dd_{(\bar u_0,\bar A_0)}^0,\qquad
     \eta_0\in\Om^0(\Sigma,\g_P),
$$
and
$$
     (u_0,A_0)\in\Tilde\Mcal_{B,\Sigma}^0(c_0),\qquad
     \zeta_0=(\xi_0,\alpha_0)\in T_{(u_0,A_0)}\Bb,
$$
satisfying~(\ref{eq:tub1}), (\ref{eq:tub2}),
and~(\ref{eq:tub3}), with $c$ replaced by $c_1$.
Hence 
$$
     \left\|\zeta_0\right\|_{1,p,\eps;(u_0,A_0)}
     \le c_1\left\|\bar\zeta\right\|_{1,p,\eps;(\bar u_0,\bar A_0)}
     \le c_1\delta\eps^{2/p+1/2}
     \le \delta_2\eps^{2/p+1/2}.
$$
This estimate together with~(\ref{eq:tub2})
shows that $(u_0,A_0)$ 
and $\zeta_0$ satisfy the hypotheses
of Corollary~\ref{cor:uniqbaby}.  
Hence, by~(\ref{eq:tub1}),
$$
     g^*(u,A)
     = (\exp_{u_0}(\xi_0),A_0+\alpha_0)
     = \Tilde{\Tt}^\eps(u_0,A_0),\qquad
     g := e^{\eta_0}.
$$
Moreover, again by~(\ref{eq:tub1}), 
$$
     (u_0,A_0) = \Ff^0_{(\bar u_0,\bar A_0)}(\bar\xi_0)
$$
and, by~(\ref{eq:tub3}), 
$$
     \left\|\bar\zeta_0\right\|_{W^{1,p}}
     + \left\|\eta_0\right\|_{2,p,\eps;\bar A_0}
     \le c_1\left\|\bar\zeta
         \right\|_{1,p,\eps;(\bar u_0,\bar A_0)}.
$$
This proves the theorem.
\end{proof}


\section{A priori estimates}\label{sec:apriori}

In this section we assume that $J\in\Jj_\G(M,\om)$ is 
independent of $z\in\Sigma$ (or in local coordinates
is independent of $s$ and $t$).
Let $\Om\subset\C$ be a bounded open set, 
$K\subset\Om$ be a compact subset, 
and $\lambda:\bar\Om\to(0,\infty)$ be a smooth function. 
Given $u:\Om\to M$ and $\Phi,\Psi:\Om\to\g$
we define $v_s,v_t:\Om\to u^*TM$ and $\kappa:\Om\to\g$
by 
$$
     v_s := \p_su+X_\Phi(u),\quad 
     v_t := \p_tu+X_\Psi(u),\quad
     \kappa := \p_s\Psi-\p_t\Phi+[\Phi,\Psi].
$$
Moreover, as in Remark~\ref{rmk:local}, we use the notation
$$
\begin{array}{rclcrcl}
     \Nabla{A,s}\xi
     &:= & \Nabla{s}\xi+\Nabla\xi X_{\Phi}(u), &\quad &
     \Nabla{A,t}\xi
     &:= & \Nabla{t}\xi+\Nabla\xi X_{\Psi}(u), \\
     \Nabla{A,s}\eta &:= & \p_s\eta+[\Phi,\eta], &\quad &
     \Nabla{A,t}\eta &:= & \p_t\eta+[\Psi,\eta], 
\end{array}
$$
for $\xi:\Om\to u^*TM$ and $\eta:\Om\to\g$.  Then
$$
     \Nabla{A,s}\mu(u) = d\mu(u)v_s = -L_u^*Jv_s,\qquad
     \Nabla{A,t}\mu(u) = d\mu(u)v_t = -L_u^*Jv_t.
$$
Moreover,
$$
     \Nabla{A,s}v_t-\Nabla{A,t}v_s = L_u\kappa,
$$
by Lemma~\ref{le:Lu-phi}, 
$$
     \Nabla{A,s}L_u\eta - L_u\Nabla{A,s}\eta 
     = \Nabla{v_s}X_\eta(u),\qquad
     \Nabla{A,t}L_u\eta - L_u\Nabla{A,t}\eta 
     = \Nabla{v_t}X_\eta(u),
$$
and, by Lemma~\ref{le:curv}, 
$$
     \Nabla{A,s}\Nabla{A,t}\xi-\Nabla{A,t}\Nabla{A,s}\xi
     = R(v_s,v_t)\xi + \Nabla{\xi}X_\kappa(u).
$$
In local coordinates equations~(\ref{eq:eps}) have the form
\begin{equation}\label{eq:eps-loc}
     v_s + Jv_t = 0,\qquad 
     \lambda^{-2}\kappa + \eps^{-2}\mu(u) = 0.
\end{equation}
If~(\ref{eq:eps-loc}) holds then 
$$
    (\Nabla{v_s}J)v_t=(\Nabla{v_t}J)v_s,\qquad
     \Nabla{A,s}v_s + \Nabla{A,t}v_t = - JL_u\kappa.
$$
Given a constant $C>0$ we denote by $M^C\subset M$ 
the compact subset of all $x\in M$ that satisfy
$$
     |\mu(x)|\le C,\qquad
     |\eta|\le C|L_x\eta|
$$
for every $\eta\in\g$. 

\begin{lemma}\label{le:mu}
For every $C>0$ and every triple
$(\Om,K,\lambda)$ as above there exist constants 
$\eps_0=\eps_0(C,\Om,K,\lambda)>0$ and $c=c(C,\Om,K,\lambda)>0$, 
such that the following holds for every $\eps\in(0,\eps_0]$. 
If $u$, $\Phi$, and $\Psi$ satisfy~(\ref{eq:eps-loc}) and
$$
     \|v_s\|_{L^\infty(\Om)}\le C,\qquad
     u(\Om)\subset M^C
$$
then 
$$
     \sup_K|\mu(u)|^2\le c\eps^3,\qquad
     \int_K|\mu(u)|^2 \le c\eps^4.
$$
\end{lemma}

\begin{proof}
Let $\Delta=\p^2/\p s^2+\p^2/\p t^2$ denote the standard
Laplacian.  For $r>0$ denote 
$B_r(z_0):=\{z\in\C\,|\,|z-z_0|<r\}$,
$B_r:=B_r(0)$, and
$$
     \Om_r := \{z\in\C\,|\,B_r(z)\subset\Om\} \subset \Om.
$$
By~(\ref{eq:eps-loc}), we have
$$
     d\mu(u)v_s=-L_u^*Jv_s=-L_u^*v_t,\qquad
     d\mu(u)v_t=-L_u^*Jv_t=L_u^*v_s,
$$
and hence
\begin{eqnarray*}
     \Nabla{A,s}\Nabla{A,s}\mu(u)
     + \Nabla{A,t}\Nabla{A,t}\mu(u)  
&= &
     \Nabla{A,s}d\mu(u)v_s
     + \Nabla{A,t}\d\mu(u)v_t  \\
&= &
     \Nabla{A,t}L_u^*v_s - \Nabla{A,s}L_u^*v_t \\
&= &
     L_u^*(\Nabla{A,t}v_s - \Nabla{A,s}v_t)
     - 2\rho(v_s,v_t) \\
&= &
     - L_u^*L_u\kappa - 2\rho(v_s,v_t) \\
&= &
     (\lambda/\eps)^2L_u^*L_u\mu(u) - 2\rho(v_s,v_t).
\end{eqnarray*}
Here $\rho\in\Om^2(M,\g)$ is as in Lemma~\ref{le:rho}.
Thus
\begin{eqnarray*}
     \Delta|\mu(u)|^2
&= &
     2|\Nabla{A,s}\mu(u)|^2 + 2|\Nabla{A,t}\mu(u)|^2 \\
&&
     +\, 2\inner{\mu(u)}{\Nabla{A,s}\Nabla{A,s}\mu(u)
     + \Nabla{A,t}\Nabla{A,t}\mu(u)}  \\
&= &
     2|\Nabla{A,s}\mu(u)|^2 + 2|\Nabla{A,t}\mu(u)|^2  \\
&&
     +\, 2(\lambda/\eps)^2|L_u\mu(u)|^2
     - 4\inner{\mu(u)}{\rho(v_s,v_t)}.
\end{eqnarray*}
Now choose a constant $c>0$ such that 
$$
     \min_\Om\lambda^2\ge c^{-1},\qquad
     2|\rho(\xi_1,\xi_2)|\le c\left|\xi_1\right|\left|\xi_2\right|,
$$
for all $x\in M^C$ and $\xi_1,\xi_2\in T_xM$.  Then 
$$
     \Delta\left|\mu(u)\right|^2
\ge 
     \frac{2}{cC^2\eps^2}\left|\mu(u)\right|^2
     - 2cC^2\left|\mu(u)\right|  
\ge 
     \frac{1}{cC^2\eps^2}\left|\mu(u)\right|^2
     - c^3C^6\eps^2
$$
and hence
\begin{equation}\label{eq:apriori1}
     \delta\left|\mu(u)\right|^2
     \le a\eps^4 + \eps^2\Delta\left|\mu(u)\right|^2,
\end{equation}
where $\delta:=1/cC^2$ and $a:=c^3C^6$. 
Fix a constant $r>0$ such that 
$$
     z\in K\qquad\IMP\qquad
     B_{3r}(z)\subset\Om.
$$
Then, by~(\ref{eq:apriori1}) and 
Lemma~\ref{le:apriori1} below, we have
$$
     \delta\int_{B_{2r}(z)}|\mu(u)|^2
     \le 9\pi r^2a\eps^4
        + \frac{4\eps^2}{r^2}\int_{B_{3r}(z)}
           \left|\mu(u)\right|^2
     \le 9\pi r^2a\eps^4 + 36\pi C^2\eps^2
$$
for every $z\in\Om_{3r}$.
Applying Lemma~\ref{le:apriori1} again we obtain
$$
     \delta\int_{B_r(z)}|\mu(u)|^2
     \le 4\pi r^2a\eps^4
     + \frac{4\eps^2}{r^2}\int_{B_{2r}(z)}
     \left|\mu(u)\right|^2
     \le c'\eps^4
$$
for every $z\in\Om_{3r}$, where 
$
     c':=4\pi(ar^2 + 9a\eps^2/\delta
         + 36 C^2/\delta r^2).
$
The $L^2$-estimate now follows by
taking the sum over finitely many balls 
of radius $r$ that cover the compact set $\Om_{3r}$. 
Moreover, by~(\ref{eq:apriori1}), the function 
$z\mapsto |\mu(u(z))|^2+4^{-1}a\eps^2|z-z_0|^2$
is subharmonic in $\Om$ for every $z_0\in\C$. Hence, 
by the mean value inequality,
$$
     |\mu(u(z))|^2
     \le \frac{a\eps^2}{4\pi s^2}\int_{B_s}|z|^2
         + \frac{1}{\pi s^2}\int_{B_s(z)}|\mu(u)|^2
     \le \frac{a\eps^2s^2}{8} 
         + \frac{c'\eps^4}{\delta\pi s^2}
$$
for $z\in\Om_{3r}$ and $0\le s\le r$. Assume $\eps\le r^2$.
Then we can choose $s:=\sqrt{\eps}$ and
this proves $L^\infty$ estimate.
\end{proof}

\begin{lemma}\label{le:apriori1}
Let $u:B_{R+r}\to\R$ be a $C^2$-function and
$f,g:B_{R+r}\to\R$ be continuous such that
$$
    f\le g + \Delta u,\qquad
     u\ge 0,\qquad f\ge 0,\qquad g\ge 0.
$$
Then 
$$
     \int_{B_R}f \le \int_{B_{R+r}}g 
     + \frac{4}{r^2}\int_{B_{R+r}\setminus B_R}u.
$$
Moreover, if $g=cu$ then
$$ 
     \frac{\pi}{2}\sup_{B_R}u
     \le \left(c+\frac{4}{r^2}\right)\int_{B_{R+r}}u.
$$
\end{lemma}

\begin{proof}
For $0\le s\le r$ we have
$$
     \int_{B_R}f-\int_{B_{R+r}}g
     \le \int_{B_{R+s}}(f-g) 
     \le \int_{B_{R+s}}\Delta u  
     = \int_{\p B_{R+s}}\frac{\p u}{\p\nu} 
$$
and hence 
$$
     \frac{d}{ds}\int_{\p B_{R+s}}u
     = \int_{\p B_{R+s}}\frac{\p u}{\p\nu} 
     + \frac{1}{R+s}\int_{\p B_{R+s}}u
     \ge \int_{B_R}f-\int_{B_{R+r}}g.
$$
Integrate this inequality over the interval
$0\le s\le t$ to obtain
$$
     \int_{B_R}f-\int_{B_{R+r}}g
     \le \frac{1}{t}\int_{\p B_{R+t}}u
     \le \frac{2}{r}\int_{\p B_{R+t}}u
$$
for $r/2\le t\le r$. The first inequality follows 
by integrating this inequality
over the interval $r/2\le t\le r$.
The second inequality was proved in~\cite[Lemma~7.3]{DS2}.
\end{proof}

\begin{lemma}\label{le:apriori}
For every triple $(\Om,K,\lambda)$ as above 
and every $C_0>0$ there exist constants 
$\eps_0=\eps_0(\Om,K,\lambda)>0$
and $c=c(\Om,K,\lambda,C_0)>0$ such that the following holds 
for every $\eps\in(0,\eps_0]$. If $u$, $\Phi$, and $\Psi$ 
satisfy~(\ref{eq:eps-loc}) and
$$
     \|v_s\|_{L^\infty(\Om)}
     + \eps^{-1}\|\mu(u)\|_{L^\infty(\Om)}
     \le C_0
$$
then 
\begin{eqnarray}\label{eq:apriori-p}
&&
     \eps^{-1}\|\mu(u)\|_{L^p(K)}
     + \|L_u^*v_s\|_{L^p(K)}
     + \|L_u^*Jv_s\|_{L^p(K)} 
     \nonumber \\
&&
     +\, \eps\|\Nabla{A,s}v_s\|_{L^p(K)}
     + \eps\|\Nabla{A,t}v_s\|_{L^p(K)}  \\
&&
     \le c\eps^{2/p}\left(
     \|v_s\|_{L^2(\Om)} + \eps^{-1}\|\mu(u)\|_{L^2(\Om)}
     \right)
     \nonumber
\end{eqnarray}
for $2\le p\le\infty$. 
\end{lemma}

\begin{proof} 
Consider the functions $u_0,v_0:\Om\to\R$ given by
$$
     u_0 := \frac{1}{2} \left(
     \left|v_s\right|^2 + \frac{\lambda^2}{\eps^2}
     \left|\mu(u)\right|^2
     \right),
$$
$$
     v_0
:= 
     \frac12\Biggl(
     \left|\Nabla{A,s}v_s\right|^2
     + \left|\Nabla{A,t}v_s\right|^2
     + \frac{\lambda^4}{\eps^4}|L_u\mu(u)|^2
     +\, \frac{\lambda^2}{\eps^2}\left|L_u^*v_s\right|^2
     + \frac{\lambda^2}{\eps^2}\left|L_u^*Jv_s\right|^2 
     \Biggr).
$$
We prove that there exists a constant $c_0>0$ such that
\begin{equation}\label{eq:u0}
     \Delta u_0 \ge v_0 - c_0u_0.
\end{equation}
To see this, recall from the proof of Lemma~\ref{le:mu}
that
$$
     \frac12\Delta|\mu(u)|^2
     = |L_u^*v_s|^2 + |L_u^*Jv_s|^2
     + \frac{\lambda^2}{\eps^2}|L_u\mu(u)|^2
     - 2\inner{\mu(u)}{\rho(v_s,v_t)},
$$
and hence 
\begin{eqnarray*}
&&
     \frac{1}{2\eps^2}\Delta(\lambda^2|\mu(u)|^2)   \\
&&=
     \frac{\Delta\lambda^2}{2\eps^2}|\mu(u)|^2
     + \frac{\lambda^2}{2\eps^2}\Delta|\mu(u)|^2  
     + \frac{\p_s\lambda^2}{\eps^2}\p_s|\mu(u)|^2 
     + \frac{\p_t\lambda^2}{\eps^2}\p_t|\mu(u)|^2 \\
&&= 
     \frac{\Delta\lambda^2}{2\eps^2}|\mu(u)|^2
     + \frac{2\p_t\lambda^2}{\eps^2}\inner{\mu(u)}{L_u^*v_s}
     - \frac{2\p_s\lambda^2}{\eps^2}\inner{\mu(u)}{L_u^*Jv_s}  \\
&&\quad
     +\,\frac{\lambda^2}{\eps^2}|L_u^*v_s|^2 
     + \frac{\lambda^2}{\eps^2}|L_u^*Jv_s|^2
     + \frac{\lambda^4}{\eps^4}|L_u\mu(u)|^2 
     - \frac{2\lambda^2}{\eps^2}\inner{\mu(u)}{\rho(v_s,v_t)}.
\end{eqnarray*}
Moreover, by Lemma~\ref{le:curv} and Lemma~\ref{le:Lu-phi},
\begin{eqnarray}\label{eq:deltavs}
&&
     \left(\Nabla{A,s}\Nabla{A,s}
     + \Nabla{A,t}\Nabla{A,t}\right)v_s   
     \nonumber \\
&&=
     (\Nabla{A,t}\Nabla{A,s}-\Nabla{A,s}\Nabla{A,t})v_t  
     + \Nabla{A,s}(\Nabla{A,s}v_s+\Nabla{A,t}v_t) 
     - \Nabla{A,t}(\Nabla{A,s}v_t-\Nabla{A,t}v_s)
     \nonumber \\
&&=
     -\, R(v_s,v_t)v_t - \Nabla{v_t}X_\kappa(u) 
     - \Nabla{A,s}(JL_u\kappa) - \Nabla{A,t}(L_u\kappa)
     \nonumber \\
&&= 
     -\, R(v_s,v_t)v_t 
     + \frac{\lambda^2}{\eps^2}J\Nabla{v_s}X_{\mu(u)}(u)
     + \frac{\lambda^2}{\eps^2}(\Nabla{v_s}J)L_u\mu(u)  
     + \frac{2\lambda^2}{\eps^2}\Nabla{v_t}X_{\mu(u)}(u)
     \nonumber \\
&&\quad
     +\, \frac{\lambda^2}{\eps^2}L_uL_u^*v_s 
     - \frac{\lambda^2}{\eps^2}JL_uL_u^*Jv_s
     + \frac{\p_s\lambda^2}{\eps^2}JL_u\mu(u) 
     + \frac{\p_t\lambda^2}{\eps^2}L_u\mu(u).
\end{eqnarray}
Hence
\begin{eqnarray*}
    \frac12 \Delta |v_s|^2
&=&
    |\Nabla{A,s}v_s|^2+|\Nabla{A,t}v_s|^2
    + \inner{v_s}{\left(\Nabla{A,s}\Nabla{A,s}
       + \Nabla{A,t}\Nabla{A,t}\right)v_s}  \\
&=&
     |\Nabla{A,s}v_s|^2  
     + |\Nabla{A,t}v_s|^2
     + \frac{\lambda^2}{\eps^2} |L_u^*v_s|^2
     + \frac{\lambda^2}{\eps^2} |L_u^*J v_s|^2  \\
&&
     -\,\frac{3\lambda^2}{\eps^2}
     \inner{\mu(u)}{\rho(v_s,v_t)}
     +\frac{\lambda^2}{\eps^2}
     \inner{v_s}{(\Nabla{v_s}J)L_u\mu(u)} \\
&&
     -\, \inner{v_s}{R(v_s,v_t)v_t}
     - \frac{\p_s\lambda^2}{\eps^2}
     \inner{L_u^* Jv_s}{\mu(u)}
     + \frac{\p_t\lambda^2}{\eps^2}
     \inner{L_u^* v_s}{\mu(u)}.
\end{eqnarray*}
Combining this with the formula for 
$\Delta (\lambda^2 |\mu(u)|^2)/2\eps^2$ 
we obtain
\begin{eqnarray}\label{eq:Delta-u0}
     \Delta u_0
&=&
     |\Nabla{A,s}v_s|^2+|\Nabla{A,t}v_s|^2
     +\frac{2\lambda^2}{\eps^2} |L_u^*v_s|^2
     +\frac{2\lambda^2}{\eps^2} |L_u^*Jv_s|^2
     +\frac{\lambda^4}{\eps^4}|L_u\mu(u)|^2 
     \nonumber \\
&&
     -\,\frac{5\lambda^2}{\eps^2}
     \inner{\mu(u)}{\rho(v_s,v_t)}
     + \frac{3\p_t\lambda^2}{\eps^2}
     \inner{\mu(u)}{L_u^*v_s}
     - \frac{3\p_s\lambda^2}{\eps^2}
     \inner{\mu(u)}{L_u^*Jv_s} 
     \nonumber \\
&&
     +\,\frac{\Delta\lambda^2}{2\eps^2} |\mu(u)|^2
     +\frac{\lambda^2}{\eps^2}    
     \inner{v_s}{(\Nabla{v_s}J)L_u\mu(u)}
     -\inner{v_s}{R(v_s,v_t)v_t}.
\end{eqnarray}
The first row on the right is bounded below by $2v_0$.
Moreover, by assumption, the image of $u$ is contained
in the compact set $\{|\mu(x)|\le\eps C_0\}$. 
Hence the last six terms can be estimated from below by 
$v_0-c_0u_0$ for some constant $c_0$ whenever $\eps$
is sufficiently small.  Thus we have proved the 
inequality~(\ref{eq:u0}).  Hence, by Lemma~\ref{le:apriori1},
there exist constants $\eps_0>0$ and $c_0'>0$ such that 
$$
     \sup_Ku_0 + \int_Kv_0\le c_0'\int_\Om u_0
$$
for every $\eps\in(0,\eps_0]$.  Since $|\mu(u)|\le C_0\eps_0$ 
and zero is a regular value of $\mu$ there is an inequality 
$|L_u\eta|\ge\delta|\eta|$ whenever $\eps_0$ is sufficiently small. 
Thus we have proved~(\ref{eq:apriori-p}) for $p=2$
as well as
\begin{equation}\label{eq:apriori-inf}
     \left\|v_s\right\|_{L^\infty(K)}
     + \eps^{-1}\left\|\mu(u)\right\|_{L^\infty(K)}
     \le c_0'\left(\left\|v_s\right\|_{L^2(\Om)}
     + \eps^{-1}\left\|\mu(u)\right\|_{L^2(\Om)}\right).
\end{equation}
Now let us define $u_1:\Om\to\R$ by 
$$
     u_1 := \frac{1}{2}\left|\Nabla{A,s}v_s\right|^2.
$$
We shall prove that there exist positive constants 
$\delta_1$, $c_1$, and $\eps_0$ such that 
\begin{equation}\label{eq:u1}
     \Delta(u_0+\eps^2u_1) \ge -c_1u_0
\end{equation}
for $0<\eps\le\eps_0$. 
We consider the equation 
$$
   \Delta u_1
   = |\Nabla{A,s}\Nabla{A,s}v_s|^2
     + |\Nabla{A,t}\Nabla{A,s}v_s|^2 
     + \inner{\left(\Nabla{A,s}\Nabla{A,s}
       + \Nabla{A,t}\Nabla{A,t}\right)\Nabla{A,s}v_s}
       {\Nabla{A,s}v_s}
$$
and use the formula
$$
     \left(\Nabla{A,s}\Nabla{A,s}
     + \Nabla{A,t}\Nabla{A,t}\right)\Nabla{A,s}v_s
     = I + II + III,
$$
where
\begin{eqnarray*}
     I 
&:= &
     \left(\Nabla{A,t}\Nabla{A,s}
       - \Nabla{A,s}\Nabla{A,t}\right)\Nabla{A,t}v_s \\
&= &
     - R(v_s,v_t)\Nabla{A,t}v_s
     + \frac{\lambda^2}{\eps^2}\Nabla{\Nabla{A,t}v_s}X_{\mu}(u), \\
     II 
&:= &
     \Nabla{A,t}\left(\Nabla{A,t}\Nabla{A,s}
     - \Nabla{A,s}\Nabla{A,t}\right)v_s \\
&= &
     - \Nabla{A,t}R(v_s,v_t)v_s 
     + \frac{\lambda^2}{\eps^2}\Nabla{A,t}(\Nabla{v_s}X_{\mu}(u))
     +\frac{\p_t\lambda^2}{\eps^2}\Nabla{v_s}X_{\mu}(u), \\
     III 
&:= &  
     \Nabla{A,s}\left(\Nabla{A,s}\Nabla{A,s}
     + \Nabla{A,t}\Nabla{A,t}\right)v_s  \\
&= &
     \Nabla{A,s}\Biggl(
     - R(v_s,v_t)v_t
     + \frac{\lambda^2}{\eps^2}L_uL_u^*v_s
     - \frac{\lambda^2}{\eps^2}JL_uL_u^*Jv_s
     + \frac{\lambda^2}{\eps^2}(\Nabla{v_s}J)X_\mu(u) \\
&&
     + \frac{\lambda^2}{\eps^2}J\Nabla{v_s}X_{\mu}(u) 
     + \frac{2\lambda^2}{\eps^2}\Nabla{v_t}X_{\mu}(u) 
     + \frac{\p_s\lambda^2}{\eps^2}JX_\mu(u)
     + \frac{\p_t\lambda^2}{\eps^2}X_\mu(u)
     \Biggr).
\end{eqnarray*}
Here we abbreviate $X_\mu(u)=X_{\mu(u)}(u)=L_u\mu(u)$.
The last equality for $III$ follows from~(\ref{eq:deltavs}).
Now consider the tensors $\nabla^2J$ and $\nabla^2X_\eta$
defined by 
\begin{eqnarray*}
    \nabla^2J(X,Y,Z)
&:= &
    \Nabla{X}((\Nabla{Y}J)Z) - (\Nabla{\Nabla{X}Y}J)Z
    - (\Nabla{Y}J)\Nabla{X}Z,\\ 
    \nabla^2X_\eta(Y,Z)
&:= &
    \Nabla{Y}(\Nabla{Z}X_\eta) - \Nabla{\Nabla{Y}Z}X_\eta
\end{eqnarray*}
for $\eta\in\g$ and $X,Y,Z\in\Vect(M)$. Then 
\begin{eqnarray*}
    \Nabla{A,s}((\Nabla{v_s}J)L_u\mu(u))
&= &
    \nabla^2J(v_s,v_s,L_u\mu(u))
    + (\Nabla{\Nabla{A,s}v_s}J)L_u\mu(u) \\
&&
    +\, (\Nabla{v_s}J)\Nabla{v_s}X_{\mu(u)}(u) 
    - (\Nabla{v_s}J)L_uL_u^*Jv_s, \\
    \Nabla{A,t}(\Nabla{v_s}X_{\mu(u)}(u))
&= &
    \nabla^2X_{\mu(u)}(v_t,v_s) 
    + \Nabla{\Nabla{A,t}v_s}X_{\mu(u)}(u)
    \Nabla{v_s}X_{L_u^*v_s}(u), \\
    \Nabla{A,s}(R(v_s,v_t)v_t) 
&= &
    \nabla R(v_s,v_s,v_t,v_t) 
    + R(\Nabla{A,s}v_s,v_t)v_t \\
&&
    + R(v_s,\Nabla{A,s}(Jv_s))v_t)
    + R(v_s,v_t)\Nabla{A,s}(Jv_s).   
\end{eqnarray*}
Hence, by a term by term inspection, we obtain
an inequality
$$
    \eps^2\inner{\Nabla{A,s}v_s}{I+II+III}
    \ge - cu_0 - v_0
$$
for $\eps>0$ sufficiently small. 
Note, in particular, that the term
$\eps^2\inner{\Nabla{A,s}v_s}{III}$
contains the two positive summands
$\lambda^2|L_u^*\Nabla{A,s}v_s|^2$ and 
$\lambda^2|L_u^*J\Nabla{A,s}v_s|^2$. 
Since $\Delta u_0\ge v_0-c_0u_0$ the
last inequality implies~(\ref{eq:u1})
with $c_1:=c+c_0$. Now it follows 
from~(\ref{eq:u1}), (\ref{eq:apriori-inf}), 
Lemma~\ref{le:apriori1}, and the formula 
$$
    \Nabla{A,t}v_s = - J\Nabla{A,s}v_s - (\Nabla{v_s}J)v_s
    + \frac{\lambda^2}{\eps^2}L_u\mu(u)
$$
that~(\ref{eq:apriori-p}) holds for $p=\infty$. 
The estimate~(\ref{eq:apriori-p}) for $2<p<\infty$ 
follows by interpolation. 
\end{proof}


\section{Proof of Theorem~D}\label{sec:D}

\begin{theorem}\label{thm:D}
Assume~$(H1)$ and~$(H4)$, let $\bar B\in H_2(\bar M;\Z)$
be a nontorsion homology class, and denote 
$B:=\kappa(\bar B)\in H_2(M_\G;\Z)$. 
Then, for every $C>0$, there exist positive constants 
$\eps_0$ and $c_0$ such that for every 
$\eps\in (0,\eps_0]$ the following holds.
If $(u,A)\in\Tilde\Mcal_{B,\Sigma}^\eps$
satisfies
$
     \left\|d_Au\right\|_{L^\infty}\le C
$
and
$
     u(P)\subset M^C
$
then 
$
     (u,A)\in\Tilde{\Tcal^\eps}
     (\Tilde\Mcal_{B,\Sigma}^0(c_0)).
$
\end{theorem}

\begin{proof}
Suppose the assertion is false.
Then there exist a constant $C>0$
and sequences $\eps_i\to0$ and 
$
   (u_i,A_i)\in\Tilde{\Mm}^{\eps_i}_{B,\Sigma}
$
such that
$$
    \left\|d_{A_i}u_i\right\|_{L^\infty}\le C,\qquad
    u_i(P)\subset M^C,\qquad
    (u_i,A_i)\notin
    \Tilde{\Tt}^{\eps_i}(\Mcal_{B,\Sigma}^0(i)).
$$
Here $\eps_i$ is chosen smaller than 
the number $\eps_0(i)$ required for the definition
of the map $\Tt^{\eps_i}$. 
We prove in four steps that there exist 
an integer $i_0\in\NN$, positive constants 
$c$ and $c_0$, and sequences 
$$
    (u_{i0},A_{i0})
    \in\Tilde{\Mm}_{B,\Sigma}^0(c_0-1),\qquad
    \zeta_{i0}=(\xi_{i0},\alpha_{i0})\in T_{(u_{i0},A_{i0})}\Bb
$$
such that  
\begin{equation}\label{eq:C}
    (u_i,A_i)=(\exp_{u_{i0}}(\xi_{i0}),A_{i0}+\alpha_{i0}),\qquad
    \left\|\zeta_{i0}\right\|_{1,p,\eps_i}
    \le c\eps_i^{2/p+1}
\end{equation}
for every $i\ge i_0$. For $i$ sufficiently large it 
then follows from Theorem~\ref{thm:loc-onto} that
$
    (u_i,A_i)
    \in\Tilde{\Tt}^{\eps_i}(\Tilde{\Mm}^0_{B,\Sigma}(c_0)),
$
in contradiction to our assumption. 

\medskip
\noindent{\bf Step~1.}
{\it There exist constants $c>0$ and 
$i_0\in\NN$ such that
$$
     \eps_i^{-1}\left\|\mu(u_i)\right\|_{L^p}
     + \left\|L_{u_i}^*d_{A_i}u_i\right\|_{L^p}
     + \left\|L_{u_i}^*Jd_{A_i}u_i\right\|_{L^p} 
     + \eps_i\left\|{\Nabla{A_i}}^*d_{A_i}u_i\right\|_{L^p}
     \le c\eps_i^{2/p}
$$
for $i\ge i_0$ and $2\le p\le\infty$. 
}

\medskip
\noindent
By the graph construction in Appendix~\ref{app:graph},
it suffices to establish the estimate
under the hypothesis that $J$ is independent of
$z\in\Sigma$. Namely, 
$$
     L_{\tilde u_i}^*d_{A_i}\tilde u_i=L_{u_i}^*d_{A_i}u_i,\qquad
     L_{\tilde u_i}^*\tilde Jd_{A_i}\tilde u_i=L_{u_i}^*Jd_{A_i}u_i,
$$
$$
     {\Tabla{A_i}}^*d_{A_i}\tilde u_i
     = ({\Nabla{\Sigma}}^*\id,{\Nabla{A_i}}^*d_{A_i}u_i),
$$
where $\tilde u_i=(\pi,u_i):P\to\tilde M=\Sigma\times M$,
$\Tabla{A}$ is the connection induced by 
$A$ on $\tilde u_i^*T\tilde M/\G$, and $\id\in\Om^1(\Sigma,T\Sigma)$.
Hence we can use the results of Section~\ref{sec:apriori}. 
Since $\left\|d_{A_i}u_i\right\|_{L^\infty}\le C$ 
and $u_i(\Om)\subset M^C$, the pair
$(\tilde u_i,A_i)$ satisfies the hypotheses of 
Lemma~\ref{le:mu} and so the sequence
${\eps_i}^{-3/2}\tilde\mu(\tilde u_i)$ is uniformly bounded.
Hence there exists a constant $c_0>0$ such that 
\begin{equation}\label{eq:linfty-i}
     \left\|d_{A_i}u_i\right\|_{L^\infty}
     + {\eps_i}^{-3/2}\left\|\mu(u_i)\right\|_{L^\infty}
     \le c_0
\end{equation}
for every $i$.  
This implies that, in local holomorphic coordinates on $\Sigma$, 
the pair $(\tilde u_i,A_i)$ satisfies the hypotheses 
of Lemma~\ref{le:apriori} for $i$ sufficiently large.  
Hence the estimate holds in local holomorphic 
coordinates on $\Sigma$ with $u_i$ replaced by $\tilde u_i$. 
Hence, by a partition of unity argument, it holds globally.

\medskip
\noindent{\bf Step~2.}
{\it There exists an integer $i_0\in\NN$
and a constant $c>0$ such that,
for every $i\ge i_0$, there exists 
a unique $\eta_i\in\Om^0(\Sigma,\g_P)$ such that
$$
     \mu(\exp_{u_i}(JL_{u_i}\eta_i)) = 0,\qquad 
     \left\|\eta_i\right\|_{L^\infty}
     \le c\left\|\mu(u_i)\right\|_{L^\infty}.
$$
Define $u_i':P\to M$ and $A_i'\in\Aa(P)$ by 
$$
     u_i' := \exp_{u_i}(JL_{u_i}\eta_i),\qquad
     L_{u_i'}^*d_{A_i'}u'_i = 0,
$$
so that $d_{A_i'}u'_i\in\Om^1(\Sigma,H_{u'_i})$,
and let $\zeta_i:=(JL_{u_i}\eta_i,A'_i-A_i)$.
Then there exists a constant $c'>0$ such that 
$$
     \left\|\zeta_i\right\|_{1,p,\eps_i}\le c'\eps_i^{1+2/p},\qquad
     \left\|\bar\p_{J,A_i'}(u'_i)\right\|_{L^p}\le c'\eps_i^{1+2/p},\qquad
     \left\|d_{A_i'}u'_i\right\|_{L^\infty}\le c'
$$
for every $i\ge i_0$.}

\medskip
\noindent
The existence of $\eta_i$ for large $i$ follows
from the implicit function theorem for the map
$\eta\mapsto\mu(\exp_{u_i(p)}(JL_{u_i(p)}\eta))$. 
This sequence satisfies an estimate of the form
$$
     \left\|\eta_i\right\|_{L^p} 
     \le c_1\left\|\mu(u_i)\right\|_{L^p}
     \le c_2\eps^{1+2/p}
$$
for every $i\ge i_0$ and every $p\in[2,\infty]$.
Here the constants $c_1$ and $c_2$ are independent of
$i$ and $p$, and the second inequality follows from Step~1.
For $p=\infty$ there is actually a better estimate
(by $\eps^{3/2}$ instead of $\eps$), but we shall
not use this here. In the following we suppress the 
subscript~$i$ and write $u,u',A,A',\eps$
instead of $u_i,u'_i,A_i,A'_i,\eps_i$, respectively.
We establish the required estimates in 
local holomorphic coordinates on $\Sigma$.
As in Remark~\ref{rmk:local}, we write 
$A'=\Phi'\,ds+\Psi'\,dt$ for some 
Lie algebra valued functiona $\Phi'$ 
and $\Psi'$, and denote
$$
     v_s'
     :=\p_su'+X_{\Phi'}(u'),
     \qquad
     v_t'
     :=\p_tu'+X_{\Psi'}(u').
$$
Then
$
     L_{u'}^*v_s'
     =
     L_{u'}^*v_t'
     =0.
$ 
We assume that the functions $u,u',\Phi,\Psi,\Phi',\Psi'$
are defined on an open set $\Om\subset\C$ and fix 
any compact subset $K\subset\Om$. 
We must prove the estimates
$$
     \left\|(\xi,\phi,\psi)\right\|_{1,p,\eps}
     \le c\eps^{1+2/p},\quad
     \|v'_s+Jv'_t\|_{L^p}
     \le c \eps^{1+2/p},\quad
     \left\|v'_s\right\|_{L^\infty}
     + \left\|v'_t\right\|_{L^\infty}\le c   
$$
on the subset $K$, where 
$$
     \xi:=JL_u\eta,\qquad 
     \phi:=\Phi'-\Phi,\qquad
     \psi:=\Psi'-\Psi.
$$
Abbreviate $E_i:=E_i(u,JL_u\eta)$, $i=1,2$.
Then
\begin{equation}\label{eq:nablaxi}
     \Nabla{A,t}\xi
     = JL_u\Nabla{A,t}\eta
       + (\Nabla{v_t}J+\p_tJ)L_u\eta
       + J\Nabla{v_t}X_\eta(u).
\end{equation}
Hence, by Lemma~\ref{le:Lu-phi},
\begin{equation}\label{eq:v'-v}
     v'_t - E_1v_t
     = L_{u'}\psi + E_2\Nabla{A,t}\xi
     = L_{u'}\psi + E_2JL_u\Nabla{A,t}\eta + R_t\eta,
\end{equation}
where 
$$
     R_t\eta := E_2(\Nabla{v_t}J+\p_tJ)L_u\eta 
        + E_2J\Nabla{v_t}X_\eta(u).
$$
Hence
$$
     d\mu(u')(v'_t-E_1v_t)
     = d\mu(u')E_2JL_u\Nabla{A,t}\eta + d\mu(u')R_t\eta.
$$
Since $L_u^*=d\mu(u)J$ we have
\begin{eqnarray*}
     L_u^*L_u\Nabla{A,t}\eta
&=&
     (d\mu(u)-d\mu(u')E_2)JL_u\Nabla {A,t}\eta \\
&&
     +\,d\mu(u')(v'_t-E_1v_t) - d\mu(u')R_t\eta,
\end{eqnarray*}
and, since $d\mu(u')v'_t=0$,
\begin{eqnarray*}
     L_u^*L_u\Nabla{A,t}\eta
&=&
     (d\mu(u)-d\mu(u')E_2)JL_u\Nabla{A,t}\eta  \\
&&
     +\, (d\mu(u)-d\mu(u')E_1)v_t + L_u^*Jv_t
     - d\mu(u')R_t\eta.
\end{eqnarray*}
It follows that
\begin{eqnarray*}
     \left\|\Nabla{A,t}\eta\right\|_{L^p}
&\le &
     c_2\left(\left\|\eta\right\|_{L^\infty}
      \left\|\Nabla{A,t}\eta\right\|_{L^p}
     + \left\|L_u^*Jv_t\right\|_{L^p} 
     + \left\|\eta\right\|_{L^p}
     \right) \\
&\le &
     c_3\left(\eps^{1+2/p}
     \left\|\Nabla{A,t}\eta\right\|_{L^p}
     + \eps^{2/p}\right).
\end{eqnarray*}
If $\eps$ is sufficiently small this gives
$$
     \left\|\Nabla{A,t}\eta\right\|_{L^p}\le c_4\eps^{2/p},\qquad
     \left\|\Nabla{A,t}\eta\right\|_{L^\infty}\le c_4.
$$
Here the second inequality follows from a similar argument as the 
first.  Combining these inequalities with~(\ref{eq:nablaxi})
we obtain
$$
     \left\|\Nabla{A,t}\xi\right\|_{L^p} \le c_5\eps^{2/p},\qquad
     \left\|\Nabla{A,t}\xi\right\|_{L^\infty} \le c_5.
$$
In order to estimate $\psi$ we apply 
the operator $L_{u'}^*$ to~(\ref{eq:v'-v}) 
and use the formula $L_{u'}^*v'_t=0$ to obtain
$$
     L_{u'}^*L_{u'}\psi
     = (L_u^*-L_{u'}^*E_1)v_t - L_u^*v_t
       - L_{u'}^*E_2JL_u\Nabla{A,t}\eta
       - L_{u'}^*R_t\eta.
$$
Combining this with Step~1 and the estimate for 
$\Nabla{A,t}\eta$ we obtain
$$
     \left\|\psi\right\|_{L^p} \le c_6\eps^{2/p},\qquad
     \left\|\psi\right\|_{L^\infty} \le c_6.
$$
Hence, by~(\ref{eq:v'-v}),
$$
     \left\|v'_t-E_1v_t\right\|_{L^p}\le c_7\eps^{2/p},\qquad
     \left\|v'_t-E_1v_t\right\|_{L^\infty}\le c_7.
$$
Similarly,
$$
     \left\|\Nabla{A,s}\eta\right\|_{L^p} 
     + \left\|\phi\right\|_{L^p} 
     + \left\|v_s^\prime-E_1v_s\right\|_{L^p}
     \le c_7\eps^{2/p}
$$
and
$$
     \left\|\Nabla{A,s}\eta\right\|_{L^\infty} 
      + \left\|\phi\right\|_{L^\infty} 
      + \left\|v_s^\prime-E_1v_s\right\|_{L^\infty}
     \le c_7.
$$
Now use~(\ref{eq:v'-v}) again to obtain
\begin{eqnarray*}
     v_s'+J(u')v_t'
&= &
     L_{u'}(\phi-\Nabla{A,t}\eta)
     + JL_{u'}(\psi+\Nabla{A,s}\eta)  \\
&& 
     +\,(E_2JL_u-JL_{u'})\Nabla{A,s}\eta
     + J(E_2JL_u-JL_{u'})\Nabla{A,t}\eta \\ 
&&
     +\, E_1v_s + JE_1v_t + (R_s+JR_t)\eta  \\
&= &
     \pi_{u'}\bigl(
       (E_2JL_u-JL_{u'})\Nabla{A,s}\eta
       + J(E_2JL_u-JL_{u'})\Nabla{A,t}\eta
     \bigr) \\ 
&&
     +\,\pi_{u'}\bigl(
        (JE_1-E_1J)v_t + (R_s+JR_t)\eta 
     \bigr).
\end{eqnarray*}
The second equality uses the fact that
$v_s+Jv_t=0$ and that the 1-form $\bar\p_{J,A'}(u')$ 
takes values in $H_{u'}$. It follows that 
$$
     \|v_s'+J(u')v_t'\|_{L^p} \le c_8\eps^{1+2/p}.
$$
It remains to show that
$$
     \left\|\Nabla{A,s}\phi\right\|_{L^p} 
     + \left\|\Nabla{A,t}\phi\right\|_{L^p}
     + \left\|\Nabla{A,s}\psi\right\|_{L^p} 
     + \left\|\Nabla{A,t}\psi\right\|_{L^p}
     \le c\eps^{2/p-1}.
$$
To estimate the term $\Nabla{A,t}\psi$ 
differentiate~(\ref{eq:v'-v}) with respect to $t$. 
Then apply the operator $d\mu(u')$ to the resulting
expression to eliminate $\Nabla{A,t}\psi$ 
and obtain an estimate of the form
$$
     \left\|\Nabla{A,t}\Nabla{A,t}\eta\right\|_{L^p}
     \le c_9\eps^{2/p-1}.
$$
Then apply the operator $L_{u'}^*$ to the equation 
obtained from differentiating~(\ref{eq:v'-v}), 
and estimate $\Nabla{A,t}\psi$ using the 
upper bound found for $\Nabla{A,t}\Nabla{A,t}\eta$.  
The estimate for $\Nabla{A,s}\psi$ 
is obtained in a similar manner. 
To estimate $\Nabla {A,s}\phi$ and $\Nabla{A,t}\phi$, 
we begin with the identity
$$
     v_s' - E_1v_s
     = L_{u'}\phi + E_2JL_u\Nabla{A,s}\eta + R_s\eta
$$
instead of~(\ref{eq:v'-v}) and then follow the same procedure.

\medskip
\noindent{\bf Step~3.}
{\it There exist an integer $i_0\in\NN$, 
a constant $c>0$, and a sequence 
$(u_i'',A_i'')\in\Tilde{\Mm}^0_{B,\Sigma}$
such that
$$
     u_i'' = \exp_{u_i'}(\xi_i'),\qquad
     \xi_i'\in\Om^0(\Sigma,H_{u_i'}),
$$
and
$$
     \left\|\xi'_i\right\|_{W^{1,p}}
     + \left\|A_i''-A_i'\right\|_{L^p}
     \le c\eps_i^{1+2/p},\qquad
     \left\|d_{A_i''}u_i''\right\|_{L^\infty}\le c,
$$
for $i\ge i_0$.
}

\medskip
\noindent
By Step~2, 
$$
     \sup_i \left\|d_{A'_i}u'_i\right\|_{L^\infty} < \infty.
$$
Hence the induced maps $\bar u_i':\Sigma\to\bar M$ 
form a sequence of approximate $\bar J$-holomorphic curves 
which satisfy a uniform $L^\infty$-bound on their 
first derivatives. Hence, by~$(H4)$ and Theorem~\ref{thm:jhol},
there is nearby a true $\bar J$-holomorphic curve
$\bar u_i'':\Sigma\to\bar M$ whose $W^{1,p}$-distance to
$\bar u_i'$ is controlled by the $L^p$-norm of
$\bar\p_{\bar J}(\bar u_i')$.  Now this $\bar J$-holomorphic
curve has a unique lift $u_i'':P\to\mu^{-1}(0)$ 
of the form 
$$
     u_i''=\exp_{u_i'}(\xi_i'),\qquad
     \xi_i'\in\Om^0(\Sigma,H_{u_i'}). 
$$
Let $A_i''\in\Aa(P)$ be the connection determined by $u_i''$ via 
$L_{u_i''}^*d_{A_i''}u_i''=0$.  Then 
$$
     \left\|\xi'_i\right\|_{W^{1,p}}
     \le c_1\left\|\bar\p_{J,A_i'}(u_i')\right\|_{L^p}
     \le c_2{\eps_i}^{1+2/p}.
$$
Here the last inequality follows from Step~2. 
Since $L_{u_i'}^*d_{A_i'}u_i'=L_{u_i''}^*d_{A_i''}u_i''=0$
we obtain
$$
     \left\|A_i''-A_i'\right\|_{L^p}\le c_3\eps_i^{1+2/p}.
$$
In particular, these inequalities together give a uniform 
$W^{1,p}$-bound on the $\bar J$-holomorphic curves 
$\bar u_i'':\Sigma\to\bar M$.  Hence, by the elliptic bootstrapping 
techniques for $J$-holomorphic curves, the sequence 
$\bar u_i''$ satisfies a uniform $L^\infty$ bound on the first 
derivatives.  This proves Step~3.  

Unfortunately, the estimate on $A_i''-A_i'$ in Step~3
is only in the $L^p$-norm and not in the $W^{1,p}$-norm.  
A further modification of the pair $(u_i'',A_i'')$ is 
required to improve this estimate. 

\medskip
\noindent{\bf Step~4.}
{\it There exist an integer $i_0\in\NN$, 
a constant $c>0$, and a sequence 
of gauge transformations $g_i\in\Gg(P)$
such that the sequence
$$
     (u_{i0},A_{i0})
     := g_i^*(u_i'',A_i'')
     \in\Tilde{\Mm}^0_{B,\Sigma}
$$
satisfies the following. For $i\ge i_0$
the original sequence $(u_i,A_i)$ has the form
$$
     (u_i,A_i) = (\exp_{u_{i0}}(\xi_i),A_{i0}+\alpha_i)
$$
where $\zeta_i:=(\xi_i,\alpha_i)\in T_{(u_{i0},A_{i0})}\Bb$
satisfies~(\ref{eq:C}).
}

\medskip
\noindent
The idea is to choose $g_i$ for large $i$ such that 
$$
     u_i = \exp_{u_{i0}}(\xi_i),\qquad
     L_{u_{i0}}^*\xi_i = 0,\qquad
     u_{i0} := g_i^{-1}u_i''.
$$
This can be done by using pointwise, for every $p\in P$, 
the implicit function theorem to obtain the local 
slice condition. This suffices 
to obtain the missing estimates for the 
first derivatives of $g_i^*A_i''-A_i$.
We sketch a proof of this estimate below.

By~(\ref{eq:linfty-i}) and Step~2, the distance 
between $u_i$ and $u_i'$ is unformly bounded 
by a constant times $\eps_i^{3/2}$ 
while the distance in the $W^{1,p}$-norm 
is bounded by a constant times $\eps_i^{2/p}$.  
By Step~3, the distance between $u_i'$ and $u_i''$ 
is bounded in the $W^{1,p}$-norm by a constant 
times $\eps_i^{1+2/p}$. 
Hence there exists a sequence of smooth sections 
$\xi_i\in\Om^0(\Sigma,{u_i}^*TM/\G)$
and a constant $c>0$ such that
\begin{equation}\label{eq:''1}
     u_i'' = \exp_{u_i}(\xi_i),\qquad
     \left\|\xi_i\right\|_{L^\infty}\le c_1\eps_i,\qquad
     \left\|\Nabla{A_i}\xi_i\right\|_{L^p}\le c_1\eps_i^{2/p}.
\end{equation}
Moreover, the sequence $d_{A_i''}u_i''$ is uniformly bounded in 
the $L^\infty$-norm and
\begin{equation}\label{eq:''2}
     \left\|A_i''-A_i\right\|_{L^p}\le c_1\eps^{2/p},\qquad
     \left\|d_{A_i''}u_i'' - E_1(u_i,\xi_i)d_{A_i}u_i\right\|_{L^p}
     \le c_1\eps^{2/p}.
\end{equation}
Here the last inequality follows from the identity
\begin{eqnarray*}
     d_{A_i''}u_i'' 
&= &
     E_1(u_i,\xi_i)d_{A_i}u_i 
     + E_2(u_i,\xi_i)\Nabla{A_i}\xi_i  \\
&&
     +\, E_1(u_i,\xi_i)X_{A_i''-A_i}(u_i) 
     + E_2(u_i,\xi_i)\Nabla{\xi_i}X_{A_i''-A_i}(u_i),
\end{eqnarray*}
which in turn follows from Lemma~\ref{le:group}.
Now, by the inverse function theorem for the map 
$
      \G\times \ker\,L_x^*\to M:(g,\xi)\mapsto g^{-1}\exp_x(\xi),
$
there exists a constant $c_2>0$ and (unique) sequences $g_i\in\Gg(P)$ 
and $\xi_i''\in\Om^0(\Sigma,{u_i''}^*TM/\G)$
such that 
$$
     u_i = g_i^{-1}\exp_{u_i''}(\xi_i''),\qquad
     L_{u_i''}^*\xi_i'' = 0,
$$
and
$$
     \left\|\xi_i''\right\|_{L^\infty}\le c_2\eps_i,\qquad
     \left\|g_i-\one\right\|_{L^\infty}\le c_2\eps_i.
$$
Define
$$
     u_{i0} := g_i^{-1}u_i'',\qquad
     A_{i0} := {g_i}^*A_i'',\qquad
     \xi_{i0} := g_i^{-1}\xi_i'',\qquad
     \alpha_{i0} := A_i-A_{i0}.
$$
We shall prove that the pair $(\xi_{i0},\alpha_{i0})$ 
satisfies~(\ref{eq:C}). To see this note first that
$$
     u_i=\exp_{u_{i0}}(\xi_{i0}),\qquad
     A_i=A_{i0}+\alpha_{i0}.
$$
The endomorphism 
$
     E_1(u_i,\xi_i)g_i^{-1}E_1(u_i'',\xi_i'')
$
of ${u_i''}^*TM$ is $\eps_i$-close to the identity,
$d_{A_{i0}}u_{i0}=g_i^{-1}d_{A_i''}u_i''$, and
\begin{eqnarray*}
&&
     E_1(u_i,\xi_i)(d_{A_i}u_i-E_1(u_{i0},\xi_{i0})d_{A_{i0}}u_{i0}) \\
&&= 
     E_1(u_i,\xi_i)d_{A_i}u_i - d_{A_i''}u_i''
     + (\one-E_1(u_i,\xi_i)g^{-1}E_1(u_i'',\xi_i''))d_{A_i''}u_i''.
\end{eqnarray*}
Hence, by~(\ref{eq:''2}), there is an estimate
$$
    \left\|d_{A_i}u_i-E_1(u_{i0},\xi_{i0})d_{A_{i0}}u_{i0}
    \right\|_{L^p}\le c_3\eps_i^{2/p},\qquad
    \left\|d_{A_{i0}}u_{i0}\right\|_{L^\infty}\le c_3
$$
for all $i$. Hence, by Corollary~\ref{cor:8.12}, 
there exists a constant $c_4>0$ such that 
$$
     \left\|\alpha_{i0}\right\|_{L^p}\le c_4\eps_i^{2/p},\qquad
     \left\|\alpha_{i0}\right\|_{L^\infty}\le c_4.
$$
Next observe that, by Lemma~\ref{le:group}, 
\begin{eqnarray*}
     d_{A_i}u_i
&= &
     E_1(u_{i0},\xi_{i0})d_{A_{i0}}u_{i0} 
     + E_2(u_{i0},\xi_{i0})\Nabla{A_{i0}}\xi_{i0}  \\
&&
     +\, E_1(u_{i0},\xi_{i0})X_{\alpha_{i0}}(u_{i0}) 
     + E_2(u_{i0},\xi_{i0})\Nabla{\xi_{i0}}X_{\alpha_{i0}}(u_{i0})   
\end{eqnarray*}
Hence there exists a constant $c_5>0$ such that
$$
     \left\|\Nabla{A_{i0}}\xi_{i0}\right\|_{L^p}\le c_5\eps_i^{2/p},\qquad
     \left\|\Nabla{A_{i0}}\xi_{i0}\right\|_{L^\infty}\le c_5
$$
for all $i$. Thus we have proved that 
\begin{equation}\label{eq:pinfty}
\begin{array}{rcl}
     \left\|\xi_{i0}\right\|_{L^p}
     + \eps_i\left\|\alpha_{i0}\right\|_{L^p}
     + \eps_i\left\|\Nabla{A_{i0}}\xi_{i0}\right\|_{L^p}
&\le &
     c_6\eps_i^{1+2/p}, \\
     \left\|\xi_{i0}\right\|_{L^\infty}
     + \eps_i\left\|\alpha_{i0}\right\|_{L^\infty}
     + \eps_i\left\|\Nabla{A_{i0}}\xi_{i0}\right\|_{L^\infty}
&\le &
     c_6\eps_i
\end{array}
\end{equation}
for all $i$.  It remains to estimate the $L^p$-norm
of the first derivatives of $\alpha_{i0}$. 
For this we shall drop the subscript $i$ 
and write $u,u_0,A,A_0,\xi_0,\alpha_0$
instead of $u_i,u_{i0},A_i,A_{i0},\xi_{i0},\alpha_{i0}$.
Moreover, we use local coordinates on $\Sigma$ as in Step~2 
and write
$$
\begin{array}{rclrcl}
     A_0 &:= & \Phi_0\,ds + \Psi_0\,dt,&
     A &:= & \Phi\,ds + \Psi\,dt, \\
     v_{0s} &:= & \p_su_0 + L_{u_0}\Phi_0,&
     v_{s} &:= & \p_su + L_{u}\Phi, \\
     v_{0t} &:= & \p_tu_0 + L_{u_0}\Psi_0,&
     v_{t} &:= & \p_tu + L_{u}\Psi,
\end{array}
$$
and  $\phi_0:=\Phi-\Phi_0$ and $\psi_0:=\Psi-\Psi_0$ . 
Consider the formula
$$
     \Nabla{A,s}v_s+J\Nabla{A,t}v_s
     = -(\Nabla{v_t}J)v_s + \frac{\lambda^2}{\eps^2}JL_u\mu(u).
$$
By Step~1, we have
$$
    \left\|\Nabla{A,s}v_s+J\Nabla{A,t}v_s\right\|_{L^p}
    \le c_7\eps^{2/p-1}
$$
and hence, by elliptic regularity 
for the Cauchy-Riemann operator,
\begin{equation}\label{eq:dvs}
    \left\|\Nabla{A,s}v_s\right\|_{L^p}
    + \left\|\Nabla{A,t}v_s\right\|_{L^p}
    \le c_8\eps^{2/p-1}.
\end{equation}
Moreover, since  $L_{u_0}^*\xi_0=0$, 
it follows from Lemma~\ref{le:rho} that
\begin{eqnarray*}\label{eq:Nabla-rho}
     L_{u_0}^*\Nabla{A_0,t}\Nabla{A_0,s}\xi_0
&=&
     \Nabla{A_0,t}(L_{u_0}^*\Nabla{A_0,s}\xi_0)
     -\rho(v_{0t},\Nabla{A_0,s}\xi_0)  \\
&=&
     -\Nabla{A_0,t}\rho(v_{0s},\xi_0)
     -\rho(v_{0t},\Nabla{A_0,s}\xi_0).
\end{eqnarray*}
and hence 
\begin{equation}\label{eq:ddxi0}
     \left\|L_{u_0}^*\Nabla{A_0,t}\Nabla{A_0,s}\xi_0\right\|
     \le c_{10}\eps^{2/p}.
\end{equation}
Here we use the fact that, 
by elliptic bootstrapping for $\bar J$-holomorphic 
curves, there is a uniform $L^p$-bound on 
$\Nabla{A_0,t}v_{0s}$.  
Now consider the pointwise inequality
\begin{eqnarray*} 
     |\Nabla{A,t}\phi_0| 
&\le & 
     c_{11}|L_u^*L_u\Nabla{A,t}\phi_0|  \\
&\le &
     c_{11}|(L_u^*-L_{u_0}^*E_2^{-1})L_u\Nabla{A,t}\phi_0|
     + c_{11}|L_{u_0}^*E_2^{-1}L_u\Nabla{A,t}\phi_0|
\end{eqnarray*}
Since the operator $(L_u^*-L_{u_0}^*E_2^{-1})L_u$ is small,
we obtain 
$$
     \left\|\Nabla{A,t}\phi_0\right\|_{L^p}
     \le c_{12}\left\|L_{u_0}^*E_2^{-1}L_u\Nabla{A,t}\phi_0\right\|_{L^p}.
$$
Now use Lemma~\ref{le:1} and the estimates~(\ref{eq:pinfty}),
(\ref{eq:dvs}), and~(\ref{eq:ddxi0}) to obtain 
$$
     \left\|\Nabla{A,t}\phi_0\right\|_{L^p}\le c_{13}\eps^{2/p-1}. 
$$
The terms $\|\Nabla{A,s}\phi_0\|_{L^p}$,
$\|\Nabla{A,t}\psi_0\|_{L^p}$,
and $\|\Nabla{A,s}\phi_0\|_{L^p}$
are estimated similarly.  This proves Step~4.

It follows from Step~4 and Theorem~\ref{thm:loc-onto}
that $(u_i,A_i)\in\Tilde{\Tt}^\eps_i(\Tilde{\Mm}^0_{B,\Sigma}(c_0)$
for some constant $c_0$ and $i$ sufficiently large.
This contradicts our assumption and hence proves the 
theorem. 
\end{proof}


\section{Vortices}\label{sec:vortex}

In this section we examine the finite energy solutions
of~(\ref{eq:1}) over the complex plane $\Sigma=\C$.
The equations have the form
\begin{equation}\label{eq:vortex}
\begin{array}{rcl}
     \p_su+X_\Phi(u) + J(\p_tu+X_\Psi(u)) &= &0,\\
     \p_s\Psi - \p_t\Phi + [\Phi,\Psi] + \mu(u) &= &0,
\end{array}
\end{equation}
where $u:\C\to M$ and $\Phi,\Psi:\C\to\g$. 
The energy of the triple $(u,\Phi,\Psi)$ is given by
$$
     E(u,\Phi,\Psi)
     := \int_\C\left(
        \left|\p_su+X_\Phi(u)\right|^2
        + \left|\mu(u)\right|^2
        \right)\,ds dt.
$$
The vortex equations~(\ref{eq:vortex}) and the energy 
are invariant under the action of the gauge group 
$$
     \Gg:=\Cinf(\C,\G)
$$ 
by 
$$
     g^*(u,\Phi,\Psi)
     := (g^{-1}u,g^{-1}\p_sg+g^{-1}\Phi g,
         g^{-1}\p_tg+g^{-1}\Psi g).
$$
A solution of~(\ref{eq:vortex}) is said to be 
in {\bf radial gauge} if 
$$
     \cos\theta\,\Phi(re^{i\theta}) 
     + \sin\theta\,\Psi(re^{i\theta}) = 0
$$
for every $\theta\in\R$ and every sufficiently large $r\ge 0$. 
It is said to be {\bf bounded} if $\sup_\C|\mu(u)|<\infty$. 

\begin{proposition}\label{prop:vortex}
Assume~$(H1)$ and~$(H2)$.
Suppose that $(u,\Phi,\Psi)$ is a 
smooth bounded finite energy solution of~(\ref{eq:vortex}) 
in radial gauge. Then there exists a $W^{1,2}$-function
$x:\R/2\pi\Z\to M$ and an $L^2$-function $\eta:\R/2\pi\Z\to\g$ 
such that 
\begin{equation}\label{eq:xeta}
     \dot x + X_\eta(x) = 0,\qquad \mu(x)=0
\end{equation}
and 
\begin{equation}\label{eq:lim}
     \lim_{r\to\infty}\sup_{\theta\in\R}
     d(u(re^{i\theta}),x(\theta)) = 0,\qquad
     \lim_{r\to\infty}\int_0^{2\pi}
     \left|\eta(\theta)-\eta_r(\theta)\right|^2\,d\theta
     = 0,
\end{equation}
where 
$
     \eta_r(\theta) := r\cos(\theta)\Psi(re^{i\theta}) 
     - r\sin\theta\Phi(re^{i\theta}).
$
Moreover, there exists a constant $\delta>0$ such that
$$
     \lim_{r\to\infty}\sup_{\theta\in\R}
     r^{2+\delta}\Bigl(\left|\p_su+X_\Phi(u)\right|^2
       + \left|\mu(u)\right|^2
      \Bigr)
     = 0,
$$
where $s+it=:re^{i\theta}$, 
$$
     E(u,\Phi,\Psi) = \int_\C u^*\om,
$$
and $\sup_\C(f\circ u)\le c$, where $c$ is as in
hypothesis~$(H2)$ and $\mu^{-1}(0)\subset f^{-1}([0,c])$. 
If~$(H3)$ holds then $E(u,\Phi,\Psi)$ is an integer multiple
of $\hbar=\tau N$. 
\end{proposition}

Note that the removable singularity theorem for $J$-holomorphic 
curves is a corollary of Proposition~\ref{prop:vortex} (consider 
the special case $\G=\{\one\}$ and $M=\bar M$).
Before entering into the proof we introduce the notion of the 
local equi\-variant symplectic action.  The definition of this 
local action functional relies on the following lemma. 
We identify $S^1\cong\R/2\pi\Z$.

\begin{lemma}\label{le:action1}
Assume~$(H1)$. Then there exist positive 
constants $\delta$ and $c$ such that, for every pair 
of smooth loops $x:S^1\to M$ and $\eta:S^1\to\g$ 
such that 
$$
     \sup_{S^1}\left|\mu(x)\right| < \delta,
$$
there exists a point $x_0\in\mu^{-1}(0)$ and a smooth 
loop $g_0:S^1\to\G$ such that 
$$
     c^{-1}\sup_{S^1}\left|\eta+\dot g_0{g_0}^{-1}\right|
     \le \ell(x,\eta) := \int_0^{2\pi}
     \left|\dot x+X_\eta(x)\right|\,d\theta,
$$
and
$$
     d(x(\theta),g_0(\theta)x_0)
     \le c\left(\left|\mu(\theta)\right| + \ell(x,\eta)\right)
$$
for every $\theta\in S^1$. 
\end{lemma}

\begin{proof}
Fix a $\G$-invariant and $\om$-compatible almost complex structure
$J\in\Jj_\G(M,\om)$. If $\delta$ is sufficiently small
then there exist unique loops $x_0:S^1\to\mu^{-1}(0)$ 
and $\eta_0:S^1\to\g$ such that 
$$
     x(\theta) = \exp_{x_0(\theta)}(\xi_0(\theta)),\qquad
     \xi_0 := JL_{x_0}\eta_0,\qquad |\xi_0|\le c_1|\mu(x)|.
$$
By Lemmata~\ref{le:group} and~\ref{le:rho}, we have
$$
     \dot x + X_\eta(x)
     = E_1(x_0,\xi_0)(\dot x_0+X_\eta(x_0))
       + E_2(x_0,\xi_0)(\nabla\xi_0+\Nabla{\xi_0}X_\eta(x_0)),
$$
$$
     \nabla\xi_0+\Nabla{\xi_0}X_\eta(x_0)   
     = JL_{x_0}(\dot\eta_0+[\eta,\eta_0])
       + \Nabla{\dot x_0+X_\eta(x_0)}(JX_{\eta_0})(x_0).  
$$
Since the image of $JL_{x_0}$ is the orthogonal complement
of the kernel of $d\mu(x_0)$ we deduce that
there exists a constant $c_2>0$ such that 
$$
     \left|\dot x_0+X_\eta(x_0)\right|
     + \left|\nabla\xi_0+\Nabla{\xi_0}X_\eta(x_0)\right|
     \le c_2\left|\dot x+X_\eta(x)\right|,
$$
pointwise for every $\theta\in S^1$.  Define 
$g:\R\to\G$ and $y_0:\R\to\mu^{-1}(0)$ by 
$$
     \dot g+\eta g=0,\qquad g(0)=\one,\qquad 
     y_0(\theta):=g(\theta)^{-1}x_0(\theta).
$$
Then
$$
     g(\theta+2\pi)=g(\theta)g(2\pi),\qquad
     \dot y_0 = g^{-1}(\dot x_0+X_\eta(x_0)).
$$
Hence 
$
     d(y_0(2\pi),x_0(0)) \le 2\pi c_2\ell(x,\eta)
$
and so
$
     d(g(2\pi),\one) \le c_3\ell(x,\eta).
$
This implies that there exists a path $h:\R\to\G$ such that 
$$
     h(\theta+2\pi)=h(\theta)g(2\pi),\qquad
     \sup_{[0,2\pi]}\left|h^{-1}\dot h\right|
     \le c_4\ell(x,\eta).
$$
Hence $g_0:=gh^{-1}$ is a loop and 
$$
     \sup_{S^1}\left|\eta+\dot g_0{g_0}^{-1}\right| 
     = \sup_{[0,2\pi]}\left|h^{-1}\dot h\right|
     \le c_4\ell(x,\eta).
$$
Moreover, with $x_0:=x_0(0)$, we obtain
\begin{eqnarray*}
     d(x(\theta),g_0(\theta)x_0)
&\le &
     d(x(\theta),x_0(\theta))
     + d(x_0(\theta),g(\theta)h(\theta)^{-1}x_0)  \\
&\le &
     d(x(\theta),x_0(\theta))
     + d(y_0(\theta),x_0)
     + d(h(\theta)x_0,x_0) \\
&\le &
     c_1|\mu(x(\theta))| + c_5\ell(x,\eta).
\end{eqnarray*}
This proves the lemma. 
\end{proof}

Let $\delta$ be as in Lemma~\ref{le:action1} 
and $(x,\eta):S^1\to M\times\g$ be a loop such that 
$\sup_{S^1}|\mu(x)|<\delta$.  Then the 
{\bf local equivariant symplectic action} of the pair 
$(x,\eta)$ is defined by 
$$
     \Aa(x,\eta) := - \int u^*\om 
     + \int_0^{2\pi}\inner{\mu(x(\theta))}{\eta(\theta)}\,d\theta,
$$
where $x_0\in\mu^{-1}(0)$ and $g_0:S^1\to\G$ are as in 
Lemma~\ref{le:action1}, $\xi_0(\theta)\in T_{g_0(\theta)x_0}M$ 
is the unique small tangent vector such that
$$
     x(\theta) = \exp_{g_0(\theta)x_0}(\xi_0(\theta)),
$$
and $u:[0,1]\times S^1\to M$ is defined by
$$
     u(\tau,\theta) := \exp_{g_0(\theta)x_0}(\tau\xi_0(\theta)).
$$
The local action is independent of the choice of $x_0$
and $g_0$ so long as the distance between $x(\theta)$ and 
$g_0(\theta)x_0$ remains sufficiently small. 

\begin{lemma}\label{le:action2}
Assume~$(H1)$.  There exist positive constants 
$\delta$ and $c$ such that the following holds. 
If $(x,\eta):S^1\to M$ is a smooth loop such that 
$
      \sup_{S^1}\left|\mu(x)\right| < \delta
$
then 
$$
      \left|\Aa(x,\eta)\right|
      \le c\int_0^{2\pi}\left(\left|\dot x+X_\eta(x)\right|^2
      + \left|\mu(x)\right|^2\right)\,d\theta.
$$
\end{lemma}

\begin{proof}  Let $\xi_0\in\Cinf(S^1,x_0^*TM)$
and $u:[0,1]\times S^1\to M$ be as above. Then the 
local equivariant symplectic action can be expressed in the form
$$
      \Aa(x,\eta) = \int_0^1\int_0^{2\pi}
      \om(\p_\tau u,\p_\theta u+X_\eta(u))\,d\theta d\tau.
$$
By Lemma~\ref{le:action1}, we have the pointwise inequality
$$
      |\p_\tau u|
      = |\xi_0|
      = d(x,g_0x_0)
      \le c_1\left(|\mu(x)|+\ell(x,\eta)\right).
$$
Moreover, by Lemma~\ref{le:group},
$$
      \p_\theta u+X_\eta(u)
      = E_1L_{g_0x_0}(\eta+\dot g_0{g_0}^{-1})
        + \tau E_2
          \left(\nabla{\xi_0}+\Nabla{\xi_0}X_\eta(g_0x_0)\right),
$$
where $E_i:=E_i(g_0x_0,\tau\xi_0)$ for $i=1,2$. 
With $\tau=1$ we obtain, by Lemma~\ref{le:action1},
$$
      \left|\nabla{\xi_0}+\Nabla{\xi_0}X_\eta(g_0x_0)\right|
      \le c_2\left(\left|\dot x+X_\eta(x)\right|
      + \ell(x,\eta)\right)
$$
and this implies
$$
      \left|\p_\theta u+X_\eta(u)\right|
      \le c_3\left(\left|\dot x+X_\eta(x)\right|
      + \ell(x,\eta)\right).
$$
Hence 
\begin{eqnarray*}
      \left|\Aa(x,\eta)\right|
&\le &
      c_1c_3\int_0^{2\pi}
      \bigl(|\mu(x)|+\ell(x,\eta)\bigr)
      \bigl(\left|\dot x+X_\eta(x)\right| + \ell(x,\eta)\bigr)
      \,d\theta \\
&\le &
      c_4\int_0^{2\pi}\left(\left|\dot x+X_\eta(x)\right|^2
      + \left|\mu(x)\right|^2\right)\,d\theta.
\end{eqnarray*}
This proves the lemma.
\end{proof}

\begin{proof}[Proof of Proposition~\ref{prop:vortex}]
Let $(f,J)$ be as in~$(H2)$ and let $(u,\Phi,\Psi)$ 
be a finite energy solution of~(\ref{eq:vortex}) in radial gauge.
We prove in seven steps that $(u,\Phi,\Psi)$ has the properties
asserted in the proposition. 

\medskip
\noindent{\bf Step~1.}
{\it 
$
     \lim_{r\to\infty}
     r^2\Bigl(\left|\p_su+X_\Phi(u)\right|^2
       + \left|\mu(u)\right|^2
      \Bigr)
     = 0
$
uniformly in~$\theta$.}

\medskip
\noindent
Abbreviate $v_s:=\p_su+L_u\Phi$ and $v_t:=\p_tu+L_u\Psi$
as in Section~\ref{sec:apriori}.  Let 
$$
     e := \frac12(|v_s|^2+|\mu(u)|^2).
$$
Then the formula~(\ref{eq:Delta-u0}) 
with $\lambda=\eps=1$ has the form 
\begin{eqnarray*}
     \Delta e
&=&
     |\Nabla{A,s}v_s|^2 + |\Nabla{A,t}v_s|^2
     + 2|L_u^*v_s|^2 + 2|L_u^*Jv_s|^2 + |L_u\mu(u)|^2 \\
&&
     -\,5\inner{\mu(u)}{\rho(v_s,v_t)}
     + \inner{v_s}{(\Nabla{v_s}J)L_u\mu(u)}
     - \inner{v_s}{R(v_s,v_t)v_t}.
\end{eqnarray*}
Since $u(\C)$ is contained in a compact subset 
of $M$ this gives an inequality
$
     \Delta e\ge -c_1e^2.
$
Namely, choose $\delta>0$ such that $L_x$ is injective
whenever $|\mu(x)|^2<\delta$.  Then the first term in 
the secomd row can be estimated from below by 
$
     - |L_u\mu(u)|^2/2 - c|v_s|^4
$
whenever $e\le\delta$. In case $e\ge\delta$ we can use
the inequalities $\Delta e\ge -c(e+e^2)$ and $e\le e^2/\delta$. 
Now it follows from~\cite[Lemma~4.3.2]{MS1} that
there is a constant $c_2>0$ such that
$$
     |z|-r\ge c_2\qquad\IMP\qquad
     e(z) \le \frac{8}{\pi r^2}\int_{B_r(z)}e.
$$
With $r:=|z|/2$ this implies 
$
     \lim_{|z|\to\infty}|z|^2e(z) = 0. 
$

\medskip
\noindent{\bf Step~2.}
{\it  For $R>0$ sufficiently large, we have 
$$
     E(u,\Phi,\Psi;\C\setminus B_R) = \Aa(x_R,\eta_R),
$$
where
$
     x_R(\theta) := u(R{e^i\theta})
$
and
$
     \eta_R(\theta) := R\cos\theta\,\Psi(Re^{i\theta})
     - R\sin\theta\,\Phi(Re^{i\theta}).
$}

\medskip
\noindent
The energy identity on $B_R=\{|z|\le R\}$ has the form
$$
     E(u,\Phi,\Psi;B_R) = \int_{B_R}u^*\om 
     - \int_0^{2\pi}\inner{\mu(x_R(\theta))}{\eta_R(\theta)}\,d\theta.
$$
For $R$ sufficiently large denote by $u_R:[0,1]\times S^1\to M$
the function used in the definition of the local symplectic 
action of $(x_R,\eta_R)$.  Then $u_R(1,\theta)=x_R(\theta)$
and $u_R(0,\theta)=g_R(\theta)x_{R0}$ for some point
$x_{R0}\in\mu^{-1}(0)$ and some loop $g_R:S^1\to\G$. 
The homotopy class of the connected sum 
$v_R:=u|_{B_R}\#(-u_R)$ (the orientation
of $u_R$ is reversed) is independent of $R$.  
Hence the number 
$$
     \int {v_R}^*\om = E(u,\Phi,\Psi;B_R) + \Aa(x_R,\eta_R)
$$
is independent of $R$. Since $\Aa(x_R,\eta_R)$ tends to
zero as $R\to\infty$ it follows that
$$
     E(u,\Phi,\Psi;B_R) + \Aa(x_R,\eta_R) = E(u,\Phi,\Psi)
$$
for every sufficiently large number $R$.  
This proves Step~2.

\medskip
\noindent{\bf Step~3.}
{\it $\sup_\C(f\circ u)\le c$.}

\medskip
\noindent
Suppose, by contradiction, that $\sup_\C(f\circ u)>c$.
Then there exists a regular value $a$ of $f\circ u$ 
such that
$
     c<a<\sup_\C(f\circ u).
$
Hence the set
$
     U := \left\{z\in\C\,|\,f(u(z))\ge a\right\}
$
is a smooth submanifold of $\C$ with boundary.
Since $\mu^{-1}(0)\subset f^{-1}([0,c])$ it follows
from Step~1 that there exists a number $R>0$ such that
$$
     \sup_{\C\setminus B_R}(f\circ u)
     < a < \sup_\C(f\circ u).
$$
Hence $U$ is compact and has a nonempty boundary.
By~$(H2)$, $\Delta(f\circ u)\ge 0$ in $U$ (see~\cite{CGMS}).  
Hence 
$$
     0 \le \int_U\Delta(f\circ u)
       = \int_{\p U}\frac{\p(f\circ u)}{\p\nu} 
       < 0.
$$
This contradiction proves Step~3.

\medskip
\noindent{\bf Step~4.}
{\it Consider equation~(\ref{eq:vortex}) in 
polar coordinates $s+it=e^{\tau+i\theta}$. 
Define $\tilde u:\R\times S^1\to M$
and $\tilde\Phi,\tilde\Psi:\R\times S^1\to\g$ by 
\begin{eqnarray*}
     \tilde u(\tau,\theta)
&:= &
     u(e^{\tau+i\theta}),    \\
     \tilde\Phi(\tau,\theta)
&:= &
     e^\tau\cos\theta\,\Phi(e^{\tau+i\theta})
     + e^\tau\sin\theta\,\Psi(e^{\tau+i\theta}),  \\
     \tilde\Psi(\tau,\theta)
&:= &
     e^\tau\cos\theta\,\Psi(e^{\tau+i\theta})
     - e^\tau\sin\theta\,\Phi(e^{\tau+i\theta}).
\end{eqnarray*}
Then $ds\wedge dt=e^{2\tau}d\tau\wedge d\theta$, 
$\Phi\,ds+\Psi\,dt=\tilde\Phi\,d\tau +\tilde\Psi\,d\theta$,
and~(\ref{eq:vortex}) is equivalent to
$$
     \tilde v_\tau + J\tilde v_\theta = 0,\qquad
     \tilde\kappa + e^{2\tau}\mu(\tilde u) = 0,
$$
where 
$$
     \tilde v_\tau:=\p_\tau\tilde u+L_{\tilde u}\tilde\Phi,\qquad
     \tilde v_\theta:=\p_\theta\tilde u+L_{\tilde u}\tilde\Psi,\qquad
     \tilde\kappa:=
     \p_\tau\tilde\Psi - \p_\theta\tilde\Phi + [\tilde\Phi,\tilde\Psi].
$$
The radial gauge condition has the form $\tilde\Phi(\tau,\theta)=0$
for large $\tau$. The energy of the triple 
$(\tilde u,\tilde\Phi,\tilde\Psi)$ is given by}
$$
     E(\tilde u,\tilde\Phi,\tilde\Psi)
     = \int_{-\infty}^\infty\int_0^{2\pi}
       \left(\left|\tilde v_\tau\right|^2
        + e^{2\tau}\left|\mu(\tilde u)\right|^2\right)
       \,d\theta d\tau.
$$

\medskip
\noindent{\bf Step~5.}
{\it There exist positive constants $c$ and $\delta$ such that,
for every $\tau_0\ge 0$,}
$$
     \eps(\tau_0) := \int_{\tau_0}^\infty\int_0^{2\pi}
     \left(\left|\tilde v_\tau\right|^2
     + e^{2\tau}\left|\mu(\tilde u)\right|^2\right)
     \,d\theta d\tau
     \le ce^{-\delta\tau_0}.
$$

\medskip
\noindent
By Step~2 and Lemma~\ref{le:action2}, we have 
\begin{eqnarray*}
     \eps(\tau) 
&= &
     \Aa(x_{e^\tau},\eta_{e^\tau})  \\
&\le &
     \delta^{-1}\int_0^{2\pi}
     \left(\left|\tilde v_\tau(\tau,\theta)\right|^2
     + e^{2\tau}\left|\mu(\tilde u(\tau,\theta))\right|^2\right)
     \,d\theta  \\
&= &
     -\delta^{-1}\eps'(\tau)
\end{eqnarray*}
for some constant $\delta>0$ and every sufficiently large
real number $\tau$.  Hence there exists a real number $\tau_0$
such that 
$$
     \tau\ge\tau_0\qquad\IMP\qquad
     \eps(\tau)\le e^{-\delta(\tau-\tau_0)}\eps(\tau_0).
$$

\medskip
\noindent{\bf Step~6.}
{\it There exist positive constants $c$ and $\delta$ such that,
for every $\tau_0\ge 0$,}
$$
     \sup_{\theta\in\R}\left(
     \left|\tilde v_\tau(\tau_0,\theta)\right|^2 
     + e^{2\tau_0}\left|\mu(\tilde u(\tau_0,\theta))\right|
     \right)
     + \int_{\tau_0}^\infty\int_0^{2\pi}
       e^{4\tau}\left|\mu(\tilde u(\tau,\theta)\right|^2\,
       d\theta d\tau
     \le ce^{-\delta\tau_0}.
$$

\medskip
\noindent
The $L^\infty$ estimate for $\tilde v_\tau$ and 
$e^{2\tau}\mu(\tilde u)$ follows from Step~5 
and~(\ref{eq:apriori-inf}).  The $L^2$-estimate for 
$\mu(\tilde u)$ follows from Step~5 and Lemma~\ref{le:apriori}.

\medskip
\noindent{\bf Step~7.}
{\it There exists a $W^{1,2}$-function $x:S^1\to\mu^{-1}(0)$
and an $L^2$-function $\eta:S^1\to\g$ such that
$$
     \lim_{\tau\to\infty}\sup_{\theta\in\R}
     d(\tilde u(\tau,\theta),x(\theta))=0,\qquad
     \lim_{\tau\to\infty}\int_0^{2\pi}
     \left|\tilde\Psi(\tau,\theta)-\eta(\theta)\right|^2\,d\theta=0,
$$
and $\dot x+X_\eta(x)=0$. Moreover, 
$E(u,\Phi,\Psi)=\int u^*\om$, and if~$(H3)$ holds
then $E(u,\Phi,\Psi)\in\Z\hbar$.}

\medskip
\noindent
By H\"older's inequality and the radial gauge assumption, 
we have, for $\tau_1>\tau_0$,
\begin{eqnarray*}
     \left|\tilde\Psi(\tau_1,\theta)
     -\tilde\Psi(\tau_0,\theta)\right|^2
&\le &
     \left(\int_{\tau_0}^{\tau_1}e^{2\tau}
     \left|\mu(\tilde u(\tau,\theta))\right|\,d\tau
     \right)^2  \\
&\le &
     (\tau_1-\tau_0)\int_{\tau_0}^{\tau_1}
     e^{4\tau}\left|\mu(\tilde u(\tau,\theta))\right|^2\,d\tau.
\end{eqnarray*}
Hence the existence of the $L^2$-limit 
of $\tilde\Psi$ follows from Step~6. 
That $\tilde u(\tau,\theta)$ converges uniformly as $\tau$ tends
to infinity follows from the exponential decay of 
$\tilde v_\tau=\p_\tau\tilde u$ in Step~6.  That the limit 
is a $W^{1,2}$-function and satisfies $\dot x+X_\eta(x)=0$
follows from the fact that 
$
    \tilde v_\theta=\p_\theta u+X_{\tilde\Psi}(\tilde u)
$
converges (exponentially) to zero as $\tau$ tends to infinity. 
That $E(u,\Phi,\Psi)=\int u^*\om$ follows from the energy identity
in  the proof of Step~2 and the $L^2$-convergence of $\tilde\Psi$.
That $E(u,\Phi,\Psi)$ is an integer multiple of $\hbar$ (when~$(H3)$ holds) 
follows from the proof of Step~2. 
\end{proof}

Every map $u:\C\to M$ that satisfies~(\ref{eq:lim})
and~(\ref{eq:xeta}) determines an equivariant homology 
class $B=[u]\in H_2(M_\G;\Z)$ as follows.  
Homotop $u$ to a map $v:D\to M$ such that
$
     v(e^{i\theta})=g(e^{i\theta})x_0.
$
Now define a principal bundle $P\to S^2\cong(\Z_2\times D)/\sim$
by 
$$
    P := (\Z_2\times D\times\G)/\sim,\qquad
    (0,e^{i\theta},h)\sim(1,e^{i\theta},g(e^{i\theta})h).
$$
Then $v$ determines a $\G$-equivariant map $w:P\to M$ by 
$$
    w(0,z,h) := h^{-1}x_0,\qquad
    w(1,z,h) := h^{-1}v(z).
$$
The equivariant homology class of $u$ is defined to
be the equivariant homology class of $w$.


\section{Proof of Theorem~A}\label{sec:bubbling}

We begin by constructing a $\Gg$-equivariant smooth function
from a suitable open subset of 
$\Bb\times P=\Cinf_\G(P,M)\times\Aa(P)\times P$ 
into a suitable finite dimensional approximation of $\EG$.   
For positive constants $\delta$ and $r$ denote
$$
     M^\delta := \left\{x\in M\,|\,|\mu(x)|\le\delta\right\},
$$
$$
     \Bb^{\delta,r}
     := \left\{(u,A)\in\Cinf_\G(P,M)\times\Aa(P)\,\Big|\,
        \exists z\in\Sigma\;\;
        \sup_{B_r(z)}|\mu\circ u|<\delta
        \right\}.
$$
Throughout we assume that $\G$ is a Lie subgroup of $\U(k)$.
Then, for every integer $n\ge k$, a finite dimensional approximation 
of the classifying space of $\G$ is the quotient 
$$
     \BG^n := \EG^n/\G,\qquad
     \EG^n := \Ff(k,n) 
     := \left\{\Theta\in\C^{n\times k}\,|\,\Theta^*\Theta=\one\right\}.
$$
Let $\delta>0$ be so small that $\G$ acts freely on $M^\delta$ 
and choose a smooth $\G$-equivariant classifying map 
$\theta:M^\delta\to\EG^m$ for some integer $m\ge k$. 

\begin{proposition}\label{prop:classify}
Let $\theta:M^\delta\to\EG^m$ be as above.
Then, for every $r>0$, there exist an integer $n\ge m$ 
and a map $\Theta:\Bb^{\delta,r}\times P\to\EG^n$
with the following properties.
\begin{description}
\item[(i)]
For $g\in\Gg$, $h\in\G$, and $(u,A,p)\in\Bb^{\delta,r}\times P$,
\begin{equation}\label{eq:theta0}
     \Theta(g^{-1}u,g^*A,pg(p)^{-1}) 
     = \Theta(u,A,p)
     = h\Theta(u,A,ph).
\end{equation}
\item[(ii)]
$\Theta$ is smooth with respect to the $C^0$ Banach 
manifold structure on (the completion of) $\Bb^{\delta,r}\times P$.  
\item[(iii)]
Let $\i:\EG^m\to\EG^n$ be the obvious inclusion. Then
$$
     |\mu(u(p))|<\delta\qquad\IMP\qquad
     \Theta(u,A,p) = \i\circ\theta(u(p)).
$$
\end{description}
\end{proposition}

\begin{proof}
Cover $\Sigma$ by finitely many distinct balls $B_r(z_i)$, 
$i=1,\dots,\ell$, and choose points $p_1,\dots,p_\ell\in P$ 
such that $\pi(p_i)=z_i$.  Then, for every $(u,A)\in\Bb^{\delta,r}$, 
there exists an $i\in\{1,\dots,\ell\}$ such that $u(p_i)\in M^\delta$.
Thus the open set $\Bb^{\delta,r}\times P$ is contained in the 
finite union of the following open sets $\Uu_{ij}$, 
$i,j=1,\dots,\ell$. Choose $\eps>0$ so small that $\theta$ 
extends to an equivariant function (still denoted by $\theta$) 
from $M^{\delta+\eps}$ to $\EG^m$ and define
\begin{eqnarray*}
     \Uu_0 
&:= &
     \left\{(u,A,p)
      \,|\,|\mu(u(p))|<\delta+\eps\right\},  \\
     \Uu_{ij}
&:= &
     \left\{(u,A,p)
     \,|\,\pi(p)\in B_r(z_i),\,|\mu(u(p_j))|<\delta+\eps\right\}.
\end{eqnarray*}
For every smooth path $\gamma:[0,1]\to\Sigma$ 
and any two points $p_0\in\pi^{-1}(\gamma(0))$
and $p_1\in\pi^{-1}(\gamma(1))$
the holonomy $\rho_A(p_1,\gamma,p_0)\in\G$ of the connection 
$A\in\Aa(P)$ is defined by 
$
     p_1\rho_A(p_1,\gamma,p_0) := \tilde\gamma(1),
$
where $\tilde\gamma:[0,1]\to P$ is the unique horizontal
lift of $\gamma$ with $\tilde\gamma(0)=p_0$. 
It satisfies
\begin{eqnarray*}
     \rho_{A}(p_1g_1,\gamma,p_0g_0) 
     &= & {g_1}^{-1}\rho_A(p_1,\gamma,p_0)g_0,  \\
     \rho_{g^*A}(p_1,\gamma,p_0) 
     &= & g(p_1)^{-1}\rho_A(p_1,\gamma,p_0)g(p_0)
\end{eqnarray*}
for $g_0,g_1\in\G$ and $g\in\Gg$. Hence the map 
$
     \Bb\times\pi^{-1}(\gamma(0))\to M:
     (u,A,p_0)\mapsto\rho_A(p_1,\gamma,p_0)^{-1}u(p_1)
$
is $\Gg$-invariant and $\G$-equivariant.  
Choose a finite sequence of smooth functions 
$\gamma_{ij}:[0,1]\times B_r(z_i)\to\Sigma$ such that
$$
     \gamma_{ij}(0,z)=z,\qquad \gamma_{ij}(1,z)=z_j.
$$
Then the functions $\Theta_0:\Uu_0\to\EG^m$ and
$\Theta_{ij}:\Uu_{ij}\to\EG^m$, defined by 
$$
     \Theta_0(u,A,p) := \theta(u(p)),\qquad
     \Theta_{ij}(u,A,p) 
     := \theta\Bigl(
        \rho_A\bigl(p_j,\gamma_{ij}(\cdot,\pi(p)),p\bigr)^{-1}
        u(p_j)\Bigr)
$$
for $i,j=1,\dots,\ell$, satisfy~(\ref{eq:theta0}).  
Now choose a $\G$-invariant smooth function $\beta:M\to[0,1]$
such that $\beta(x)=1$ for $x\in M^\delta$ and 
$\beta(x)=0$ for $x\in M\setminus M^{\delta+\eps/2}$.
Define $\rho_j:\bigcup_{i=1}^\ell\Uu_{ij}\to[0,1]$ 
and $\rho_0:\Uu_0\to[0,1]$ by 
\begin{eqnarray*}
     \rho_0(u,A,p) &:= & \frac{\beta(u(p))}
     {\sqrt{\beta(u(p))^2+(1-\beta(u(p)))^2}}, \\
     \rho_j(u,A,p) &:= & \frac{\beta(u(p_j))(1-\beta(u(p)))}
     {\sqrt{(\beta(u(p))^2+(1-\beta(u(p)))^2)
      \sum_{k=1}^\ell\beta(u(p_k))^2}}
\end{eqnarray*}
for $j=1,\dots,\ell$.  Then $\rho_j$ is smooth with respect 
to the $C^0$ Banach manifold structure on 
(the completion of) $\Bb\times P$.  Moreover, 
$$
     \sum_{j=0}^\ell \rho_j(u,A,p)^2 = 1
$$
and $\rho_0(u,A,p)=1$ whenever $|\mu(u(p))|\le\delta$.
Now choose a finite sequence of smooth functions 
$\sigma_i:\Sigma\to[0,1]$, $i=1,\dots,\ell$, 
such that $\supp\,\sigma_i\subset B_r(z_i)$
and 
$
     \sum_{j=1}^\ell{\sigma_i}^2 = 1.
$
Then the function $\Theta:\Bb^{\delta,r}\to\EG^{(\ell^2+1)m}$, 
defined by
$$
     \Theta(u,A) := \left(\begin{array}{c}
     \rho_0(u,A,p)\Theta_0(u,A,p) \\ 
     \vdots \\
     \sigma_i(\pi(p))\rho_j(u,A,p)\Theta_{ij}(u,A,p)\\
     \vdots \\
     \end{array}\right),
$$
is the required classifying map.
\end{proof}

The integer $n=(\ell+1)m$ in the proof of 
Proposition~\ref{prop:classify} diverges to infinity 
as $r$ tends to zero.  In general, there is no $\Gg$-equivariant
map from $\Bb^\delta:=\{(u,A)\in\Bb\,|\,\min|\mu\circ u|<\delta\}$ 
to any finite dimensional approximation of~$\EG$. 

\begin{lemma}\label{le:Meps}
Assume~$(H1-3)$.  Let $B\in H_2(M_\G;\Z)$, 
$(\Sigma,j_\Sigma,\dvol_\Sigma)$ be a compact Riemann surface,
$\pi:P\to\Sigma$ be the principal $\G$-bundle determined by $B$,
and $\Sigma\to\Jj_\G(M,\om)$ be a smooth family 
of $\G$ invariant $\om$-compatible almost complex structures 
on $M$ such that each $J_z$ agrees with the 
almost complex structure of~$(H2)$ outside of a sufficiently
large compact subset of~$M$. Then for every $\delta>0$ there 
exist positive constants $r$ and $\eps_0$ such that 
$$
      \Tilde{\Mm}_{B,\Sigma}^\eps\subset\Bb^{\delta,r}
$$
for $0<\eps\le\eps_0$.
\end{lemma}

\begin{proof}
Suppose the result is false.  Then there exist
a constant $\delta>0$ and sequences 
$$
     r_i\to 0,\qquad \eps_i\to 0,\qquad 
     (u_i,A_i)\in\Tilde{\Mm}_{B,\Sigma}^{\eps_i}
$$
such that $(u_i,A_i)\notin\Bb^{\delta,r_i}$ 
for every~$i$. This means that, for every $p\in P$,
there exists a sequence $p_i\in P$ such that 
$$
      \lim_{i\to\infty}p_i=p,\qquad
      |\mu(u_i(p_i))|\ge\delta.
$$
This contradicts the bubbling argument in Step~5 of 
the proof of Theorem~A below. 
\end{proof}

Let $\eps_0>0$ be as in Lemma~\ref{le:Meps}.
For $0\le\eps\le\eps_0$ we consider the evaluation map 
$$
     \ev_\G^\eps:\Mm^\eps_{B,\Sigma}\times\Sigma\to M\times_\G\EG^n,
$$
given by
$$
     \ev_\G^\eps([u,A,p]) := [u(p),\Theta(u,A,p)]
$$
where $\delta>0$ is chosen such that $\G_x=\{\one\}$ for every
$x\in M^\delta$, $r>0$ is as in Lemma~\ref{le:Meps},
and $\Theta:\Bb^{\delta,r}\times P\to\EG^n$ is the map of 
Proposition~\ref{prop:classify}.
Recall that $\Mcal_{B,\Sigma}^0$ and $\Mcal_{B,\Sigma}^\eps$
have the same dimension.  

\begin{proposition}\label{prop:ev}
For every $c_0>0$ there exist positive constants $c$ and $\eps_0$
such that the following holds. 

\smallskip
\noindent{\bf (i)} 
For $0<\eps\le \eps_0$ the map 
$\Tcal^\eps:\Mcal^0_{B,\Sigma}(c_0)\to \Mcal_{B,\Sigma}^\eps$
is an orientation preserving embedding.

\smallskip
\noindent{\bf (ii)}
For $0<\eps\le \eps_0$, 
$$
     d_{C^1}
     (\ev_\G^0,\ev_\G^\eps\circ(\Tcal^\eps\times\id))
     \le c\eps^{1-2/p},
$$
where the $C^1$-distance is understood on the space of continuously 
differentiable maps from $\Mm^0_{B,\Sigma}(c_0)\times\Sigma$
to $M\times_\G\EG^n$. 
\end{proposition}

\begin{lemma}\label{le:dT}
Assume~$(H1)$ and~$(H4)$ and let $\bar B\in H_2(\bar M;\Z)$
be a nontorsion homology class. 
For every $p>2$ and every $c_0>0$ there exist positive constants 
$\eps_0$ and $c$ such that the following holds for every 
$\eps\in(0,\eps_0]$.  Let $I\subset\R$ be an interval and 
$$
     I\to\Tilde{\Mm}^0_{B,\Sigma}(c_0):r\mapsto(u_0(r),A_0(r))
$$
be a smooth path that satisfies~(\ref{eq:path}).
Then every smooth vector field 
$$
     r\mapsto\zeta(r)\in\im\,\left(\Dd^\eps_{(u_0(r),A_0(r))}\right)^*
$$
satisfies the inequality
$$
      \left\|\Tabla{r}\zeta\right\|_{1,p,\eps}
      \le c\left(
      \eps\left\|\Tabla{r}\Dd^\eps\zeta\right\|_{0,p,\eps}
      + \left\|\Tabla{r}\pi_{u_0}\Dd^\eps\zeta\right\|_{L^p} 
      + \eps^{-1}\left\|\zeta\right\|_{1,p,\eps}
      \right)
$$
for $r\in I$, where $\Dd^\eps:=\Dd^\eps_{(u_0(r),A_0(r))}$
and $\pi_{u_0}$ is defined by~(\ref{eq:pi-u}).
\end{lemma}

\begin{proof}
Let $r\mapsto\zeta'(r)\in\Xx'_{u_0(r)}$ be the smooth path
defined by
$
     \zeta={\Dd^\eps}^*\zeta'.  
$
Then, by Lemmata~\ref{le:onto}, \ref{le:onto1}, and~\ref{le:commute}, 
we have
\begin{eqnarray*}
      \left\|\Tabla{r}\zeta\right\|_{1,p,\eps}
&\le &
      \left\|{\Dd^\eps}^*\Tabla{r}\zeta'\right\|_{1,p,\eps}
      + \left\|\Tabla{r}{\Dd^\eps}^*\zeta'
        - {\Dd^\eps}^*\Tabla{r}\zeta'\right\|_{1,p,\eps}  \\
&\le &
      c_2\left(
      \eps\left\|\Dd^\eps{\Dd^\eps}^*\Tabla{r}\zeta'\right\|_{0,p,\eps}
      + \left\|\pi_{u_0}\Dd^\eps{\Dd^\eps}^*\Tabla{r}\zeta'\right\|_{L^p}
      \right) \\
&&
      +\, \left\|\Tabla{r}{\Dd^\eps}^*\zeta'
        - {\Dd^\eps}^*\Tabla{r}\zeta'\right\|_{1,p,\eps}  \\
&\le &
      c_2\left(
      \eps\left\|\Dd^\eps\Tabla{r}\zeta\right\|_{0,p,\eps}
      + \left\|\pi_{u_0}\Dd^\eps{\Dd^\eps}^*\Tabla{r}\zeta'\right\|_{L^p}
      \right) \\
&&
      +\, c_3\left\|\Tabla{r}{\Dd^\eps}^*\zeta'
          - {\Dd^\eps}^*\Tabla{r}\zeta'\right\|_{1,p,\eps} \\
&\le &
      c_2\left(
      \eps\left\|\Tabla{r}\Dd^\eps\zeta\right\|_{0,p,\eps}
      + \left\|\pi_{u_0}\Dd^\eps{\Dd^\eps}^*\Tabla{r}\zeta'\right\|_{L^p}
      \right) \\
&&
      +\, c_2\eps\left\|\Dd^\eps\Tabla{r}\zeta
          - \Tabla{r}\Dd^\eps\zeta\right\|_{0,p,\eps}
      + c_4\eps^{-1}\left\|\zeta'\right\|_{2,p,\eps}  \\
&\le &
      c_2\left(
      \eps\left\|\Tabla{r}\Dd^\eps\zeta\right\|_{0,p,\eps}
      + \left\|\Tabla{r}\pi_{u_0}\Dd^\eps\zeta\right\|_{L^p}
      \right) \\
&&
      +\, c_2\left\|[\pi_{u_0}\Dd^\eps{\Dd^\eps}^*,\Tabla{r}]\zeta'
          \right\|_{L^p}
      + c_5\eps^{-1}\left\|\zeta\right\|_{1,p,\eps}  \\
&\le &
      c_2\left(
      \eps\left\|\Tabla{r}\Dd^\eps\zeta\right\|_{0,p,\eps}
      + \left\|\Tabla{r}\pi_{u_0}\Dd^\eps\zeta\right\|_{L^p}
      \right) \\
&&
      +\, c_6\left(
      \left\|\Dd^\eps\zeta\right\|_{0,p,\eps}
      + \eps^{-1}\left\|\pi_{u_0}\Dd^\eps\zeta\right\|_{L^p}
      + \eps^{-1}\left\|\zeta\right\|_{1,p,\eps}
      \right).
\end{eqnarray*}
The last inequality follows as in Step~7 in the proof
of Lemma~\ref{le:onto1}. Since
$$
     \left\|\Dd^\eps\zeta\right\|_{0,p,\eps}
     + \eps^{-1}\left\|\pi_{u_0}\Dd^\eps\zeta\right\|_{L^p}
     \le c_7\eps^{-1}\left\|\zeta\right\|_{1,p,\eps},
$$
the lemma is proved.
\end{proof}

\begin{proof}[Proof of Proposition~\ref{prop:ev}]
Let $r\mapsto(u_0(r),A_0(r))$ be as in Lemma~\ref{le:dT}
and 
$
     r\mapsto\zeta_\eps(r)
     =(\xi_\eps(r),\alpha_\eps(r))
     \in\im\,(\Dd^\eps_{(u_0(r),A_0(r))})^*
$
be as in Theorem~\ref{thm:exist} so that
$$
     \Tilde{\Tt}^\eps(u_0(r),A_0(r))
     = (u_\eps(r),A_\eps(r))
     := (\exp_{u_0(r)}(\xi_\eps(r)),A_0(r)+\alpha_\eps(r)).
$$
Let $\Ff^\eps_r:\Xx_{u_0(r)}\to\Xx'_{u_0(r)}$ be defined 
by~(\ref{eq:Feps}).  Then $\Ff^\eps_r(\zeta_\eps(r))=0$ and hence 
$$
     \Dd^\eps\zeta_\eps
     = \left(\begin{array}{c} 0 \\ 0 \\ *F_{A_0(r)} \end{array}\right)
       + \Bigl(\Ff^\eps_r(\zeta_\eps(r)) - \Ff^\eps_r(0)
         - d\Ff^\eps_r(0)\zeta_\eps(r)\Bigr).
$$
Differentiating this identity with respect to $r$ we find
$$
     \left\|\Tabla{r}\Dd^\eps\zeta_\eps\right\|_{0,p,\eps}
     \le
     c\left(\eps
     +\eps^{-1-2/p}\left\|\zeta_\eps(r)\right\|_{1,p,\eps}
      \left(\left\|\zeta_\eps(r)\right\|_{1,p,\eps}
      + \left\|\Tabla{r}\zeta_\eps(r)\right\|_{1,p,\eps}\right)
     \right).
$$
(See Proposition~\ref{prop:quadra1}.)
Hence, by Lemma~\ref{le:dT},
\begin{equation}\label{eq:dT}
     \left\|\Tabla{r}\zeta_\eps\right\|_{1,p,\eps}
     \le c'\eps.
\end{equation}
Since 
$$
     (\p_ru_\eps,\p_rA_\eps)
     = \left(E_1(u_0,\xi_\eps)\p_ru_0+E_2(u_0,\xi_\eps)\Nabla{r}\xi_\eps,
       \p_rA_0+\p_r\alpha_\eps\right),
$$
this shows that $\Tt^\eps:\Mm^0_{B,\Sigma}(c_0)\to\Mm^\eps_{B,\Sigma}$
is an orientation preserving embedding. Indeed, it follows 
that the restriction of $\Tt^\eps$ to every ball of radius $\delta$
is an embedding for $\delta$ and $\eps$ sufficiently small and hence, 
by Theorem~\ref{thm:exist}, $\Tt^\eps$ is an embedding for 
$\eps$ sufficiently small.  For $\eps\ge 0$ denote
$$
      \Pp^\eps_{B,\Sigma}(c_0):=\Tilde{\Mm}^\eps_{B,\Sigma}(c_0)\times_{\Gg}P
$$
and consider the map $\ev^\eps:\Pp^\eps_{B,\Sigma}(c_0)\to M$ defined by 
$$
      \ev^\eps([u_\eps,A_\eps,p]) := u_\eps(p).
$$
Then it follows from~(\ref{eq:dT}) and the inequality
$
     \left\|\zeta_\eps\right\|_{2,p,\eps}
     \le c'\eps^2 
$
of Theorem~\ref{thm:exist} that
$$
      d_{C^1}(\ev^0,\ev^\eps\circ(\Tt^\eps\times\id))
      \le c''\eps^{1-2/p}.
$$
For $\eps$ sufficiently small we have $\ev^\eps_\G=\theta^\delta\circ\ev^\eps$, 
where $\theta^\delta:M^\delta\to M\times_\G\EG^n$ is given by  
$\theta^\delta(x):=[x,\theta(x)]$.  This proves the proposition.
\end{proof}

\begin{proof}[Proof of Theorem~A]
The result is obvious when $\bar B=0$.
Moreover, both moduli spaces are empty when $\bar B$ is a nonzero
torsion class.  Hence assume that $\bar B\in H_2(\bar M;\Z)$
is a nontorsion homology class, denote by $B\in H_2(M_\G;\Z)$ 
the corresponding equivariant homology class,
fix a compact Riemann surface $(\Sigma,j_\Sigma,\dvol_\Sigma)$, 
and let $\pi:P\to\Sigma$ be a principal $\G$-bundle 
whose characteristic class $b\in H_2(\BG;\Z)$
is the pushforward of $B$.  In the course of the proof 
it will be necessary to also consider other bundles
$\pi':P'\to\Sigma$ with corresponding equivariant 
homology clsses $B'\in H_2(M_\G;\Z)$. 
By~$(H2)$, there exists a constant $c>0$ such that
$
     u(P) \subset M^c := \{x\in M\,|\,|\mu(x)|\le c\}
$ 
for every solution $(u,A)$ of~(\ref{eq:eps}) over any 
Riemann surface.  Note that $c$ can be chosen 
to be a regular value of the function $M\to\R:x\mapsto|\mu(x)|$.

\smallbreak

Let $\delta>0$ such that $\G_x=\{\one\}$
for every $x\in M^\delta$ and let $r$ and $\eps_0$ 
be as in Lemma~\ref{le:Meps}. 
Fix $k$ points $p_1,\dots,p_k\in P$ such that the points
$z_i:=\pi(p_i)\in\Sigma$ are pairwise distinct.
Choose an integer $n$, a $\G$-equivariant smooth map
$\theta:M^\delta\to\EG^n$, and $k$ smooth classifying maps 
$
     \Theta_i:\Bb^{\delta,r}\to\EG^n,
$ 
defined by $\Theta_i(u,A):=\Theta(u,A,p_i)$, where
$\Theta$ is as in Proposition~\ref{prop:classify}.  Then
$$
     \Theta_i(g^{-1}u,g^*A) = g(p_i)^{-1}\Theta_i(u,A)
$$
and
$$
     |\mu(u(p_i))|\le\delta\qquad\IMP\qquad
     \Theta_i(u,A) = \theta(u(p_i))
$$
for $i=1,\dots,k$. For $0\le\eps\le\eps_0$ consider the evaluation maps
$$
     \ev^\eps_{B,i}:\Mm^\eps_{B,\Sigma}\to M^c_\G := M^c\times_\G\EG^n
$$
given by 
$
     \ev^\eps_{B,i}([u,A]) := [u(p_i),\Theta_i(u,A)].
$
Let $\ev^\eps_B:\Mm^\eps_{B,\Sigma}\to(M^c_\G)^k$
denote the product map defined by 
$$
     \ev^\eps_B([u,A]) := \left(\ev^\eps_{B,1}([u,A]),\dots,
     \ev^\eps_{B,k}([u,A])\right).
$$
For any subset $I=\{i_1,\dots,i_j\}\subset\{1,\dots,k\}$
such that $i_1<\cdots<i_j$ and any class $B'\in H_2(M_\G;\Z)$ 
that descends to $H_2(\bar M;\Z)$ we consider the evaluation map
$
     \ev^0_{B',I}:\Mm^0_{B',\Sigma}\to(M^c_\G)^{|I|}
$
given by
$$
     \ev^0_{B',I}([u,A]) 
     := \left([u(p_{i_1}),\theta(u(p_{i_1}))],\dots,
              [u(p_{i_j}),\theta(u(p_{i_j}))]\right).
$$
Now fix equivariant cohomology classes 
$\alpha_1,\dots,\alpha_k\in H^*(M_\G;\Z)$ such that 
$$
     m_i := \deg(\alpha_i)<2N,\qquad
     \sum_{i=1}^k\deg(\alpha_i)=\dim\,\Mm^0_{B,\Sigma}.
$$
There is a natural embedding $M^c_\G\to M_\G$
and we denote by $\alpha_i^c\in H^{m_i}(M^c_\G;\Z)$
the pullback of $\alpha_i$ under this embedding. 
Note that $M^c_\G$ is a compact manifold with boundary.
Replacing $\alpha_i$ by some integer multiple of $\alpha_i$,
if necessary, we may assume without loss of generality that,
for every $i$, there exists a compact oriented manifold with boundary 
$Y_i$ of dimension
$$
     \dim\,Y_i=\dim\,M^c_\G-m_i
$$
and a smooth map
$$
     f_i:(Y_i,\p Y_i)\to (M^c_\G,\p M^c_\G)
$$
such that the homology class in $H_*(M^c_\G,\p M^c_\G)$
represented by $f_i$ is Poincar\'e dual to $\alpha_i^c$. 
For $I=\{i_1,\dots,i_j\}\subset\{1,\dots,k\}$
such that $i_1<\cdots<i_j$ we denote the corresponding 
product map by
$$
     Y_I:=Y_{i_1}\times\cdots\times Y_{i_j},\qquad
     f_I:=f_{i_1}\times\cdots\times f_{i_j}:Y_I\to(M^c_\G)^{|I|}.
$$
For $I=\{1,\dots,k\}$ we abbreviate
$
     Y:=Y_{\{1,\dots,k\}}
$
and
$
     f:=f_{\{1,\dots,k\}}.
$
The functions $f_1,\dots,f_k$ can be chosen such 
that the following holds.
\begin{description}
\item[(H5)]
{\it $f_i$ is transverse to $\mu^{-1}(0)\times_\G\EG^n$
for every $i$ and $f_I$ is transverse to $\ev^0_{B',I}$
for every subset $I\subset\{1,\dots,\ell\}$ and every
equivariant homology class $B'\in H_2(M_\G;\Z)$.}
\end{description}
Now the notation has been set up and we shall prove Theorem~A
in five steps.  For $0\le\eps\le\eps_0$ and $B\in H_2(M_\G;\Z)$
consider the set
$$
      \Mm^\eps_{B,\Sigma;f}
      := \left\{([u,A],y_1,\dots,y_k)\in\Mm^\eps_{B,\Sigma}
         \times Y\,|\,
         \ev^\eps_{B,i}([u,A])=f_i(y_i)
         \right\}.
$$

\medskip
\noindent{\bf Step~1.}
{\it The map $\ev^0_B:\Mm^0_{B,\Sigma}\to(\mu^{-1}(0)\times_\G\EG^n)^k$
is a pseudo-cycle.}

\medskip
\noindent
The map $\ev^0_B$ is the composition 
$\ev^0_B=\i^k\circ\bev_B$, where the evaluation map
$
     \bev_B:\Mm^0_{B,\Sigma}\to\bar M^k
$ 
is given by 
$
     \bev_B([u,A]) := ([u(p_1)],\dots,[u(p_k)])
$
and the embedding
$
     \i:\bar M\to\mu^{-1}(0)\times_\G\EG^n
$
is given by 
$
     \i([x]) := [x,\theta(x)].
$
That $\bev_B$ is a pseudo-cycle was proven in~\cite{MS1}.
Hence $\ev^0_B$ is a pseudo-cycle.
(see~\cite{MS1} for the definitions). 

\medskip
\noindent{\bf Step~2.}
{\it
$\Mm^0_{B,\Sigma;f}$ is a finite set and
the number of elements of $\Mm^0_{B,\Sigma;f}$,
counted with appropriate signs, is the Gromov--Witten invariant:
$$
      \GW_{\bar B,\Sigma}(\bar\alpha_1,\dots,\bar\alpha_k)
= 
      \ev^0_B\cdot f 
= 
      \sum_{([u_0,A_0],y)\,\in\,\Mm^0_{B,\Sigma;f}}
      \nu^0([u_0,A_0],y).
$$
Here the function $\nu^0:\Mm^0_{B,\Sigma;f}\to\{\pm1\}$ 
denotes the intersection index of the maps $\ev^0_B$ 
and $f$.}

\medskip
\noindent
Consider the functions $\phi_i:X_i\to\mu^{-1}(0)\times_\G\EG^n$
given by 
$$
     X_i := f_i^{-1}(\mu^{-1}(0)\times_\G\EG^n)\subset Y_i,\qquad
     \phi_i:=f_i|_{X_i}. 
$$
Since $f_i$ is transverse to $\mu^{-1}(0)\times_\G\EG^n$,
$X_i$ is a smooth submanifold of $Y_i$ and $\phi_i$ is dual
to the cohomology class $\alpha_i^0\in H^*(\mu^{-1}(0)\times_\G\EG^n;\Z)$
obtained from $\alpha_i^c$ by pullback under the obvious 
inclusion $\mu^{-1}(0)\times_\G\EG^n\to M^c_\G$.  The class
$\alpha_i^0$ agrees with the image of the class 
$\bar\alpha_i\in H^*(\bar M;\Z)$ under the homomorphism
$
     H^*(\bar M;\Z)
     \cong H^*(\mu^{-1}(0)\times_\G\EG;\Z)
     \to   H^*(\mu^{-1}(0)\times_\G\EG^n;\Z):
$
$$
\xymatrix
{
\alpha_i^c \;\; \in \hspace{-20pt} & 
{\rm H}^*(M^c\times_\G\EG^n;\Z)\ar[d]  &
{\rm H}^*(M_\G;\Z) \ar[l]\ar[d]^\kappa & 
\hspace{-20pt} \ni \;\;\alpha_i  \\
\alpha^0_i\;\;\in \hspace{-20pt} & 
{\rm H}^*(\mu^{-1}(0)\times_\G\EG^n;\Z) & 
\ar[l] {\rm H}^*(\bar M;\Z) &
\hspace{-20pt} \ni\;\; \bar\alpha_i 
}
$$
Hence another representative of the class $\alpha_i^0$ 
can be obtained as follows.  Let $\bar\psi_i:Z_i\to \bar M$
be a smooth function, defined on a compact manifold $\bar Z_i$ 
that is dual to $\bar\alpha_i$ (replace $\bar\alpha_i$ 
by an integer multiple of $\bar\alpha_i$, if necessary).  
Lift $\bar\psi_i$ to a $\G$-equivariant map 
$\tilde\psi_i:Q_i\to\mu^{-1}(0)$, 
defined on the total space of a principal $\G$-bundle 
$Q_i\to Z_i$, and consider the induced map 
$$
     \psi_i:Q_i\times_\G\EG^n\to\mu^{-1}(0)\times_\G\EG^n.
$$
It is homologous to $\phi_i$. 
Let $\phi:=\phi_1\times\cdots\times\phi_k$
and $\psi:=\psi_1\times\cdots\times\psi_k$.
Then
$$
     \ev^0_B\cdot f 
     = \ev^0_B\cdot\phi 
     = \ev^0_B\cdot\psi
     = \bev_B\cdot\bar\psi
     = \GW_{\bar B,\Sigma}(\bar\alpha_1,\dots,\bar\alpha_k).
$$
The first equality follows from the definition of $\phi$,
the second from the fact that $\ev^0_B$ is a pseudo-cycle (Step~1)
and $\phi$ is homologous to $\psi$, and the third equality
follows from the definition of the Gromov--Witten invariants
(see~\cite{MS1} for example). 

\medskip
\noindent{\bf Step~3.}
{\it The invariant $\Phi_{B,\Sigma}$ can be expressed 
as the intersection number
$$
      \Phi_{B,\Sigma}(\alpha_1\smile\cdots\smile\alpha_k)
      = \ev^\eps_B\cdot f
$$
for $\eps>0$ sufficiently small.}

\medskip
\noindent
The map $f:Y\to(M^c_\G)^k$ is dual to the 
class $\pi_1^*\alpha_1\smile\cdots\smile\pi_k^*\alpha_k$,
where $\pi_i:(M^c_\G)^k\to M^c_\G$ denotes the projection onto the
$i$th factor.  Moreover, $\ev^\eps_{B,i}=\pi_i\circ\ev^\eps_B$.
Hence
\begin{eqnarray*}
     \ev_B^\eps\cdot f
&= &
     \int_{\Mm^\eps_{B,\Sigma}}
     (\ev_B^\eps)^*\left(\pi_1^*\alpha_1^c\smile\cdots\smile
     \pi_k^*\alpha_k^c\right)  \\
&= &
     \int_{\Mm^\eps_{B,\Sigma}}
     (\ev_{B,1}^\eps)^*\alpha_1^c\smile\cdots\smile
     (\ev_{B,k}^\eps)^*\alpha_k^c \\
&= &
     \int_{\Mm^\eps_{B,\Sigma}}
     (\i^c\circ\ev_{B,1}^\eps)^*\alpha_1\smile\cdots\smile
     (\i^c\circ\ev_{B,k}^\eps)^*\alpha_k \\
&= &
     \int_{\Mm^\eps_{B,\Sigma}}
     \ev_\G^*\left(\alpha_1\smile\cdots\smile\alpha_k\right).
\end{eqnarray*}
Here $\i^c:M_\G^c\to M_\G$ denotes the obvious inclusion.
The last equality follows from the fact that 
$\i^c\circ\ev^\eps_{B,i}:\Mm^\eps_{B,\Sigma}\to M_\G$
is homotopic to the evaluation map $\ev_\G$ in the definition
of $\Phi_{B,\Sigma}$.

\medskip
\noindent{\bf Step~4.}
{\it For $\eps>0$ sufficiently small there is an injective map
$$
      \Tt^\eps_{B,\Sigma;f}:\Mm^0_{B,\Sigma;f}\to\Mm^\eps_{B,\Sigma;f}
$$
such that 
$$
      \Tt^\eps_{B,\Sigma;f}([u_0,A_0],y_{0,1},\dots,y_{0,k})
      = ([u_\eps,A_\eps],y_{\eps,1},\dots,y_{\eps,k})
$$
satisfies
$$
     (u_\eps,A_\eps) = (\exp_{u_0}(\xi_\eps),A_0+\alpha_\eps),\qquad
     \left\|(\xi_\eps,\alpha_\eps)\right\|_{2,p,\eps}
     \le c\eps^{2-2/p},
$$
$$
     \nu^\eps([u_\eps,A_\eps],y_{\eps,1},\dots,y_{\eps,k})
     = \nu^0([u_0,A_0],y_{0,1},\dots,y_{0,k}).
$$
Here $\nu^\eps:\Mm^\eps_{B,\Sigma;f}\to\{\pm1\}$ 
denotes the intersection index of the maps $\ev^\eps_B$ and $f$ 
(in the transverse case).}

\medskip
\noindent
Choose $c_0>0$ such that 
$
     \Mm^0_{B,\Sigma;f} \subset \Mm^0_{B,\Sigma}(c_0)
$
and consider the map 
$$
     (\ev_B^\eps \circ \Tcal^\eps)\times f:
     \Mcal_{B,\Sigma}^0(c_0)\times Y\to M_\G^c\times M_\G^c.
$$
By Proposition~\ref{prop:ev}~(ii), this map converges 
to $\ev^0_B\times f$ in the $C^1$-topology as $\eps$ tends to zero.
By~$(H5)$ the map $\ev^0_B\times f$ is transverse to the 
diagonal $\Delta\subset M_\G^c\times M_\G^c$.  Hence 
$(\ev_B^\eps \circ \Tcal^\eps)\times f$ is transverse to
$\Delta$ for $\eps$ sufficiently small. 
Moreover, by Theorem~\ref{thm:exist}, the image of 
$\Mcal_{B,\Sigma;f}^0$ under $(\ev_B^\eps \circ \Tcal^\eps)\times f$ 
is $\eps^{2-2/p}$-close to $\Delta$.  
Hence, by the implicit function theorem, 
there is, for $\eps$ sufficiently small, 
a unique injective map 
$$
     \Mm^0_{B,\Sigma;f}
     \to ((\ev_B^\eps \circ \Tcal^\eps)\times f)^{-1}(\Delta)
     \subset \Mm^0_{B,\Sigma}(c_0)\times Y
$$
such that the distance between each point and its image 
is bounded above by a constant times $\eps^{2-2/p}$. 
Composing this map with 
$$
     \Tt^\eps\times\id:\Mm^0_{B,\Sigma}\times Y
     \to\Mm^\eps_{B,\Sigma}\times Y
$$
we obtain the required map $\Tt^\eps_{B,\Sigma;f}$. 
By Proposition~\ref{prop:ev}~(i), the map $\Tt^\eps_{B,\Sigma;f}$
identifies the two intersection indices. 

\medskip
\noindent{\bf Step~5.}
{\it Assume $\Sigma=S^2$.  Then there exists a constant
$\eps_0>0$ such that the map 
$\Tcal_{B,\Sigma;f}^\eps:\Mm^0_{B,\Sigma;f}\to\Mm^\eps_{B,\Sigma;f}$ 
of Step~4 is surjective for $0<\eps\le\eps_0$.}

\medskip
\noindent
Suppose, by contradiction, that there exist sequences
$\eps_\nu\to 0$ and 
$$
    ([u_\nu,A_\nu],y_{1\nu},\cdots,y_{k\nu})
    \in\Mm_{B,S^2;f}^{\eps_\nu}
$$
such that 
$$
     ([u_\nu,A_\nu],y_{1\nu},\cdots,y_{k\nu})
     \notin\im\,\Tcal_{B,S^2;f}^{\eps_\nu}.
$$
Consider the sequence 
$$
     C_\nu:=\sup_P \left(|d_{A_\nu}u_\nu|
     +\eps_\nu^{-1}|\mu(u_\nu)|\right).
$$
We prove that $C_\nu$ diverges to $\infty$.
Assume otherwise that $C_\nu$ is bounded.
Then, by Theorem~D, there exists a constant~$c_1>0$
such that $[u_\nu,A_\nu]$ belongs to the image of the map
$
  \Tcal_{B,S^2}^{\eps_\nu}:
  \Mcal_{B,S^2}^0(c_1)\to \Mcal_{B,S^2}^{\eps_\nu}
$
for $\nu$ sufficiently large.   Write
$$
     (u_\nu,A_\nu)=\Tilde\Tt^{\eps_\nu}(u_{0\nu},A_{0\nu}),\qquad
     (u_{0\nu},A_{0\nu})\in\Tilde{\Mm}^0_{B,S^2}(c_1).
$$
Since $\Mm^0_{B,S^2}(c_1)$ is compact 
we may assume that the limit
$$
     (u_0,A_0)
     =\lim_{\nu\to\infty}(u_{0\nu},A_{0\nu})
       \in\Tilde{\Mm}^0_{B,S^2}(c_1)
$$  
exists. Moreover, since $Y$ is compact, we may assume, by passing to a 
further subsequence if necessary, that the limit
$$
     (y_1,\dots,y_k)=
     \lim_{\nu\to\infty}(y_{1\nu},\dots,y_{k\nu})
$$
exists. Since $\ev_B^{\eps_\nu}\circ\Tt^{\eps_\nu}$ converges 
to $\ev_B^0$ in the $C^1$-topology, and 
$$
     \ev_B^{\eps_\nu}\circ\Tt^{\eps_\nu}([u_{0\nu},A_{0\nu}])
     = f(y_{1\nu},\cdots,y_{k\nu})
$$
we deduce that $([u_0,A_0],y_1,\dots,y_k)\in\Mm^0_{B,S^2;f}$ 
and, for $\nu$ sufficiently large,
$$
     (\Tt^{\eps_\nu}([u_{0\nu},A_{0\nu}]),y_{1\nu},\cdots,y_{k\nu})
     = \Tt^{\eps_\nu}_{B,S^2;f}([u_0,A_0],y_1,\dots,y_k).
$$
The last assertion follows from 
the uniqueness part of the implicit function theorem used
in the definition of the maps $\Tt^{\eps_\nu}_{B,S^2;f}$. 
This contradicts our assumption.  Thus we have proved 
that $C_\nu$ diverges to $\infty$ as claimed.

Now choose a sequence $p_\nu\in P$ such that 
$$
     c_\nu := |d_{A_\nu}u_\nu(p_\nu)| 
     + \eps_\nu^{-1}|\mu(u_\nu(p_\nu))| \to \infty.
$$
Passing to a subsequence, if necessary, we may assume that
$p_\nu$ converges.  Denote
$$
    \w := \lim_{\nu\to\infty}\pi(p_\nu).
$$
Moreover, by applying Hofer's trick (see~\cite[Lemma~4.5.3]{MS1} 
for example) we may assume that 
$$
    \sup_{B_{r_\nu}(\pi(p_\nu))}\left(|d_{A_\nu}u_\nu| 
    + \eps_\nu^{-1}|\mu(u_\nu)|\right)
    \le 2c_\nu,\qquad
    r_\nu c_\nu\to\infty.
$$
We distinguish three cases.
\begin{description}
\item[Case~1:]
$c_\nu\eps_\nu\to\infty$.
\item[Case~2:]
There exists a $\delta>0$ such that
$\delta\le c_\nu\eps_\nu\le\delta^{-1}$
for all~$\nu$.
\item[Case~3:]
$c_\nu\eps_\nu\to0$.
\end{description}
We shall prove that in Case~1 a nonconstant 
$J$-holomorphic sphere in $M$ bubbles off 
at the point $\w$, in Case~2 a nontrivial solution of 
the vortex equations~(\ref{eq:vortex}) bubbles off,
and in Case~3 a nonconstant $\bar J$-holomorphic sphere 
in $\bar M$ bubbles off.  
To see this, we choose a local holomorphic coordinate
chart $s+it$ on $\Sigma$ that maps $\w$ to zero, 
identifies a neighbourhood of $\w$ with the ball $B_{2r}$, 
and identifies the volume form $\dvol_\Sigma$ with the form
$\lambda^2ds\wedge dt$, where $\lambda(0)=1$.  
Moreover, we choose a local frame of the bundle
$P$ along this coordinate chart.  
We use the notation of Remark~\ref{rmk:local}.
Then the sequences $u_\nu:B_{2r}\to M$ and 
$\Phi_\nu,\Psi_\nu:B_{2r}\to\g$ satisfy
$$
    v_{\nu s} + Jv_{\nu t} = 0,\qquad
    \lambda^{-2}\kappa_\nu + \eps_\nu^{-2}\mu(u_\nu) = 0,
$$ 
$$
    v_{\nu s}:=\p_su_\nu+X_{\Phi_\nu}(u_\nu),\qquad
    v_{\nu t}:=\p_tu_\nu+X_{\Psi_\nu}(u_\nu),
$$
$$
    \kappa_\nu := \p_s\Psi_\nu - \p_t\Phi_\nu + [\Phi_\nu,\Psi_\nu].
$$
Moreover, there is a sequence $\w_\nu:=(s_\nu,t_\nu)\to 0$ 
such that
$$
    c_\nu
    = \lambda(\w_\nu)^{-1}\left|v_{\nu s}(\w_\nu)\right|
      + \eps_\nu^{-1}\left|\mu(u(\w_\nu))\right| 
    \ge \frac12\sup_{B_{r_\nu}(\w_\nu)}
      \left(\lambda^{-1}\left|v_{\nu s}\right|
      + \eps_\nu^{-1}\left|\mu(u_\nu)\right|
      \right).
$$
Let us define $\tilde u_\nu:B_{r_\nu c_\nu}\to M$
and $\tilde\Phi_\nu,\tilde\Psi_\nu:B_{r_\nu c_\nu}\to\g$ by 
$$
       \tilde u_\nu(z) := u(\w_\nu+c_\nu^{-1}z), 
$$
$$
    \tilde\Phi_\nu(z) := c_\nu^{-1}\Phi_\nu(\w_\nu+c_\nu^{-1}z),\qquad
    \tilde\Psi_\nu(z) := c_\nu^{-1}\Psi_\nu(\w_\nu+c_\nu^{-1}z),
$$
and $\tilde\lambda_\nu:B_{r_\nu c_\nu}\to(0,\infty)$ and
$\tilde J_\nu:B_{r_\nu c_\nu}\to\Jj_\G(M,\om)$ by
$$
    \tilde\lambda_\nu(z) := \lambda(\w_\nu+c_\nu^{-1}z),\qquad
    \tilde J_\nu(z) = J_{\w_\nu+c_\nu^{-1}z}.
$$
Then $\tilde\lambda_\nu$ converges to $1$
in the $\Cinf$-topology and $\tilde J_\nu$ converges
to $J_0$ in the $\Cinf$-topology.  Moreover,
$$
    \tilde v_{\nu s} + \tilde J_\nu\tilde v_{\nu t} = 0,\qquad
    \tilde\lambda_\nu^{-2}\tilde\kappa_\nu 
    + (c_\nu\eps_\nu)^{-2}\mu(\tilde u_\nu) = 0,
$$
$$
    \sup_{B_{r_\nu c_\nu}}\left(\frac{1}{\tilde\lambda_\nu}
    \left|\tilde v_{\nu s}\right|
    + \frac{1}{c_\nu\eps_\nu}\left|\mu(\tilde u_\nu)\right|
    \right)  
    \le  2\left(\frac{1}{\tilde\lambda_\nu(0)}
    \left|\tilde v_{\nu s}(0)\right|
    + \frac{1}{c_\nu\eps_\nu}\left|\mu(\tilde u_\nu(0))\right|\right)
    = 2.
$$
\noindent
{\bf Case~1:} Suppose that $c_\nu\eps_i$ diverges to infinity.
Then, by hypothesis~(H3), the curvature 
$\tilde\kappa_\nu$ converges uniformly to zero.
Hence, by Uhlenbeck's weak compactness theorem~\cite{Uh,We}, 
we may assume that $\tilde\Phi_\nu$ and $\tilde\Psi_\nu$ converge
in the sup-norm and weakly in $W^{1,p}$.  
This implies that the sequence $\tilde u_i$
is bounded in $W^{1,p}$.  Hence, by the usual elliptic 
bootstrapping argument for pseudoholomorphic curves,
it is bounded in $W^{2,p}$ (the lower order terms in the 
equation have the form $X_{\tilde\Phi_i}(\tilde u_i)$
and hence are bounded in $W^{1,p}$). 
Hence there exists a subsequence, 
still denoted by $\tilde u_i$,
that converges strongly in $W^{1,p}$ 
to a $J_0$-holomorphic curve
$\tilde u:\C\to M$ with finite energy.  
Since the sequence $\mu(\tilde u_i(0))$ 
is bounded it follows that
$
     |\p_s\tilde u(0)|
     = \lim_{\nu\to\infty}|\tilde v_{\nu s}(0)|
     =1,
$
and hence $\tilde u$ extends to a nonconstant 
holomorphic sphere in $M$. This contradicts~$(H2)$.

\medskip
\noindent
{\bf Case~2:}
Suppose that the sequence $c_\nu\eps_\nu$
is bounded and does not converge to zero.
Let us assume, without loss of generality, that 
$
    \lim_{\nu\to\infty} c_\nu\eps_\nu = 1.
$
Then we can use the compactness result 
of~\cite{CGMS} to deduce that, after a suitable gauge
transformation and after passing to a further
subsequence, the triple $(\tilde u_\nu,\tilde\Phi_\nu,\tilde\Psi_\nu)$
converges to a solution $(\tilde u,\tilde\Phi,\tilde\Psi)$
of the vortex equations~(\ref{eq:vortex}) with finite
energy.  Moreover, 
$$
     |\p_s\tilde u(0)+X_{\tilde\Phi(0)}(\tilde u(0))|
     + |\mu(\tilde u(0))| = 1
$$
and hence the energy is nonzero.  Hence, by 
Proposition~\ref{prop:vortex}, we have
$$
     E(\tilde u,\tilde\Phi,\tilde\Psi) \ge\hbar.
$$
{\bf Case~3:} Suppose that 
$
    \lim_{\nu\to\infty} c_\nu\eps_\nu = 0.
$
Then, by Lemma~\ref{le:mu}, 
$$
    \sup_\nu(c_\nu\eps_\nu)^{-3/2}\|\mu(\tilde u_\nu)\|_{L^\infty(K)}
    + \sup_\nu(c_\nu\eps_\nu)^{-2}\|\mu(\tilde u_\nu)\|_{L^2(K)}
    < \infty
$$
for every compact set $K\subset\C$. It follows that the 
sequence $\kappa_i$ is uniformly bounded in $L^2$.
Hence, by Uhlenbeck's weak compactness theorem, 
we may assume that $\tilde\Phi_\nu$ and $\tilde\Psi_\nu$ 
converge weakly in $W^{1,2}$ and strongly in $L^p$,
on every compact subset of $\C$.  Here $p$ is any fixed 
real number, say $p>4$. Hence, the sequence $\tilde u_\nu$ 
is bounded in $W^{1,p}$.  Now it follows again from the 
elliptic bootstrapping analysis for pseudoholomorphic
curves that $\tilde u_\nu$ is bounded in $W^{2,2}$
and hence has a subsequence that converges strongly
in $W^{1,p}$ on every compact subset of $\C$. 
The limit $(\tilde u,\tilde\Phi,\tilde\Psi)$ 
is a finite energy solution of~(\ref{eq:jhol})
on $\C$.  This solution represents a $\bar J$-holomorphic 
sphere in the quotient $\bar M$. Moreover, since 
$(c_\nu\eps_\nu)^{-1}|\mu(\tilde u_\nu(0))|\to 0$ it follows
that 
$$
     |\p_s\tilde u(0)+X_{\tilde\Phi(0)}\tilde u(0)|
     = 1
$$
and hence the resulting holomorphic sphere in $\bar M$ is 
nonconstant.  Hence 
$$
     E(\tilde u,\tilde\Phi,\tilde\Psi) \ge \hbar.
$$
Thus we have proved in all three cases that 
$$
     \lim_{\nu\to\infty}E_{B_r(\w)}(u_\nu,A_\nu) \ge \hbar
$$
for every $r>0$.  

This shows that, after passing to 
a suitable subsequence, bubbling can only take place
at finitely many points $\w_1,\dots,\w_\ell\in\Sigma$. 
On every compact subset of $\Sigma\setminus\{\w_1,\dots,\w_\ell\}$ 
the sequence $|d_{A_\nu}u_\nu|+\eps_\nu^{-1}|\mu(u_\nu)|$
is uniformly bounded.  (As an aside: this is used in the proof of 
Lemma~\ref{le:Meps}.)  Hence it follows as in Case~3,
that a suitable subsequence in a suitable gauge 
converges on this complement to a finite energy
solution of~(\ref{eq:jhol}). The limit $(u,A)$ 
descends to a holomorphic curve 
$$
     \bar u:\Sigma\setminus\{\w_1,\dots,\w_\ell\}\to \bar M
$$
with finite energy.  Hence, by the removable singularity 
theorem for $\bar J$-holomorphic curves, it extends
to a holomorphic curve on all of $\Sigma$, 
still denoted by $\bar u$.  The energy of this
$\bar J$-holomorphic curve satisfies 
$$
     E(\bar u) \le \inner{[\bar\om]}{\bar B}-\ell\hbar.
$$
By hypothesis~(H3), the dimension of the moduli space reduces 
by at least $2N$ at each bubble.  Thus the limit $[u,A]$ 
belongs to a moduli space $\Mm^0_{B',S^2}$ of dimension
$$
     \dim\,\Mm^0_{B',S^2} \le \dim\,\Mm^0_{B,S^2} -2N\ell.
$$
If $\{\w_1,\dots,\w_\ell\}\cap\{z_1,\dots,z_k\}=\emptyset$
then the limit curve $(u,A)$ still satisfies 
$\ev_i([u,A])\in f_i(Y_i)$ for every $i$
and hence cannot exist, by the transversality 
condition~$(H5)$.  In general, denote
$$
     I := \left\{i\in\{1,\dots,k\}\,|\,
     z_i\notin \{\w_1,\dots,\w_\ell\}\right\}.
$$
Then the limit $[u,A]$ satisfies
$$
     i\in I\qquad\IMP\qquad
     \ev_i([u,A])\in f_i(Y_i).
$$
Since the points $z_i\in\Sigma$ are pairwise distinct we have
$$
     \ell\ge k-|I|
$$
and so
$$
     \dim\,\Mm^0_{B',S^2} 
     \le \dim\,\Mm^0_{B,S^2} - 2N(k-|I|) 
     < \sum_{i\in I}\deg(\alpha_i).
$$
Here we have used the fact that $\deg(\alpha_i)<2N$ 
for each $i\in\{1,\dots,k\}\setminus I$.
It follows again from~$(H5)$ that such a limit 
curve cannot exist. Hence our assumption that the map
$\Tt^{\eps_i}_{B,S^2;f}$ were not surjective for every $i$
must have been wrong. This proves the theorem.
\end{proof}

\begin{remark}\label{rmk:index}\rm
A more subtle argument, as in Gromov compactness 
for pseudoholomorphic spheres, shows that in the higher 
genus case the limit curve $\bar u$ also satisfies
$$
     \inner{c_1(T\bar M)}{[\bar u]}
     \le\inner{c_1(T\bar M)}{\bar B}-N\ell,
$$
where $\ell$ denotes the number of points near which
bubbling occurs.  Here one needs to prove that no energy gets 
lost and one obtains convergence to a bubble
tree that represents the homology class $B$.
With this refined compactness argument one can extend 
Theorem~A to the higher genus case. 
\end{remark}

\begin{remark}\label{rmk:minenergy}\rm
Assume~$(H1)$, $(H2)$, and~$(H4)$, but not the monotonicity
hypothesis~$(H3)$.  Suppose that the number 
$
     \hbar > 0
$
is a lower bound for the energy of the nonconstant
$\bar J$-holomorphic spheres in $\bar M$ as well 
as for the energy of the nontrivial (that is positive
energy) solutions of the vortex equations~(\ref{eq:vortex}).
Let $(\Sigma,\dvol_\Sigma,j_\Sigma)$ be a compact
Riemann surface of genus $g>0$ and 
suppose that $\bar B\in H_2(\bar M;\Z)$ 
satisfies
$$
    0 \le \inner{[\bar\om]}{\bar B} < \hbar.
$$
Then the moduli space $\Mm^0_{B,\Sigma}$ is compact
and the bubbling argument in the proof of Theorem~A 
together with Proposition~\ref{prop:ev} shows that the map 
$\Tt^\eps:\Mm^0_{B,\Sigma}\to\Mm^\eps_{B,\Sigma}$
of Theorem~\ref{thm:exist} is a diffeomorphism for $\eps>0$
sufficiently small.  Hence in this case the invariants
$\Phi_{B,\Sigma}$ agree with the Gromov--Witten invariants
$\GW_{\bar B,\Sigma}$. 
\end{remark}

 
\appendix


\section{The graph construction}\label{app:graph}

Let $\G$ be a compact Lie group whose Lie algebra 
$\g=\Lie(\G)$ is equipped with an invariant
inner product and $(M,\om)$ be a symplectic manifold 
with a Hamiltonian $\G$-action
generated by a moment map $\mu:M\to\g$. 
We denote by $\g\to\Vect(M):\eta\mapsto X_\eta$ the
infinitesimal action, by $\Cinf_\G(M)$ the space 
of $\G$-invariant smooth functions on $M$,
and by $\Jj_\G(M,\om)$ the space of $\G$-invariant 
and $\om$-compatible almost complex structures on $M$.  
We fix a Riemann surface $(\Sigma,\dvol_\Sigma,j_\Sigma)$ 
and a principal $\G$-bundle $P\to\Sigma$.   
Given a family of almost complex structures 
$
     \Sigma\to\Jj_\G(M,\om):z\mapsto J_z
$
and a 1-form
$
     T\Sigma\to\Cinf_\G(M):\hat z\mapsto H_{\hat z}
$
we consider the perturbed equations
\begin{equation}\label{eq:eps-ham}
     \bar\p_{J,H,A}(u) = 0,\qquad
     *F_A + \eps^{-2}\mu(u) = 0,
\end{equation}
where 
$$
     \bar\p_{J,H,A}(u) := \bar\p_{J,A}(u) + X_H(u)^{0,1}.
$$
Here the $(0,1)$-form $\bar\p_{J,A}(u)\in\Om^{0,1}(\Sigma,u^*TM/\G)$
is understood with respect to the family of almost 
complex structures $J_z$, parametrized by $z\in\Sigma$.  
Moreover, the Hamiltonian perturbation is defined as follows.  
Associated to $H\in\Om^1(\Sigma,\Cinf_\G(M))$
is the $1$-form $X_H\in\Om^1(\Sigma,\Vect_\G(M,\om))$
which assigns to every $\hat z\in T_z\Sigma$ the 
Hamiltonian vector field $X_{H,\hat z}$ associated
to the Hamiltonian function $H_{\hat z}:M\to\R$. 
Thus 
$
     \i(X_{H,\hat z})\om = dH_{\hat z}.
$
The 1-form $X_H(u)\in\Om^1(\Sigma,u^*TM/\G)$
lifts to an equivariant and horizontal $1$-form 
on $P$ with values in $u^*TM$, also denoted 
by $X_H(u)$ and defined by
$$
     (X_H(u))_p(v) := X_{H,d\pi(p)v}(u(p)).
$$
The complex anti-linear part of this 1-form
is the Hamiltonian term in the definition
of $\bar\p_{J,H,A}(u)$. In this section we show how
to reduce the perturbed equations~(\ref{eq:eps-ham})
to~(\ref{eq:eps}) via Gromov's graph construction~\cite{Gr}.

Let us denote by $\alpha_H\in\Om^1(\Sigma\times M)$ 
the 1-form associated to $H$. Thus $\alpha_H$ assigns
to every pair of tangent vactors 
$(\hat z,\hat x)\in T_z\Sigma\times T_xM$
the real number $H_{\hat z}(x)$. Denote
$$
    \tilde M := \Sigma\times M.
$$
The 2-form
$$
    \tilde\om := \om - d\alpha_H + c\,\dvol_\Sigma
$$
is a symplectic form on $\tilde M$ whenever the 
constant $c$ is sufficiently large. Here we have abused 
notation and denoted by $\om$ the pullback of the 
2-form $\om$ on $M$ under the obvious
projection $\Sigma\times M\to M$ and likewise for $\dvol_\Sigma$.
To see that $\tilde\om$ is symplectic for large $c$, note
first that $\tilde\om$ is a connection form: it is closed
and its restriction to each fibre $\{z\}\times M$ is symplectic. 
The curvature of this connection form is the $2$-form
$$
    \Om_H\,\dvol_\Sigma
    := dH +\frac12\{H\wedge H\}
    \in\Om^2(\Sigma,\Cinf_\G(M)).
$$
This identity defines the function $\Om_H:\Sigma\times M\to\R$. 
Now the top exterior power of $\tilde\om$ is given by 
$$
    \frac{\tilde\om^{n+1}}{(n+1)!} 
    = (c-\Om_H)\frac{\om^n}{n!}\wedge\dvol_\Sigma,
$$
where $\dim M=2n$.  Hence $\tilde\om$ is nondegenerate
whenever $c>\max\Om_H$. 
Now consider the almost complex structure $\tilde J$
on $\tilde M$ given by 
$$
    \tilde J(z,x) 
    := \left(
    \begin{array}{cc}
     j_\Sigma(z) & 0 \\
     J(z,x)\circ X_H(z,x)-X_H(z,x)\circ j_\Sigma(z) & J(z,x)
    \end{array}\right).
$$
Here $J(z,x):=J_z(x)$ and we denote by 
$
    X_H(z,x):T_z\Sigma\to T_xM
$
the linear map $\hat z\mapsto X_{H,\hat z}(x)$.
Lemma~\ref{le:Jtilde} below shows that 
$\tilde J$ is compatible with $\tilde\om$. 

\begin{lemma}\label{le:eps-tilde}
Let $(u,A)\in\Cinf_\G(P,M)\times\Aa(P)$ and 
define $\tilde u:P\to\tilde M$ by 
$\tilde u(p):=(\pi(p),u(p))$.  Then 
$u$ and $A$ satisfy~(\ref{eq:eps-ham})
if and only if $\tilde u$ and $A$ satisfy
$$
     \bar\p_{\tilde J,A}(\tilde u) = 0,\qquad
     *F_A + \eps^{-2}\tilde\mu(\tilde u) = 0.
$$
Here $\tilde\mu:\tilde M\to\g$ is defined by
$\tilde\mu(z,x):=\mu(x)$. 
\end{lemma}

\begin{proof}
By definition of $\tilde J$ we have 
$\bar\p_{\tilde J,A}(\tilde u) = (0,\bar\p_{J,H,A}(u))$.
Alternatively, we can compute in local holomorphic
coordinates $s+it$ on $\Sigma$. In such coordinates
the Hamiltonian perturbation, the connection $A$,
and the volume form on $\Sigma$ have the form 
$$
     H = F\,ds + G\,dt,\qquad 
     A = \Phi\,ds + \Psi\,dt,\qquad
     \dvol_\Sigma = \lambda^2\,ds\wedge dt
$$
and the equations~(\ref{eq:eps-ham}) have the form
\begin{eqnarray*}
     \p_su+L_u\Phi+X_F(u) + J(\p_tu+L_u\Psi+X_G(u)) &= &0, \\
     \p_s\Psi - \p_t\Phi + [\Phi,\Psi] 
     + (\lambda/\eps)^2\mu(u) &= &0.
\end{eqnarray*}
Moreover, the almost complex structure $\tilde J$ 
is given by 
$$
    \tilde J
    := \left(
    \begin{array}{ccc}
     0 & -1  & 0 \\
     1 &  0  & 0 \\
     JX_F-X_G & JX_G+X_F & J
    \end{array}\right).
$$
This proves the lemma.
\end{proof}

\begin{lemma}\label{le:Jtilde}
Let $\hat z_i\in T_z\Sigma$ and $\hat x_i\in T_xM$
for $i=1,2$.  Then 
$$
    \tilde\om((\hat z_1,\hat x_1),(\hat z_2,\hat x_2))
    = \INNER{\hat x_1+X_{H,\hat z_1}(x)}{\hat x_2+X_{H,\hat z_2}(x)}_z
    + (c-\Om_H)\inner{\hat z_1}{\hat z_2},
$$
where $\INNER{\cdot}{\cdot}_z:=\om(\cdot,J_z\cdot)$.
\end{lemma}

\begin{proof}
Continue the notation of the proof of Lemma~\ref{le:eps-tilde}.
Then the curvature $\Om_H$ is given by 
$$
    \lambda^2\Om_H
    = \p_sG - \p_tF + \{F,G\},
$$
where $\{F,G\}:=\om(X_F,X_G)$ denotes the Poisson bracket on $M$,
and
$$
    \tilde\om = \om - dF\wedge ds - dG\wedge dt 
    + (\p_tF-\p_sG + c\lambda^2)ds\wedge dt
$$
where $dF$ and $dG$ denote the differential on $M$. 
Abbreviate $\zeta_i:=(\hat s_i,\hat t_i,\hat x_i)$
and $\xi_i:=\hat x_i+\hat s_iX_F+\hat t_iX_G$
for $i=1,2$. Then
\begin{eqnarray*}
   \tilde\om(\zeta_1,\tilde J\zeta_2)
&= &
   \tilde\om(\zeta_1,
   (-\hat t_2,\hat s_2,J\xi_2+\hat t_2X_F-\hat s_2X_G)) \\
&= &
   \om(\hat x_1,J\xi_2+\hat t_2X_F-\hat s_2X_G)  \\
&&
   +\, \hat t_2dF(\hat x_1) 
   + \hat s_1dF(J\xi_2+\hat t_2X_F-\hat s_2X_G)  \\
&&
   -\,\hat s_2dG(\hat x_1)
   + \hat t_1dG(J\xi_2+\hat t_2X_F-\hat s_2X_G)  \\
&&
   +\,(\p_tF-\p_sG+c\lambda^2)(\hat s_1\hat s_2+\hat t_1\hat t_2) \\
&= &
   \om(\xi_1,J\xi_2)
   + (c\lambda^2-\p_sG+\p_tF-\{F,G\})
   (\hat s_1\hat s_2+\hat t_1\hat t_2).
\end{eqnarray*}
The last identity follows from the fact that
$\{F,G\}=dF(X_G)=-dG(X_F)$ and $dF(J\xi_2)=\om(X_F,J\xi_2)$.
\end{proof}


\section{Cauchy--Riemann operators}\label{app:CR}

Fix a compact Lie group $\G$,
an invariant inner product on the Lie algebra $\g=\Lie(\G)$, 
a symplectic manifold $(M,\om)$, a Hamiltonian $\G$-action 
on $M$ generated by a moment map $\mu:M\to\g$,
a compact Riemann surface $(\Sigma,j_\Sigma,\dvol_\Sigma)$,
a principal $\G$-bundle $P\to\Sigma$, 
and a family of $\G$-invariant and $\om$-compatible 
almost complex structures 
$\Sigma\to\Jj_\G(M,\om):z\mapsto J_z$.
Each almost complex structure determines a Riemannian metric
$
      \INNER{\cdot}{\cdot}_z := \om(\cdot,J_z\cdot)
$
on $M$ and hence a Levi-Civita connection $\nabla=\nabla_z$.
The value of $z$ will usually be clear from the context
and we shall omit the subscript $z$. 
Let $u:P\to M$ be an equivariant smooth map 
and $A$ be a connection on $P$. Then $A$ and $\nabla$
determine a connection $\Nabla{A}$ on $u^*TM/\G$ given by
$$
      \Nabla{A}\xi := \nabla{\xi} + \Nabla{\xi}X_A(u)
$$
for $\xi\in\Om^0(\Sigma,u^*TM/\G)$. 
More precisely, we think of $A$ as an equivariant 
1-form on $P$ with values in the Lie algebra $\g$
which identifies the vertical tangent bundle with $\g$. 
A section $\xi$ of $u^*TM/\G$ lifts to an equivariant
section of the bundle $u^*TM\to P$ (also denoted by $\xi$)
and a 1-form $\theta\in\Om^1(\Sigma,u^*TM/\G)$ lifts to
an equivariant and horizontal 1-form on $P$ with values in $u^*TM$ 
(also denoted by $\theta$). In this notation the
1-form $\Nabla{A}\xi$ is given by 
$$
      (\Nabla{A}\xi)_p(v)
      := \Nabla{v}\xi(p) + \Nabla{\xi(p)}X_{A_p(v)}(u(p))
$$
for $v\in T_pP$. In general, $\Nabla{A}$ 
preserves neither the inner product nor the complex structure
on $u^*TM/\G$.  More precisely, let $J_u\in\Om^0(P,\End(u^*TM))$ 
be given by $J_u(p):=J_{\pi(p)}(u(p))\in\End(T_{u(p)}M)$.  
This section is equivariant and hence descend to a complex structure,
also denoted by $J_u$, of the bundle $\End(u^*TM/\G)\to\Sigma$.

\begin{lemma}\label{le:nablaJ}
The covariant derivative of $J_u$ is given by 
$$
     \Nabla{A}J_u = \Nabla{d_Au}J(u) + \dot J(u).
$$
where $\dot J(u)\in\Om^1(\Sigma,\End(u^*TM/\G))$
is defined by 
$$
     \dot J(u)_p(v) := \left.\frac{d}{dt}\right|_{t=0} 
     J_{\gamma(t)}(u(p))
$$
for $v\in T_pP$ and a smooth path $\gamma:\R\to P$ 
such that $\gamma(0)=p$ and $\dot\gamma(0)=v$. 
\end{lemma}

\begin{proof}
Since $J$ is $\G$-invariant we have 
$\Ll_{X_\eta}J=0$ for every $\eta\in\g$.
This formula can be expressed in the form 
\begin{equation}\label{eq:LJ}
     (\Nabla{X_\eta}J)\xi 
     = \Nabla{J\xi}X_\eta - J\Nabla{\xi}X_\eta.
\end{equation}
Using this formula we obtain
\begin{eqnarray*}
    (\Nabla{A}J_u)\xi
&= &
    \Nabla{A}(J_u\xi) - J_u(\Nabla{A}\xi)  \\
&= &
    \nabla(J_u\xi) - J_u\nabla\xi  
    + \Nabla{J_u\xi}X_A(u) 
    - J_u\Nabla{\xi}X_A(u) \\
&= &
    \left(\Nabla{du}J(u)+\dot J(u)\right)\xi
    + \left(\Nabla{X_A(u)}J(u)\right)\xi  \\
&= &
    \left(\Nabla{d_Au}J(u)+\dot J(u)\right)\xi
\end{eqnarray*}
as claimed.
\end{proof}

It follows from Lemma~\ref{le:nablaJ} 
that the complex linear part of the 
connection $\Nabla{A}$ is the connection 
$\Tabla{A}$ on $u^*TM/\G$ given by 
\begin{eqnarray}\label{eq:Tabla}
     \Tabla{A}\xi 
&:= &
     \Nabla{A}\xi - \frac12J_u(\Nabla{A}J_u)\xi \\
&= &
     \nabla\xi + \Nabla{\xi}X_A(u)
     - \frac12J_u\left(\Nabla{d_Au}J(u)+\dot J(u)\right)\xi.
     \nonumber
\end{eqnarray}

\begin{lemma}\label{le:nabla}
$\Tabla{A}$ is a Hermitian connection on $u^*TM/\G$.
\end{lemma}

\begin{proof}
We shall use the identity
$
     \inner{(\Nabla{\xi_1}J)\xi_2}{\xi_3}
     + \mbox{cyclic}
     = 0.
$
By~(\ref{eq:LJ}), we have
\begin{eqnarray*}
&&
     \inner{\xi_1}{\Nabla{\xi_2}X_\eta}
     + \inner{\Nabla{\xi_1}X_\eta}{\xi_2}  \\
&&= 
     \inner{J\xi_1}{J\Nabla{\xi_2}X_\eta}  
     + \inner{J\Nabla{\xi_1}X_\eta}{J\xi_2}  \\
&&= 
     \inner{J\xi_1}{\Nabla{\xi_2}(JX_\eta)-(\Nabla{\xi_2}J)X_\eta} 
     + \inner{\Nabla{J\xi_1}X_\eta-(\Nabla{X_\eta}J)\xi_1}{J\xi_2} \\
&&= 
     \inner{J\xi_1}{\Nabla{\xi_2}(JX_\eta)}  
     + \inner{\Nabla{J\xi_1}X_\eta}{J\xi_2} 
     - \inner{J\xi_1}{(\Nabla{\xi_2}J)X_\eta}
     - \inner{\xi_2}{(\Nabla{X_\eta}J)J\xi_1} \\
&&=  
     \inner{J\xi_1}{\Nabla{\xi_2}(JX_\eta)}
     + \inner{\Nabla{J\xi_1}X_\eta}{J\xi_2}
     + \inner{X_\eta}{(\Nabla{J\xi_1}J)\xi_2} \\
&&=  
     \inner{\Nabla{J\xi_1}(JX_\eta)}{\xi_2}
     - \inner{J(\Nabla{J\xi_1}X_\eta)}{\xi_2}
     - \inner{(\Nabla{J\xi_1}J)X_\eta}{\xi_2} \\
&&= 
     0
\end{eqnarray*}
for $\xi_1,\xi_2\in T_xM$ and $\eta\in\g$. 
Here the penultimate equality follows from the 
fact that $JX_\eta$ is a gradient vector field
and that $\Nabla{J\xi_1}J$ is skew-adjoint. 
This shows that $\nabla X_A(u)$ is a $1$-form on $\Sigma$
with values in the bundle of skew-Hermitian endomorphisms
of $u^*TM/\G$, and so is $J(\Nabla{d_Au}J)$.  Moreover,
since
$$
       d\inner{\xi_1}{\xi_2}
       = \inner{\nabla\xi_1}{\xi_2} + \inner{\xi_1}{\nabla\xi_2}
         - \inner{\xi_1}{J\dot J\xi_2},
$$
the operator $\xi\mapsto\nabla\xi-\frac12J\dot J\xi$
is a Riemannian connection.  By~(\ref{eq:Tabla}), this shows
that $\Tabla{A}$ is a Riemannian connection.  It follows 
directly from the definition that $\Tabla{A}$ preserves the 
complex structure $J_u$. 
\end{proof}

\begin{lemma}\label{le:gauge1}
For every gauge transformation $g\in\Gg(P)$ and every
section $\xi\in\Om^0(\Sigma,u^*TM/\G)$ we have
$$
     \Nabla{g^*A}(g^{-1}\xi) = g^{-1}\Nabla{A}\xi,\qquad
     \Tabla{g^*A}(g^{-1}\xi) = g^{-1}\Tabla{A}\xi.
$$
\end{lemma}

\begin{proof}
Since the metric $\INNER{\cdot}{\cdot}_z$ is $\G$-invariant
for every $z\in\Sigma$ we have 
$$
     \nabla(g^{-1}\xi) 
     = g^{-1}\nabla\xi
       -\Nabla{g^{-1}\xi}X_{g^{-1}dg}(g^{-1}u).
$$
Hence the first identity follows from the fact that
$g^*A=g^{-1}dg+g^{-1}Ag$ and that 
$
     \Nabla{g^{-1}\xi}X_{g^{-1}\eta g}(g^{-1}x)
     = g^{-1}\Nabla{\xi}X_\eta(x).
$
The second identity follows from the first and the 
fact that $J_z$ is $\G$-invariant for every $z$. 
\end{proof}

\begin{lemma}\label{le:curv}
Suppose that $J$ is independent of $z\in\Sigma$. 
Then the curvature of the connection $\Nabla{A}$ is the 
equivariant and horizontal 2-form 
$
     F^{\Nabla{A}}\in\Om^2(P,\End(u^*TM))
$
given by
$$
    F^{\Nabla{A}}(v_1,v_2)\xi
    = R(d_Au(p)v_1,d_Au(p)v_2)\xi
      + \Nabla{\xi}X_{F_A(v_1,v_2)}(u(p))
$$
for $v_1,v_2\in T_pP$ and $\xi\in T_{u(p)}M$,
where $R\in\Om^2(M,\End(TM))$ is the Riemann
curvature tensor of the metric 
$\inner{\cdot}{\cdot}=\om(\cdot,J\cdot)$. 
This 2-form descends to a 2-form 
on $\Sigma$ with values in $\End(u^*TM/\G)$, 
also denoted by $F^{\Nabla{A}}$.
\end{lemma}

\begin{proof}
Given a map $\R^2\to M:(s,t)\mapsto u(s,t)$,
a vector field $\xi(s,t)\in T_{u(s,t)}M$ along $u$,
and a $\G$-connection 
$
     A = \Phi\,ds + \Psi\,dt,
$
where $\Phi,\Psi:\R^2\to\g$, we denote
$$
\begin{array}{rclcrcl}
     v_s &:=& \p_su+X_\Phi(u), &\qquad &
     v_t &:=& \p_tu+X_\Psi(u),  \\
     \Nabla{A,s}\xi &:= & \Nabla{s}\xi+\Nabla{\xi}X_\Phi(u),&\qquad&
     \Nabla{A,t}\xi &:= & \Nabla{t}\xi+\Nabla{\xi}X_\Psi(u).
\end{array}
$$
Then the assertion can be restated in the form
$$
     \Nabla{A,s}\Nabla{A,t}\xi-\Nabla{A,t}\Nabla{A,s}\xi
     = R(v_s,v_t)\xi 
     + \Nabla{\xi}X_{\p_s\Psi-\p_t\Phi+[\Phi,\Psi]}(u).
$$
To prove this we use the identities 
\begin{eqnarray*}
     \Nabla{A,s}\Nabla{A,t}\xi-\Nabla{A,t}\Nabla{A,s}\xi
 &= &
     \Nabla{s}\Nabla{t}\xi-\Nabla{t}\Nabla{s}\xi
    + \Nabla{s}\Nabla{\xi}X_\Psi(u)
     - \Nabla{\Nabla{s}\xi}X_\Psi(u) \\
 &&
     -\, \Nabla{t}\Nabla{\xi}X_\Phi(u)
     + \Nabla{\Nabla{t}\xi}X_\Phi(u) \\
 &&    
     + \Nabla{\Nabla{\xi}X_\Psi}X_\Phi(u)
     - \Nabla{\Nabla{\xi}X_\Phi}X_\Psi(u), \\
     R(\p_su,\p_tu)\xi 
 &= & \Nabla{s}\Nabla{t}\xi-\Nabla{t}\Nabla{s}\xi, \\
     R(\p_su,X_\Psi(u))\xi
 &= & \Nabla{s}\Nabla{\xi}X_\Psi(u)
       - \Nabla{\Nabla{s}\xi}X_\Psi(u)
       - \Nabla{\xi}X_{\p_s\Psi}(u),  \\
     R(X_\Phi(u),X_\Psi(u))\xi
 &= &
     -\Nabla{\xi}X_{[\Phi,\Psi]}(u)
       + \Nabla{\Nabla{\xi}X_\Psi}X_\Phi(u)
       - \Nabla{\Nabla{\xi}X_\Phi}X_\Psi(u).
\end{eqnarray*}
The first and second identities are the definition 
of the connection $\Nabla A$ and the curvature
tensor~$R$. The other identities use the equations
$$
    \Nabla{X_\eta}Z-\Nabla{Z}X_\eta
    = [Z,X_\eta] = 0,\qquad
    \Nabla{Z}[X_{\eta_1},X_{\eta_2}]
    = [\Nabla{Z}X_{\eta_1},X_{\eta_2}],
$$
for every $\G$-invariant vector field $Z\in\Vect_\G(M)$ 
and all $\eta,\eta_1,\eta_2\in\g$. 
\end{proof}

Now consider the Cauchy--Riemann operator 
$$
     D := D_{(u,A)}: 
     \Om^0(\Sigma,u^*TM/\G)\to\Om^{0,1}(\Sigma,u^*TM/\G)
$$
given by 
\begin{equation}\label{eq:Du}
     D\xi
     := \left(\Nabla{A}\xi\right)^{0,1}
       - \frac12J(\Nabla{\xi}J)\p_{J,A}(u).
\end{equation}
In the case $\bar\p_{J,A}(u)=0$ this operator 
is the vertical differential of the section
$u\mapsto\bar\p_{J,A}(u)$ of the infinite dimensional
vector bundle over the space $\Cinf_\G(P,M)$
with fibre $\Om^{0,1}(\Sigma,u^*TM/\G)$ over $u$. 

In the following we denote the Nijenhuis tensor
of $J$ by $N\in\Om^2(TM,TM)$.  It is given by 
\begin{eqnarray*}
     N(\xi_1,\xi_2)
&= &
     [\xi_1,\xi_2] + J[J\xi_1,\xi_2]
     + J[\xi_1,J\xi_2] - [J\xi_1,J\xi_2] \\
&= &
     2J(\Nabla{\xi_2}J)\xi_1 
     - 2J(\Nabla{\xi_1}J)\xi_2.
\end{eqnarray*}

\begin{lemma}\label{le:D}
The complex linear part of $D$ is the operator
$\xi\mapsto(\Tabla{A}\xi)^{0,1}$.  Moreover,
$$
     D\xi = (\Tabla{A}\xi)^{0,1} 
     + \frac14 N(\xi,\p_{J,A}(u))
     + \frac12 (J\dot J\xi)^{0,1}.  
$$
\end{lemma}

\begin{proof}
By definition of $\Tabla{A}$, we have 
$$
    D\xi = (\Tabla{A}\xi)^{0,1} 
      + \frac12J\left(\Nabla{\p_{J,A}(u)}J(u)+\dot J(u)^{0,1}\right)\xi
      - \frac12J(\Nabla{\xi}J)\p_{J,A}(u).
$$
Hence the formula for $D\xi$ follows from the 
relation between the Nijenhuis tensor and 
$\nabla J$. Now this equation shows that the 
operator $\xi\mapsto D\xi - (\Tabla{A}\xi)^{0,1}$
is complex anti-linear. 
\end{proof}


\section{Invariant metrics}\label{app:invariant}
 
Let $M$ be a (complete) Riemannian $m$-manifold.
For $v\in T_xM$ and $i,j\in\{1,2\}$ there exist
linear maps
$$
    E_i(x,v):T_xM\to T_{\exp_x(v)}M,\qquad
    E_{ij}(x,v):T_xM\oplus T_xM\to T_{\exp_x(v)}M
$$
characterized by the following conditions. 
If $x:\R\to M$ is a smooth curve and $v,w:\R\to x^*TM$ 
are vector fields along $x$ then
\begin{eqnarray*}
     \frac{d}{dt} \exp_{x(t)}(v(t))
&= & 
     E_1(x,v)\dot x + E_2(x,v)\Nabla{t}v,  \\
     \Nabla{t}(E_1(x,v)w)
&= &
     E_{11}(x,v)(w,\dot x) + E_{12}(x,v)(w,\Nabla{t}v)
      + E_1(x,v)\Nabla{t}w,  \\
     \Nabla{t}(E_2(x,v)w)
&= &
     E_{21}(x,v)(w,\dot x) + E_{22}(x,v)(w,\Nabla{t}v)
      + E_2(x,v)\Nabla{t}w.
\end{eqnarray*}
Note that the map $E_{11}(x,v)(w,w')$ 
is not symmetric in $w$ and $w'$.  It satisfies
$$
     E_{11}(x,v)(w,w') - E_{11}(x,v)(w',w)
     = E_2(x,v)R(w,w')v,
$$
where $R\in\Om^2(M,\End(TM))$ denotes the 
curvature tensor. However, 
$$
     E_{12}(x,v)(w,w')=E_{21}(x,v)(w',w),
$$
and $E_{22}(x,v)(w,w')$ is symmetric in $w$ and $w'$.
(See~\cite{Ga} for more details.)
Now let 
$
   \G\times M\to M:(g,x)\mapsto gx
$
be a smooth action of a compact Lie group $\G$
with infinitesimal action
$
    \g\to\Vect(M):\eta\mapsto X_\eta.
$
Assume that $M$ is equipped with a $\G$-invariant
Riemannian metric.

\begin{lemma}\label{le:group}
\begin{eqnarray*}
     X_\eta(\exp_x(v))
&= & 
     E_1(x,v)X_\eta(x) + E_2(x,v)\Nabla{v}X_\eta(x),  \\
     \Nabla{E_i(x,v)w}X_\eta(\exp_x(v))  
&= &
     E_{i1}(x,v)(w,X_\eta(x)) 
     + E_{i2}(x,v)(w,\Nabla{v}X_\eta(x)) \\
&&
     +\, E_i(x,v)\Nabla{w}X_\eta(x).  
\end{eqnarray*}
\end{lemma}
 
\begin{proof}  
Since the group action preserves 
geodesics,
$
     g\exp_x(v) = \exp_{gx}(gv).
$
Differentiate this identity with respect to $g$ 
to obtain the first identity. To prove the second differentiate 
the first identity covariantly and use the definition
of $E_i$ and $E_{ij}$.  
For more details see~\cite{Ga}.
\end{proof}
 
For each $x\in M$ denote by 
$
     L_x:\g\to T_xM
$
the infinitesimal action, i.e.
$
     L_x\eta = X_\eta(x).
$
Given a map $u:\R^2\to M$,
a vector field $\xi:\R^2\to u^*TM$ along $u$,
a function $\eta:\R^2\to\g$, 
and a $\G$-connection 
$
     A = \Phi\,ds + \Psi\,dt,
$
where $\Phi,\Psi:\R^2\to\g$, we denote
$$
\begin{array}{rclrcl}
     v_s &= &\p_su+X_\Phi(u),&    
     v_t &= &\p_tu+X_\Psi(u),  \\
     \Nabla{A,s}\xi &= & \Nabla{s}\xi+\Nabla{\xi}X_\Phi(u),&
     \Nabla{A,t}\xi &= & \Nabla{t}\xi+\Nabla{\xi}X_\Psi(u),  \\
     \Nabla{A,s}\eta &=& \p_s\eta+[\Phi,\eta],&
     \Nabla{A,t}\eta &=& \p_t\eta+[\Psi,\eta].
\end{array}
$$
Define $\rho\in\Om^2(M,\G)$ by
$$
    \inner{\eta}{\rho(\xi,\xi')} 
    := \inner{\Nabla{\xi}X_\eta(x)}{\xi'}
     = -\inner{\Nabla{\xi'}X_\eta(x)}{\xi}
$$
for $\xi,\xi'\in T_xM$ and $\eta\in\g$.

\begin{lemma}\label{le:rho}
With the above notation we have
$$
\begin{array}{rclrcl}
   \Nabla{A,s}L_u\eta - L_u\Nabla{A,s}\eta
   &= &\Nabla{v_s}X_\eta(u),&
   \Nabla{A,t}L_u\eta - L_u\Nabla{A,t}\eta
   &= &\Nabla{v_t}X_\eta(u),\\
   \Nabla{A,s}L_u^*\xi - L_u^*\Nabla{A,s}\xi
   &= &\rho(v_s,\xi),&
   \Nabla{A,t}L_u^*\xi - L_u^*\Nabla{A,t}\xi
   &= &\rho(v_t,\xi)
\end{array}
$$
\end{lemma}
 
\begin{proof} See~\cite{Ga}. 
\end{proof}
 
Now let $M^*$ denote the subset of all points $x\in M$
with finite isotropy subgroup 
$
     \G_x := \{g\in\G\,|\,gx=x\}.
$ 
Thus $x\in M^*$ if and only if the linear map
$L_x:\g\to T_xM$ is injective. 
Hence, for every map 
$
     u_0:\R^2\to M^*
$
there exists a unique $\G$-connection 
$A_0=\Phi_0\,ds+\Psi_0\,dt$ such that 
$$
     L_{u_0}^*v_{0s} = L_{u_0}^*v_{0t} = 0,
$$
where 
$$
     v_{0s} := \p_su_0+L_{u_0}\Phi_0,\qquad 
     v_{0t} := \p_tu_0+L_{u_0}\Psi_0.
$$
Let 
$
     \xi_0:\R^2\to u_0^*TM
$ 
be a vector field along 
$u_0$, consider the map 
$$
     u(s,t) := \exp_{u_0(s,t)}(\xi_0(s,t)).
$$
and abbreviate 
$$
     \phi_0 := \Phi-\Phi_0,\qquad \psi_0 := \Psi-\Psi_0.
$$
 
\begin{lemma}\label{le:Lu-phi}
$$
\begin{array}{rcl}
     L_u\phi_0
     &= & v_s-E_1(u_0,\xi_0)v_{0s}-E_2(u_0,\xi_0)\Nabla{A_0,s}\xi_0, \\
     L_u\psi_0
     &= & v_t-E_1(u_0,\xi_0)v_{0t}-E_2(u_0,\xi_0)\Nabla{A_0,t}\xi_0.
\end{array}
$$
\end{lemma}
 
\begin{proof}  We compute
\begin{eqnarray*}
     L_u(\Phi-\Phi_0)
 &= &
     X_\Phi(u) - E_1(u_0,\xi_0)X_{\Phi_0}(u_0)
     - E_2(u_0,\xi_0)\Nabla{\xi_0}X_{\Phi_0}(u_0)  \\
 &= &
     \p_su - E_1(u_0,\xi_0)\p_su_0 - E_2(u_0,\xi_0)\Nabla{s}\xi_0 \\
 &&
     + X_\Phi(u) - E_1(u_0,\xi_0)X_{\Phi_0}(u_0)
     - E_2(u_0,\xi_0)\Nabla{\xi_0}X_{\Phi_0}(u_0)  \\
 &= &
     v_s-E_1(u_0,\xi_0)v_{0s}-E_2(u_0,\xi_0)\Nabla{A_0,s}\xi_0.
\end{eqnarray*}
Here the first equation follows from Lemma~\ref{le:group},
the second from the definition of $E_i$, 
and the last from the definitions of $v_s$ and $v_{0s}$.
\end{proof}
 
In the proof we did not use the fact that
$L_{u_0}^*v_{0s}=L_{u_0}^*v_{0t}=0$.
Now suppose
$
    L_{u_0}^*\xi_0 = 0.
$
Then, by Lemma~\ref{le:rho}, 
$$
    L_{u_0}^*\Nabla{A_0,s}\xi_0 = -\rho(v_{0s},\xi_0),\qquad
    L_{u_0}^*\Nabla{A_0,t}\xi_0 = -\rho(v_{0t},\xi_0).
$$
Abbreviating $E_1=E_1(u_0,\xi_0)$ and 
$E_2=E_2(u_0,\xi_0)$ we obtain the following.

\begin{corollary}\label{cor:8.12}
If $L_{u_0}^*\xi_0=0$ then
$$
\begin{array}{rcl}
    L_{u_0}^*{E_2}^{-1}L_u(\Phi-\Phi_0)
    &= &L_{u_0}^*{E_2}^{-1}(v_s-E_1v_{0s})
      +\rho(v_{0s},\xi_0),\\
    L_{u_0}^*{E_2}^{-1}L_u(\Psi-\Psi_0)
    &= &L_{u_0}^*{E_2}^{-1}(v_t-E_1v_{0t})
      +\rho(v_{0t},\xi_0).
\end{array}
$$
\end{corollary}
 
\begin{proof}  Lemma~\ref{le:Lu-phi}.
\end{proof}
 
\begin{lemma}\label{le:1}
Assume $L_{u_0}^*\xi_0=0$ and abbreviate 
$E_i:=E_i(u_0,\xi_0)$ and $E_{ij}:=E_{ij}(u_0,\xi_0)$.
Then 
\begin{eqnarray*}
     L_u\Nabla{A,t}\phi_0
 &= &
     \Nabla{A,t}v_s
     + \Nabla{X_{\phi_0}}X_{\psi_0}(u) 
     - \Nabla{v_t}X_{\phi_0}(u) 
     - \Nabla{v_s}X_{\psi_0}(u) \\
 &&
     - E_{11}(v_{0s},v_{0t})
     - E_{12}(v_{0s},\Nabla{A_0,t}\xi_0)   \\
 &&
     - E_{21}(\Nabla{A_0,s}\xi_0,v_{0t}) 
     - E_{22}(\Nabla{A_0,s}\xi_0,\Nabla{A_0,t}\xi_0) \\
 &&
     - E_1\Nabla{A_0,t}v_{0s} 
     - E_2\Nabla{A_0,t}\Nabla{A_0,s}\xi_0,
\end{eqnarray*}
\begin{eqnarray*}
     L_u\Nabla{A,s}\phi_0
&= &
     \Nabla{A,s}v_s
     + \Nabla{X_{\phi_0}}X_{\phi_0}(u) 
     - 2\Nabla{v_s}X_{\phi_0}(u)  \\
&&
     - E_{11}(v_{0s},v_{0s})
     - E_{12}(v_{0s},\Nabla{A_0,s}\xi_0)   \\
&&
     - E_{21}(\Nabla{A_0,s}\xi_0,v_{0s}) 
     - E_{22}(\Nabla{A_0,s}\xi_0,\Nabla{A_0,s}\xi_0) \\
&&
     - E_1\Nabla{A_0,s}v_{0s} 
     - E_2\Nabla{A_0,s}\Nabla{A_0,s}\xi_0.
\end{eqnarray*}
\end{lemma}
 
\begin{proof}  
We only prove the first identity.  The proof of the 
second is similar.  By Lemmata~\ref{le:rho} 
and~\ref{le:Lu-phi},
\begin{eqnarray*}
     L_u\Nabla{A,t}\phi_0
&= &
     \Nabla{A,t}L_u\phi_0
     - \Nabla{v_t}X_{\phi^0}(u) \\
&= &
     \Nabla{A,t}v_s
     - \Nabla{A,t}(E_1v_{0s}) 
     - \Nabla{A,t}(E_2\Nabla{A_0,s}\xi_0)
     - \Nabla{v_t}X_{\phi_0}(u).
\end{eqnarray*}
Hence, by the definition of $E_{ij}$ and 
Lemma~\ref{le:group},
\begin{eqnarray*}
     \Nabla{A,t}(E_1v_{0s}) 
&= &
     \Nabla{t}(E_1v_{0s}) + \Nabla{E_1v_{0s}}X_\Psi(u) \\
&= &
     E_{11}(v_{0s},\p_tu_0)
     + E_{12}(v_{0s},\Nabla{t}\xi_0)  \\
&&
     +\, E_{11}(v_{0s},X_{\Psi_0}(u_0))
     + E_{12}(v_{0s},\Nabla{\xi_0}X_{\Psi_0}(u_0)) \\
&&
     +\, E_1\Nabla{t}v_{0s}  
     + \Nabla{E_1v_{0s}}X_{\psi_0}(u) 
     + E_1\Nabla{v_{0s}}X_{\Psi_0}(u_0)  \\
&= &
     \Nabla{E_1v_{0s}}X_{\psi_0}(u)
     + E_{11}(v_{0s},v_{0t})
     + E_{12}(v_{0s},\Nabla{A_0,t}\xi_0)
     + E_1\Nabla{A_0,t}v_{0s}
\end{eqnarray*}
and
\begin{eqnarray*}
     \Nabla{A,t}(E_2\Nabla{A_0,s}\xi_0) 
&= &
     \Nabla{t}(E_2\Nabla{A_0,s}\xi_0) 
     + \Nabla{E_2\Nabla{A_0,s}\xi_0}X_\Psi(u)  \\
&= &
     E_{21}(\Nabla{A_0,s}\xi_0,\p_tu_0)
     + E_{22}(\Nabla{A_0,s}\xi_0,\Nabla{t}\xi_0)
     + E_2\Nabla{t}\Nabla{A_0,s}\xi_0 \\
&&
     +\, \Nabla{E_2\Nabla{A_0,s}\xi_0}X_{\psi_0}(u) 
     + E_{21}(\Nabla{A_0,s}\xi_0,X_{\Psi_0}(u_0)) \\
&&
     +\, E_{22}(\Nabla{A_0,s}\xi_0,\Nabla{\xi_0}X_{\Psi_0}(u_0)) 
     + E_2\Nabla{\Nabla{A_0,s}\xi_0}X_{\Psi_0}(u_0)  \\
&= &
     \Nabla{E_2\Nabla{A_0,s}\xi_0}X_{\psi_0}(u) 
     + E_{21}(\Nabla{A_0,s}\xi_0,v_{0t})  \\
&&
     +\, E_{22}(\Nabla{A_0,s}\xi_0,\Nabla{A_0,t}\xi_0)
     + E_2\Nabla{A_0,t}\Nabla{A_0,s}\xi_0.
\end{eqnarray*}
In\-ser\-ting these two i\-den\-ti\-ties into the pre\-vious
for\-mula we obtain
\begin{eqnarray*}
     L_u\Nabla{A,t}\phi_0
&= &
     \Nabla{A,t}v_s
     - \Nabla{v_t}X_{\phi_0}(u) \\
&&
     - \Nabla{E_1v_{0s}}X_{\psi_0}(u)
     - \Nabla{E_2\Nabla{A_0,s}\xi_0}X_{\psi_0}(u) \\
&&
     - E_{11}(v_{0s},v_{0t})
     - E_{12}(v_{0s},\Nabla{A_0,t}\xi_0)   \\
&&
     - E_{21}(\Nabla{A_0,s}\xi_0,v_{0t}) 
     - E_{22}(\Nabla{A_0,s}\xi_0,\Nabla{A_0,t}\xi_0) \\
&&
     - E_1\Nabla{A_0,t}v_{0s} 
     - E_2\Nabla{A_0,t}\Nabla{A_0,s}\xi_0.
\end{eqnarray*}
Now the result follows from Lemma~\ref{le:Lu-phi}.
\end{proof}
 

%
 


\begin{thebibliography}{99}
\addcontentsline{toc}{section}{References}
 
\small
 
\bibitem[AB]{AB} M.F.~Atiyah and R.~Bott,
      The Yang--Mills equations over Riemann surfaces,
      {\it Phil. Trans. R. Soc. Lond. A} {\bf 308} (1982), 523--615.
 
\bibitem[CGS]{CGS}
      C.~Cieliebak, A.R.~Gaio, D.A.~Salamon, 
      $J$-homolorphic curves, moment maps and invariants of
      Hamiltonian group actions,
      {\it Int. Math. Res. Notes} {\bf 10} (2000), 831--882.

\bibitem[CGMS]{CGMS}
      C.~Cieliebak, R.~Gaio, I.~Mundet, D.A.~Salamon, 
      Invariants of Hamiltonian group actions, in preparation.

\bibitem[DK]{DK}
      S.K.~Donaldson and P.B.~Kronheimer,
      {\it The Geometry of Four-Manifolds},
      Oxford University Press, 1990.
 
\bibitem[DS1]{DS1}  S.~Dostoglou and D.A.~Salamon,
      Cauchy-Riemann operators, self-duality, and the spectral flow, 
      in {\it First European Congress of Mathematics,
      Volume I, Invited Lectures (Part 1)},
      edited by A.~Joseph, F.~Mignot, F.~Murat,
      B.~Prum, R.~Rentschler,  Birkh\"auser Verlag,
      Progress in Mathematics, {\bf Vol. 119},
      1994, pp. 511--545.
 
\bibitem[DS2]{DS2}  S.~Dostoglou and D.A.~Salamon,
      Self-dual instantons and holomorphic curves,
      {\it Annals of Mathematics} {\bf 139} (1994), 581--640.

\bibitem[Ga]{Ga} A.~R.~Gaio,
      $J$-holomorphic curves and moment maps,
      PhD thesis, University of Warwick, November 1999.

\bibitem[GT]{GT} D.~Gilbarg and N.S.~Trudinger,
      {\it Elliptic Partial Differential Equations of 
      the Second Order}, Springer 1983. 
 
\bibitem[Gr]{Gr} M.Gromov,
      Pseudo holomorphic curves in symplectic manifolds,
      {\it Invent. Math.} {\bf 82} (1985), 307--347.

\bibitem[K]{Kw} F.~Kirwan,
      {\it Cohomology of Quotients in Symplectic and Algebraic
      Geometry},
      Princeton University Press, 1984.
 
\bibitem[MS1]{MS1}D.~McDuff and D.A.~Salamon,
      {\it $J$-holomorphic Curves and Quantum Cohomology},
      University Lecture Series, {\bf 6},
      American Mathematical Society, Providence, RI, 1994;
      second edition, 2001.
 
\bibitem[MS2]{MS2}D.~McDuff and D.A.~Salamon,
      {\it Introduction to Symplectic Topology},
      Oxford University Press, 1995, 2nd edition 1998.
 
\bibitem[M]{Mu} I.~Mundet,
      Yang-Mills-Higgs theory for symplectic fibrations,
      PhD thesis, Madrid, April 1999. 

\bibitem[RT]{RT} Y.~Ruan and G.~Tian,
      A mathematical theory of quantum cohomology, 
      {\it J. Diff. Geom.} {\bf 42} (1995), 259--367.

\bibitem[S]{Sa}  D.A.~Salamon,
      Quantum products for mapping tori and the Atiyah-Floer conjecture,
      {\it Amer. Math. Soc. Transl.} {\bf 196} (1999), 199--235. 
      Revised in December 2000,
      http://www.math.ethz.ch/~salamon
 
\bibitem[U]{Uh}  K.~Uhlenbeck,
      Connections with $L^p$-Bounds on Curvature,
      {\it Commun. Math. Physics.} {\bf 83} (1982), 31-42.
 
\bibitem[W]{We} K.~Wehrheim,
      {\it Uhlenbeck Compactness},
      to appear in Birkh\"auser Verlag, Basel. 

\end{thebibliography}
\end{document}